\documentclass[reqno, 12pt]{amsart}

\usepackage{amsmath}
\usepackage{amsfonts}
\usepackage{amssymb}
\usepackage{amsthm}
\usepackage{mathrsfs}
\usepackage{bbm}
\usepackage{graphicx}
\usepackage{tikz}
\usetikzlibrary{calc,decorations.markings}
\usetikzlibrary{arrows.meta}
\usetikzlibrary{positioning}
\usetikzlibrary{cd, intersections, calc, decorations.pathmorphing, arrows, decorations.pathreplacing}
\usepackage{color}
\usepackage{comment}
\usepackage[normalem]{ulem}
\usepackage{soul}

\usepackage{epsfig,amssymb,amsmath,version}
\usepackage{amssymb,version,graphicx,fancybox,mathrsfs}
\usepackage{subfigure}
\usepackage{color}
\usepackage{stmaryrd}
\usepackage{multirow}
\usepackage{booktabs,siunitx}
\usepackage{booktabs}
\usepackage{multicol,epstopdf}
\usepackage{xypic}
\usepackage[]{graphicx}
\usepackage[all]{xy}
\usepackage{tikz,ifthen}
\usetikzlibrary{positioning, math, decorations.markings, arrows.meta}
\usetikzlibrary{calc,shapes,cd}
\usetikzlibrary{knots} %, cd,patterns,hobby} %hobby draws ugly paths
\usepgfmodule{decorations}
\allowdisplaybreaks
\usepackage{bm}

\usepackage{physics2}
\usephysicsmodule{ab}   
\usephysicsmodule{ab.braket}   
\usephysicsmodule{diagmat}  
\usephysicsmodule[showleft=2,showtop=2]{xmat}

\usepackage{adjustbox}

\usepackage{bmpsize}

\usepackage{mathtools}

\usepackage{enumerate}

\usepackage{subfiles}
\usepackage[colorlinks, linkcolor=blue, citecolor=red, urlcolor=cyan,pagebackref]{hyperref}

\allowdisplaybreaks

\numberwithin{equation}{section}

\usepackage{cleveref}
\crefname{theorem}{Theorem}{Theorems}
\crefname{cor}{Corollary}{Corollaries}
\crefname{lemma}{Lemma}{Lemmas}
\crefname{sublem}{Sublemma}{Sublemmas}
\crefname{prop}{Proposition}{Propositions}
\crefname{definition}{Definition}{Definitions}
\crefname{example}{Example}{Examples}
\crefname{claim}{Claim}{Claims}
\crefname{conj}{Conjecture}{Conjectures}
\crefname{conv}{Notation}{Notations}
\crefname{remark}{Remark}{Remarks}
\crefname{rmk}{Remark}{Remarks}
\crefname{prob}{Problem}{Problems}
\crefname{figure}{Figure}{Figures}
\crefname{table}{Table}{Tables}
\crefname{section}{Section}{Sections}
\crefname{subsection}{Section}{Sections}
\crefname{appendix}{Appendix}{Appendices}
\crefname{introthm}{Theorem}{Theorems}
\crefname{introcor}{Corollary}{Corollaries}
\crefname{introconj}{Conjecture}{Conjectures}

\newtheorem{theorem}{Theorem}[section]
\newtheorem{prop}[theorem]{Proposition}
\newtheorem{cor}[theorem]{Corollary}
\newtheorem{lemma}[theorem]{Lemma}

\newtheorem{introthm}{Theorem}
\newtheorem{introcor}[introthm]{Corollary}
\newtheorem{introconj}[introthm]{Conjecture}

\theoremstyle{definition}
\newtheorem{definition}[theorem]{Definition}

\newtheorem{conj}[theorem]{Conjecture}

\theoremstyle{remark}

\newtheorem{remark}[theorem]{Remark}

\newcommand{\beq}{\begin{equation}}
\newcommand{\eeq}{\end{equation}}

\tikzset{
    squigarrow/.style={-{Classical TikZ Rightarrow[length=4pt]}, decorate, decoration={snake, amplitude=1.8pt, pre length=2pt, post length=3pt}}
}

\newcommand{\bC}{\mathbb C}

\newcommand{\bN}{\mathbb N}

\newcommand{\bQ}{\mathbb Q}
\newcommand{\bR}{\mathbb R}

\newcommand{\bZ}{\mathbb Z}

\newcommand{\cA}{\mathcal A}
\newcommand{\cB}{\mathcal B}
\newcommand{\cC}{\mathcal C}

\newcommand{\cI}{\mathcal I}

\newcommand{\cN}{\mathcal N}

\newcommand{\cP}{\mathcal P}

\newcommand{\cT}{\mathcal T}

\newcommand{\End}{\operatorname{End}}

\DeclareMathOperator{\Aut}{{\rm Aut}}

\DeclareMathOperator{\tangle}{\mathsf{Tangle}}
\newcommand{\Jin}{J_{\rm in}}
\newcommand{\Jout}{J_{\rm out}}
\newcommand{\id}{\operatorname{id}}
\DeclareMathOperator{\Ker}{\operatorname{Ker}}

\DeclareMathOperator{\Gr}{\operatorname{Gr}}
\newcommand{\rot}{\operatorname{rot}}

\newcommand{\plim}[1][]{\mathop{\varprojlim}\limits_{#1}}

\newcommand{\abar}{\overline{a}}

\title[Graded quotients of HOMFLYPT skein algebras and a genus $0$ formality]{Quantized necklace Lie algebras, graded quotients of HOMFLYPT skein algebras, and formality for genus $0$ surfaces}

\author[Hiroaki Karuo]{Hiroaki Karuo}
\address{Hiroaki Karuo, Department of Mathematics, Gakushuin University, Mejiro, Toshima-ku, Tokyo, Japan.}
\email{hiroaki.karuo@gakushuin.ac.jp}

\author[Shunsuke Tsuji]{Shunsuke Tsuji}
\address{Shunsuke Tsuji, Department of Mathematical and Statistical Sciences Based on Modeling and Analysis, Meiji University, 4-21-1 Nakano, Nakano-ku, Tokyo, 164-8525, Japan}
\email{641x6700417@gmail.com}

\date{}

\begin{document}
\maketitle

We would like to extend our warmest congratulations to Professors Nariya Kawazumi and Tomotada Ohtsuki on the occasion of their 61st and 60th birthdays, respectively.

\begin{abstract}
In the classical setting, it is known that a graded quotient of the Goldman--Turaev Lie bialgebra is isomorphic to the necklace Lie bialgebra for any compact oriented surface with boundary. 
In this paper, we mainly deal with the HOMFLYPT skein algebra, which is a bi-quantization of the Goldman--Turaev Lie bialgebra.  
As their quantizations, and its graded quotient following Turaev's bi-quantization $h\neq \hbar$. 
For any compact oriented surface with boundary, our first main result is to give an explicit algebra isomorphism between the associated graded quotient and the quantized necklace algebra introduced by Schedler under bi-quantization $h\neq \hbar$ \cite{Sch05}. 
This implies that a basis of the graded quotient is given from the the quantized necklace Lie module. 
In particular, if we ignore the framing, we show that Turaev's Lie bialgebra structure of the HOMFLYPT skein algebra recovers that of the quantized necklace Lie bialgebra. 
As an application, for any surface of genus $0$, we show the formality of the HOMFLYPT skein algebra with respect to the multiplication and the height involution using the Kontsevich integral. This formality can be regarded as a quantization of the result of \cite{BN+25}. More precisely, we use a HOMFLYPT quotient obtained by imposing additional relations on the Conway quotient appearing in \cite{BN+25}, and our results include the corresponding version of the results of \cite{BN+25}. 
%i.e. this completely covers their result. 
Another application is to give simpler description and proof of the Dehn twist formula in the HOMFLYPT skein algebra given by \cite{Tsu18a}.
\end{abstract}

\setcounter{tocdepth}{2}
\tableofcontents

\section{Introduction}
\subsection{Background}
\paragraph{\textbf{The classical setting}} The Kashiwara--Vergne problem \cite{KV78} was formulated as a conjecture concerning the Baker--Campbell--Hausdorff (BCH) series in the free Lie algebra on two generators and was solved affirmatively by Alekseev--Meinrenken \cite{AM06}. It can be viewed as the problem of finding solutions to a system of equations. Subsequently, Alekseev--Torossian clarified its relation to Drinfeld associators \cite{AT12}, and Alekseev--Enriquez--Torossian gave an explicit construction of solutions to the Kashiwara--Vergne equations from associators \cite{AET10}.
%This construction is equivariant with respect to the embeddings of the Grothendieck--Teichm\"uller groups into the Kashiwara--Vergne groups,
%$\mathrm{GT}_1\hookrightarrow\mathrm{KV}$ and its graded counterpart $\mathrm{GRT}_1\hookrightarrow\mathrm{KRV}_2$ \cite{AET10}.
Through associators and Grothendieck--Teichm\"uller theory, the problem is also closely related to arithmetic objects such as the absolute Galois group and multiple zeta values.

While the above formulation is purely algebraic, the formality of the Goldman--Turaev Lie bialgebra provides a topological interpretation of the problem, as follows.
Let $\Sigma_{0,3}$ be a compact oriented surface of genus zero with three boundary components. Its fundamental group is the free group $F_2$ of rank two. The Goldman Lie bracket \cite{Gol86} and the Turaev cobracket \cite{Tur91}, defined using an intersection form, give rise to the Goldman--Turaev Lie bialgebra of $\Sigma_{0,3}$, denoted by $\bQ|\pi_1(\Sigma_{0,3})|$.
Here, $|\pi_1(\Sigma_{0,3})|$ denotes the cyclic quotient of $\pi_1(\Sigma_{0,3})$. The augmentation map $\varepsilon$ defines a descending filtration $\{|\ker\varepsilon|^n\}_n$ on the Goldman--Turaev Lie bialgebra. Both the completion
$\plim[n]\bQ|\pi_1(\Sigma_{0,3})|/|\ker\varepsilon|^n$ and the completed associated graded space $\widehat{\Gr}\bQ|\pi_1(\Sigma_{0,3})|:=\prod_{n=0}^{\infty}|\ker\varepsilon|^n/|\ker\varepsilon|^{n+1}$ carry Lie bialgebra structures induced by that on $\bQ|\pi_1(\Sigma_{0,3})|$.
A key point is that $\widehat{\Gr}\bQ|\pi_1(\Sigma_{0,3})|$ is isomorphic, as a Lie bialgebra, to a necklace Lie bialgebra with a simpler algebraic structure \cite{KK14,KK15, MT14}.
This necklace Lie bialgebra can be described in terms of $H_1(\Sigma;\bQ)$ and admits an explicit basis, making it amenable to computation.
The formality problem for the Goldman--Turaev Lie bialgebra asks for a filtered Lie bialgebra isomorphism from its completion to its completed associated graded space that induces the identity on the associated graded.
Alekseev--Kawazumi--Kuno--Naef interpreted the conditions KV1 and KV2 of the Kashiwara--Vergne equations in terms of formality of the Goldman Lie bracket and of the Goldman--Turaev Lie bialgebra structure, respectively. They established a one-to-one correspondence between the set of Lie bialgebra isomorphisms realizing formality, called homomorphic expansions, and the set of solutions to the Kashiwara--Vergne problem \cite{AKKN18}.
The Goldman--Turaev Lie bialgebra is also defined for surfaces of higher genus. Its formality thus provides a topological approach to extending the original Kashiwara--Vergne problem to higher genus \cite{AKKN18,AKKN20}.

\medskip
\paragraph{\textbf{The quantum setting}}
The Kontsevich integral is a powerful invariant of framed links and framed tangles \cite{Kon93}.
Drinfeld associators play an essential role in its algebraic formulation for parenthesized tangles. Massuyeau related loops on punctured disks to braids and used the Kontsevich integral to construct homomorphic expansions of the Goldman--Turaev Lie bialgebra \cite{Mas18}. Bar-Natan--Dancso--Hogan--Liu--Scherich gave a three-dimensional interpretation of the Turaev cobracket using the height involution, which reverses the height direction of tangles, and obtained another proof of formality for genus-zero surfaces via the Kontsevich integral \cite{BN+25}. These works realize the formality of the Goldman--Turaev Lie bialgebra by methods of quantum topology.

These results naturally lead to the question of whether the formality of the classical Lie bialgebra extends to the setting of quantum algebras. In this paper, we study this question using HOMFLYPT skein modules and algebras. These are generated by framed oriented links or tangles, subject to local skein relations. Our HOMFLYPT skein modules and algebras involve two quantum parameters $h,\hbar$ and an indeterminate $\rho$ (Definition~\ref{def:HOMFLYPT_skein}).

Turaev's biquantization provides the foundation for this approach \cite{Tur91}. By introducing skein relations that distinguish crossings between distinct components from self-crossings, Turaev constructed a skein algebra with a bialgebra structure encoding both the Goldman bracket and the Turaev cobracket. In a biquantization, the failure of commutativity of the product and the failure of cocommutativity of the coproduct are controlled by independent parameters $h,\hbar$; the Lie bracket and Lie cobracket are recovered from the corresponding first-order terms.

Other directions in the study of HOMFLYPT skein algebras include the following:
\begin{enumerate}
    \item The study of Torelli groups and homology cylinders, and the construction of invariants of $3$-manifolds \cite{Tsu18a,Tsu18b,Tsu20}.
    
    \item Relations between the HOMFLYPT skein algebra of the torus, the elliptic Hall algebra, and double affine Hecke algebras \cite{MS17,MS21}.
    
    \item Categorical formulations using skein categories \cite{Gar25}.
\end{enumerate} 
These developments demonstrate the importance of skein algebras as a link between low-dimensional topology and representation theory.

On the other hand, Schedler constructed a quantization of the necklace Lie bialgebra associated with a quiver as a Hopf algebra \cite{Sch05}. Ginzburg--Schedler gave a Moyal--Weyl description of this quantization and proved that it is isomorphic to Schedler's construction \cite{GS06}.
Chen--Huang--Liu--Zhang subsequently constructed a biquantization, in Turaev's sense, of the necklace Lie bialgebra of a doubled quiver \cite{CHLZ26}. Their construction recovers Schedler's quantization under a suitable specialization of parameters and treats the quantizations of the bracket and the cobracket independently using two parameters.

Ginzburg--Schedler constructed a Moyal--Weyl quantization of the necklace Lie bialgebra and proved that it is isomorphic to Schedler's quantization, thereby giving a natural interpretation of the latter \cite{GS06}.
More recently, a refinement of Schedler's quantization was given in \cite{CHLZ26}.

In view of these results, one expects that a suitable filtration on the HOMFLYPT skein algebra has an associated graded algebra isomorphic to the quantized necklace Lie bialgebra. This is a basic step toward extending classical formality to the quantum setting. In the quantum setting, however, changes in the number of components under crossing resolutions and the contributions of framing must be treated consistently with the degrees of the parameters $h,\hbar$. Moreover, the filtration must be compatible with both multiplication and comultiplication and must recover the Goldman--Turaev Lie bialgebra and its associated graded in the classical limit. Constructing a filtration satisfying these requirements and determining its associated graded algebra are therefore central problems of this paper.

\subsection{Main results}
Following the philosophy of Turaev's biquantization \cite{Tur91}, we study HOMFLYPT skein modules and algebras with $h$ and $\hbar$. 
As stated below, $h$ quantizes the Goldman bracket (Kawazumi--Kuno action) and $\hbar$ quantizes the Turaev cobracket (Kawazumi--Kuno coaction). 
For $h\neq\hbar$, we introduce a finite-type filtration on the HOMFLYPT skein modules and algebras of surfaces of higher genus and determine the algebraic structure of their associated graded objects. As an application, we prove formality for the HOMFLYPT skein modules and algebras of genus-zero surfaces. We also give a simpler proof of the Dehn twist formula for general surfaces.
To incorporate the action of the Goldman--Turaev Lie bialgebra on the fundamental group, we also consider a finite marking $J\subset\partial\Sigma$, $J$-tangles with endpoints in $J$, and the HOMFLYPT skein module generated by these tangles.

\medskip
\paragraph{\textbf{The algebraic structure of the associated graded of the HOMFLYPT skein algebra}}
We introduce a finite-type filtration $\{F^n\}_{n\in\bN}$ on the HOMFLYPT skein modules and algebras and consider the associated graded algebra $\Gr\widetilde{\cA}(\Sigma)=\bigoplus_n\Gr^n\widetilde{\cA}(\Sigma)$.
This filtration provides a quantization of the filtration $\{|\ker\varepsilon|^n\}$ on the Goldman--Turaev Lie bialgebra considered, for example, in \cite{AKKN18, AKKN17}.

Our main theorem is the following.
\begin{introthm}[Theorem~\ref{thm:isom}]\label{introthm:isomorphism_higher_genera}
Let $\Sigma$ be a compact oriented surface with non-empty boundary.
The map $\Gr\Phi\colon\widetilde{\cN}(\Sigma)\to\Gr\widetilde{\cA}(\Sigma)$ is a well-defined $\bQ[\rho,h,\hbar]$-algebra isomorphism. 
In other words, the quantized necklace algebra $\widetilde{\cN}(\Sigma)$ and the associated graded algebra $\Gr\widetilde{\cA}(\Sigma)$ of the HOMFLYPT skein algebra are isomorphic as $\bQ[\rho,h,\hbar]$-algebras.
\end{introthm}

Since the quantized necklace algebra $\widetilde{\cN}(\Sigma)$ has no zero divisors, we obtain the following corollary.
\begin{introcor}
Let $\Sigma$ be a compact oriented surface with non-empty boundary.
Then the completion $\widehat{\widetilde{\cA}(\Sigma)}:=\plim[n]\widetilde{\cA}(\Sigma)/F^n\widetilde{\cA}(\Sigma)$ has no zero divisors.
\end{introcor}

Our results lead to the following conjecture.
\begin{introconj}
Let $\Sigma$ be a compact oriented surface with non-empty boundary.
\begin{enumerate}
    \item We have $\bigcap_nF^n\widetilde{\cA}(\Sigma)=\{0\}$.
    \item More weakly, $\widetilde{\cA}(\Sigma)$ has no zero divisors.
\end{enumerate}
\end{introconj}

The key ingredients in the proof of Theorem~\ref{introthm:isomorphism_higher_genera} are an explicit description of the Lie action on the associated graded and the following theorem.
\begin{introthm}[Theorem~\ref{thm:h_times}, cf.\cite{Tur91} for (3)]\label{introthm:h_times}
Let $\Sigma=\Sigma_{g,b}$ with $g\geq0$ and $b\geq1$. Then the following hold.
\begin{enumerate}
    \item The module $\Gr\widetilde{\cA}(\Sigma)=\bigoplus_n\Gr^n\widetilde{\cA}(\Sigma)$ is free over $\bQ[\rho,h,\hbar]$.
    
    \item The $\bQ[\rho]$-module homomorphisms between the associated graded modules
    \begin{align*}
    \Gr^n\widetilde{\cA}(\Sigma)\to \Gr^{n+2}\widetilde{\cA}(\Sigma),\quad x\mapsto hx\\ 
    \Gr^n\widetilde{\cA}(\Sigma)\to \Gr^{n+2}\widetilde{\cA}(\Sigma),\quad x\mapsto \hbar x
    \end{align*}
    are injective.

    \item The $\bQ[\rho]$-module homomorphism
    \begin{align*}
    \Gr^n\widetilde{\cA}(\Sigma)/\eta\Gr^{n-2}\widetilde{\cA}(\Sigma)\to\Gr^{n+2}\widetilde{\cA}(\Sigma)/\eta\Gr^n\widetilde{\cA}(\Sigma),\quad x\mapsto\eta'x\ (\eta\in\{h,\hbar\})
    \end{align*}
    is injective, where $\{\eta,\eta'\}=\{h,\hbar\}$. We set $\Gr^r\widetilde{\cA}(\Sigma)=0$ for $r<0$.
\end{enumerate}
\end{introthm}

As illustrated by \cite{Tur91, DL16, HM06, MS17}, finding a basis of a HOMFLYPT skein module is difficult not only for general $3$-manifolds but even for trivial $I$-bundles over simple surfaces. The constructions in \cite{DL16, HM06, MS17} use arguments specific to the $3$-manifolds or surfaces under consideration. Since bases have played an essential role in the development of skein modules and algebras, understanding a basis of the associated graded of the HOMFLYPT skein algebra should provide a useful foundation for further study.
Theorem~\ref{introthm:isomorphism_higher_genera} allows us to obtain a basis of $\Gr\widetilde{\cA}(\Sigma)$ from a basis of the quantized necklace algebra.

\medskip
\paragraph{\textbf{Coalgebra structures on the HOMFLYPT skein algebra and the quantized necklace Lie bialgebra}}

Following Turaev, we specialize to $\rho=0$, so that the HOMFLYPT skein algebra carries Turaev's coalgebra structure in the sense of \cite{Tur91}.
In \cite{Gar25}, the Turaev coproduct is defined using skein categories, and it is expected to induce the coproduct on Schedler's quantized necklace algebra \cite{Sch05} under a suitable filtration.
We give a partial answer to this question. More precisely, for $\rho=0$, we show that the Turaev coproduct induces the coproduct of the quantized necklace algebra in \cite{Sch05} upon passage to the associated graded.

\begin{introthm}[Theorems~\ref{thm:coproduct_filtration} and~\ref{thm:coproduct_gr}]\label{introthm:coproduct_filtration}
Suppose that $\rho=0$. Then the Turaev coproduct preserves the finite-type filtration $\{F^n\widetilde{\cA}(\Sigma)\}_{n\in\bN}$. Namely, for $x\in F^n\widetilde{\cA}(\Sigma)$, we have $\Delta(x)\in\sum_{i+j=n}F^i\widetilde{\cA}(\Sigma)\otimes F^j\widetilde{\cA}(\Sigma)$.
Moreover, the induced coproduct on the associated graded admits an explicit description.
\end{introthm}
For $\rho=0$ and $h=1$, the coproduct in Theorem~\ref{introthm:coproduct_filtration} agrees with that in \cite{Sch05} (Remark~\ref{rmk:bialg}). Here, $\rho=0$ is required in \cite{Tur91}, whereas $h=1$ is the specialization used in \cite{Sch05}.

Furthermore, for $\rho=0$ and $\Sigma=\Sigma_{g,1}$, the Hopf algebra structure of the quantized necklace Lie bialgebra induces a Hopf algebra structure on $\Gr\widetilde{\cA}(\Sigma)$ via the bialgebra isomorphism.

\medskip
\paragraph{\textbf{Formality of HOMFLYPT skein modules and algebras and of the Goldman--Turaev Lie bialgebra for genus-zero oriented surfaces}}
For an oriented genus-zero surface $\Sigma$ with $\partial\Sigma\neq\emptyset$, we prove that the HOMFLYPT skein algebra, viewed as a quantization of the Goldman--Turaev Lie bialgebra, has formality with respect to several operations.
The relations of the quantized necklace algebra of $\Sigma$ are defined using the intersection form on $H_1(\Sigma,\bZ)$, so this algebra is commutative in genus zero.
Consequently, for genus-zero surfaces, the associated graded algebra with respect to $\{F^n\}_{n\in\bN}$ is commutative.
It therefore does not retain the noncommutativity of the original HOMFLYPT skein algebra. To study formality, we modify the filtration used in Theorem~\ref{introthm:isomorphism_higher_genera} and introduce the modified finite-type filtration $\{I^n\}_{n\in\bN}$. Its associated graded algebra $\Gr_I\widetilde{\cA}(\Sigma)$ retains this noncommutativity.
We also introduce the modified quantized necklace algebra $\widetilde{\cN}_I(\Sigma)$ corresponding to this filtration and prove that it is isomorphic to $\Gr_I\widetilde{\cA}(\Sigma)$ as an algebra, giving a modified version of Theorem~\ref{introthm:isomorphism_higher_genera}.
In particular, the modified quantized necklace algebra is a quantization of the necklace Lie bialgebra arising as the associated graded of the Goldman--Turaev Lie bialgebra of a genus-zero surface with respect to the augmentation-ideal filtration.

The following commutative diagram is central to this paper. We prove formality for HOMFLYPT skein modules and algebras by establishing property $(\ast)$ below.
\begin{equation}
\begin{tikzcd}
\widehat{\widetilde{\cA}(\Sigma,J)} 
\arrow[r] \ar[rd, dashed, "\text{expansion}"'] 
&\mathcal{D}_{\widetilde{\cA}}^J \ar[d]\ar[rr] &&\mathcal{D}_{\widetilde{\cN}}^J \ar[d]\\
&\widehat{\Gr}_I\widetilde{\cA}(\Sigma,J) 
&\widehat{\widetilde{\cN}_I(\Sigma,J)} \ar[l] \ar[r] 
&\widehat{\widetilde{\cN_a}(\Sigma,J)}
\end{tikzcd}\label{eq:commutative_intro}
\end{equation}
\paragraph{\textbf{Property $(\ast)$}(Section~\eqref{aaa})}
\begin{enumerate}
    \item If $|J|=0$, each module has an algebra structure, and the maps in the commutative diagram are isomorphisms.
    
    \item Let $J_1$ and $J_2$ be markings with $|J_1|=|J_2|=2$. If the incoming point of $J_1$ and the outgoing point of $J_2$ lie on the same boundary component, a composition is defined on the corresponding modules, and the maps in the commutative diagram preserve composition.
    
    \item The HOMFLYPT skein module $\widetilde{\cA}(\Sigma,J)$ is a $\widetilde{\cA}(\Sigma)$-bimodule.
    The other modules in the commutative diagram also carry bimodule structures, and the maps preserve the left and right actions.
\end{enumerate}

The modules (algebras if $J=\emptyset$) $\mathcal{D}_{\widetilde{\cA}}^J$ and $\mathcal{D}_{\widetilde{\cN}}^J$ in the diagram are quotients of spaces of chord diagrams with endpoints by specified relations, and $Z_{\widetilde{\cA}}^J$ is a homomorphism induced by the Kontsevich integral.
In \cite{BN+25}, formality is studied in the framework of $\mathcal{D}_{\widetilde{\cA}}^J$. A key feature of our approach is the use of $\mathcal{D}_{\widetilde{\cN}}^J$, or equivalently, the perspective of quantized necklace Lie bialgebras.

\begin{introthm}[Theorem~\ref{thm:isom_modified}]\label{introthm:isom_modified}
Let $\Sigma=\Sigma_{0,p+1}$, and let $J\subset\partial\Sigma_{0,p+1}$ with $|J|=0,2$.
For the modified quantized necklace Lie module $\widetilde{\cN_I}(\Sigma,J)$ and the modified associated graded module $\Gr_I\widetilde{\cA}(\Sigma,J)$ of the HOMFLYPT skein module,
the map $\Gr_I\Phi\colon\widehat{\widetilde{\cN_I}(\Sigma,J)}\to\widehat{\Gr_I}\widetilde{\cA}(\Sigma,J)$ is a well-defined $\bQ\llbracket\rho\hbar,h\hbar\rrbracket[\rho,h,\hbar]$-module isomorphism.
Thus the completed modified quantized necklace module $\widehat{\widetilde{\cN_I}(\Sigma,J)}$ and the completed modified associated graded module $\widehat{\Gr_I}\widetilde{\cA}(\Sigma,J)$ of the HOMFLYPT skein module are isomorphic as $\bQ[\rho,h,\hbar]$-modules. Moreover, this isomorphism satisfies property $(\ast)$.
\end{introthm}

\begin{introcor}[Corollary~\ref{cor:h_times_genus0}]\label{introcor:h_times_genus0}
Let $(\Sigma,J)$ be a marked surface with $\Sigma=\Sigma_{0,p+1}$ and $p\in\bN$. Then the following hold.
\begin{enumerate}
    \item The associated graded module $\Gr_I\widetilde{\cA}(\Sigma,J)$ and the modified quantized necklace module $\widetilde{\cN_I}(\Sigma,J)$ are isomorphic as $\bQ[\rho,h,\hbar]$-modules by a degree-preserving isomorphism. Moreover, this isomorphism satisfies property $(\ast)$.
    
    \item The $\bQ[\rho]$-module homomorphisms between the modified associated graded modules %multiplication by $h$ and by $\hbar$
\begin{align*}
    \Gr_{I}^n\widetilde{\cA}(\Sigma,J)\to \Gr^{n+1}_{I}\widetilde{\cA}(\Sigma,J),\quad x\mapsto hx\\ 
    \Gr_{I}^n\widetilde{\cA}(\Sigma,J)\to \Gr_{I}^{n+1}\widetilde{\cA}(\Sigma,J),\quad x\mapsto \hbar x
\end{align*}
are injective.
\end{enumerate}
\end{introcor}

\begin{introthm}[Theorem~\ref{thm:formality}]\label{introthm:formality}
With respect to the modified finite-type filtration $\{I^n\}_{n\in\bN}$ on $\widetilde{\cA}(\Sigma,J)$,
there is an isomorphism from the completion $\widehat{\widetilde{\cA}(\Sigma,J)}$ to $\widehat{\Gr_I}\widetilde{\cA}(\Sigma,J)$ realizing formality with respect to multiplication, composition, and the left and right actions, whenever these operations are defined. In particular, it is also compatible with the Lie action $\sigma$.
\end{introthm}

Since the classical limit of $\sigma$ recovers the Goldman Lie bracket, the formality with respect to $\sigma$ in Theorem~\ref{introthm:formality} can be viewed as a quantization of KV1.

Theorem~\ref{introthm:isom_modified} implies the following corollary.
\begin{introcor}[Corollary~\ref{cor:isom_completions}]\label{introcor:isom_completions}
For $\Sigma=\Sigma_{0,p+1}$ with $p\in\bN$, the completed quantized necklace algebra $\widehat{\widetilde{\cN}_I(\Sigma)}$, the completed associated graded algebra $\widehat{\Gr_I}\widetilde{\cA}(\Sigma)$ of the HOMFLYPT skein algebra, and its completion $\widehat{\widetilde{\cA}(\Sigma)}$ are isomorphic as $\bQ[\rho]\llbracket h,\hbar\rrbracket$-algebras.
\end{introcor}

The formality of HOMFLYPT skein modules and algebras for genus-zero surfaces shows that a basis of $\Gr_I\widetilde{\cA}(\Sigma)$ gives a topological basis of its completion $\widehat{\Gr_I\widetilde{\cA}(\Sigma)}$. Through the formality isomorphism, it also yields a topological basis of $\widehat{\widetilde{\cA}(\Sigma)}$.

From the perspective of HOMFLYPT skein modules and algebras, we recover the results of \cite{AKKN18} and \cite{BN+25}.
\begin{introthm}[Theorem~\ref{thm:formality_Goldman_Turaev}, cf. {\cite[Theorem~5.2]{BN+25}}]\label{introthm:formality_Goldman_Turaev}
Let $\Sigma=\Sigma_{0,p+1}$ and set $h=\hbar=0$. Under \eqref{eq:commutative_intro}, the map $Z_{\widetilde{\cA}}$ induces a module homomorphism from the Goldman--Turaev Lie bialgebra to the necklace Lie bialgebra.
The induced homomorphism is a homomorphic expansion preserving both the Goldman bracket and the Turaev cobracket. It is also a homomorphic expansion with respect to the Kawazumi--Kuno action, which is the classical limit of our action $\sigma$.
\end{introthm}
The proof uses the height involution on HOMFLYPT skein modules and algebras. Generalizing \cite{BN+25}, this can be viewed as a quantization of the Kawazumi--Kuno coaction in terms of the height involution.
Since the coaction is described by the height involution, it is invariant under the mapping class group. This may be useful in describing obstructions to the surjectivity of Johnson homomorphisms.

In \eqref{eq:difference_mu+}, we give an explicit description of the Kawazumi--Kuno coaction under the formality isomorphism. In particular, \eqref{eq:difference_mu+} agrees with Eq.~(47) of \cite{AKKN18}.
The Duflo function appears naturally in the coefficients of this description. It is known that the Goldman--Turaev Lie bialgebra does not have formality with respect to the Kawazumi--Kuno coaction \cite[Theorem~7.6 (i)]{AKKN18}.
We describe this obstruction explicitly using associators, by a method different from that of \cite{AKKN18}.

In view of the identification in \cite{AKKN18} of formality of the Goldman--Turaev Lie bialgebra structure with KV1 and KV2, Theorem~\ref{introthm:formality_Goldman_Turaev} can be regarded as a quantization of KV1 and KV2.

\medskip
\paragraph{\textbf{The Dehn twist formula}}
Using our HOMFLYPT quotient and the formality of the HOMFLYPT skein algebra, we give a simpler proof of the existence of $\sigma(L_{\cA})$, which represents the logarithm of a Dehn twist in \cite{Tsu18a} (Lemma~\ref{lem:Dehn_twist_annulus}).
An explicit formula for $L_{\cA}$ is given in \cite{Tsu18a}. Although this formula does not follow immediately from formality for genus-zero surfaces, we show that it can be recovered from uniqueness (Theorem~\ref{thm:Dehn_twist_general}).
In terms of the modified quantized necklace algebra, the element corresponding to $L_{\cA}$ is represented by a necklace with two beads.
We expect that this alternative proof using formality will lead to refinements of S.T.'s results based on the Dehn twist formula. More specifically, we are working toward an explicit connection between the LMO functor and skein algebras.

\subsection{Further directions}
We consider the following three directions for further research.

\paragraph{\textbf{Formality for higher-genus surfaces}}
In ongoing work, building on Theorem~\ref{introthm:isomorphism_higher_genera} for higher-genus surfaces, we seek to establish formality by constructing an isomorphism between the completion of the HOMFLYPT skein module or algebra and that of the quantized necklace algebra. This can be viewed as a quantization of part of the results of \cite{AKKN17}. We have already obtained partial results in this direction.
This approach uses the isomorphism of Theorem~\ref{introthm:isomorphism_higher_genera}, which is also essential to our genus-zero formality theorem, but the method for proving formality is different.
We therefore expect an interpretation of the higher-genus Kashiwara--Vergne problem of \cite{AKKN17} in terms of HOMFLYPT skein modules and algebras.

\medskip
\paragraph{\textbf{Invariants of homology cylinders}}
In a forthcoming paper \cite{KT26+}, we construct an invariant of homology cylinders with values in the HOMFLYPT skein algebra. This allows us to construct quantizations of Johnson homomorphisms by the following procedure.
Johnson homomorphisms can be interpreted as composition of an embedding of the Torelli group into the Goldman Lie algebra with a projection from the Goldman Lie algebra to its associated graded, the necklace Lie algebra \cite{MT13,MT14,KK14,KK15}. This gives a topological interpretation of the algebraically defined Johnson homomorphisms.
Since the monoid of homology cylinders contains the Torelli group as a submonoid, we obtain an invariant of the Torelli group with values in the HOMFLYPT skein algebra. Composing this invariant with the isomorphism of Theorem~\ref{introthm:isomorphism_higher_genera} defines a quantization of Johnson homomorphisms. Its classical limit recovers the embedding into the Goldman Lie algebra, namely, the Johnson homomorphism.

The explicit description of the associated graded structure allows a detailed analysis of the quantized Johnson homomorphisms and a finer understanding of where their classical images lie in the filtration induced by the augmentation ideal. For example, we show that the Casson core can be described as a map to the coefficient ring, corresponding to the trivial-necklace part. Thus understanding the algebraic structures of $\widetilde{\cA}(\Sigma)$ and $\Gr\widetilde{\cA}(\Sigma)$ is an important step toward extracting further information from these invariants.

\medskip
\paragraph{\textbf{The formality problem for HOMFLYPT skein modules and algebras}}
In Section~\ref{sec:formality_genus0}, we start with the Kontsevich integral associated with a Drinfeld associator
and pass through the HOMFLYPT quotient to prove formality for genus-zero HOMFLYPT skein modules and algebras.
This isomorphism recovers the formality of the Goldman--Turaev Lie bialgebra.
This raises the question of which information in an associator is needed to establish formality in the HOMFLYPT skein setting.
In the genus-zero setting above, expansions realizing the formality of the Goldman--Turaev Lie bialgebra
are in one-to-one correspondence with solutions to the Kashiwara--Vergne problem \cite{AKKN18}.
Kuno weakened the associator equations by focusing on a first-order approximation, called emergent braids, obtained by quotienting the group algebra of pure braids on a punctured disk by the square of the two-sided ideal generated by differences across crossing changes of tangles \cite{Kun25}. This approach is closely related to the viewpoint of HOMFLYPT skein modules and algebras.
It is also important to use solutions of the weakened equations to identify obstructions to formality of HOMFLYPT skein modules and algebras and necessary conditions on associators that realize formality.

%In Section~\ref{sec:formality_genus0}, we establish formality for HOMFLYPT skein modules and algebras for solutions induced by Drinfeld associators and recover the formality of the Goldman--Turaev Lie bialgebra as a corollary.
%In view of \cite{AKKN18}, this can be interpreted as another proof of the existence of solutions to the Kashiwara--Vergne problem.
%A further problem is to impose the HOMFLYPT relations connecting positive crossings, negative crossings, and smoothings before taking the first-order approximation with respect to crossing changes, and to determine which conditions of the pentagon and hexagon equations survive in the HOMFLYPT quotient.

\subsection{Organization of the paper}
In Section~\ref{sec;HOMFLYPT}, we review HOMFLYPT skein modules and algebras for surfaces of arbitrary genus and describe quantized necklace algebras in terms of surface homology.
We also introduce the finite-type filtration, the finite-type degree, and the component degree.
In Sections~\ref{sec;isom} and~\ref{sec:proof_mainthm}, we prove that the associated graded algebra of the HOMFLYPT skein algebra with respect to the finite-type filtration is isomorphic to the quantized necklace algebra. In particular, an explicit description of the Lie action $\sigma$ is used to prove injectivity of multiplication by $h$ and $\hbar$ and injectivity of the map defining the algebra isomorphism.
In Section~\ref{sec:coalgebra}, we show that Turaev's coalgebra structure on the HOMFLYPT skein algebra recovers Schedler's coalgebra structure through the above isomorphism, and explain its interpretation in terms of the component degree.
In Section~\ref{sec:formality_genus0}, we modify the finite-type filtration for genus-zero surfaces and prove that the resulting associated graded modules and algebras of the HOMFLYPT skein modules and algebras are isomorphic to the corresponding quantized necklace modules and algebras. We then combine this isomorphism with the Kontsevich integral to prove formality for HOMFLYPT skein modules and algebras. As an application, we show that this formality recovers that of the Goldman--Turaev Lie bialgebra.
As another application, we give an alternative proof of the Dehn twist formula in \cite{Tsu18a}.

Since \cite{Tsu18a} is unpublished, we give simpler proofs of some of the results in \cite{Tsu18a} to make this paper as self-contained as possible.
Although the construction is technically involved, it is self-contained within this paper (Remark~\ref{rmk:tautological}).

\medskip
\paragraph{\textbf{Use of AI}}
We used ChatGPT to strengthen the introduction, to check the Japanese version of the manuscript for typos and to help translate it into English. 
We carefully proofread the English version ourselves and used AI for no other purpose.

\medskip
\paragraph{\textbf{Acknowledgement}} 
We are grateful to Hisatoshi Kodani for leading us to the formality problem of the HOMFLYPT skein algebra in the genus 0 case through reading \cite{BN+25} together. 
H.K. was supported by JSPS KAKENHI Grant Number JP26K06786.
S.T. was supported by JSPS KAKENHI Grant Number 22K13918.

%%%%%%%%%%%%%%%%%%%%%%%%%%%%%%%%%%%%%%%%%%%%%%%%%%%%%%%%%
\section{HOMFLYPT skein algebras and quantized necklace Lie algebras}\label{sec;HOMFLYPT}
\paragraph{\textbf{Notation}}
In this paper, let $\bN$ denote the set of non-negative integers. 

\subsection{HOMFLYPT skein modules and algebras}

Let $M$ be an oriented $3$-manifold with boundary, and consider an embedding $\bigsqcup[0,1]\to\partial M$.
Denote its image by $J$, with each component oriented from $0$ to $1$.
We call the pair $(M,J)$ a \textbf{marked $3$-manifold}.
In particular, if $\Sigma$ is a compact oriented surface with boundary and $J\subset\partial\Sigma$ is a finite subset, then $(\Sigma\times[0,1],J\times[0,1])$ is a marked $3$-manifold.
We call $(\Sigma,J)$ a \textbf{marked surface}.
Each component of $J\times[0,1]$ is oriented in the direction of increasing $[0,1]$-coordinate.

Assume that $J=\Jout\sqcup\Jin$, that $m:=|\pi_0(\Jout)|=|\pi_0(\Jin)|<\infty$, and fix a bijection $\pi_0(\Jout)\to\pi_0(\Jin)$. For convenience, number the connected components of $\Jout$ by $1,2,\dots,m$ and those of $\Jin$ by $m+1,m+2,\dots,2m$, and identify each component with its index. We choose the numbering so that the fixed bijection sends $i$ to $m+i$ for $i\in\{1,2,\dots,m\}$.
Throughout this paper, we consider only markings $J=\Jout\sqcup\Jin$ equipped with these data.

A compact connected $1$-dimensional submanifold of a $3$-manifold is homeomorphic to either $S^1$ or $[0,1]$.
A component of a $1$-manifold is called a \textbf{closed component} if it is homeomorphic to $S^1$, and an \textbf{arc component} if it is homeomorphic to $[0,1]$.
\begin{definition}
A \textbf{$J$-tangle} $\alpha$ is a properly embedded, framed, oriented compact $1$-dimensional submanifold of $M$ satisfying the following conditions:
\begin{enumerate}
    \item $\partial\alpha\subset\operatorname{Int}J$, and $|\alpha\cap j|=1$ for every component $j$ of $J$,
    
    \item Each arc component of $\alpha$ joins a point of $\Jout$ to a point of $\Jin$ and is oriented from $\Jout$ to $\Jin$,
    
    \item At the point $\alpha\cap j$ on each component $j$ of $J$, the framing vector is tangent to $j$ and agrees with its orientation.
\end{enumerate}
A $J$-tangle has exactly $m$ arc components and may also have closed components.
When $J=\emptyset$, we call a $J$-tangle an \textbf{oriented link}.

For a $J$-tangle $T$, join its endpoints on components $i$ and $m+i$ of $J$ for each $1\leq i\leq m$. The resulting abstract $1$-manifold is called the \textbf{virtual closure} of $T$ and is denoted by $\check T$. Set $\#\check T:=|\pi_0(\check T)|$.
In general, $\#\check T\neq|\pi_0(T)|$.
\begin{figure}
	\centering\includegraphics[scale=0.15]{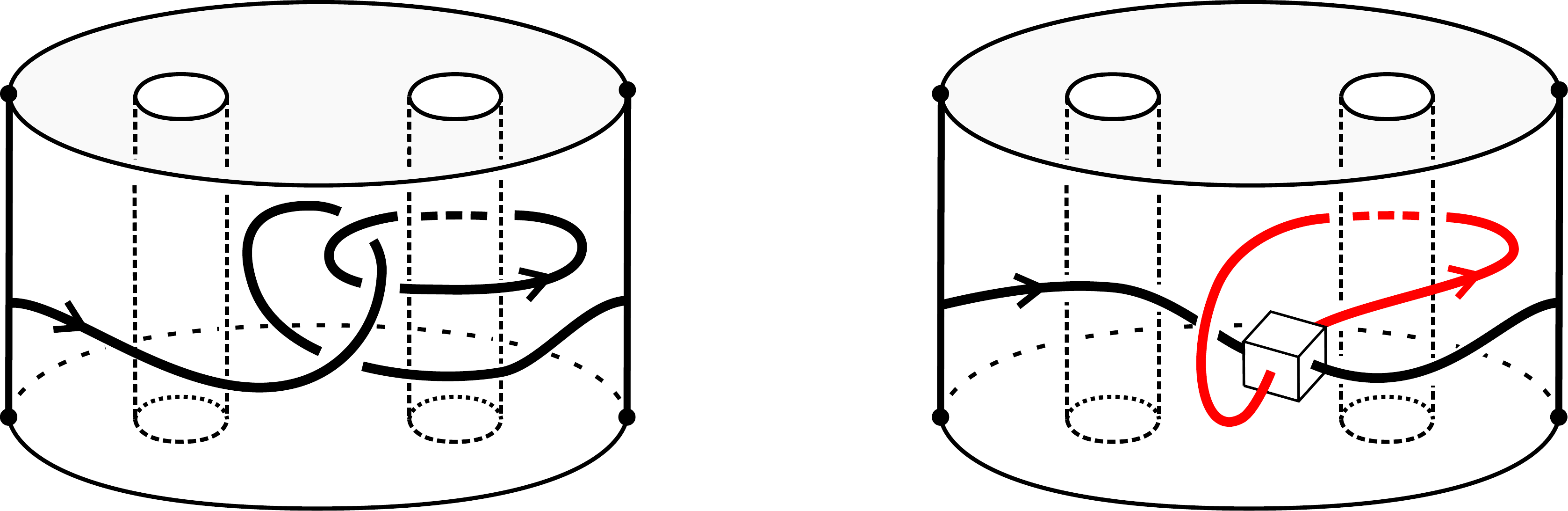}
    \caption{Left: a $J$-tangle for $\Sigma=\Sigma_{0,3}$. Right: a $J$-tangle with boxes for $\Sigma=\Sigma_{0,3}$; see Section~\ref{subsec:finite_type_filtration}. The tangles are equipped with a vertical framing.}
	\label{fig:tangle}   
\end{figure}

A \textbf{$J$-relative isotopy} is an ambient isotopy preserving $J$ setwise. Two $J$-tangles are \textbf{$J$-relative isotopic} if they are isotopic in the class of $J$-tangles.
When $M=\Sigma\times[0,1]$, every $J$-tangle is $J$-relative isotopic to one with a vertical framing, where a $J$-tangle $\alpha$ has a \textbf{vertical framing} if its framing vector at every point points in the direction of increasing $[0,1]$-coordinate.

Let $\tangle(M,J)$ be the set of $J$-relative isotopy classes of $J$-tangles in $(M,J)$.
When $M=\Sigma\times[0,1]$ and $J=J_\Sigma\times[0,1]\subset\partial\Sigma\times[0,1]$, where $J_\Sigma$ is a finite subset of $\partial\Sigma$, we also use $\tangle(\Sigma,J_\Sigma)$ for $\tangle(M,J)$. Set $\cT(M,J):=\bZ\tangle(M,J)$.

Suppose that $M=\Sigma\times[0,1]$ and $J=J_\Sigma\times[0,1]\subset\partial\Sigma\times[0,1]$, where $J_\Sigma$ is a finite subset of $\partial\Sigma$. Isotope a $J$-tangle $\alpha$ so that it has a vertical framing, and consider its image under the projection $\Sigma\times[0,1]\to\Sigma\times\{1/2\}$. 
After a further isotopy if necessary, assume that the projection has only transverse double points. 
At each double point, assign crossing information according to the heights in $[0,1]$. 
The resulting diagram is called a \textbf{$J$-tangle diagram} of $\alpha$. When $J=\emptyset$, it is called a \textbf{link diagram}.

Two $J$-tangles are $J$-relative isotopic if and only if their diagrams are related by a finite sequence of the framed Reidemeister moves in Figure~\ref{fig:Reidemeister} and $\partial$-relative isotopies of $(\Sigma,J_\Sigma)$.
\begin{figure}
	\centering\includegraphics[width=14cm]{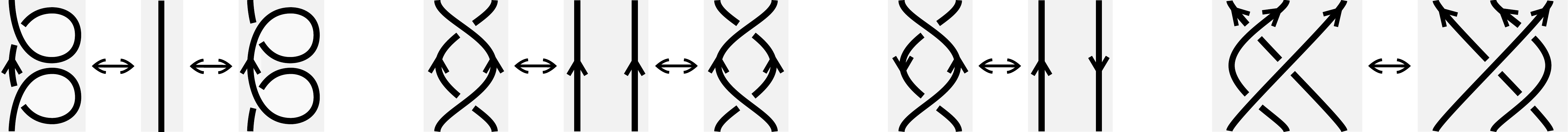}
    \caption{Framed oriented Reidemeister moves. Moves with all other possible orientations follow from combinations of those shown. The left move, the two middle moves, and the right move are called framed Reidemeister moves of types I, II, and III, respectively.}
	\label{fig:Reidemeister} 
\end{figure}
\end{definition}

\begin{remark}
Since each component of $J$ contains at most one endpoint of a $J$-tangle, our tangles are more restricted than those used in stated skein algebras or Muller skein algebras.
\end{remark}

Let $\bQ\llbracket\rho\hbar,h\hbar\rrbracket$ be the completion of the polynomial ring over $\bQ$ generated by $\rho\hbar$ and $h\hbar$. 
Denote $\bQ\llbracket\rho\hbar,h\hbar\rrbracket[\rho,h,\hbar]$ for the polynomial ring in $\rho,h,\hbar$ over $\bQ\llbracket\rho\hbar,h\hbar\rrbracket$.
\begin{definition}\label{def:HOMFLYPT_skein}
For a marked $3$-manifold $(M,J)$, consider the quotient of the free $\bQ\llbracket\rho\hbar,h\hbar\rrbracket[\rho,h,\hbar]$-module generated by isotopy classes of oriented framed $J$-tangles in $M$ by the relations (A)--(D).
This module is called the \textbf{HOMFLYPT skein module} denoted by $\widetilde{\cA}(M,J)$.

For a surface $\Sigma$, 
$\widetilde{\cA}(\Sigma\times[0,1],\emptyset)$ has an algebraic structure with the multiplication defined by stacking with respect to $[0,1]$. The algebra $\widetilde{\cA}(\Sigma):=\widetilde{\cA}(\Sigma\times[0,1],\emptyset)$ is called the \textbf{HOMFLYPT skein algebra}. 
\end{definition}
\begin{align*}
&{\rm (A)} 
\begin{array}{c}\includegraphics[scale=0.20]{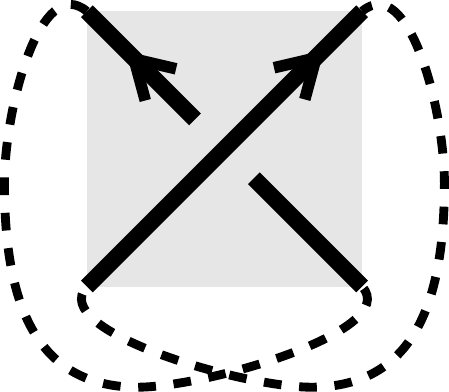}\end{array}-\begin{array}{c}\includegraphics[scale=0.20]{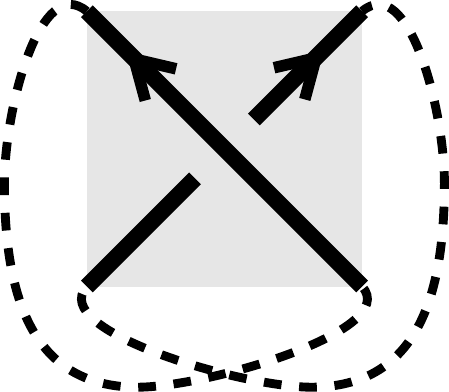}\end{array}=h\begin{array}{c}\includegraphics[scale=0.20]{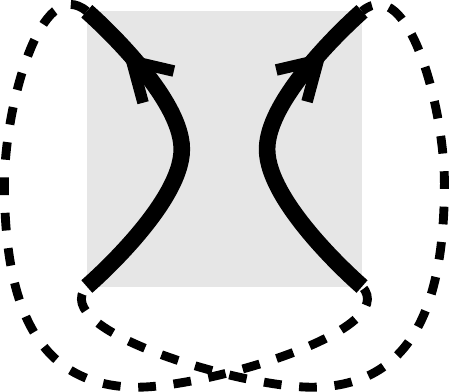}\end{array},\quad
&&{\rm (B)} 
\begin{array}{c}\includegraphics[scale=0.20]{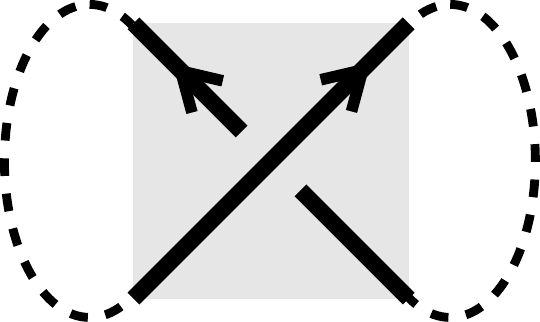}\end{array}-\begin{array}{c}\includegraphics[scale=0.20]{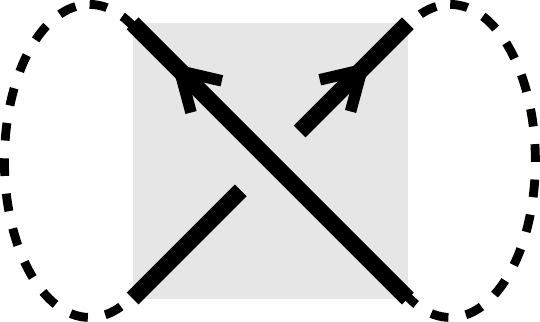}\end{array}=\hbar\begin{array}{c}\includegraphics[scale=0.20]{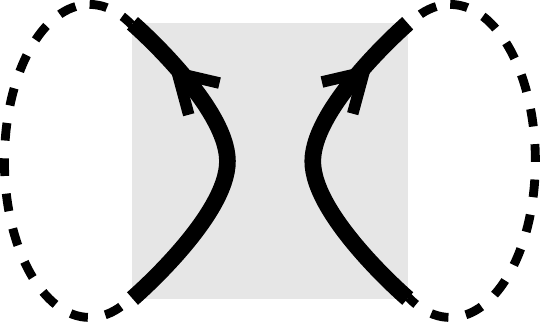}\end{array},\\
&{\rm (C)} 
\begin{array}{c}\includegraphics[scale=0.18]{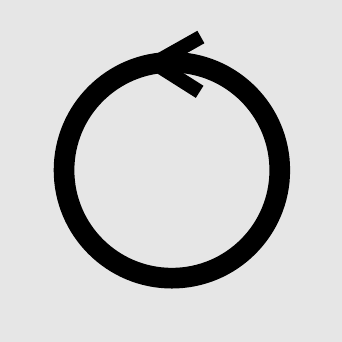}\end{array}=\dfrac{2\sinh{(\rho \hbar)}}{\hbar}\begin{array}{c}\includegraphics[scale=0.18]{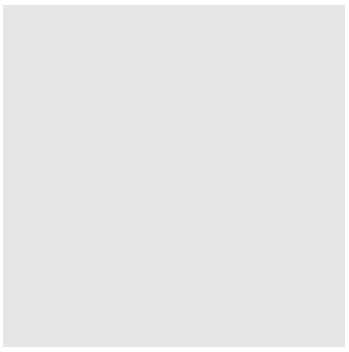}\end{array},\quad
&&{\rm (D)} 
\begin{array}{c}\includegraphics[scale=0.18]{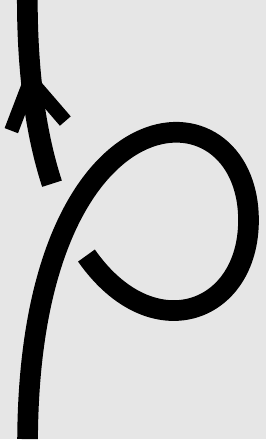}\end{array}=\exp(\rho \hbar)\begin{array}{c}\includegraphics[scale=0.18]{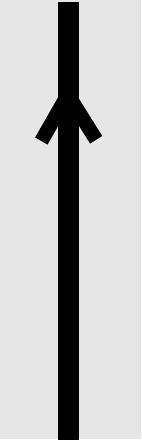}\end{array}. 
\end{align*}
The dotted lines indicate which endpoints are connected in the virtual closure $\check T$. Here, $\sinh x=(e^x-e^{-x})/2$, and
$\dfrac{2\sinh(\rho\hbar)}{\hbar}$ is regarded as a formal power series containing no negative powers of $\hbar$.

\begin{remark}\label{rem;HOMFLYPT_properties}
\begin{enumerate}
    \item Our definition is a specialization of the more general HOMFLYPT skein module: the coefficients on the LHS's of (A) and (B) are specialized to be $\pm1$.
    
    \item In \cite{Tur91}, the notation $\cA$ is used for $\widetilde{\cA}$.
    In this paper, we use $\cA$ for the specialization $h=\hbar$ and use $\widetilde{\cA}$ to distinguish the two-parameter version.
    
    \item 
    For $J$-tangles, we have
    \begin{align*}
    \#\begin{array}{c}\includegraphics[scale=0.20]{draws/crossing_posi_2.pdf}\end{array}&=\#\begin{array}{c}\includegraphics[scale=0.20]{draws/crossing_nega_2.pdf}\end{array}=\#\begin{array}{c}\includegraphics[scale=0.20]{draws/resolution_1.pdf}\end{array}+1, \\
    \#\begin{array}{c}\includegraphics[scale=0.20]{draws/crossing_posi_1.pdf}\end{array}&=\#\begin{array}{c}\includegraphics[scale=0.20]{draws/crossing_nega_1.pdf}\end{array}=\#\begin{array}{c}\includegraphics[scale=0.20]{draws/resolution_2.pdf}\end{array}-1. 
    \end{align*}
    Here, the dotted lines indicate the connections in the virtual closure of the $J$-tangle, and $\#$ denotes the number of components of that closure.
    
    \item When $J$-tangle diagrams on a marked surface $(\Sigma,J)$ are regarded as elements of $\widetilde{\cA}(\Sigma,J)$, the framed oriented Reidemeister moves in Figure~\ref{fig:Reidemeister} hold.
\end{enumerate}
\end{remark}

For $\ell\in\bZ$, let $\widetilde{\cA}^{(\ell)}(M,J)$ be the $\bQ\llbracket\rho\hbar,h\hbar\rrbracket$-submodule of $\widetilde{\cA}(M,J)$ generated by
$$
\{h^i \hbar^j \rho^k T \mid T\in \tangle(M,J) \text{ and }
\# \check{T} + i-j+k=\ell\}.
$$

The quantity
$\#\check T+i-j+k$ is preserved by relations (A)--(D). We call it the \textbf{component degree},
and call $\widetilde{\cA}^{(\ell)}(M,J)$ the \textbf{component degree $\ell$-part} of $\widetilde{\cA}(M,J)$.
\begin{remark}
The component degree $\ell$-part of $\widetilde{\cA}(M,J)$ corresponds to the pullback under $\omega$ of the $y^\ell$-part, or degree $\ell$ part, in \cite[Section~5.3]{Tur91}.
\end{remark}

The following proposition extends the statement in \cite[Section~5.3]{Tur91} to tangles.
\begin{prop}As $\bQ\llbracket\rho\hbar,h\hbar\rrbracket$-modules, we have
$\widetilde{\cA}(M,J)=\bigoplus_{\ell\in\bZ}\widetilde{\cA}^{(\ell)}(M,J)$.
\end{prop}
\begin{proof}
It suffices to show that the defining relations are homogeneous with respect to the component degree.
From Remark~\ref{rem;HOMFLYPT_properties}, the relations (A) and (B) preserve the component degree.
In (C), the RHS has one fewer component, and its coefficient is $2\rho$ times a power series in $\rho\hbar$. Thus, (C) preserves the component degree.
The coefficient on the RHS of (D) is a power series in $\rho\hbar$, so (D) also preserves the component degree.
\end{proof}

\subsection{A filtration on $\widetilde{\cA}(M,J)$}\label{subsec:finite_type_filtration}
We now define a finite-type filtration on $\widetilde{\cA}(M,J)$. For this, we first introduce tangles with boxes using proper immersions of $1$-manifolds into $3$-manifolds.

Let $L$ be the image of a framed oriented immersion of a disjoint union of circles $S^1$ and intervals $[0,1]$ into a $3$-manifold $M$. Each circle is designated as either a \textbf{$d$-component} or an \textbf{$s$-component}, while each interval is an $s$-component. We use the same terms for their images.
We require that $L$ contains at least one $s$-component.
The allowed singularities are transverse double points.
We further impose the following conditions:
\begin{enumerate}
    \item Every double point is formed by a $d$-component and an $s$-component,
    \item Each $d$-component contains exactly one double point.
\end{enumerate}
For each double point $p$ of $L$, identify an open $3$-ball in $M$ containing $p$ with $B_2=\{(x,y,z)\in\bR^3\mid -2<x,y,z<2\}$, taking $p$ to $(0,0,0)$. Assume that the framing on $L\cap B_1$ points in the positive $z$-direction, where $B_1$ is defined as follows.
Inside $B_2$, consider the box $B_1=\{(x,y,z)\in\bR^3\mid -1\leq x,y,z\leq1\}$.
Each of the four faces of $B_1$ that are not parallel to the $xy$-plane is required to meet $L$ transversely in exactly one point. Opposite faces meet components of the same type.
Replacing every double point by such a box gives a \textbf{$J$-tangle with boxes}.
By convention, an ordinary $J$-tangle is also regarded as a $J$-tangle with boxes, with no boxes.
Every $d$-component $\alpha$ of a $J$-tangle with boxes meets a unique box, which we call the \textbf{box of $\alpha$}; see the right-hand diagram in Figure~\ref{fig:tangle}.

The definition of $J$-relative isotopy extends to $J$-tangles with boxes by requiring the isotopy to remain within this class.
Let $\mathsf{Box}(M,J)$ be the set of $J$-relative isotopy classes of tangles with boxes in the marked $3$-manifold $(M,J)$. When $J=\emptyset$, write $\mathsf{Box}(M):=\mathsf{Box}(M,J)$.

Suppose that $M=\Sigma\times[0,1]$.
A $J$-tangle with boxes is called \textbf{vertical} if it satisfies the following conditions:
\begin{enumerate}
    \item Its framing vector at each point points in the direction of increasing $[0,1]$-coordinate.
    \item Under the projection $\Sigma\times[0,1]\to\Sigma$, each box $B_1$ projects to the square $\{(x,y,0)\in\bR^3\mid -1\leq x,y\leq1\}$.
    \item For each point $(x,y,0)$ of this square, the fiber of the projection restricted to $B_1$ is the interval $\{(x,y,z)\in\bR^3\mid -1\leq z\leq1\}$ with $x$ and $y$ fixed.
\end{enumerate}
Under the projection $\Sigma\times[0,1]\to\Sigma$, a neighborhood of a box is represented on $\Sigma$ by
$\begin{array}{c}\includegraphics[scale=0.20]{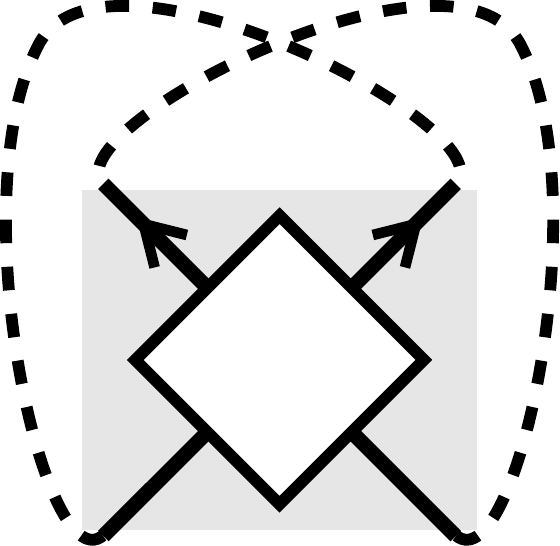}\end{array}$. The shaded region depicts $B_1\subset\bR^3$ viewed from the positive $z$-direction, and the tangle has a vertical framing pointing toward the reader.
To distinguish the two types of components, we draw $d$-components in red. For a general $M$, a vertical framing can still be considered locally.

For a diagram, replace each box $\Box$ with the linear combination obtained by subtracting the second of the following tangles from the first:
\begin{enumerate}
    \item[\textbf{detour}] The tangle with boxes obtained by replacing $\Box$ with the unique smoothing compatible with the orientations and with trivial framing, whose framing vectors point in the positive $z$-direction in $B_2$, then decreasing the number of double points.
    \item[\textbf{straight}] The tangle with boxes obtained by removing the $d$-component meeting $\Box$.
\end{enumerate}
We call these the detour and straight terms at $\Box$, respectively, and call their difference the \textbf{ds-linear combination} at $\Box$.
For example, on the LHS of \eqref{eq:ds-linear}, the vertical component is a $d$-component and the horizontal component is an $s$-component. The box represents the linear combination on the right:
\begin{align}
\begin{array}{c}\includegraphics[scale=0.23]{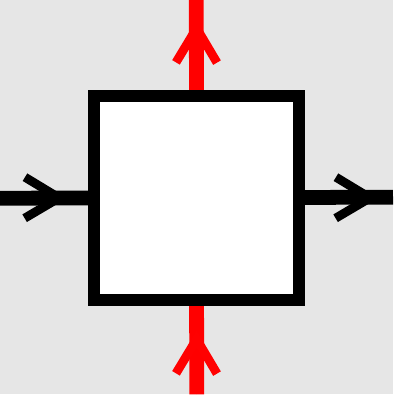}\end{array}=
\begin{array}{c}\includegraphics[scale=0.23]{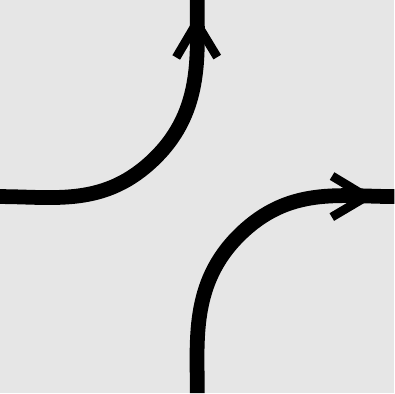}\end{array}-\begin{array}{c}\includegraphics[scale=0.23]{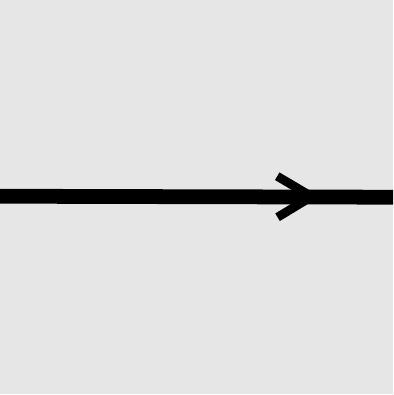}\end{array}. \label{eq:ds-linear}
\end{align}
The drawn strands on the RHS are all parts of $s$-components.
By definition, each $d$-component meets exactly one box. For a general $M$, the ds-linear combination is defined by choosing a local projection and using a vertical framing.

Replacing all boxes by their ds-linear combinations allows us to regard a tangle with boxes as an element of $\widetilde{\cA}(M,J)$. Thus, a box is only a convenient notation for a ds-linear combination in $\widetilde{\cA}(M,J)$. 
A $d$-component is a virtual component, so the skein relations cannot be applied directly to it before the boxes are resolved.

To simplify later arguments, we also introduce the following type of box:
\begin{align}
\begin{array}{c}\includegraphics[scale=0.23]{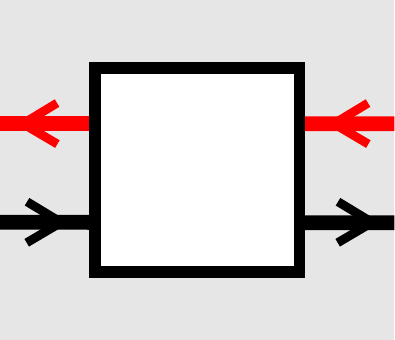}\end{array}:=
\begin{array}{c}\includegraphics[scale=0.23]{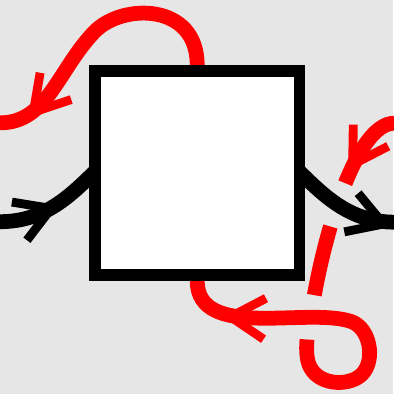}\end{array}
=\begin{array}{c}\includegraphics[scale=0.23]{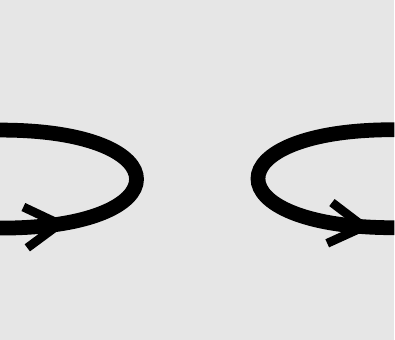}\end{array}-
\begin{array}{c}\includegraphics[scale=0.23]{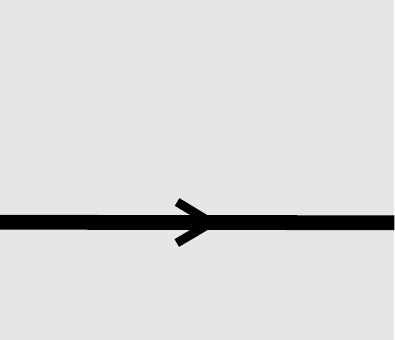}\end{array}
\end{align}
We call this a \textbf{parallel box}, and call the previous type a \textbf{transverse box} when we should distinguish them. Both types will simply be called boxes.

\begin{definition}\label{def:finite-type_filtration} 
Let $(M,J)$ be a marked oriented $3$-manifold, and consider its HOMFLYPT skein module $\widetilde{\cA}(M,J)$.
Let $F^n\widetilde{\cA}(M,J)$ be the $\bQ\llbracket\rho\hbar,h\hbar\rrbracket[\rho]$-submodule generated by
$$\{h^i\hbar^jL\mid L\in\mathsf{Box}(M,J),\ \text{(the number of boxes of $L$)}+2i+2j\geq n\}.$$
The descending filtration $\{F^n\widetilde{\cA}(M,J)\}_{n\in\bN}$ is called the \textbf{finite-type filtration}.
On the component degree $\ell$-part, set $F^n\widetilde{\cA}^{(\ell)}(M,J):=F^n\widetilde{\cA}(M,J)\cap\widetilde{\cA}^{(\ell)}(M,J)$, regarded as a $\bQ\llbracket\rho\hbar,h\hbar\rrbracket$-submodule.
\end{definition}

\begin{remark}
\begin{enumerate}
    \item To prove one of the main theorems (Theorem~\ref{thm:isom}) of this paper, it is not necessary to distinguish $h$ and $\hbar$ in the coefficient ring $\bQ\llbracket\rho\hbar,h\hbar\rrbracket[h,\hbar,\rho]$. In the sequel \cite{KT26+}, however, we determine the component degree parts containing the values of a homology cylinder invariant with values in $\Gr\widetilde{\cA}(\Sigma)$. For this purpose, the two parameters must be kept distinct.
    
    \item In the classical specialization $h=\hbar=0$, the filtration is determined solely by the number of boxes. It corresponds to the filtration $\{(\ker\varepsilon)^i\}_{i\in\bN}$ associated with the augmentation map $\varepsilon$ of the group ring of $\pi_1(\Sigma)$.
    
    \item In Section~\ref{sec:formality_genus0}, we define a filtration using $\text{(the number of boxes of $L$)}+i+j\geq n$ instead of $\text{(the number of boxes of $L$)}+2i+2j\geq n$. Although this only changes the relative weights assigned to boxes and to the powers of $h,\hbar$, and the two filtrations essentially have the same topology, the corresponding statements and proofs do not carry over unchanged. See Section~\ref{sec:formality_genus0} for details.
\end{enumerate}
\end{remark}
The finite-type filtration gives the graded quotients
\begin{align*}
    \Gr^{n}\widetilde{\cA}(M,J):=F^{n}\widetilde{\cA}(M,J)/F^{n+1}\widetilde{\cA}(M,J),\qquad \Gr^{n}\widetilde{\cA}^{(\ell)}(M,J):=F^{n}\widetilde{\cA}^{(\ell)}(M,J)/F^{n+1}\widetilde{\cA}^{(\ell)}(M,J)
\end{align*}

We will prove key properties of $\Gr^n\widetilde{\cA}(M,J)$ in Section~\ref{sec;isom}.

%%%%%%%%%%%%%%%%%%%%%%%%%%%%%%%%%%%%%%%%%%%%%%%%%%%%%%%%%
\subsection{Schedler's quantized necklace algebra}\label{subsec:quantized_necklace}
Necklace Lie algebras were introduced by Ginzburg \cite{Gin01} and Bocklandt--Le Bruyn \cite{BLB02}.
Following Turaev's biquantization, Schedler constructed a quantization of the necklace Lie bialgebra as a Hopf algebra \cite{Sch05}.
In this paper, we interpret Schedler's quantization for a bouquet quiver in terms of the homology of $\Sigma_{g,1}$, and use this viewpoint to define quantized necklace algebras for general surfaces.
This is based on the observation that the arrows of the double of a bouquet quiver can be identified with a symplectic basis of $H_1(\Sigma_{g,1})$.

A \textbf{necklace} is an oriented circle $S^1$ with finitely many (possibly none) marked points. The circle itself is called the \textbf{wire}, and each marked point is called a \textbf{bead}.
Thus, by convention, a necklace is allowed to have no beads.
Necklaces are abstract objects and are not embedded in $\Sigma$.
In this paper, we also regard the empty set as a necklace.

Assigning an element of $H_1(\Sigma):=H_1(\Sigma;\bQ)$ to each bead gives an \textbf{$H_1(\Sigma)$-necklace}.
In diagrams, each bead is drawn as a white circle containing its label in $H_1(\Sigma)$. Unless otherwise stated, wires are oriented counterclockwise.

Consider a disjoint union of $H_1(\Sigma)$-necklaces, together with a total order on the set of its beads. If there are $m$ beads, number them $1,\dots,m$ according to this order. We call this number the \textbf{height} of a bead, and say that bead $i$ is lower than bead $j$ if $i<j$.
Such a disjoint union, equipped with heights on its beads, is called a \textbf{quantum $H_1(\Sigma)$-necklace ensemble} for $\Sigma$.

\begin{definition}\label{def:necklace_Lie_algebra}
Consider the $\bQ[\rho,h,\hbar]$-module generated by necklace ensembles, subject to the relation (H).
For necklace ensembles $\alpha_i$ with $m_i$ beads, $i=1,2$, define the product $\alpha_1\alpha_2$ to be the disjoint union $\alpha_1\sqcup\alpha_2$, keeping the heights of the beads of $\alpha_2$ unchanged and replacing each height $k$ in $\alpha_1$ with $m_2+k$. The resulting $\bQ[\rho,h,\hbar]$-algebra is called the
\textbf{universal quantized necklace algebra} of $\Sigma$ and is denoted by $\widetilde{\mathrm N}(\Sigma)$.
Similarly, imposing (H) and (A')--(C') and using the same multiplication gives the \textbf{quantized necklace algebra} of $\Sigma$, denoted by $\widetilde{\cN}(\Sigma)$.
\end{definition}
For $x,y\in H_1(\Sigma)$ and $a,b\in\bQ$, the relations are
\begin{align*}
&{\rm (H)} 
\begin{array}{c}\includegraphics[scale=0.26]{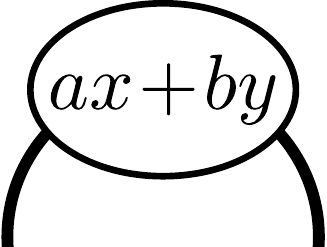}\end{array}=a\begin{array}{c}\includegraphics[scale=0.26]{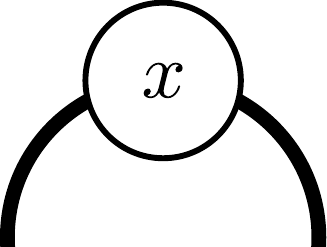}\end{array}+b\begin{array}{c}\includegraphics[scale=0.26]{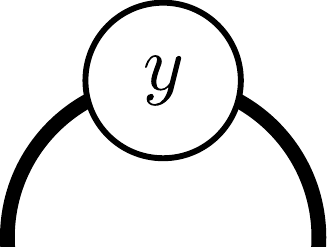}\end{array},\\
&{\rm (A')} 
\begin{array}{c}\includegraphics[scale=0.26]{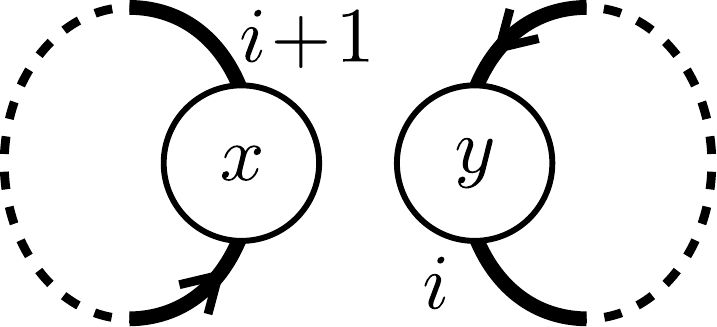}\end{array}-\begin{array}{c}\includegraphics[scale=0.26]{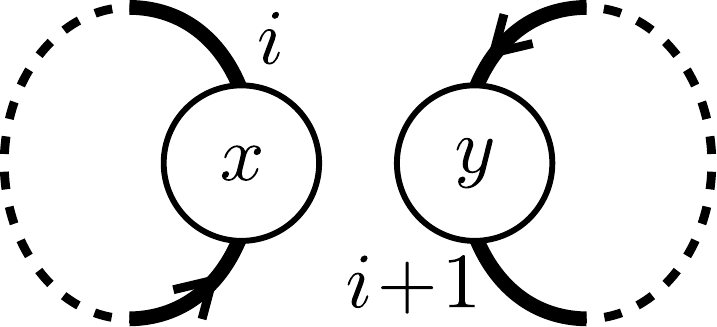}\end{array}=\mu(x,y) h\begin{array}{c}\includegraphics[scale=0.26]{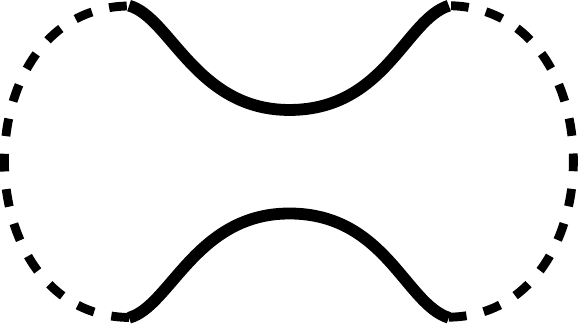}\end{array},\\
&{\rm (B')} 
\begin{array}{c}\includegraphics[scale=0.26]{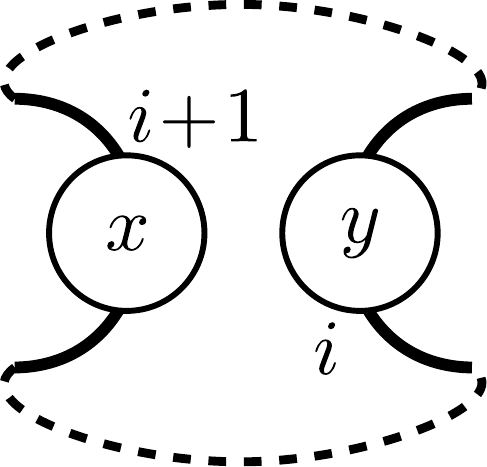}\end{array}- \begin{array}{c}\includegraphics[scale=0.26]{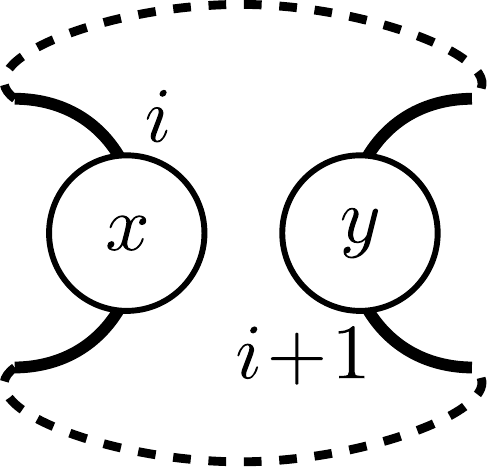}\end{array}=\mu(x,y) \hbar\begin{array}{c}\includegraphics[scale=0.26]{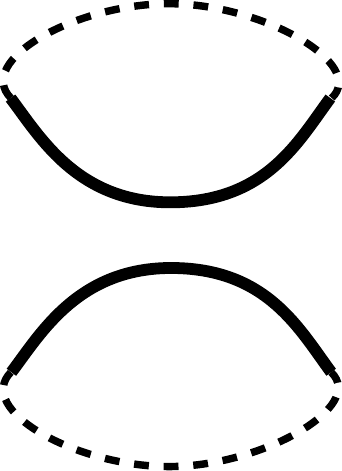}\end{array},\\
&{\rm (C')} 
\begin{array}{c}\includegraphics[scale=0.25]{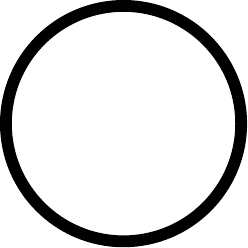}\end{array}=2\rho.
\end{align*}
In (A'), the reconnection is compatible with the orientations. the LHS of (C') is a necklace with no beads, and $\mu\colon H_1(\Sigma)\times H_1(\Sigma)\to\bQ$ denotes the intersection form.
These relations do not change the relative order of the heights of beads not shown in the diagrams.

\begin{remark}\label{rmk:HOMFLYPT}
\begin{enumerate}
    \item Relations (A') and (B') correspond to Eqs.~(3.2) and~(3.3) of \cite{Sch05}, respectively. In Eq.~(3.2) of \cite{Sch05}, our parameter $h$ is specialized to $1$. We impose (C') for compatibility with the HOMFLYPT skein relations.
    
    \item For $\Sigma_{g,1}$, the classical limit $h=\hbar=0$ of the HOMFLYPT skein algebra with $\rho=0$ is known to be isomorphic, as a bi-Poisson bialgebra, to the symmetric algebra of the Goldman--Turaev Lie bialgebra modulo the trivial loop \cite{Tur91}. Moreover, the associated graded Lie bialgebra of the Goldman--Turaev Lie bialgebra, with respect to the filtration induced by the augmentation ideal, is isomorphic to the necklace Lie bialgebra modulo the necklace with no beads \cite{MT13,KK14,KK15}. Extending this isomorphism to the symmetric algebras shows that the associated graded of the classical limit of the HOMFLYPT skein algebra, with its induced filtration, is isomorphic as a bi-Poisson bialgebra to the classical limit of the quantized necklace algebra at the same specialization.
    
    \item If $\Sigma$ has genus $0$, its intersection form vanishes. Thus, (A') and (B') simply allow heights to be interchanged. This loss of information is a serious issue for formality. To deal with it, we modify the definition of the quantized necklace algebra, as in the classical setting. This is discussed in Section~\ref{sec:formality_genus0}.
    \item Table~\ref{table:correspondence} gives the correspondence between our quantized necklace algebra and Schedler's construction \cite{Sch05}.
    Here, $A_{\overline Q}/{\rm cyclic}$ denotes the cyclic quotient of $A_{\overline Q}$.
    Our quantized necklace algebra is obtained from this correspondence by replacing $H_1(\Sigma_{g,1})$ with $H_1(\Sigma)$ for an arbitrary compact oriented surface $\Sigma$ with non-empty boundary, replacing Schedler's $h$ with $\hbar$, introducing a new quantum parameter $h$, and imposing the additional relation (C'). Following \cite{Tur91}, a coalgebra structure is also defined in \cite{Sch05}; we discuss its relation to our construction in Remark~\ref{rmk:bialg}.
\end{enumerate}
\begin{table}[h]
    \centering
    \caption{Correspondence with Schedler's notation.}
    \label{table:correspondence}
    \begin{tabular}{ccc}
        \hline
        This paper & Schedler's notation \\ \hline \hline
        A symplectic basis of $H_1(\Sigma_{g,1})$ & The arrows of the double of the $g$-bouquet quiver $Q$\\ \hline
        A necklace with no beads & A vertex idempotent of the quiver \\ \hline
        An $H_1(\Sigma_{g,1})$-necklace & A product of quiver paths in $A_{\overline Q}/{\rm cyclic}$ \\ \hline
        A quantum $H_1(\Sigma_{g,1})$-necklace ensemble & A disjoint union of products of paths with heights \\ \hline
        Quantum parameters $\hbar,\ h$ & Quantum parameters $h,\ 1$ \\ \hline
        $2\rho=$\,\text{trivial necklace} & trivial necklace\\ \hline
    \end{tabular}
\end{table}
Thus, our parameter $h$ is set to $1$ in \cite{Sch05}. A refinement of Schedler's construction as a biquantization, without specializing this parameter to $1$, is considered in \cite{CHLZ26}.
\end{remark}

As for $\widetilde{\cA}^{(\ell)}(\Sigma)$, let $\widetilde{\cN}^{(\ell)}(\Sigma)$ be the $\bQ$-submodule of $\widetilde{\cN}(\Sigma)$ generated by
\begin{align*}
\{h^i\hbar^j\rho^kN\in\widetilde{\cN}(\Sigma)\mid\text{$N$ is an $H_1(\Sigma)$-necklace ensemble},\ \#\{\text{necklaces in $N$}\}+i-j+k=\ell\}.
\end{align*}
We call $\#\{\text{necklaces in $N$}\}+i-j+k$ the \textbf{component degree}.
The quantity $\#\{\text{necklaces in $N$}\}+i-j+k$ is preserved by (H), (A'), (B'), and (C').
Hence, we have a direct sum decomposition $\displaystyle\widetilde{\cN}(\Sigma)=\bigoplus_\ell\widetilde{\cN}^{(\ell)}(\Sigma)$ as $\bQ$-modules.

Let $\Gr^n\widetilde{\cN}(\Sigma)$ denote the $\bQ[\rho]$-submodule of $\widetilde{\cN}(\Sigma)$ generated by
\begin{align}
\{h^i\hbar^j N\mid\text{$N$ is a necklace ensemble},\ \#\{\text{beads}\}+2i+2j=n, k\in\bN\}.\label{eq:degree_necklace}
\end{align}
Set $\Gr^n\widetilde{\cN}^{(\ell)}(\Sigma):=\widetilde{\cN}^{(\ell)}(\Sigma)\cap \Gr^n\widetilde{\cN}(\Sigma)$. 
We call $\#\{\text{beads}\}+2i+2j$ the \textbf{finite-type degree}.
The notation $\Gr^n\widetilde{\cN}(\Sigma)$ does not refer to an associated graded quotient of a filtration. It is used to emphasize the correspondence with $\Gr^n\widetilde{\cA}(\Sigma)$.

%%%%%%%%%%%%%%%%%%%%%%%%%%%%%%%%%%%%%%%%%%%%%%%%%%%%%%%%%%%%%%%%%%%%%%%%%%%%%%%%%
\section{The associated graded algebra of the HOMFLYPT skein algebra and the quantized necklace algebra}\label{sec;isom}
\subsection{An algebra homomorphism from the quantized necklace algebra to the associated graded algebra of the HOMFLYPT skein algebra}\label{subsec:map_GrPhi}
In this subsection, we define a map $\Gr^n\Phi\colon \Gr^n\widetilde{\cN}(\Sigma)\to \Gr^n\widetilde{\cA}(\Sigma)=F^n\widetilde{\cA}(\Sigma)/F^{n+1}\widetilde{\cA}(\Sigma)$ and show that it is an isomorphism of $\bQ[\rho]$-modules.
These maps induce a $\bQ[\rho,h,\hbar]$-algebra isomorphism $\Gr\Phi:=\bigoplus_n\Gr^n\Phi\colon \widetilde{\cN}(\Sigma)\to \bigoplus_n \Gr^n\widetilde{\cA}(\Sigma)$.

We define $\Gr^n\Phi\colon \Gr^n\widetilde{\cN}(\Sigma)\to \Gr^n\widetilde{\cA}(\Sigma)=F^n\widetilde{\cA}(\Sigma)/F^{n+1}\widetilde{\cA}(\Sigma)$ as follows. %The map is defined by assigning a link diagram with boxes.
Let $N$ be an $H_1(\Sigma)$-necklace ensemble with $k$ wires. Denote $\{W_1,\dots,W_k\}$ for its connected components and $\{b_1^{(i)},\dots,b_{m_i}^{(i)}\}$ for the beads on $W_i$.
\begin{enumerate}
    \item \ul{Draw an oriented link diagram on $\Sigma$} with as many components as there are necklaces, or equivalently wires. Let $D_i$ be the $S^1$-component corresponding to $W_i$, and let $p_j^{(i)}$ be the point on $D_i$ corresponding to the bead $b_j^{(i)}$. We require that the cyclic order of the beads along the orientation of $W_i$ agrees with that of the points $p_j^{(i)}$ along the orientation of $D_i$. Independence of the choice of the link diagram representing the wires is proved in Theorem~\ref{thm:welldefined}.
    
    \item \ul{Represent the element of $H_1(\Sigma;\bQ)$ assigned to each bead $b_j^{(i)}$ as a rational multiple of an element of $H_1(\Sigma;\bZ)$}. \ul{Then, choose a lift of this integral homology class to $\pi_1(\Sigma,p_j^{(i)})$}, and draw it as a loop on $\Sigma$ based at $p_j^{(i)}$. We assume that the only singularities are transverse double points. Replace $p_j^{(i)}$ with a box to obtain a diagram on $\Sigma$ in which the image $D_i$ of the wire $W_i$ is an $s$-component and the representative of the element of $\pi_1(\Sigma)$ is a $d$-component.
    Independence of the expression as a rational multiple and of the choice of the lift to $\pi_1(\Sigma)$ is proved in Theorem~\ref{thm:welldefined}.

    \item At each double point between distinct $d$-components, assign over/under information according to the heights of the corresponding beads. Namely, if the height of a bead $b'$ is greater than that of a bead $b$, the representative $\gamma'$ corresponding to $b'$ passes over the representative $\gamma$ corresponding to $b$ at every point of $\gamma\cap\gamma'$.
    \ul{The crossing information at all other double points may be chosen arbitrarily}. Independence of these choices is proved in Theorem~\ref{thm:welldefined}.
\end{enumerate}
This construction corresponds an $H_1(\Sigma)$-necklace ensemble with a link diagram with boxes on $\Sigma$, with one box for each bead.
Regarded as an element of $\widetilde{\cA}(\Sigma)$, this diagram determines the image of the $H_1(\Sigma)$-necklace ensemble. If the label at a bead $b$ is written as $q_bx_b$, where $q_b\in\bQ$ and $x_b\in H_1(\Sigma;\bZ)$, the resulting diagram is multiplied by $\prod_b q_b$.
More generally, if $N$ has $r$ beads and $n=r+2i+2j$, multiply its image by $h^i\hbar^j$ to define the image of $h^i\hbar^jN$ in $F^n\widetilde{\cA}(\Sigma)/F^{n+1}\widetilde{\cA}(\Sigma)$. 
Extending $\bQ[\rho]$-linearly gives the map $\Gr^n\Phi\colon \Gr^n\widetilde{\cN}(\Sigma)\to\Gr^n\widetilde{\cA}(\Sigma)$.

\begin{theorem}\label{thm:welldefined}
For a marked 3-manifold $(M,J)$, the map $\Gr^n\Phi$ is a well-defined $\bQ[\rho]$-module homomorphism.
\end{theorem}
Theorem~\ref{thm:welldefined} implies the following.
\begin{cor}\label{cor:welldefined}
The map
\begin{align*}
\Gr\Phi:=\bigoplus_n\Gr^n\Phi\colon \widetilde{\cN}(\Sigma)=\bigoplus_n \Gr^n\widetilde{\cN}(\Sigma)\longrightarrow \bigoplus_n \Gr^{n}\widetilde{\cA}(\Sigma)
\end{align*}
is a $\bQ[\rho,h,\hbar]$-module homomorphism. 
\end{cor}

We first prepare a lemma for the proof of Theorem~\ref{thm:welldefined}.
\begin{lemma}\label{lem:transparency}
The following hold in $\Gr^n\widetilde{\cA}(M,J)$.

\begin{enumerate}
    \item A crossing between (possibly the same) $s$-components can be switched. In other words, switching such a crossing does not change the represented element of $\Gr^n\widetilde{\cA}(M,J)$.
    
    \item A crossing between a $d$-component and an $s$-component can be switched without changing the represented element of $\Gr^n\widetilde{\cA}(M,J)$.
    
    \item A self-crossing of a $d$-component can be switched without changing the represented element of $\Gr^n\widetilde{\cA}(M,J)$.
    
    \item Framing information disappears, i.e., adding or removing a kink does not change the represented element of $\Gr^n\widetilde{\cA}(M,J)$.
\end{enumerate}
\end{lemma}
\begin{proof}
\begin{enumerate}
    \item Apply the relation (A) or (B) of the HOMFLYPT skein module at a crossing $c$ between $s$-components. The result is a linear combination of the tangle with boxes obtained by switching $c$ and the one obtained by smoothing $c$.
    The smoothing term contains one additional factor of $h$ or $\hbar$, and Hence, vanishes in $\Gr^n\widetilde{\cA}(M,J)$. Thus, the crossing $c$ may be switched.
    
    \item For a crossing $c$ between a $d$-component and an $s$-component, first replace the box on the $d$-component with its ds-linear combination. In the first term of \eqref{eq:ds-linear}, $c$ becomes a crossing between $s$-components. In the second term, the $d$-component has been removed, so $c$ is absent. Therefore, we can switch $c$ in the first term as in (1). Although one box has been removed, the smoothing term has a factor of $h$ or $\hbar$, and still vanishes in $\Gr^n\widetilde{\cA}(M,J)$. The resulting linear combination is exactly the ds-linear combination obtained by first switching the original crossing between the $d$-component and the $s$-component. This proves (2).
    
    \item For a self-crossing $c$ of a $d$-component, replace its box with the ds-linear combination and apply the relation (B) to the first term. This allows us to switch $c$.
    
    \item The same argument as in (3), with the relation (D) in place of (B), proves the claim. For example, in $\Gr^n\widetilde{\cA}(M,J)$, we have \begin{align*}\begin{array}{c}\includegraphics[scale=0.17]{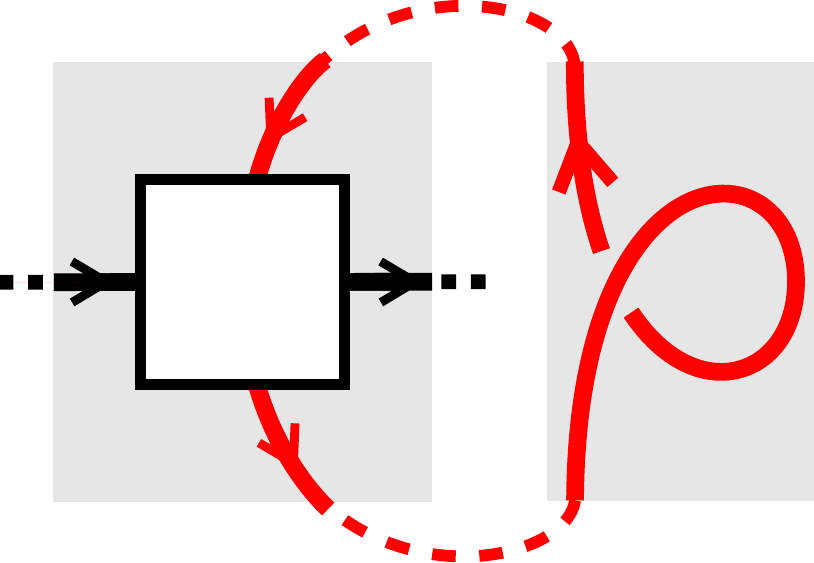}\end{array}-\begin{array}{c}\includegraphics[scale=0.17]{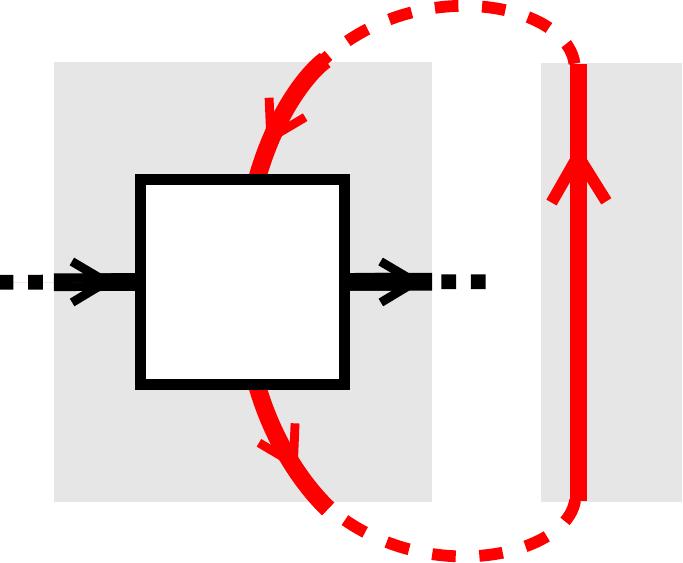}\end{array}=\begin{array}{c}\includegraphics[scale=0.17]{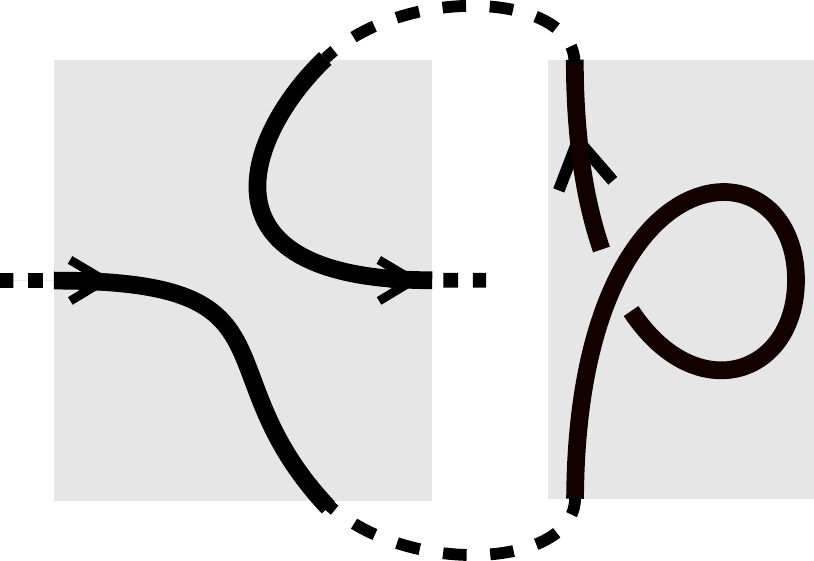}\end{array}-\begin{array}{c}\includegraphics[scale=0.17]{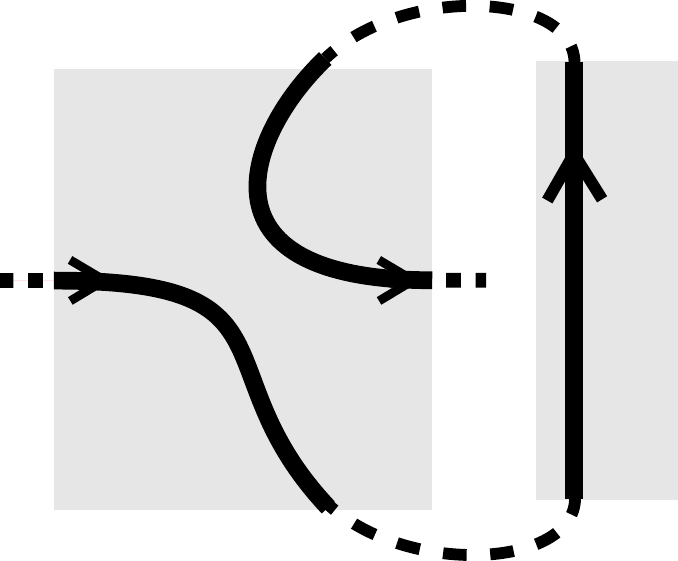}\end{array}=(\exp(\rho\hbar)-1)\begin{array}{c}\includegraphics[scale=0.17]{draws/resolv_no_dd_self.pdf}\end{array}=0.
\end{align*}
\end{enumerate}
\end{proof}

\begin{remark}
    Switching crossings between distinct $d$-components corresponds to the relations (A') and (B') of the quantized necklace algebra.
\end{remark}

\begin{proof}[\textbf{Proof of Theorem~$\ref{thm:welldefined}$}]
To prove well-definedness, it suffices to check independence of the underlined choices in the construction and compatibility of $\Gr^n\Phi$ with the defining relations (H), (A'), (B'), and (C') of the quantized necklace algebra. More explicitly, the choices to be checked are the following:
\begin{enumerate}
    \item[(1)] The link diagram representing the wires,
    
    \item[(2-1)] The representation of a label as a rational multiple of an element of $H_1(\Sigma;\bZ)$,
    
    \item[(2-2)] The lift of a homology class to $\pi_1(\Sigma)$,
    
    \item[(3)] The crossing information at double points other than those between distinct $d$-components.
\end{enumerate}
Independence of the crossing information in (3) follows from Lemma~\ref{lem:transparency}.

Next, we prove (1). In $\Gr^n\widetilde{\cA}(\Sigma)$, we have
\begin{align}
\begin{array}{c}\includegraphics[scale=0.3]{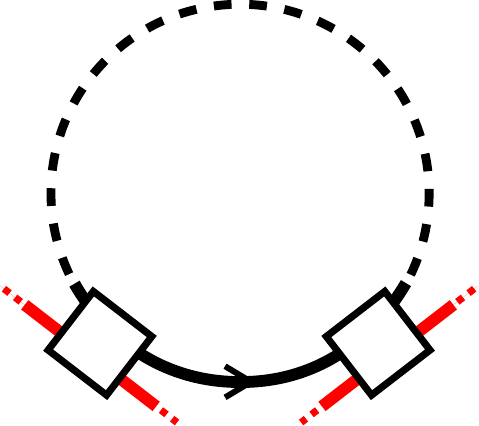}\end{array}=-\begin{array}{c}\includegraphics[scale=0.3]{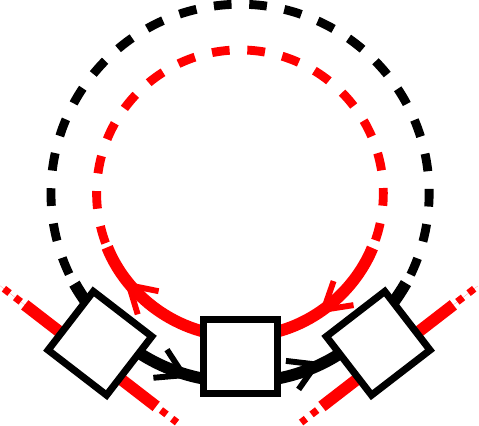}\end{array}+\begin{array}{c}\includegraphics[scale=0.3]{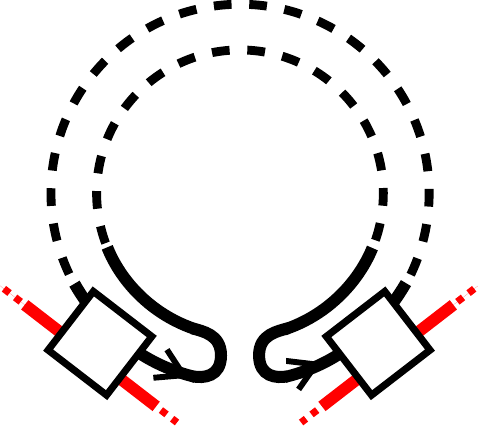}\end{array}
=\begin{array}{c}\includegraphics[scale=0.3]{draws/wire_parallel_d.pdf}\end{array},\label{eq:s-comp_trivial}
\end{align}
where the $s$-component, including its dotted part, and the $d$-component are parallel, i.e., they bound an embedded annulus. Other boxes may occur along the dotted part of the $s$-component.
This allows us to replace any link diagram representing the wires with a trivial one, and conversely. Thus, the construction is independent of the choice of the link diagram.

Next, we prove (2-2).
By the definition of a box, we have the following identity in $\Gr^n\widetilde{\cA}(\Sigma)$:
\begin{align}
\begin{array}{c}\includegraphics[scale=0.3]{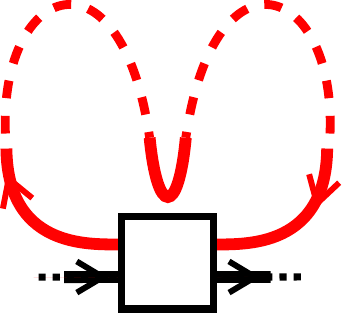}\end{array}
-\begin{array}{c}\includegraphics[scale=0.3]{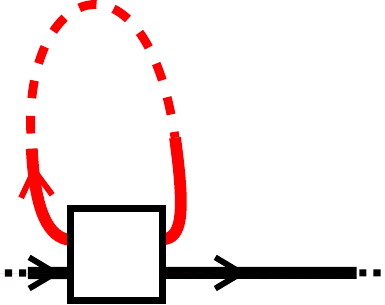}\end{array}
-\begin{array}{c}\includegraphics[scale=0.3]{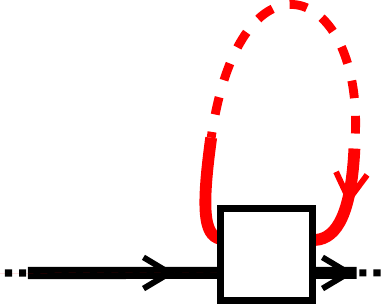}\end{array}=
\begin{array}{c}\includegraphics[scale=0.3]{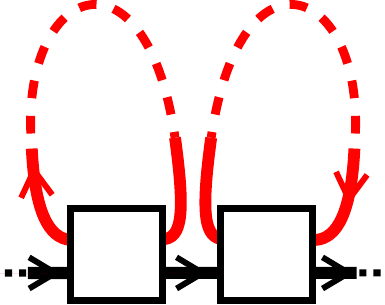}\end{array}=0.\label{eq:homologous}
\end{align}

The equation is interpreted as follows: if each term on the left represents an element of $\Gr^n\widetilde{\cA}(\Sigma)$, the term on the right vanishes since it has an additional box.
Let $\gamma_1$ and $\gamma_2$ be the $d$-components in the second and third terms on the left of \eqref{eq:homologous}, respectively. Then, the loop in the first term is $\gamma_1\gamma_2$, where products of paths are read from left to right.
Let $B_{\gamma_1}$, $B_{\gamma_2}$, and $B_{\gamma_1\gamma_2}$ denote the corresponding link diagrams with boxes.
Then
\eqref{eq:homologous} gives $B_{\gamma_1}+B_{\gamma_2}=B_{\gamma_1\gamma_2}$. Similarly, $B_{\gamma_1}+B_{\gamma_2}=B_{\gamma_2\gamma_1}$, and Hence, $B_{\gamma_1\gamma_2}=B_{\gamma_2\gamma_1}$.
Taking $\gamma_1=\alpha\beta$ and $\gamma_2=\beta^{-1}$, we obtain
$B_\alpha=B_{\gamma_1\gamma_2}=B_{\gamma_2\gamma_1}=B_{\beta^{-1}\alpha\beta}$, which proves invariance under conjugation. 
More generally, the additivity $B_{\gamma_1\gamma_2}=B_{\gamma_1}+B_{\gamma_2}$ shows that the map factors through the abelianization of $\pi_1(\Sigma)$. Thus, it depends only on the integral homology class.

We next prove (2-1). Let $\gamma\in\pi_1(\Sigma)$ be a lift of $x\in H_1(\Sigma;\bZ)$. For any $m\in\bZ$, the element $\gamma^m$ lifts $mx$, and \eqref{eq:homologous} shows that its image is $m$ times the image obtained from $\gamma$. Thus, the image is independent of the expression as a rational multiple of an integral homology class. 
In the following, we can assume that every bead is labelled by an element of $H_1(\Sigma;\bZ)$.

So far, we have proven that the $\bQ[\rho]$-module homomorphism $\Gr^n\widetilde{\mathrm{N}}(\Sigma)\to\Gr^n\widetilde{\cA}(\Sigma)$ is well-defined.

Next, we show that it factors through a map $\Gr^n\widetilde{\mathcal{N}}(\Sigma)\to\Gr^n\widetilde{\cA}(\Sigma)$. It suffices to check compatibility with the relations (H), (A'), (B'), and (C').

In the relation (H), take $a,b\in\bZ$ and denote the $H_1(\Sigma)$-necklace ensembles in the three terms, from left to right, by $N_{ax+by}$, $N_x$, and $N_y$. Then, the relation~\eqref{eq:homologous} gives $\Gr^n\Phi(N_{ax+by})=a\Gr^n\Phi(N_x)+b\Gr^n\Phi(N_y)$.

Next, we check compatibility of $\Gr\Phi$ with (A') and (B').
Let $D$ be a diagram on $\Sigma$ of the tangle with boxes obtained by applying $\Gr\Phi$ to the $H_1(\Sigma)$-necklace ensemble in the first term on the LHS of (A') or (B'). Let $\gamma_1$ and $\gamma_2$ be the $d$-components corresponding to the two beads in the relation. Replace both boxes with their ds-linear combinations, and let $D'$ be the diagram of the term in which both $\gamma_1$ and $\gamma_2$ remain, namely the first term in each ds-linear combination.
Let $p_1,p_2,\dots,p_k$ be the intersections of $\gamma_1$ and $\gamma_2$ in $D'$, and let $\varepsilon_i$ be the sign at $p_i$ according to Figure~\ref{fig:crossings_signs}.
For $\varepsilon_i\in\{-1,1\}$, $i=1,\dots,k$, let $D(\varepsilon_1,\dots,\varepsilon_k)$ (resp. $D'(\varepsilon_1,\dots,\varepsilon_k)$) denote the diagram obtained from $D$ (resp. $D'$) by replacing a neighborhood of each $p_i$ with the local diagram indicated by $\varepsilon_i$, respecting orientations.
The image of the LHS of (A') or (B') is then
\begin{align*}
&D(\varepsilon_1,\varepsilon_2,\dots,\varepsilon_k)
-D(-\varepsilon_1,-\varepsilon_2,\dots,-\varepsilon_k)
=
D'(\varepsilon_1,\varepsilon_2,\dots,\varepsilon_k)
-D'(-\varepsilon_1,-\varepsilon_2,\dots,-\varepsilon_k)\\
&=
D'(\varepsilon_1,\varepsilon_2,\dots,\varepsilon_k)-D'(-\varepsilon_1,\varepsilon_2,\dots,\varepsilon_k)
+D'(-\varepsilon_1,\varepsilon_2,\dots,\varepsilon_k)
-D'(-\varepsilon_1,-\varepsilon_2,\dots,-\varepsilon_k)\\
&=\varepsilon_1\eta D'(0,\varepsilon_2,\dots,\varepsilon_k)
+D'(-\varepsilon_1,\varepsilon_2,\dots,\varepsilon_k)
-D'(-\varepsilon_1,-\varepsilon_2,\dots,-\varepsilon_k)\\
&\ \vdots\\
&=\varepsilon_1\eta D'(0,\varepsilon_2,\dots,\varepsilon_k)
+\varepsilon_2\eta D'(-\varepsilon_1,0,\dots,\varepsilon_k)+\dots+
\varepsilon_k\eta D'(-\varepsilon_1,\dots,-\varepsilon_{k-1},0),
\end{align*}
where $\eta=h$ if the boxes of $\gamma_1$ and $\gamma_2$ lie on distinct $s$-components, and $\eta=\hbar$ if they lie on the same $s$-component.
If $D(\varepsilon_1,\dots,\varepsilon_k)\in\Gr^n\widetilde{\cA}(\Sigma)$, each term in the last expression also belongs to $\Gr^n\widetilde{\cA}(\Sigma)$.
%Using \eqref{eq:s-comp_trivial},
By the independence of the link diagram representing the wires, proved in (1), we have
$\eta D'(0,\varepsilon_2,\dots,\varepsilon_k)
=\eta D'(-\varepsilon_1,0,\dots,\varepsilon_k)=\eta D'(-\varepsilon_1,\dots,-\varepsilon_{k-1},0).$
Thus, the last expression is
$$\left(\sum_{i=1}^k \varepsilon_i \right)\eta D'(0,\varepsilon_2,\dots,\varepsilon_k)$$
Let $[\gamma_i]$ denote the homology class of $\gamma_i$. 
Then, we have $\sum_{i=1}^k\varepsilon_i=\mu([\gamma_1],[\gamma_2])$.
This is exactly the image under $\Gr\Phi$ of the RHS of (A') or (B').

Finally, we check compatibility with (C'). Let $N$ be the element on its LHS of (C'). It suffices to consider $N\in\Gr^0\widetilde{\cN}(\Sigma)$. By the definition of $\Gr^0\Phi$, its image $\Gr^0\Phi(N)$ belongs to $\Gr^0\widetilde{\cA}(\Sigma)$. In this quotient, $\dfrac{2\sinh(\rho\hbar)}{\hbar}=2\rho+(\text{higher-order terms in $\hbar$})=2\rho$. 
Hence, (C') is preserved.
\end{proof}

\begin{lemma}\label{lem:alghom}
The map $\Gr\Phi$ is a $\bQ[\rho,h,\hbar]$-algebra homomorphism. In particular, multiplication by $h$ or $\hbar$ increases the finite-type degree by $2$. 
\end{lemma}
\begin{proof}
By Corollary~\ref{cor:welldefined}, the map $\Gr\Phi$ is a $\bQ[\rho,h,\hbar]$-module homomorphism. Thus, it suffices to show that it preserves multiplication.
Let $N_1$ and $N_2$ be quantum $H_1(\Sigma)$-necklace ensembles.
In the product $N_1N_2$, every bead of $N_1$ is higher than every bead of $N_2$.
For the product $\Gr\Phi(N_1)\Gr\Phi(N_2)$, we can choose representatives with $\Gr\Phi(N_1)\subset\Sigma\times(1/2,1)$ and $\Gr\Phi(N_2)\subset\Sigma\times(0,1/2)$. For these representatives, $\Gr\Phi(N_1)\Gr\Phi(N_2)=\Gr\Phi(N_1)\cup\Gr\Phi(N_2)$.
The heights in the quantum $H_1(\Sigma)$-necklace ensemble $N_1N_2$ place the beads of $N_1$ above those of $N_2$.
Together with the construction of $\Gr\Phi$, this gives $\Gr\Phi(N_1N_2)=\Gr\Phi(N_1)\Gr\Phi(N_2)$.
\end{proof}

\begin{lemma}\label{lem:surj}
The map $\Gr\Phi$ is surjective.
\end{lemma}
\begin{proof}
Every element of $\Gr\widetilde{\cA}(\Sigma)$ can be expressed as a $\bQ[\rho,h,\hbar]$-linear combination of isotopy classes of links with boxes in $\Sigma\times[0,1]$.
Therefore, it suffices to realize the class of each such a link as the image under $\Gr\Phi$ of an element of $\widetilde{\cN}(\Sigma)$.
For $L\in\mathsf{Box}(\Sigma)$, let $s_1,\dots,s_k$ denote its $s$-components and $d_1,\dots,d_l$ denote its $d$-components.
Let $D$ be a diagram of $L$.

First suppose that every crossing between distinct $d$-components $d_i,d_j$ in $D$ satisfies
\begin{align}\label{text:monotone}
    \text{$i<j$ and $d_j$ is the overpassing arc.}
\end{align}
Let $\gamma_1,\dots,\gamma_l\in H_1(\Sigma)$ be the classes represented by $d_1,\dots,d_l$.
Since the construction of $\Gr\Phi$ is independent of the choice of the wires, take $s_i\ (i=1,\dots,k)$ as wires and place a bead labelled by $\gamma_j$ on $s_i$ wherever $d_j$ is attached to $s_i$. Give this bead height $j$. The resulting $H_1(\Sigma)$-necklace ensemble maps to $L$ in $\Gr\widetilde{\cA}(\Sigma)$.
Thus, surjectivity follows if every link diagram with boxes can be replaced, in $\Gr\widetilde{\cA}(\Sigma)$, by a linear combination of diagrams satisfying \eqref{text:monotone}.

Order the complexity lexicographically by the number of boxes and then by the number of crossings that fail \eqref{text:monotone}. 
The ds-linear combination and the relation (A) or (B) give
\begin{align}\label{eq:surj_crossing}
\begin{array}{c}\includegraphics[scale=0.20]{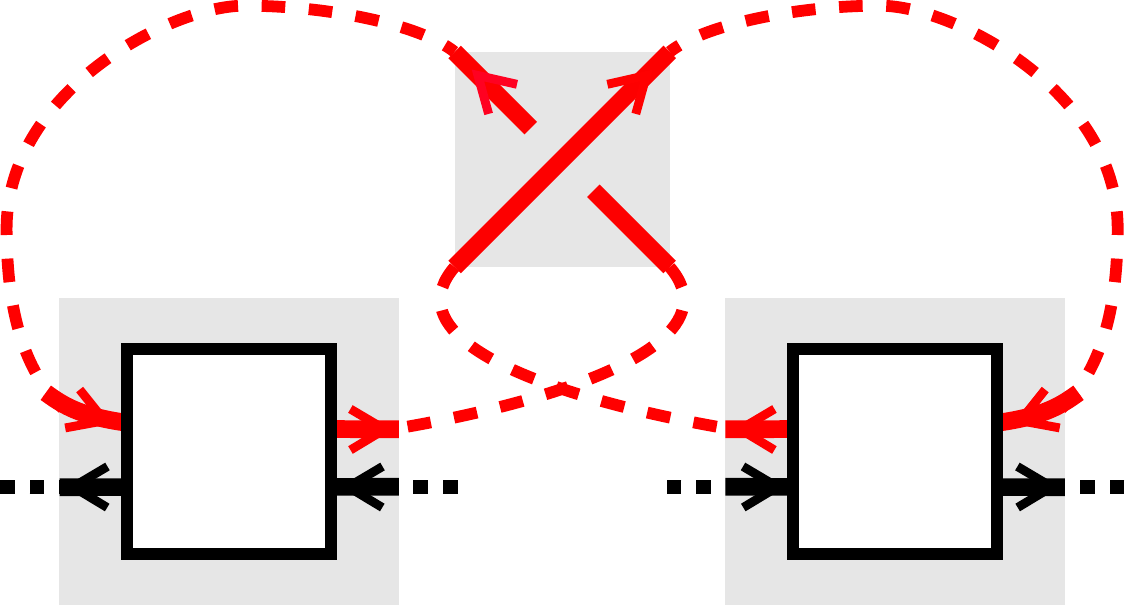}\end{array}
-\begin{array}{c}\includegraphics[scale=0.20]{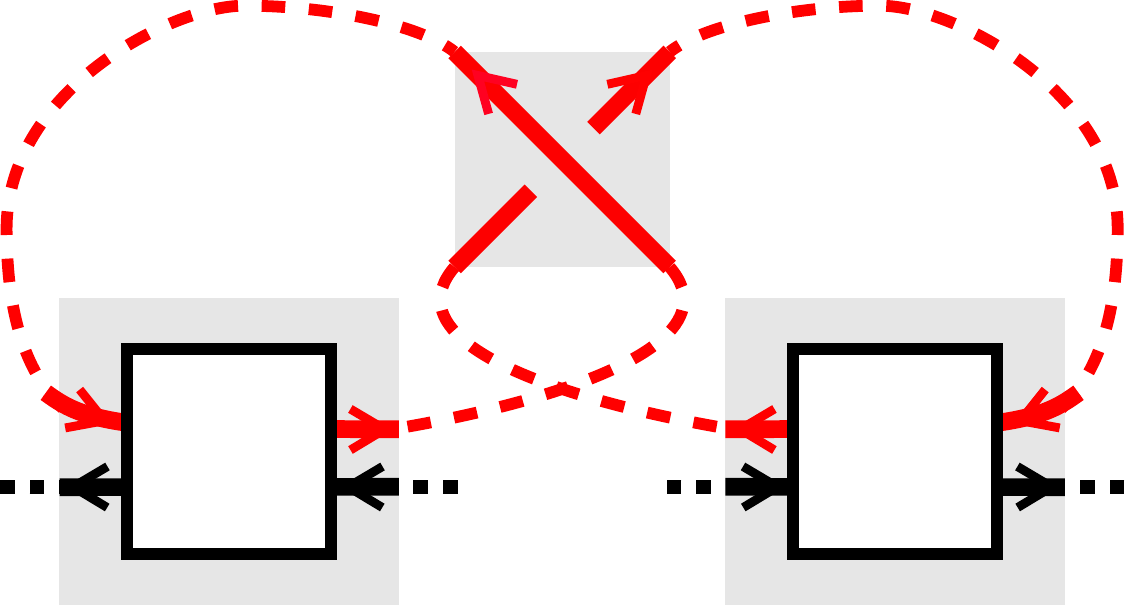}\end{array}=\eta
\begin{array}{c}\includegraphics[scale=0.20]{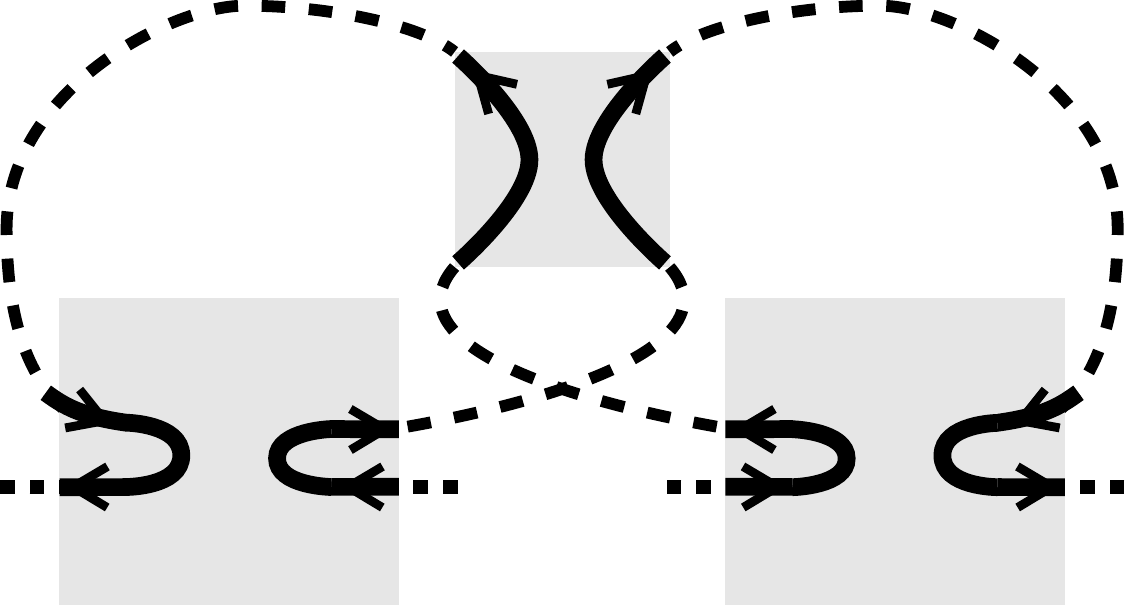}\end{array},
\end{align}
where $\eta=h$ if relation (A) is used, and $\eta=\hbar$ if relation (B) is used.
At a crossing between distinct $d$-components, this expresses the original diagram as a linear combination of the diagram with that crossing switched and a diagram with two fewer boxes. The factor $h$ or $\hbar$ on the right compensates for the loss of two boxes, so the finite-type degree is unchanged.
There are finitely many boxes and crossings between $d$-components. In the switched term, the second quantity of the complexity decreases, whereas in the smoothing term the number of boxes decreases. Thus, finitely many applications express any link diagram with boxes as a linear combination of diagrams satisfying \eqref{text:monotone}.
\end{proof}

\begin{theorem}\label{thm:inj}
Let $\Sigma$ be a compact oriented surface with non-empty boundary. Then $\Gr\Phi$ is injective.
\end{theorem}
The proof of Theorem~\ref{thm:inj} is given in Section~\ref{sec:proof_mainthm}. The assumption in Theorem~\ref{thm:inj} that $\Sigma$ has non-empty boundary is essential to the proof; the same argument cannot be applied to closed surfaces.

Combining Corollary~\ref{cor:welldefined}, Lemmas~\ref{lem:alghom} and~\ref{lem:surj}, and Theorem~\ref{thm:inj}, which establish well-definedness, multiplicativity, surjectivity, and injectivity, respectively, we obtain the following.
\begin{theorem}\label{thm:isom}
Let $\Sigma$ be a compact oriented surface with non-empty boundary.
The map $\Gr\Phi\colon\widetilde{\cN}(\Sigma)\to\Gr\widetilde{\cA}(\Sigma)$ is a well-defined $\bQ[\rho,h,\hbar]$-algebra isomorphism. 
In other words, the quantized necklace algebra $\widetilde{\cN}(\Sigma)$ and the associated graded algebra $\Gr\widetilde{\cA}(\Sigma)$ of the HOMFLYPT skein algebra are isomorphic as $\bQ[\rho,h,\hbar]$-algebras.
\end{theorem}

\begin{remark}\label{rmk:extending_Jversion}
If $|J|=2$ and both points of $J$ lie on the same component of $\partial\Sigma$, the analogue of Theorem~\ref{thm:isom} holds as an isomorphism of $\bQ[\rho,h,\hbar]$-modules. The finite-type filtration in Definition~\ref{def:finite-type_filtration} and the finite-type degree in \eqref{eq:degree_necklace} are defined in the same way. Indeed, well-definedness and surjectivity follow as above. For injectivity, attach a handle to $\partial\Sigma$ connecting the two points of $J$ to obtain a surface $\Sigma'$, and regard a $J$-tangle as a framed link in the thickening of $\Sigma'$.
If the two points of $J$ lie on distinct boundary components, the relations of the quantized necklace Lie module must be adjusted using the intersection form on $H_1(\Sigma')$.
\end{remark}

As is implicit in \cite{Sch05}, the quantized necklace algebra $\widetilde{\cN}(\Sigma)$ has no zero divisors. Hence, we obtain the following.
\begin{cor}
Let $\Sigma$ be a compact oriented surface with non-empty boundary.
Then the completion $\widehat{\widetilde{\cA}(\Sigma)}:=\plim[n]\widetilde{\cA}(\Sigma)/F^n\widetilde{\cA}(\Sigma)$ has no zero divisors.
\end{cor}

\begin{conj}
Let $\Sigma$ be a compact oriented surface with non-empty boundary.
\begin{enumerate}
    \item We have $\bigcap_nF^n\widetilde{\cA}(\Sigma)=\{0\}$.
    \item More weakly, $\widetilde{\cA}(\Sigma)$ has no zero divisors.
\end{enumerate}
\end{conj}

\subsection{The Lie algebra structure on $\widetilde{\cA}(\Sigma)$ and the Lie module structure on $\widetilde{\cA}(\Sigma,J)$}\label{subsec:Lie_structure}
To prove Theorem~\ref{thm:inj},
we define a Lie action $\sigma$ of $\widetilde{\cA}(\Sigma)$ on $\widetilde{\cA}(\Sigma,J)$ and establish its basic properties.

\begin{figure}
	\centering\includegraphics[width=60mm]{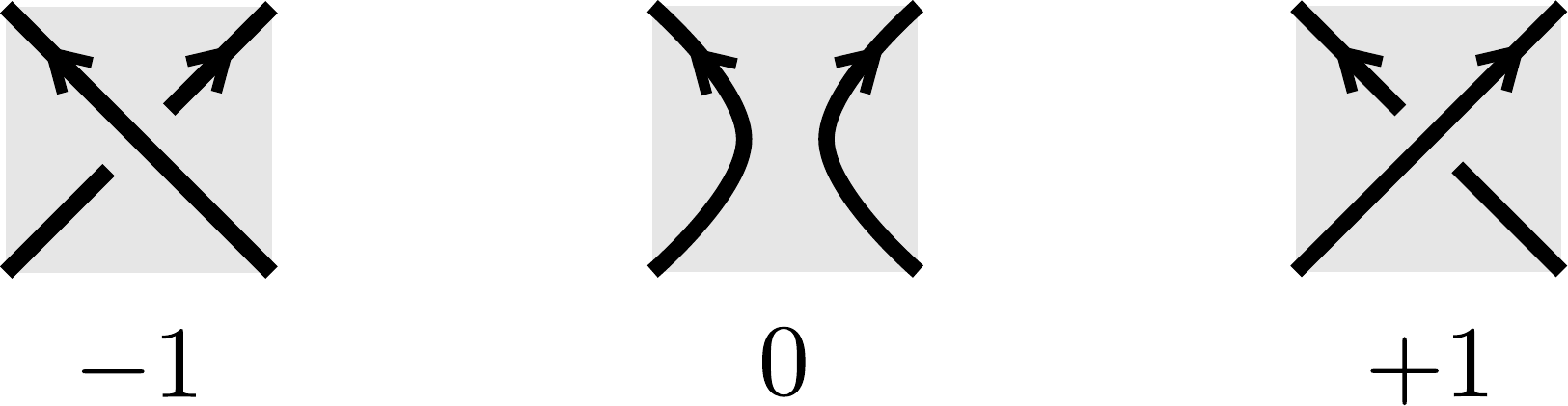}
    \caption{Signs and the corresponding local diagrams at a crossing.}
	\label{fig:crossings_signs}   
\end{figure}

Let $D_1$ be a link diagram and $D_2$ a tangle diagram. Assume that their intersections consist of $n$ transverse double points $p_1,\dots,p_n\in D_1\cap D_2$, where $n\in\bN$.
For $\varepsilon_i'\in\{-1,0,1\}$, let $(D_1\cup D_2)(\varepsilon_1',\dots,\varepsilon_n')$ denote the diagram obtained by replacing the neighborhood of each $p_i$ with the local diagram corresponding to $\varepsilon_i'$ in Figure~\ref{fig:crossings_signs}.

Let $\varepsilon_i$ be the local intersection number of $D_1$ and $D_2$ at $p_i$.
Then $(D_1\cup D_2)(\varepsilon_1,\dots,\varepsilon_n)$ is the diagram of the union in which $D_1$ lies above $D_2$.

\begin{definition}\label{def:sigma}
Let $D_1$ be a link diagram and $D_2$ be a $J$-tangle diagram. We regard them as elements of $\widetilde{\cA}(\Sigma)$ and $\widetilde{\cA}(\Sigma,J)$, respectively.
Define $\sigma\colon\widetilde{\cA}(\Sigma)\to\End(\widetilde{\cA}(\Sigma,J))$ by extending
$$\sigma(D_1)(D_2)=
\begin{cases}
\sum_{i=1}^n \varepsilon_{i}(D_1\cup D_2)(-\varepsilon_1,\dots,-\varepsilon_{i-1},0,\varepsilon_{i+1},\dots,\varepsilon_{n}) & (n\geq 1)\\
0 &(n=0)
\end{cases}$$
$\bQ\llbracket\rho\hbar,h\hbar\rrbracket[\rho,h,\hbar]$-linearly, where the intersection signs satisfy $\varepsilon_i\in\{-1,1\}$ for $1\leq i\leq n$. As shown below, Proposition~\ref{prop:sigma_welldef} and Lemma~\ref{lem:sigma_isotopy} imply that $\sigma$ is well-defined, and Lemma~\ref{lem:sigma_Lie} shows that it is a Lie action.
\end{definition}
\begin{remark}\label{rmk:Goldman_bracket_quantum}
\begin{enumerate}
    \item The map $\sigma$ is defined so that $D_1D_2-D_2D_1=h\sigma(D_1)(D_2)$. Informally, one can write $\dfrac{D_1D_2-D_2D_1}{h}=\sigma(D_1)(D_2)\medskip$. A direct computation shows that when $D_1$ and $D_2$ are single-component link diagrams, the semiclassical bracket obtained by setting $\hbar=0$ and reducing $[D_1,D_2] =\sigma(D_1)(D_2)$ modulo $h$ agrees with the Goldman bracket on their free homotopy classes. This gives a skein quantization of the Goldman Lie algebra \cite[Theorem~3.3]{Tur91}.
    The map $\sigma$ refines the extension to tangles of the formula in the proof of \cite[Theorem~3.3]{Tur91}, and can be regarded as a quantization of the Kawazumi--Kuno action \cite{KK14,KK15}.
    
    \item A similar formula appeared in the proof of Theorem~\ref{thm:welldefined}. We considered all crossings in a tangle diagram with boxes there, whereas Definition~\ref{def:sigma} only uses crossings between the link diagram $D_1$ and the tangle diagram $D_2$. Therefore, the two constructions are different.
\end{enumerate}
\end{remark}

\begin{remark}
The action $\sigma$ above is a quantization of the Kawazumi--Kuno--Goldman action \cite{KK14,KK15}, whose target is $\operatorname{Der}(\bQ\pi_1(\Sigma))$.
When $J=\emptyset$, Lemma~\ref{lem:sigma_Leibniz} shows that the target of $\sigma$ can be restricted to $\operatorname{Der}(\widetilde{\cA}(\Sigma))$. For general $J$, the action instead satisfies the Leibniz identities for the module operations stated in that lemma.
\end{remark}

\begin{prop}\label{prop:sigma_welldef}
The map $\sigma$ is independent of the numbering of the intersections.
\end{prop}
\begin{proof}
Number the points of $D_1\cap D_2$ as $p_1,\dots,p_n$.
Let $p_1',\dots,p_n'$ be the numbering obtained by interchanging $p_i$ and $p_{i+1}$, so that $p_i'=p_{i+1}$, $p_{i+1}'=p_i$, and $p_j'=p_j$ otherwise.
Since any permutation is a product of adjacent transpositions, it suffices to compare these two numberings.
The difference between the corresponding values of $\sigma$ is
\begin{align*}
&\varepsilon_i\left(\begin{array}{c}\includegraphics[scale=0.17]{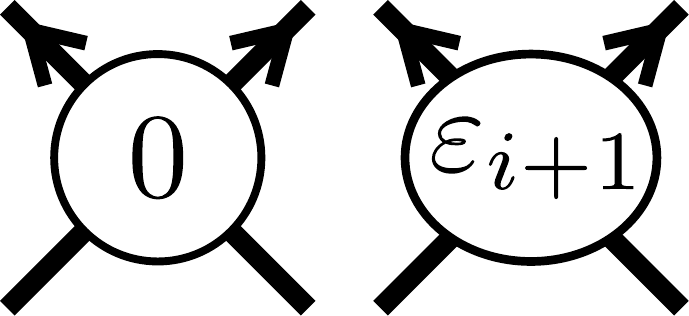}\end{array}\right)
+\varepsilon_{i+1}\left(\begin{array}{c}\includegraphics[scale=0.17]{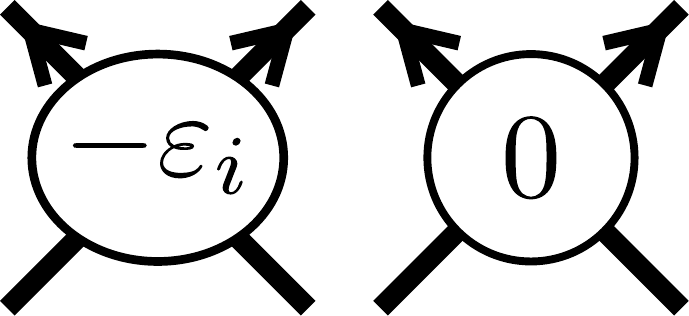}\end{array}\right)
-\left(\varepsilon_{i+1}\left(\begin{array}{c}\includegraphics[scale=0.17]{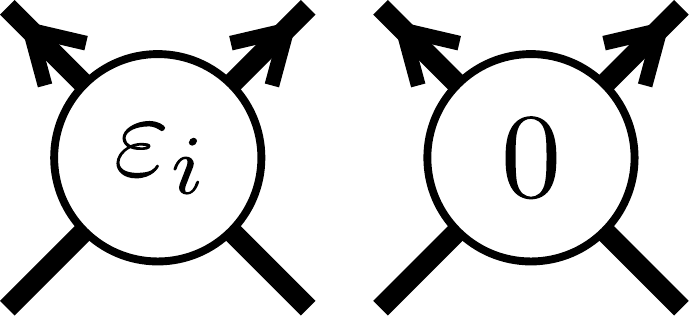}\end{array}\right)
+\varepsilon_{i}\left(\begin{array}{c}\includegraphics[scale=0.17]{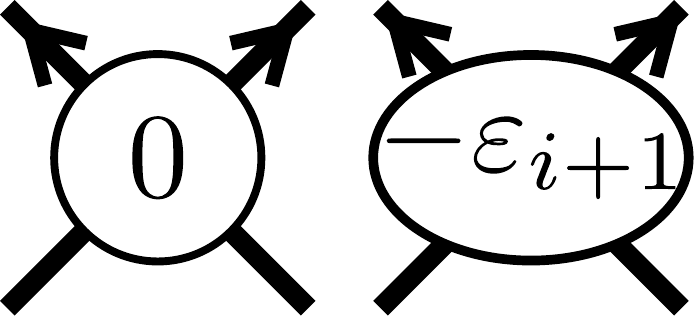}\end{array}\right)\right)\\
&= \eta\varepsilon_i\varepsilon_{i+1}\left(\begin{array}{c}\includegraphics[scale=0.17]{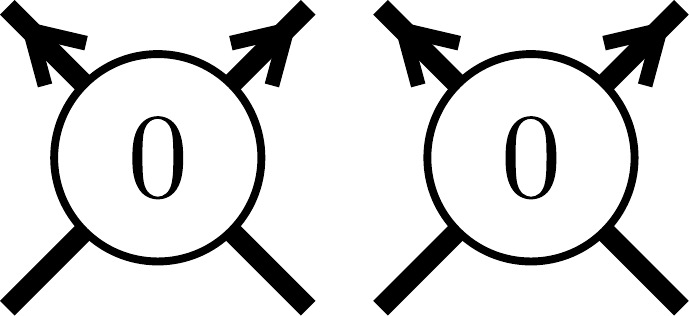}\end{array}\right)
+\eta\varepsilon_{i+1}(-\varepsilon_i)\left(\begin{array}{c}\includegraphics[scale=0.17]{draws/00.pdf}\end{array}\right)=0,
\end{align*}
where $\eta\in\{h,\hbar\}$ is determined by preservation of the component degree. The picture $\left(\begin{array}{c}\includegraphics[scale=0.17]{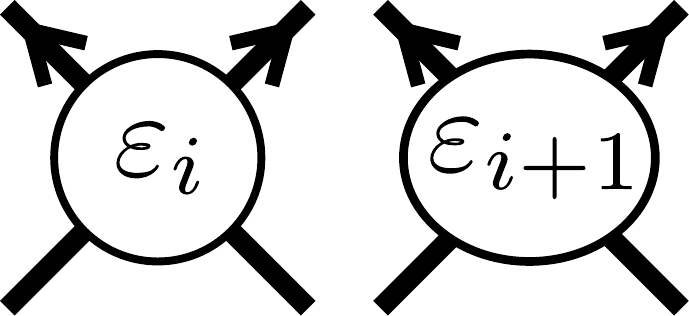}\end{array}\right)$
shows only neighborhoods of $p_i$ and $p_{i+1}$ in $(D_1\cup D_2)(\varepsilon_1,\dots,\varepsilon_n)$, with $p_i=p_{i+1}'$ on the left and $p_{i+1}=p_i'$ on the right. 
All diagrams are identical outside the shown neighborhoods.
The labels $\varepsilon_i,\varepsilon_{i+1}\in\{-1,0,1\}$ in the white circles represent the local diagrams in Figure~\ref{fig:crossings_signs}.
\end{proof}
\begin{lemma}[cf. {\cite[Lemma 3.14]{Tsu18a}}]\label{lem:sigma_isotopy}
Let $D_1$ be a link diagram and $D_2$ a $J$-tangle diagram. Then, the following hold.
\begin{enumerate}
    \item[(i)] The map $\sigma$ is compatible with the skein relations.
    
    \item[(ii)] The map $\sigma$ is invariant under Reidemeister moves of $D_1\cup D_2$, where $D_1$ lies above $D_2$ at each intersection between them.
\end{enumerate}
\end{lemma}
\begin{proof}
Let $D_1$ be a link diagram and $D_2$ be a $J$-tangle diagram. We first prove (i). A skein relation applied to $D_1$ does not affect a neighborhood of $D_1\cap D_2$. The same holds for $D_2$.

Next, we prove (ii). If a Reidemeister move involves only components of $D_1$, let $D_1'$ be the resulting diagram. Then, $\sigma(D_1)(D_2)=\sigma(D_1')(D_2)$. The same argument applies to $D_2$.

It remains to check Reidemeister moves involving intersections in $D_1\cap D_2$.
For example, suppose that $D_1$ and $D_2$ locally look like $\begin{array}{c}\includegraphics[scale=0.17]{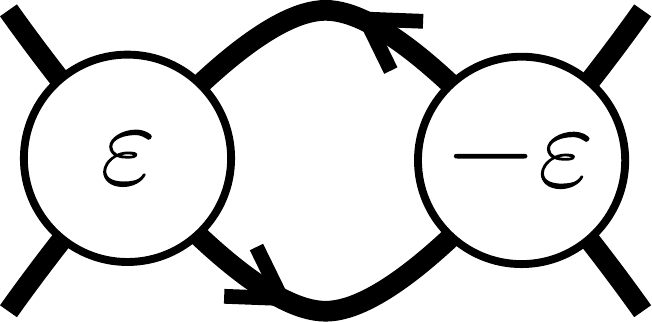}\end{array}$,where $\varepsilon\in\{-1,1\}$ indicates the crossing as in Figure~\ref{fig:crossings_signs}. Invariance under the second Reidemeister II move from the right in Figure~\ref{fig:Reidemeister} follows from Proposition~\ref{prop:sigma_welldef} and
$$
\varepsilon\begin{array}{c}\includegraphics[scale=0.17]{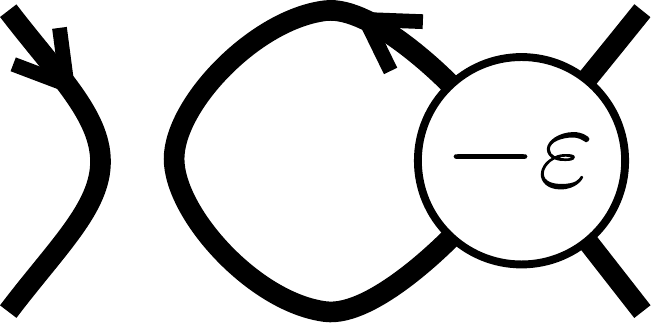}\end{array}-\varepsilon\begin{array}{c}\includegraphics[scale=0.17]{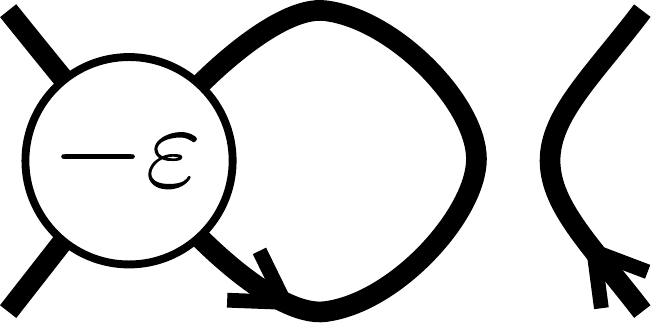}\end{array}=(\varepsilon-\varepsilon)\exp(\rho \hbar)\begin{array}{c}\includegraphics[scale=0.17]{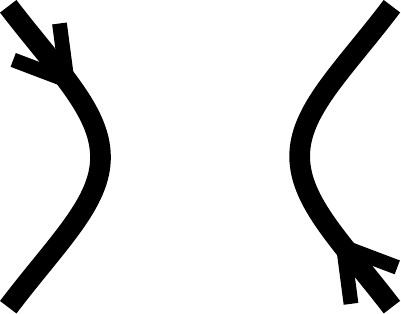}\end{array}=0$$
The second equality uses the relation (D).
The other Reidemeister moves are checked similarly.
\end{proof}

By Lemma~\ref{lem:sigma_isotopy}, for a link $L$ and a $J$-tangle $T$ represented by diagrams $D_L$ and $D_T$, the expression
\begin{align*}
\sigma(L)(T):= \sigma(D_L)(D_T)\in \widetilde{\cA}(\Sigma,J)
\end{align*}
is well-defined. Extending $\bQ\llbracket\rho\hbar,h\hbar\rrbracket[\rho,h,\hbar]$-linearly gives a well-defined action $\sigma\colon\widetilde{\cA}(\Sigma)\to\End(\widetilde{\cA}(\Sigma,J))$.

\begin{lemma}[cf. {\cite[Lemma 3.14]{Tsu18a}}]\label{lem:sigma_Lie}
\begin{enumerate}
    \item $\sigma(L_1)(L_2)=-\sigma(L_2)(L_1)$. 
    \item Set $[L_1,L_2]:=\sigma(L_1)(L_2)$. Then, $\sigma([L_1,L_2])(T)=\sigma(L_1)(\sigma(L_2)(T))-\sigma(L_2)(\sigma(L_1)(T))$.
\end{enumerate}
\end{lemma}
The case $h=\hbar$ is proved in \cite[Lemma~3.14]{Tsu18a}.
The same argument, with a suitable ordering of the crossings, applies the case when $h\neq\hbar$, so we omit the proof.

By Lemma~\ref{lem:sigma_Lie}, the bracket $[\cdot,\cdot]$ makes $\widetilde{\cA}(\Sigma)$ a Lie algebra, and $\sigma$ equips $\widetilde{\cA}(\Sigma,J)$ with the structure of a module over this Lie algebra.

Let $T_i$ be a $J_i$-tangle diagram for $i=1,2$.
Assume that $J_1,J_2\subset\partial\Sigma$ are disjoint.
For $T_i\in\widetilde{\cA}(\Sigma,J_i)$, define $T_1\boxtimes T_2\in\widetilde{\cA}(\Sigma,J_1\cup J_2)$
by placing $T_1$ in $\Sigma\times(1/2,1)$ and $T_2$ in $\Sigma\times(0,1/2)$. We call this the \textbf{tensor product} of $T_1$ and $T_2$. The bijection between $(J_1\cup J_2)_{\rm out}$ and $(J_1\cup J_2)_{\rm in}$ is inherited from those for $J_1$ and $J_2$.
We also allow the bijection to be modified by exchanging one incoming endpoint in $J_1$ with one in $J_2$, and denote the resulting operation by $T_1\mathring{\boxtimes}T_2$. This decreases the component degree by exactly $1$.
When $J_1=J_2=\emptyset$, this is the multiplication on $\widetilde{\cA}(\Sigma)$. When $J_1=\emptyset$ (resp. $J_2=\emptyset$), it gives $\widetilde{\cA}(\Sigma,J_2)$ (resp. $\widetilde{\cA}(\Sigma,J_1)$) the structure of a left (resp. right) $\widetilde{\cA}(\Sigma)$-module.
Together these actions make $\widetilde{\cA}(\Sigma,J)$ a $\widetilde{\cA}(\Sigma)$-bimodule.

\begin{figure}
\centering\includegraphics[scale=0.15]{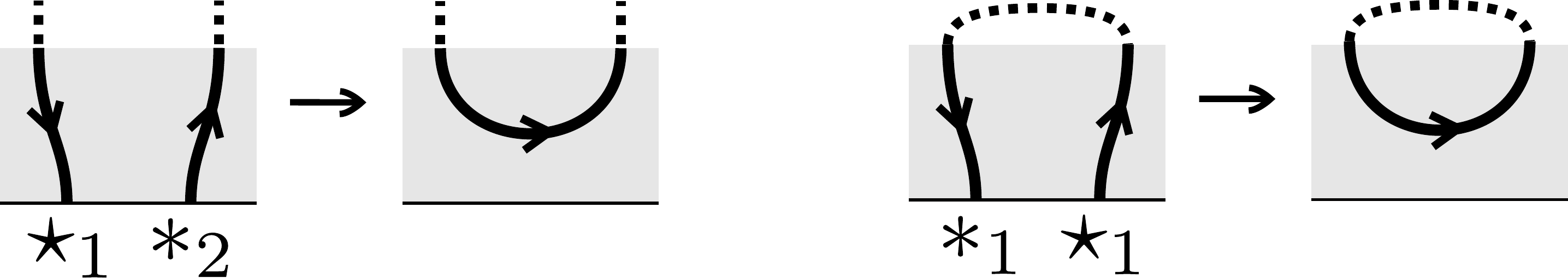}
    \caption{In each operation, there are no other endpoints between the shown ones. Left: joining endpoints for composition. Right: joining endpoints for closure.}
	\label{fig:composition_bdr}   
\end{figure}

Let $J_i=\{\ast_i,\star_i\}$ for $i=1,2$. Assume that $\ast_2$ and $\star_1$ lie on the same boundary component and either coincide or are consecutive, in that order, along the boundary orientation.
Applying the local operation in Figure~\ref{fig:composition_bdr} near the boundary to the projection of $T_1\boxtimes T_2$, and extending linearly, gives a map
$$\circ \colon \widetilde{\cA}(\Sigma,J_1)\times \widetilde{\cA}(\Sigma,J_2)\to \widetilde{\cA}(\Sigma,\{\ast_1,\star_2\})$$
Let $T_1\circ T_2$ denote the element obtained by the operation in Figure~\ref{fig:composition_bdr} and we call it the \textbf{composition} of $T_1$ and $T_2$. Our convention for composition differs from those in some earlier works; see \cite{AKKN18,AKKN20} for example.
If $\ast_1=\ast_2$ (resp. $\star_1=\star_2$), first move $\ast_2$ slightly along the boundary orientation (resp. move $\star_1$ slightly against it), and then apply the same construction. 
Suppose that $J$ lies on a single boundary component.
Let $\textbf{1}_J$ denote the $J$-tangle that is the identity for composition.
It is represented by an embedded path regular homotopic, relative to its endpoints, to a boundary interval whose endpoints appear in this order along the boundary orientation.

Suppose that $|J|=2$ and both points lie on the same boundary component of $\Sigma$. The \textbf{closure} of a $J$-tangle is the link diagram obtained by joining its endpoints as in Figure~\ref{fig:composition_bdr}. 
The closure of $\textbf{1}_J$ is a trivial loop.
Unlike composition, which involves two tangles, closure is an operation on a single tangle.

The tensor product, composition, and closure preserve the skein relations and Hence, induce the maps
\begin{align*}
\boxtimes &\colon \widetilde{\cA}(\Sigma,J_1)\times \widetilde{\cA}(\Sigma,J_2)\to \widetilde{\cA}(\Sigma,J_1\sqcup J_2),\\
\circ &\colon \widetilde{\cA}(\Sigma,J_1)\times \widetilde{\cA}(\Sigma,J_2)\to \widetilde{\cA}(\Sigma,\{\ast_1,\star_2\}),\\
\text{closure}&\colon \widetilde{\cA}(\Sigma,J)\to \widetilde{\cA}(\Sigma).
\end{align*}

These operations are also defined for tangles with boxes. Since they preserve the number of boxes, they preserve the finite-type degree.

In what follows, we omit $J$ from the term $J$-tangle when there is no need to specify the marking.

\begin{lemma}\label{lem:sigma_Leibniz}
The action $\sigma$ satisfies the Leibniz rule for the multiplication, composition, the module actions, and the tensor products. Namely,
\begin{align}
    &\sigma(L_1 L_2)(T)=\sigma(L_1)(T)L_2+L_1\sigma(L_2)(T)\label{eq:Leibniz_L}\\
    &\sigma(L)(T_1\boxtimes T_2)=(\sigma(L)(T_1))\boxtimes T_2+T_1\boxtimes(\sigma(L)(T_2))\label{eq:Leibniz_boxtimes}\\
    &\sigma(L)(T_1\circ T_2)=\sigma(L)(T_1)\circ T_2+T_1\circ \sigma(L)(T_2)\label{eq:Leibniz_circ}\\
    &\text{setting $[L,L_1]:=\sigma(L)(L_1)$, }\sigma(L)(L_1T_1)=[L,L_1]T_1+L_1\sigma(L)(T_1)\label{eq:Leibniz_bracket}
\end{align}
Moreover, closure commutes with $\sigma(x)$.
\end{lemma}

\begin{proof}
Since closure is performed near the boundary, Lemma~\ref{lem:sigma_isotopy} shows that it commutes with $\sigma(x)$.

We prove \eqref{eq:Leibniz_L}.
Represent $L_1,L_2$ by link diagrams and $T$ by a tangle diagram on $(\Sigma,J)$, with transverse intersections in $(L_1L_2)\cup T$. Number the intersections of $L_2$ with $T$ as $p_1,\dots,p_n$, and those of $L_1$ with $T$ as $p_{n+1},\dots,p_{n+m}$. Let $D$ denote the diagram representing the stacked product of $L_1L_2$ with $T$. Then
\begin{align}
    \sigma(L_1L_2)(T)&=\sum_{i=1}^{n+m} \varepsilon_i D(-\varepsilon_1,\dots, -\varepsilon_{i-1},0,\varepsilon_{i+1},\dots,\varepsilon_{n+m})\nonumber\\
    &=\sum_{i=1}^n \varepsilon_i D(-\varepsilon_1,\dots, -\varepsilon_{i-1},0,\varepsilon_{i+1},\dots,\varepsilon_{n+m})+\sum_{i=n+1}^{n+m} \varepsilon_i D(-\varepsilon_1,\dots, -\varepsilon_{i-1},0,\varepsilon_{i+1},\dots,\varepsilon_{n+m})\label{eq:eq_leibnitz}\\
    &=L_1 \sigma(L_2)(T)+\sigma(L_1)(T)L_2\nonumber
\end{align}
which proves \eqref{eq:Leibniz_L}.
Here, $\varepsilon_{n+1},\dots,\varepsilon_{n+m}$ place $L_1$ above $T$, whereas $-\varepsilon_1,\dots,-\varepsilon_n$ place $T$ above $L_2$.

Equation \eqref{eq:Leibniz_boxtimes} follows by the same argument as \eqref{eq:Leibniz_L}, numbering the intersections of $L$ with $T_1$ as $p_1,\dots,p_n$ and those of $L$ with $T_2$ as $p_{n+1},\dots,p_{n+m}$.

Since composition is performed near the boundary, Lemma~\ref{lem:sigma_isotopy} and the case $|J|=2$ of \eqref{eq:Leibniz_boxtimes} imply \eqref{eq:Leibniz_circ}.

Equation \eqref{eq:Leibniz_bracket} is the case of \eqref{eq:Leibniz_boxtimes} in which the first tangle is the link $L_1$.
\end{proof}

\begin{lemma}\label{lem:sigma_filtration}
We have $\sigma(F^n\widetilde{\cA}(\Sigma))(F^m\widetilde{\cA}(\Sigma,J))\subset F^{n+m-k}\widetilde{\cA}(\Sigma,J)$,
where $k=\min\{2,n+1,m+1,n+m\}$, or equivalently,
$$k=\begin{cases}
0&m=n=0\\
1&\text{$mn=0$ and $m+n>0$}\\
2&mn\neq 0\\
\end{cases}$$
\end{lemma}

\begin{proof}
The inclusion $\sigma(F^0\widetilde{\cA}(\Sigma))(F^0\widetilde{\cA}(\Sigma,J))\subset F^0\widetilde{\cA}(\Sigma,J)=\widetilde{\cA}(\Sigma,J)$ holds immediately.

Consider $\sigma(L)(T)$ for a link diagram with boxes $L\in F^n\widetilde{\cA}(\Sigma)$ and a $J$-tangle diagram with boxes $T\in F^m\widetilde{\cA}(\Sigma,J)$.
Replacing all boxes of $L$ and $T$ with their ds-linear combinations gives $\sum_i b_iL_i$ and $\sum_jc_jT_j$, where $b_i,c_j\in\{-1,1\}$ and $L_i,T_j$ are fixed diagrams without boxes rather than isotopy classes.  
Fix compatible orders on the intersections $L_i\cap T_j$ to define $\sigma$. Outside the boxes of $L$ (resp. $T$), the diagram $L_i$ (resp. $T_j$) agrees with $L$ (resp. $T$). We call the parts of $T_j$ that coincide with detour components of $T$ the \textbf{detour parts} of $T_j$.
Each term $\sigma(L_i)(T_j)$ in $\sigma(L)(T)=\sum_{i,j}b_ic_j\sigma(L_i)(T_j)$ is obtained by taking the signed sum of the smoothings at intersections of $L_i$ and $T_j$.
A smoothing in the definition of $\sigma$ does not affect detour parts away from that intersection. Thus, at most two detour parts, and Hence, at most two boxes, are affected. 
Therefore, the result can be written as a ds-linear combination in which at most two boxes have been removed. 
In particular, if $mn=0$ and $m+n>0$, we obtain $\sigma(F^n\widetilde{\cA}(\Sigma))(F^m\widetilde{\cA}(\Sigma,J))\subset F^{n+m-1}\widetilde{\cA}(\Sigma,J)$, regardless of whether the chosen representatives $L,T$ contain boxes.
\end{proof}

\subsection{The free $\widetilde{\cA}(\Sigma)$-module structure of $\widetilde{\cA}(\Sigma,J)$}\label{subsec:module_str_HOMFLYPT}
We now show that $\widetilde{\cA}(\Sigma,J)$ is a free $\widetilde{\cA}(\Sigma)$-module. This will be used in the proof of Theorem~\ref{thm:inj}.

Let $J=\{\ast,\star\}$, where $\star$ is obtained by moving $\ast$ slightly along the orientation of the boundary. Define
\begin{align*}
&\pi_J^+:=\{\text{paths on $\Sigma$ from $\ast$ to $\star$}\}/\text{regular homotopy},\\
&\cP_J:=\bQ\llbracket\rho\hbar\rrbracket[\rho,\hbar]\pi_J^+/(\text{the relation (D) with crossing information forgotten}).
\end{align*}
For the augmentation $\varepsilon\colon\bZ\pi_J^+\to\bZ$, $\gamma\mapsto1$ for $\gamma\in\pi_J^+$, let $F^n\cP_J$ be the $\bQ\llbracket\rho\hbar\rrbracket[\rho,\hbar]$-submodule generated by $\{\hbar^iL\mid L\in(\ker\varepsilon)^j,\ 2i+j\geq n\}$.
Regular homotopy classes of paths are regarded as morphisms in a groupoid. For the fixed endpoint set $J$ on a single boundary component, composition convention in Figure~\ref{fig:composition_bdr} equips $\pi_J^+$ with a group structure.
The family $\{F^n\cP_J\}_{n\in \bN}$ is a filtration. Let $\Gr\cP_J$ denote its associated graded module.

\begin{remark}
More generally, $\cP_{J'}$ and $F^n\cP_{J'}$ can be defined when $J'=\{\ast',\star'\}$ has its endpoints on distinct boundary components.
Choose paths $\gamma_\ast$ from $\ast$ to $\ast'$ and $\gamma_\star$ from $\star$ to $\star'$. Extending $\alpha\mapsto\gamma_\ast^{-1}\alpha\gamma_\star$, $\alpha\in\pi_J^+$, linearly gives a $\bQ\llbracket\rho\hbar\rrbracket$-module isomorphism $\cP_J\to\cP_{J'}$. Transporting the filtration along this isomorphism defines $F^n\cP_{J'}$ and $\Gr\cP_{J'}$.
\end{remark}

Let $J_i=\{\ast_i,\star_i\}$ and $J=\bigsqcup_{i=1}^nJ_i$, with $J_{\rm out}=\{\ast_1,\dots,\ast_n\}$. For $\tau\in\mathfrak{S}_n$, set $\tau(J_i)=\{\ast_i,\tau(\star_i)\}$.
Take $\bigotimes_{i=1}^n\cP_{\tau(J_i)}$ over $\bQ\llbracket\rho\hbar\rrbracket[\rho,\hbar]$.
For $(T_1,\dots,T_n)\in\prod_{i=1}^n\pi_{\tau(J_i)}^+$,
consider $T_1\otimes\dots\otimes T_n\in\bigotimes_{i=1}^n\cP_{\tau(J_i)}$. Let $\widetilde T_i\in\widetilde{\cA}(\Sigma,\tau(J_i))$ be obtained from $T_i$ by assigning crossing information so that
\begin{enumerate}
    \item The projection of $\widetilde T_i$ is $T_i$,
    
    \item The tangle $\widetilde T_i$ lies below $\widetilde T_j$ for $j=1,\dots,i-1$,
    
    \item The height decreases as we traverse $T_i$ from its initial point; that is, the lift is descending.
\end{enumerate}
Set $\cP_J:=\bigoplus_\tau\bigotimes_{i=1}^n\cP_{\tau(J_i)}$. Extending the above construction linearly defines $\widetilde T$ for every $T\in\cP_J$. Note that all tensor products are taken over $\bQ\llbracket\rho\hbar\rrbracket[\rho,\hbar]$.
\begin{align}
f_J\colon \widetilde{\cA}(\Sigma)\otimes \cP_J\to \widetilde{\cA}(\Sigma,J),\quad L\otimes T \mapsto L\widetilde{T}\label{eq:map_fJ}
\end{align}
This defines $f_J$. Since the relation in $\cP_J$ corresponds to (D), the map $f_J$ is
a homomorphism of left $\widetilde{\cA}(\Sigma)$-modules.

\begin{prop}[{cf.\cite[Propositions 3.4, 3.5]{Tsu18a}}]\label{prop:free_module}
\begin{enumerate}
    \item As a left $\widetilde{\cA}(\Sigma)$-module, $\widetilde{\cA}(\Sigma,J)$ is isomorphic to the free module $\widetilde{\cA}(\Sigma)\otimes\cP_J$.
    In particular, when $|J|=2$, a basis is obtained from a basis of $\cP_J$ over $\bQ\llbracket\rho\hbar\rrbracket[\rho,\hbar]$.
    \item 
    Equip $\widetilde{\cA}(\Sigma)\otimes\cP_J$ with the tensor product filtration $\{F^n(\widetilde{\cA}(\Sigma)\otimes\cP_J)\}_{n\in\bN}$ induced by $\{F^n\widetilde{\cA}(\Sigma)\}_{n\in\bN}$ and $\{F^n\cP_J\}_{n\in\bN}$. Then the isomorphism in (1) preserves this filtration.
\end{enumerate}
\end{prop}
Part (1) is proved in Propositions~3.4 and~3.5 of \cite{Tsu18a}. Here, we extend the complexity used in \cite[Propositions~3.4 and~3.5]{Tsu18a}
and give a proof using the diamond lemma from confluence theory.

\begin{figure}
\centering\includegraphics[width=12cm]{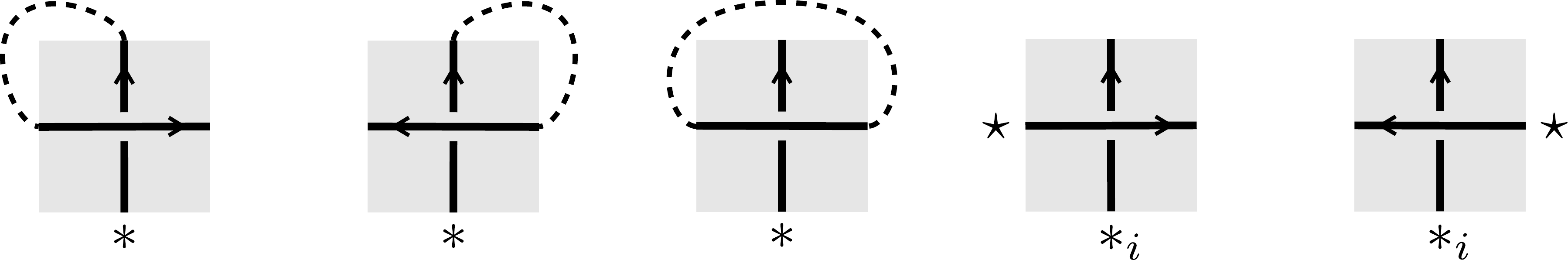}
	\caption{$\ast\in J_{\rm out}\subset J,\quad \star\in \{\ast_{i+1},\dots,\ast_{n}\}$}
	\label{fig:confluence}   
\end{figure}

\begin{proof}
Let $\Jout=\{\ast_1,\ast_2,\dots,\ast_n\}$.
Let ${\rm Cross}(D)$ be the set of crossings of $D$, and let ${\rm Cross}_k(D)$ be the subset of crossings of the types shown in Figure~\ref{fig:confluence} with $i=k$. Figure~\ref{fig:confluence} lists the crossings that prevent a diagram from being the image under $f_J$ of a tensor product of a link and an element of $\prod_{i=1}^n\pi_{\tau(J_i)}^+$ for some $\tau\in\mathfrak{S}_n$.
Let $\mathcal{D}(N,N_1,\dots,N_n)$ denote the set of $J$-tangle diagrams satisfying
$$(\#{\rm Cross}(D),\#{\rm Cross}_1(D),\dots,\#{\rm Cross}_n(D))<(N,N_1,\dots,N_n)$$ in lexicographic order.
For a $J$-tangle diagram $D$, we call the tuple
$$(\#{\rm Cross}(D), \#{\rm Cross}_1(D),\dots ,\#{\rm Cross}_n(D))$$
the \textbf{complexity} of $D$.

Consider two crossings of a $J$-tangle diagram of the types shown in Figure~\ref{fig:confluence}. We show that reducing their complexity gives the same result in either order.
The defining relation
\begin{align} 
\begin{array}{c}\includegraphics[scale=0.20]{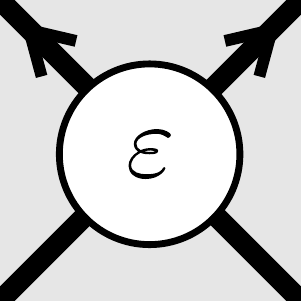}\end{array}=\begin{array}{c}\includegraphics[scale=0.20]{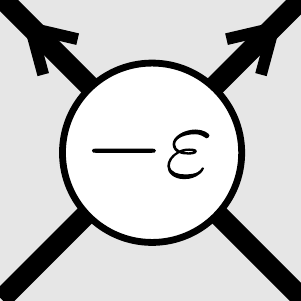}\end{array}+\varepsilon\eta\begin{array}{c}\includegraphics[scale=0.20]{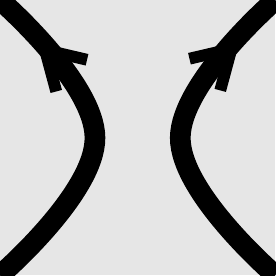}\end{array}\label{eq:skein_relations_complexity}
\end{align}preserves the component degree.
In \eqref{eq:skein_relations_complexity}, $\varepsilon\in\{-1,1\}$ indicates the crossing in Figure~\ref{fig:crossings_signs}, and $\eta$ is $h$ or $\hbar$ according to whether (A) or (B) is applied.
When both crossings are reduced, corresponding terms in the resulting linear combinations have the same component degree and involve the same number of applications of \eqref{eq:skein_relations_complexity}.
Thus, the exponents of $h$ and $\hbar$ in each term are uniquely determined. The resulting $J$-tangle diagrams may differ, but an argument similar to the proof of Proposition~\ref{prop:sigma_welldef} shows that they agree modulo the defining relations of $\widetilde{\cA}(\Sigma,J)$.
Hence, the resulting linear combination is independent of the order in which \eqref{eq:skein_relations_complexity} is applied at the two crossings. 
The trivial-loop relation also commutes with the other reductions, so the system is locally confluent. Moreover, by applying \eqref{eq:skein_relations_complexity}, the complexity of $D$ decreases.
By the diamond lemma \cite{SW07}, the linear combination obtained when the complexity reaches $(k,0,\dots,0)$ is independent of the reductions chosen, where $k$ is the total number of crossings and need not be zero.
At complexity $(k,0,\dots,0)$, the closed and open components are separated, giving a linear combination of elements of $\widetilde{\cA}(\Sigma)\otimes\cP_J$. Trivial loops may still occur, but they can be removed using the relation (C).
Compatibility with (D) follows from the defining relation of $\cP_J$.

Since the arc components have specified initial points, invariance under Reidemeister moves involving arcs follows by extending Lickorish's argument \cite[Section~1]{Lic88} to general surfaces; see also \cite[Lemma~4.3]{Tur91}. %This proves invariance under Reidemeister moves.
Invariance under Reidemeister moves involving only $S^1$-components follows since these moves do not affect the complexity.
It follows that $\widetilde{\cA}(\Sigma,J)$ is a free $\widetilde{\cA}(\Sigma)$-module.

The reduction to complexity $(k,0,\dots,0)$ proves surjectivity of the map $f_J$ in \eqref{eq:map_fJ}.
Independence of the resulting linear combination from the reduction sequence proves injectivity of $f_J$.

For (2), the degree on $\cP_J$ can be represented using boxes. Each defining relation of $\widetilde{\cA}(\Sigma,J)$ used in the local reductions in (1) does not decrease the filtration degree. This proves the claim.
\end{proof}

\begin{remark}
The relation $\begin{array}{c}\includegraphics[scale=0.18]{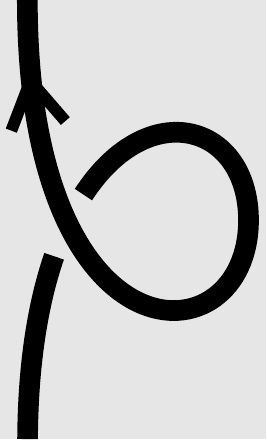}\end{array}=\exp(-\rho\hbar)\begin{array}{c}\includegraphics[scale=0.18]{draws/kink_no.pdf}\end{array}$ follows from (A)--(D). The relation (D) does not appear in Figure~\ref{fig:confluence}, and is not used in the complexity reductions.
\end{remark}

Let $\partial_1,\dots,\partial_b$ be the boundary components of $\Sigma$.
Let $\gamma_i$ be a simple path from $\partial_1$ to $\partial_i$. For $J\subset\partial_1$ with $|J|=2$, set $\gamma_1=\textbf{1}_J$.
Write $\gamma_i^{-1}$ for the same path with its orientation reversed.
Choose a generating set of $\pi_1(\Sigma)$, based on $\partial_1$, such that moving the incoming endpoint of each path slightly along the boundary orientation gives a diagram without self-crossings. Let $\cB_+$ be the set of the resulting diagrams and $\cB_-$ the set obtained by reversing their orientations.
Set
\begin{align}
\widetilde{\cP}_+:=\cB_+ \cup \{\gamma_i\mid  i=2,\dots, b\},\qquad
\widetilde{\cP}_-:=\cB_- \cup \{\gamma_i^{-1}\mid  i=2,\dots, b\},\qquad\widetilde{\cP}:=\widetilde{\cP}_+\cup\widetilde{\cP}_-\label{eq:tildeP}
\end{align}

Let $J,J'$ be markings with $|J|=|J'|=2$, and suppose that $J_{\rm out},J'_{\rm out}$ (resp. $J_{\rm in},J'_{\rm in}$) lie on the same boundary component $\partial$. Choose a neighborhood $N(\partial)$ in $\Sigma$ and a homeomorphism that is the identity on $\Sigma\setminus N(\partial)$ and sends $J_{\rm out}$ to $J'_{\rm out}$ (resp. $J_{\rm in}$ to $J'_{\rm in}$). We call this operation \textbf{sliding}.
We have the following lemma.
\begin{lemma}\label{lem:generating_tildeP}
The collection $\widetilde{\cA}(\Sigma)\sqcup\bigsqcup_{|J|=2}\widetilde{\cA}(\Sigma,J)$ is generated by $\widetilde{\cP}$ under the left action of $\widetilde{\cA}(\Sigma)$, sliding, composition, and closure.
\end{lemma}
\begin{proof}
We first show that $\widetilde{\cA}(\Sigma)$ is generated by $\widetilde{\cP}$. It suffices to obtain link diagrams by composing and closing elements of $\widetilde{\cP}$. We use mathematical induction on the number of crossings. The case of no crossings is immediate. Suppose the claim holds for diagrams with at most $k$ crossings. For a diagram with $k+1$ crossings, the relation (A) or (B) expresses it as a linear combination of the diagram with a chosen crossing switched and a diagram with $k$ crossings. The latter is covered by the induction hypothesis. Thus, we may switch crossings so that the $S^1$-components lie in disjoint height intervals in $\Sigma\times[0,1]$: there is a subdivision of $[0,1]$ such that the product of each subinterval with $\Sigma$ contains only one $S^1$-component. We may also move a point of each component near the boundary and switch its self-crossings so that it is descending when we move from that point in its orientation. Each resulting component is obtained by composing and closing elements of $\widetilde{\cP}$.

For $\bigsqcup_{|J|=2}\widetilde{\cA}(\Sigma,J)$, Proposition~\ref{prop:free_module} separates the closed components from the arc component. The latter represents an element of $\cP_J$, and, after sliding if necessary, $\widetilde{\cP}$ generates $\cP_J$.
Together with sliding and the left action of $\widetilde{\cA}(\Sigma)$, this proves that $\widetilde{\cP}$ also generates $\bigsqcup_{|J|=2}\widetilde{\cA}(\Sigma,J)$.
\end{proof}

%%%%%%%%%%%%%%%%%%%%%%%%%%%%%%%%%%%%%%%%%%%%%%%%%%%%%%%%%%%%%%%%%%%%%%%%%%%%%%%%
\section{Proof of Theorem~\ref{thm:inj}}\label{sec:proof_mainthm}
Throughout this section, unless otherwise stated, let $\Sigma$ be a compact connected oriented surface with non-empty boundary, and let $(\Sigma,J)$ be a marked surface.
For $x\in\Gr^n\widetilde{\cA}(\Sigma)$, consider the action $\sigma(x)(\,\cdot\,)$
on $\Gr^m\widetilde{\cA}(\Sigma,J)$ for each $m\geq0$. This gives a $\bQ[\rho]$-module homomorphism
\begin{align}
\sigma_{n}\colon \Gr^n\widetilde{\cA}(\Sigma)\to\displaystyle\prod_{m,J}{\rm Hom}(\Gr^m\widetilde{\cA}(\Sigma, J),\Gr^{k}\widetilde{\cA}(\Sigma, J)),  \quad x\mapsto \sigma(x),\label{eq:sigma_n}
\end{align}
where $\Gr^{-1}\widetilde{\cA}(\Sigma)=\{0\}$ and $k=n+m-\min\{2,n+1,m+1\}$ records the decrease in the number of boxes. 
Let ${\rm Hom}(\Gr^m\widetilde{\cA}(\Sigma,J),\Gr^k\widetilde{\cA}(\Sigma,J))$
denote the $\bQ[\rho]$-module of $\bQ[\rho]$-module homomorphisms from $\Gr^m\widetilde{\cA}(\Sigma,J)$ to $\Gr^k\widetilde{\cA}(\Sigma,J)$.

\begin{lemma}\label{lem:gene_groupoid}
For any $x\in\widetilde{\cA}(\Sigma)$,
the action of $\sigma(x)$ on the generating set $\widetilde{\cP}_+$ in \eqref{eq:tildeP} determines its action on
$\widetilde{\cA}(\Sigma,J')$ for every $J'$.
\end{lemma}
\begin{proof}
For $|J|=0,2$, Lemma~\ref{lem:generating_tildeP} and the Leibniz rule (Lemma~\ref{lem:sigma_Leibniz}) imply that the action of $\sigma(x)$ on $\widetilde{\cP}$ determines its action on
$\widetilde{\cA}(\Sigma,J)$ for every $J$ with $|J|=0,2$.
Moreover, the Leibniz rule shows that the action on the generating set $\widetilde{\cP}_+$ determines the action on $\widetilde{\cP}$. This proves the assertion for $|J'|=0,2$.

For general $J'$,
relations (A)--(D) allow us to express any $J'$-tangle diagram as a linear combination of diagrams of complexity $(k,0,\dots,0)$, where complexity is defined in the proof of Proposition~\ref{prop:free_module}.
A computation analogous to \eqref{eq:eq_leibnitz} shows that the action of $\sigma(x)$ on a $J'$-tangle diagram of complexity $(k,0,\dots,0)$ is determined by its action on each connected component.
Each connected component has either no endpoints or two endpoints, so the assertion follows from the cases above.
\end{proof}

The following proposition gives generators of the skein module as a filtered module.
\begin{prop}\label{prop:closure_composition}
Suppose that $|J|=|J'|=|J''|=2$ and that the conditions for composition hold with $J_1=J'$ and $J_2=J''$. Namely,
assume that $J=J'_+\sqcup J''_-$, that $J'_-$ and $J''_+$ lie on the same boundary component, and that either $J'_-=J''_+$ or the point of $J'_-$ immediately follows the point of $J''_+$ along the boundary orientation.
\begin{enumerate}
    \item If $J'_+$ and $J''_-$ lie on the same boundary component, every element of $F^{n+1}\widetilde{\cA}(\Sigma)$ is a $\bQ$-linear combination of closures of compositions of elements of $F^1\widetilde{\cA}(\Sigma,J')$ and $F^n\widetilde{\cA}(\Sigma,J'')$, and elements of $hF^{n-1}\widetilde{\cA}(\Sigma)+\hbar F^{n-1}\widetilde{\cA}(\Sigma)$.\item Every element of $F^{n+1}\widetilde{\cA}(\Sigma,J)$ is a $\bQ$-linear combination of elements obtained by composing elements of $F^1\widetilde{\cA}(\Sigma,J')$ and $F^n\widetilde{\cA}(\Sigma,J'')$, by the left action of $F^{n+1}\widetilde{\cA}(\Sigma)$, and from $hF^{n-1}\widetilde{\cA}(\Sigma,J)+\hbar F^{n-1}\widetilde{\cA}(\Sigma,J)$.
\end{enumerate}
Here, $F^{-1}\widetilde{\cA}(\Sigma):=\widetilde{\cA}(\Sigma)$.
\end{prop}
\begin{proof}
For the first assertion, it suffices to show that every $h^i\hbar^jL\in F^n\widetilde{\cA}(\Sigma)$, where $L$ is a link with boxes, is closure of an element of $F^n\widetilde{\cA}(\Sigma,J)$.
By Lemma~\ref{lem:sigma_isotopy}, Reidemeister moves and isotopies of $\Sigma$ allow us to move a portion of $L$ close to the boundary, as in the right-hand diagram of Figure~\ref{fig:composition_bdr}. Let $T$ be the $J$-tangle with boxes shown in the left-hand diagram of Figure~\ref{fig:composition_bdr}. Its closure is $L$. Since $L$ and $T$ have the same number of boxes, we have $h^i\hbar^jT\in F^n\widetilde{\cA}(\Sigma,J)$, as required.

By a similar observation, one can show the second assertion easily. 
\end{proof}

Since the operations in Propositions~\ref{prop:free_module} and~\ref{prop:closure_composition} and Lemmas~\ref{lem:generating_tildeP} and~\ref{lem:gene_groupoid} preserve the finite-type filtration, their statements pass to the associated graded modules. Therefore, we obtain the following lemma.

\begin{lemma}\label{lem:minimal_generators}
Let $J'\subset\partial_1$ with $|J'|=2$.
Under the left action of $\Gr\widetilde{\cA}(\Sigma)$, sliding, composition, and closure,
the collection $\Gr\widetilde{\cA}(\Sigma)\sqcup\bigsqcup_{|J|=2}\Gr\widetilde{\cA}(\Sigma,J)$ is generated by
$$\cB'_+:=\{b- \textbf{1}_{J'}\mid b\in \cB_+\}\subset \Gr^1\widetilde{\cA}(\Sigma,J'),\qquad \{\gamma_i^{\pm 1}\mid i=1,2,\dots ,b\}\subset \bigsqcup_{|J|=2}\Gr^0\widetilde{\cA}(\Sigma,J),$$
where $\textbf{1}_{J'}$ denotes the identity $J'$-tangle for composition, and $\gamma_1=\textbf{1}_{\partial\gamma_1}$.
\end{lemma}
\begin{proof}
By Proposition~\ref{prop:closure_composition}, the collection $\Gr\widetilde{\cA}(\Sigma)\sqcup\bigsqcup_{|J|=2}\Gr\widetilde{\cA}(\Sigma,J)$ is generated by $\bigsqcup_{|J|=2}(\Gr^0\widetilde{\cA}(\Sigma,J)\sqcup\Gr^1\widetilde{\cA}(\Sigma,J))$.

First, in degree zero, neither the parameters $h,\hbar$ nor boxes occur.
Suppose that the $J$-tangle joins $\partial_i$ to $\partial_j$.
By Lemma~\ref{lem:transparency} and relation (C), the module $\Gr^0\widetilde{\cA}(\Sigma,J)$ is generated over $\bQ[\rho]$ by $\gamma_i^{-1}\gamma_j\ (i=1,\dots,b)$.
Here, by sliding the endpoints, we may assume that the initial and terminal points of the $J$-tangle agree with those of $\gamma_i^{-1}\gamma_j$. Since $\gamma_1=\textbf{1}_{\partial\gamma_1}$, the module $\Gr^0\widetilde{\cA}(\Sigma,J)$ is spanned over $\bQ[\rho]$ by compositions of the $\gamma_i^{\pm1}$.

By sliding and composing with the $\gamma_i^{\pm1}$, we see that $\bigsqcup_{|J|=2}\Gr^1\widetilde{\cA}(\Sigma,J)$ is generated by $\Gr^1\widetilde{\cA}(\Sigma,J')$ and $\{\gamma_i^{\pm1}\mid i=1,2,\dots,b\}$.
Next, we show that $\Gr^1\widetilde{\cA}(\Sigma,J')$ is generated by $\Gr^0\widetilde{\cA}(\Sigma,J')$, knots with boxes, and paths with boxes.
The parameters $h,\hbar$ do not occur in degree one.
By Lemma~\ref{lem:transparency}, every element of $\Gr^1\widetilde{\cA}(\Sigma,J')$ is a $\bQ[\rho]$-linear combination of the following: 
$$\begin{cases}
\text{a product of a knot with a box and $\textbf{1}_{J'}$} &\text{if the box lies on an $S^1$-component},\\
\text{a path with a box} &\text{if the box lies on an arc component}.
\end{cases}$$
A knot with a box can be realized as the closure of a path with a box. Hence, $\Gr^1\widetilde{\cA}(\Sigma,J')$ is generated by paths with boxes and $\textbf{1}_{J'}$.
As in the proof of Theorem~\ref{thm:welldefined}, paths with boxes modulo $F^2\widetilde{\cA}(\Sigma,J')$ can be identified with elements of the homology group $H_1(\Sigma,\bQ[\rho])$.
In particular, $\cB'_+$ is a generating set of $H_1(\Sigma,\bQ[\rho])$.

This proves the assertion.
\end{proof}

The assumption that $\Sigma$ is a compact connected oriented surface with non-empty boundary is essential in the following lemma.
\begin{lemma}\label{lem:sigma_generators_gr}
Let $J'\subset\partial_1$ with $|J'|=2$, let $J$ be any marking,
and let $x\in\Gr^n\widetilde{\cA}(\Sigma)$.
The action of $\sigma_n(x)$
on $\Gr\widetilde{\cA}(\Sigma,J)$ is determined by its action on $\bigsqcup_J\Gr^0\widetilde{\cA}(\Sigma,J)\sqcup\Gr^1\widetilde{\cA}(\Sigma,J')$.
In particular, it is determined by the action of $\sigma_n(x)$ on
$$\cB'_+:=\{b- \textbf{1}_{J'}\mid b\in \cB_+\}\subset \Gr^1\widetilde{\cA}(\Sigma,J'),\qquad \{\gamma_i\mid i=2,\dots ,b\}\subset \bigsqcup_{|J|=2}\Gr^0\widetilde{\cA}(\Sigma,J)$$
alone.
Here $\textbf{1}_{J'}$ denotes the identity $J'$-tangle for composition, and $\gamma_1=\textbf{1}_{\partial\gamma_1}$.
\end{lemma}
\begin{proof}
For any marking $J$ and any $x\in\Gr^n\widetilde{\cA}(\Sigma)$, the action of $\sigma_n(x)$
on $\Gr\widetilde{\cA}(\Sigma,J)$
is determined by its action on $\Gr\widetilde{\cA}(\Sigma)\sqcup\bigsqcup_{|J|=2}\Gr\widetilde{\cA}(\Sigma,J)$, by passing the second part of the proof of Lemma~\ref{lem:gene_groupoid} to the associated graded modules.
Moreover, by Lemmas~\ref{lem:sigma_Leibniz} and~\ref{lem:minimal_generators}, the action of $\sigma_n(x)$ on $\Gr\widetilde{\cA}(\Sigma)\sqcup\bigsqcup_{|J|=2}\Gr\widetilde{\cA}(\Sigma,J)$ is determined by its action on
$\cB'_+\sqcup\{\gamma_i^{\pm1}\mid i=1,2,\dots,b\}$.

Since $\sigma_n(x)$ acts trivially on $\gamma_1^{\pm1}=\textbf{1}_{J'}^{\pm1}$, we may omit these elements.

Lemma~\ref{lem:sigma_Leibniz} gives $$0=\sigma_n(x)(\textbf{1}_{J'})=\sigma_n(x)(\gamma_i\gamma_i^{-1})=\sigma_n(x)(\gamma_i)\gamma_i^{-1}+\gamma_i\sigma_n(x)(\gamma_i^{-1}).$$
Thus the action of $\sigma_n(x)$ on $\gamma_i$ determines its action on $\gamma_i^{-1}$.
It follows that the action on $\{\gamma_i\mid i=2,\dots,b\}$ determines the action on $\{\gamma_i^{\pm1}\mid i=1,2,\dots,b\}$.

This proves the lemma.
\end{proof}

\begin{lemma}\label{lem:ker_sigma}
Let $R:=\bQ\llbracket\rho\hbar,h\hbar\rrbracket[\rho,h,\hbar]$.
For the map $\sigma_n$ in \eqref{eq:sigma_n}, if $x\in\Ker(\sigma_n)$, then $x=\Psi(x)\emptyset$ for some $\Psi(x)\in\Gr^nR$. Here, $\Gr^nR$ is the associated graded module with respect to the finite-type filtration on the coefficient ring, without boxes.
\end{lemma}
An analogous statement and proof for the Kauffman bracket skein algebra can be found in \cite[Lemma~3.8]{Tsu18b}.
\begin{proof}
Attach bands joining the boundary components of $\Sigma$ until the resulting oriented surface has a single boundary component. Take two copies $\widetilde{\Sigma}_i$, $i=1,2$, of this surface, and let $\widetilde{\Sigma}$ be their boundary connected sum.
For the two copies $\Sigma_i$, $i=1,2$, of $\Sigma$, there are natural embeddings $\iota_i\colon\Sigma_i\to\widetilde{\Sigma}$. These induce homomorphisms $\iota_{i\ast}\colon\widetilde{\cA}(\Sigma_i)\to\widetilde{\cA}(\widetilde{\Sigma})$.
Choose an embedding $e\colon\widetilde{\Sigma}\times I\to\Sigma\times I$ such that the image of $\iota_1(\Sigma_1)$ is isotopic to $\Sigma\times\{0\}$ and the image of $\iota_2(\Sigma_2)$ is contained in a closed ball. By Theorem~\ref{thm:Dehn_twist_general},
if $\sigma_n(x)$ acts trivially, then $\iota_{i\ast}(x)$ is invariant under Dehn twists in $\widetilde{\Sigma}$, and Hence, under the mapping class group. Applying a mapping class that exchanges the two embedded surfaces gives $\iota_{1\ast}(x)=\iota_{2\ast}(x)$. Consequently, $e_\ast\circ\iota_{1\ast}(x)=e_\ast\circ\iota_{2\ast}(x)$.
By the choice of $e$, we have $e_\ast\circ\iota_{1\ast}(x)=x$, whereas $e_\ast\circ\iota_{2\ast}(x)$ lies in the coefficient ring.
\end{proof}

\begin{remark}\label{rmk:tautological}
Lemma~\ref{lem:Dehn_twist_annulus} and Theorem~\ref{thm:Dehn_twist_general} are also proved in \cite{Tsu18a}.
The existence of the element $L_{\widetilde{\cA}}$ in the Dehn twist formula of Theorem~\ref{thm:Dehn_twist_general}, used in the proof of Lemma~\ref{lem:ker_sigma},
follows solely from the existence assertion in Lemma~\ref{lem:Dehn_twist_annulus}.
The proof of the existence assertion in Lemma~\ref{lem:Dehn_twist_annulus} is independent of this section. Thus the proof of Lemma~\ref{lem:ker_sigma} is not circular.
Likewise, although the uniqueness proof in Lemma~\ref{lem:Dehn_twist_annulus} uses Lemma~\ref{lem:ker_sigma}, this does not introduce circularity.
Consequently, the proofs of Lemma~\ref{lem:ker_sigma}, Lemma~\ref{lem:Dehn_twist_annulus}, and Theorem~\ref{thm:Dehn_twist_general} given in this paper are self-contained.
\end{remark}

\begin{theorem}[cf.\cite{Tur91} for (3)]\label{thm:h_times}
Let $\Sigma=\Sigma_{g,b}$ with $g\geq0$ and $b\geq1$. Then the following hold.
\begin{enumerate}
    \item The module $\Gr\widetilde{\cA}(\Sigma)=\bigoplus_n\Gr^n\widetilde{\cA}(\Sigma)$ is free over $\bQ[\rho,h,\hbar]$.
    
    \item The $\bQ[\rho]$-module homomorphisms between the associated graded modules
    \begin{align*}
    \Gr^n\widetilde{\cA}(\Sigma)\to \Gr^{n+2}\widetilde{\cA}(\Sigma),\quad x\mapsto hx\\ 
    \Gr^n\widetilde{\cA}(\Sigma)\to \Gr^{n+2}\widetilde{\cA}(\Sigma),\quad x\mapsto \hbar x
    \end{align*}
    are injective.

    \item The $\bQ[\rho]$-module homomorphism
    \begin{align*}
    \Gr^n\widetilde{\cA}(\Sigma)/\eta\Gr^{n-2}\widetilde{\cA}(\Sigma)\to\Gr^{n+2}\widetilde{\cA}(\Sigma)/\eta\Gr^n\widetilde{\cA}(\Sigma),\quad x\mapsto\eta'x\ (\eta\in\{h,\hbar\})
    \end{align*}
    is injective, where $\{\eta,\eta'\}=\{h,\hbar\}$. We set $\Gr^r\widetilde{\cA}(\Sigma)=0$ for $r<0$.
\end{enumerate}
\end{theorem}
\begin{proof}
Fix an embedding $\iota\colon\Sigma\times[0,1]\to D^2\times[0,1]$. It induces a $\bQ[\rho]$-module homomorphism $\iota_\ast\colon\Gr^n\widetilde{\cA}(\Sigma)\to\Gr^n\widetilde{\cA}(D^2)=\Gr^nR$.
Then, we obtain for $x\in\Gr^n\widetilde{\cA}(\Sigma)$,
an injective $\bQ[\rho]$-module homomorphism
\begin{align}
\Gr^n\widetilde{\cA}(\Sigma)\hookrightarrow\displaystyle\prod_{m\geq0,J}{\rm Hom}(\Gr^m\widetilde{\cA}(\Sigma,J),\Gr^k\widetilde{\cA}(\Sigma,J))\oplus\Gr^n\widetilde{\cA}(\Sigma_{0,1}),\quad
x\mapsto (\sigma_n(x), \iota_\ast(x))\label{eq:GrA_into_Hom}
\end{align}
where $k=n+m-\min\{2,n+1,m+1\}$ records the decrease in the number of boxes.
Indeed, if $(\sigma_n(x),\iota_*(x))=(0,0)$, then Lemma~\ref{lem:ker_sigma} gives $x=\Psi(x)\emptyset$, and $\iota_*(x)=0$ implies $\Psi(x)=0$. This proves injectivity.

By Lemma~\ref{lem:sigma_generators_gr}, the action $\sigma_n(x)$ in \eqref{eq:GrA_into_Hom} is determined by its values on $\cB'_+\sqcup\{\gamma_i\mid i=2,\dots,b\}$. Together with the component in $\Gr^n\widetilde{\cA}(\Sigma_{0,1})$, these values determine the element $x\in\Gr^n\widetilde{\cA}(\Sigma)$ by \eqref{eq:GrA_into_Hom}. By Lemma~\ref{lem:sigma_filtration}, evaluation on the degree-zero and degree-one generators lowers the degree by one. Hence, there is an injective module homomorphism
\begin{align}
\Gr^n\widetilde{\cA}(\Sigma)\hookrightarrow \bigoplus_{i}\Gr^{n-1}\widetilde{\cA}(\Sigma,J_i)\oplus \Gr^{n}\widetilde{\cA}(\Sigma_{0,1}),\label{eq:Grn_n-1_inj}
\end{align}
where $\Gr^{-1}\widetilde{\cA}(\Sigma,J_i)=\{0\}$, and the first direct sum records the values on $\cB'_+\sqcup\{\gamma_i\mid i=2,\dots,b\}$.

By Proposition~\ref{prop:free_module}, the module $\widetilde{\cA}(\Sigma,J_i)$ is free over $\widetilde{\cA}(\Sigma)$, and $\widetilde{\cA}(\Sigma,J_i)\cong\widetilde{\cA}(\Sigma)\otimes\cP_{J_i}$.
In particular, it is free over $\bQ[\rho]$, so
we obtain a $\bQ[\rho]$-module isomorphism
\begin{align}
\displaystyle\bigoplus_{i}\Gr^{n-1}\widetilde{\cA}(\Sigma, J_i)\xrightarrow{\cong}\bigoplus_{i}\bigoplus_{n'+n''=n-1}
\left(\Gr^{n'}\widetilde{\cA}(\Sigma)\otimes_{\bQ[\rho]}\Gr^{n''}\cP_{J_i}\right)/\langle \hbar x\otimes y-x\otimes \hbar y\rangle
\label{eq:Grn-1_isom}
\end{align}
 
Since the fundamental group $\pi_1(\Sigma)$ of $\Sigma=\Sigma_{g,b}$ with $b\geq1$ is free, we apply the Magnus expansion, in each degree in $\hbar$, to the filtration $\{(\ker\varepsilon)^j\}$ defined by the augmentation map $\varepsilon$ of the group algebra $\bQ[\rho]\pi_1(\Sigma)$.
The isomorphism $(\ker\varepsilon)^j/(\ker\varepsilon)^{j+1}\cong H^{\otimes j}$ then gives
\begin{align}
\Gr^k \cP_{J_i}\cong \bigoplus_{j=0}^{\lfloor k/2\rfloor} \hbar^j H^{\otimes k-2j},\label{eq:GrPj_isom}
\end{align}
where $H=H_1(\Sigma;\bQ[\rho])$, and $\lfloor\,\cdot\,\rfloor$ denotes the floor function.

By \eqref{eq:Grn_n-1_inj} and \eqref{eq:Grn-1_isom}, we may regard $\Gr^n\widetilde{\cA}(\Sigma)$ as a submodule of $$\displaystyle\bigoplus_i\bigoplus_{n'+n''=n-1}\left(\Gr^{n'}\widetilde{\cA}(\Sigma)\otimes_{\bQ[\rho]}\Gr^{n''}\cP_{J_i}\right)/\langle\hbar x\otimes y-x\otimes\hbar y\rangle\oplus\Gr^n\widetilde{\cA}(\Sigma_{0,1}).$$
Using \eqref{eq:GrPj_isom}, we further identify it with a submodule of
$$\displaystyle\bigoplus_{i}\bigoplus_{n'+n''=n-1}\left(\Gr^{n'}\widetilde{\cA}(\Sigma)\otimes_{\bQ[\rho]}  H^{\otimes n''}\right)/\langle \hbar x\otimes y-x\otimes \hbar y\rangle\oplus \Gr^{n}\widetilde{\cA}(\Sigma_{0,1}),$$
where the factors $\hbar^j$ are omitted by transferring them, through the tensor product, to the first factor.

Apply the same procedure to $\Gr^{n'}\widetilde{\cA}(\Sigma)$, using \eqref{eq:GrPj_isom}, and repeat until each direct summand of the target is isomorphic to $h^{i_1}\hbar^{i_2}\bQ[\rho]\otimes H^{\otimes n''}$.
This proves (1).

Multiplication by $h$ and by $\hbar$ is injective on the ambient free module. Hence, the maps $\Gr^n\widetilde{\cA}(\Sigma)\to\Gr^{n+2}\widetilde{\cA}(\Sigma)$ given by multiplication by $h$ and by $\hbar$ are injective. This proves (2) and (3).
\end{proof}

\begin{remark}\label{rmk:notorsion_hhbar_completion}
By Theorem~\ref{thm:h_times}, the module $\Gr\widetilde{\cA}(\Sigma)$ has no $h$-torsion or $\hbar$-torsion.
The same holds for the completed skein algebra $\widehat{\widetilde{\cA}(\Sigma)}=\plim[n]\widetilde{\cA}(\Sigma)/F^n\widetilde{\cA}(\Sigma)$. Indeed,
for every nonzero $x\in\widehat{\widetilde{\cA}(\Sigma)}$, there is an $n$ such that its leading class $\overline{x}$ in $\Gr^n\widetilde{\cA}(\Sigma)$ is nonzero. By the injectivity asserted in the theorem,
$h\overline{x}$ and $\hbar\overline{x}$ are nonzero in $\Gr^{n+2}\widetilde{\cA}(\Sigma)$. Hence, $hx$ and $\hbar x$ are nonzero in the completion.
\end{remark}

To prove Theorem~\ref{thm:inj}, it suffices to establish the following theorem.
\begin{theorem}
The map $\Gr^n\Phi$ is injective.
\end{theorem}
\begin{proof}
We first prove that the maps in degrees zero and one are isomorphisms, and then proceed by induction. In degree zero, both modules can be identified with $\bQ[\rho]$, and the map is an isomorphism.
The module $\Gr^1\widetilde{\cN}(\Sigma)$ is spanned by necklace ensembles with exactly one bead. Since $h$ and $\hbar$ do not occur, this degree is unchanged by specialization at $h=\hbar=0$.
Thus $\Gr^1\widetilde{\cN}(\Sigma)$ can be identified with the homology group $H_1(\Sigma;\bQ[\rho])$.
Similarly, $\Gr^1\widetilde{\cA}(\Sigma)$ is spanned by link diagrams with exactly one box, and neither $h$ nor $\hbar$ occurs. Hence, this degree is also unchanged by specialization at $h=\hbar=0$. It can be identified with the degree-one part of the associated graded of the symmetric algebra of the Goldman Lie algebra over $\bQ[\rho]$, for the filtration induced by the augmentation ideal of the group algebra of the fundamental group. A loop without a box equals a trivial loop and can therefore be replaced by $2\rho$.
There is exactly one loop carrying a box, and the degree-one part of the associated graded of the symmetric algebra of the Goldman Lie algebra over $\bQ[\rho]$ is $H_1(\Sigma;\bQ[\rho])$. Thus both $\Gr^1\widetilde{\cN}(\Sigma)$ and $\Gr^1\widetilde{\cA}(\Sigma)$ are identified with $H_1(\Sigma;\bQ[\rho])$. If $\gamma$ is a lift of the homology class labeling a bead, then $\Gr^1\Phi$ sends that bead to $\gamma-\textbf{1}$, namely, the corresponding ds-linear combination. Hence, $\Gr^1\Phi$ is an isomorphism.

Consider the $\bQ[\rho]$-submodule $h\Gr^n\widetilde{\cN}(\Sigma)+\hbar\Gr^n\widetilde{\cN}(\Sigma)\subset\Gr^{n+2}\widetilde{\cN}(\Sigma)$ and the quotient $$\Gr^{n+2}\widetilde{\cN}(\Sigma)|_{h=0}^{\hbar=0}=\Gr^{n+2}\widetilde{\cN}(\Sigma)/(h\Gr^n\widetilde{\cN}(\Sigma)+\hbar\Gr^n\widetilde{\cN}(\Sigma)).$$
Similarly, consider the quotient $$\Gr^{n+2}\widetilde{\cA}(\Sigma)|_{h=0}^{\hbar=0}=\Gr^{n+2}\widetilde{\cA}(\Sigma)/(h\Gr^n\widetilde{\cA}(\Sigma)+\hbar\Gr^n\widetilde{\cA}(\Sigma)).$$
These are the specializations at $h=\hbar=0$ of $\Gr^{n+2}\widetilde{\cN}(\Sigma)$ and $\Gr^{n+2}\widetilde{\cA}(\Sigma)$, respectively. As in degrees zero and one,
the classical isomorphism between the associated graded of the Goldman Lie algebra and the necklace Lie algebra
shows that these two quotients are isomorphic. Denote the isomorphism induced by $\Gr^{n+2}\Phi$ by $\Gr^{n+2}\Psi$.

Assume that $\Gr^n\Phi$ is injective.
Let $x\in\Gr^{n+2}\widetilde{\cN}(\Sigma)$ satisfy $\Gr^{n+2}\Phi(x)=0$. Its class $[x]\in\Gr^{n+2}\widetilde{\cN}(\Sigma)|_{h=0}^{\hbar=0}$ then satisfies $\Gr^{n+2}\Psi([x])=0$. Since $\Gr^{n+2}\Psi$ is injective, we have $[x]=0$.
Thus $x\in h\Gr^n\widetilde{\cN}(\Sigma)+\hbar\Gr^n\widetilde{\cN}(\Sigma)$, and then $x=hy+\hbar z$ for some $y,z\in\Gr^n\widetilde{\cN}(\Sigma)$. Since $\Gr^{n+2}\Phi(x)=0$, we have
$$0=\Gr^{n+2}\Phi(x)=h\Gr^{n}\Phi(y)+\hbar\Gr^{n}\Phi(z)\Longleftrightarrow h\Gr^{n}\Phi(y)=-\hbar\Gr^{n}\Phi(z)$$

This equality gives
\begin{align}
\hbar \Gr^n (\Phi(z))=\Gr^{n+2} \Phi( \hbar z)=0\mod h \Gr \widetilde{\cA}(\Sigma)\label{eq:injectivity_h_hbar}
\end{align}
By Theorem~\ref{thm:h_times}(3), we obtain
$\Gr^n\Phi(z)=0\mod h\Gr\widetilde{\cA}(\Sigma)$, or equivalently,
$\Gr^n\Phi(z)\in h\Gr^{n-2}\widetilde{\cA}(\Sigma)$.
Using also the surjectivity of $\Gr^{n-2}\Phi$, choose $w_z\in\Gr^{n-2}\widetilde{\cN}(\Sigma)$ such that $\Gr^n\Phi(z)=h\Gr^{n-2}\Phi(w_z)$. Then, the injectivity of $\Gr^n\Phi$ gives $z=hw_z$. 
%For $n=0,1$, the negative-degree modules are regarded to be zero.
Together with \eqref{eq:injectivity_h_hbar}, this yields
$h \Gr^{n}\Phi (y)=-\hbar h \Gr^{n-2}\Phi(w_z)$
 
By Theorem~\ref{thm:h_times}(2), we further obtain $\Gr^n\Phi(y)=-\hbar\Gr^{n-2}\Phi(w_z)$. The induction hypothesis then gives $y=-\hbar w_z$.
Consequently, $x=0$. This proves that $\Gr^{n+2}\Phi$ is injective.
\end{proof}

\begin{remark}
The injectivity of multiplication by $h$ and $\hbar$ on $\Gr\widetilde{\cA}(\Sigma)$ (Theorem~\ref{thm:h_times}), together with the injectivity of the map from the quantized necklace algebra to $\Gr\widetilde{\cA}(\Sigma)$ (Theorem~\ref{thm:inj}), implies that multiplication by $h$ and $\hbar$ is also injective on the quantized necklace algebra. This gives another proof of the injectivity implicit in \cite{Sch05}.
\end{remark}

Proposition~\ref{prop:free_module} and Theorem~\ref{thm:isom} give the following corollary. Here, virtual closure is taken as in Proposition~\ref{prop:free_module}.
\begin{cor}
If $J\neq\emptyset$,
then $\Gr\widetilde{\cA}(\Sigma,J)$ is isomorphic to
$$\bigoplus_{J'}\left(\Gr\widetilde{\cA}(\Sigma)\otimes  (TH_1)^{\otimes |J'|/2}\right)\cong 
\bigoplus_{J'}\left(\widetilde{\cN}(\Sigma)\otimes  (TH_1)^{\otimes |J'|/2}\right)$$
Here the direct sum ranges over $J'=\tau(J_{\rm out})\sqcup J_{\rm in}$ obtained from $J=J_{\rm out}\sqcup J_{\rm in}$ by permutations $\tau\in\mathfrak{S}_{|J|/2}$.
\end{cor}

%%%%%%%%%%%%%%%%%%%%%%%%%%%%%%%%%%%%%%%%%%%%%%%%%%%%%%%%%%%%%%%%%%%%%%%%%%%%%%%%
\section{Coalgebra structure of the HOMFLYPT skein algebra}\label{sec:coalgebra}
Throughout this section, we assume that $\rho=0$ and retain the notation $\widetilde{\cA}(\Sigma)$ for the HOMFLYPT skein algebra. The finite-type filtration is defined in the same way under this assumption, and we consider the associated graded algebra $\Gr(\widetilde{\cA}(\Sigma))$.

We show that Turaev's coalgebra structure on $\widetilde{\cA}(\Sigma)$ induces a coalgebra structure on $\Gr\widetilde{\cA}(\Sigma)$.
We also give an explicit description of the induced coproduct.
As explained in Remark~\ref{rmk:bialg}, this description identifies the induced coproduct, for a bouquet quiver, with the coproduct in Schedler's quantized necklace Lie bialgebra \cite{Sch05}, and also yields the biquantization of the Lie bialgebra structure in \cite{CHLZ26}.

\begin{figure}	\centering\includegraphics[width=60mm]{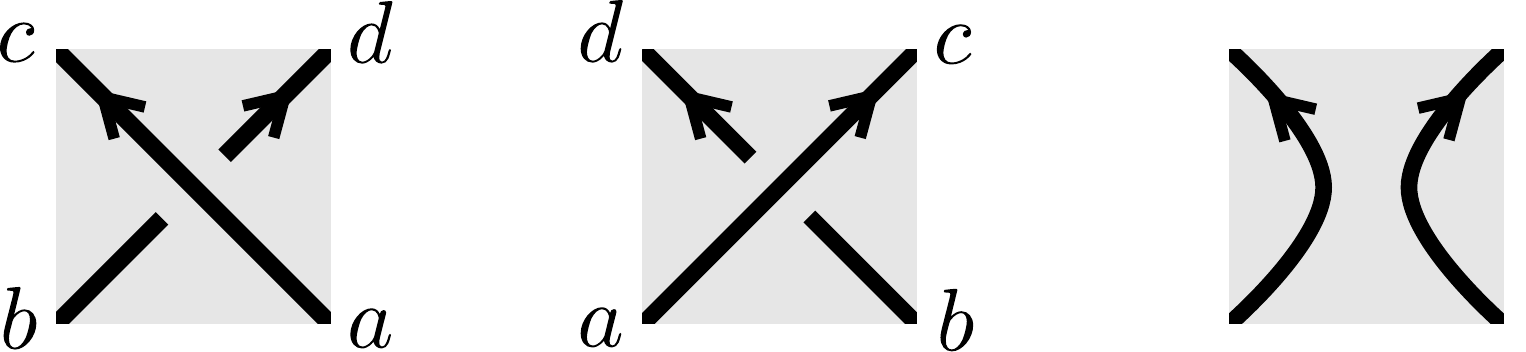}
    \caption{Left and middle: the edges adjacent to a crossing. Right: reconnecting the edges at an $f$-cutting vertex.}
	\label{fig:coloring_Turaev}  
\end{figure}

We first recall Turaev's coalgebra structure on the HOMFLYPT skein algebra \cite[Section 9]{Tur91}.
Let $D$ be a link diagram on $\Sigma$.
We call each connected component of $D\setminus{\rm Cross}(D)$ an \textbf{edge} of $D$. An edge is homeomorphic to $(0,1)$, except when it is an entire component of $D$ containing no crossings, in which case it is homeomorphic to $S^1$. 
Let $E(D)$ denote the set of edges of $D$.
A map $f\colon E(D)\to\{1,\dots,n\}$ is called an \textbf{$n$-coloring} if, at every crossing of $D$, the adjacent edges $a,b,c,d\in E(D)$ shown in Figure~\ref{fig:coloring_Turaev} satisfy either $f(a)=f(c)$ and $f(b)=f(d)$, or $f(a)=f(d)>f(b)=f(c)$.
A crossing satisfying $f(a)=f(d)>f(b)=f(c)$ is called an \textbf{$f$-cutting vertex}.
Let $|f|$ denote the number of $f$-cutting vertices of $D$, and let $|f|_+$ and $|f|_-$ denote the numbers of positive and negative $f$-cutting vertices, respectively. Thus $|f|=|f|_++|f|_-$.
Let ${\rm lb}_n(D)$ denote the set of $n$-colorings of $D$. For $f\in{\rm lb}_n(D)$ and $i=1,\dots,n$, let $D_{f,i}$ be the link diagram obtained from $f^{-1}(i)$ by reconnecting the edges near the $f$-cutting vertices.
Set $\|f\|:=\#D-\sum_{i=1}^n\#D_{f,i}$. This records the change in the number of components when the $f$-cutting vertices of $D$ are replaced as in the right-hand diagram of Figure~\ref{fig:coloring_Turaev}. Define
\begin{align}
\langle D\mid f\rangle:=(-1)^{|f|_-}h^{(|f|+\|f\|)/2}\hbar^{(|f|-\|f\|)/2}\label{eq:coproduct_coeff}
\end{align}
 
In \cite[Section 9.1]{Tur91}, the parentheses in the exponents of $h$ and $\hbar$ are missing; this is a typographical error.

We then set $$\Delta(D,f):=\langle D\mid f\rangle D_{f,1}\otimes\dots\otimes D_{f,n}.$$
 
Writing ${\rm lb}(D):={\rm lb}_2(D)$, define
\begin{align}
\Delta(D)=\sum_{f\in {\rm lb}(D)} \Delta(D,f)\in \widetilde{\cA}(\Sigma)\otimes \widetilde{\cA}(\Sigma)\label{eq:Turaev_cobracket_quantum}
\end{align}
The following theorem holds.
\begin{theorem}[\cite{Tur91}]\label{thm:coalg_Turaev}
Suppose that $\rho=0$.
Then \eqref{eq:Turaev_cobracket_quantum} extends uniquely to a $\bQ[h,\hbar]$-module homomorphism
\begin{align}
    \Delta\colon \widetilde{\cA}(\Sigma)\to \widetilde{\cA}(\Sigma)\otimes \widetilde{\cA}(\Sigma)
\end{align}
In particular, $\widetilde{\cA}(\Sigma)$ is a bialgebra with coproduct $\Delta$.
\end{theorem}
The following theorem shows that this bialgebra structure induces a bialgebra structure on $\Gr\widetilde{\cA}(\Sigma)$.
Via the isomorphism of Theorem~\ref{thm:isom}, it also defines a bialgebra structure on $\widetilde{\cN}(\Sigma)$.

\begin{theorem}\label{thm:coproduct_filtration}
Suppose that $\rho=0$. The coproduct in Theorem~\ref{thm:coalg_Turaev} preserves the filtration $\{F^n\widetilde{\cA}(\Sigma)\}_{n\in\bN}$. Namely, for $x\in F^n\widetilde{\cA}(\Sigma)$, we have $\Delta(x)\in\sum_{i+j=n}F^i\widetilde{\cA}(\Sigma)\otimes F^j\widetilde{\cA}(\Sigma)$.
\end{theorem}
\begin{proof}
By Theorem~\ref{thm:coalg_Turaev}, the coproduct $\Delta$ is compatible with multiplication. Arguments analogous to those in the proofs of Lemmas~\ref{lem:transparency} and~\ref{lem:surj} show that the module is generated by links with boxes whose distinct connected components are separated in height. Hence, the algebra is generated by connected links with boxes, and it suffices to consider these generators.

Fix a diagram $D$ of a connected link with boxes on $\Sigma$. If $D$ has no boxes, the assertion is immediate.

Suppose that $D$ has at least one box. Replace each box by its ds-linear combination to obtain $\sum_i c_iD_i$ with $c_i\in\{-1,1\}$. Here, each $D_i$ is a fixed link diagram without boxes, rather than an isotopy class, and agrees with $D$ outside the boxes.
The coproduct of Theorem~\ref{thm:coalg_Turaev} sends $\sum_i c_iD_i$ to $\sum_{f_i,i}c_i\Delta(D_i,f_i)$, where $f_i$ ranges over the $2$-colorings of $D_i$. Each summand is determined by a pair $(D_i,f_i)$. Since $D_i$ is connected, a $2$-coloring $f_i$ with cut vertices is uniquely determined by the colors near those vertices. If there are no cut vertices, there are exactly two colorings, assigning either $1$ or $2$ to every edge.

We rewrite $\sum_{f_i,i}c_i\Delta(D_i,f_i)$ in terms of link diagrams with boxes.
Let $S_v$ be a subset of the crossings of $D$, where crossings carry over/under information and are not boxes. Let $\mathcal{D}_v\subset\{D_i\}_i$ consist of the diagrams containing every $d$-component of $D$ that contains a crossing in $S_v$. Each $D_i\in\mathcal{D}_v$ contains crossings corresponding to those in $S_v$; their set is denoted by $S_v$ as well.
The pairs consisting of a diagram and a $2$-coloring whose set of cut vertices is $S_v$ are precisely
\begin{align}
\begin{cases}
   \{(D_i,f_i)\mid D_i\in\mathcal{D}_v,\ \text{$f_i$ is a $2$-coloring of $D_i$ with cut-vertex set $S_v$}\} &(S_v\neq\emptyset),\label{eq:coloring}\\
   \{(D_i,f_i^{(j)})\mid D_i\in\mathcal{D}_v,\ j=1,2,\ \text{$f_i^{(j)}$ is the constant $j$-coloring of $D_i$}\}&(S_v=\emptyset).
\end{cases}
\end{align}
Taking the disjoint union of \eqref{eq:coloring} over all $S_v$
gives a bijection with $\{(D_i,f_i)\}_{i,f_i}$.
For each $S_v$, the coefficients $\langle D_i\mid f_i\rangle$ are the same for all pairs $(D_i,f_i)$ in \eqref{eq:coloring}, independently of $i$. We write this common value as $\langle S_v\rangle:=\langle D_i\mid f_i\rangle$.
In the corresponding sum $\sum_{D_i\in\mathcal{D}_v}c_i\langle S_v\rangle D_i$, each $d$-component containing no crossing of $S_v$ occurs in a pair of terms forming a ds-linear combination, which can be grouped into a box.
By \eqref{eq:coproduct_coeff}, the finite-type degree of $\langle S_v\rangle$ is $2|S_v|$. Moreover, $(\text{number of $d$-components of $D$ meeting $S_v$})\leq2|S_v|$. Hence, the number of boxes in $\sum_i c_i\langle S_v\rangle D_i$ decreases by at most $2|S_v|$ relative to the original diagram $D$.
This proves the assertion.
\end{proof}

By Theorem~\ref{thm:coproduct_filtration}, the coproduct on $\Gr\widetilde{\cA}(\Sigma)$ induced by Theorem~\ref{thm:coalg_Turaev} is well-defined.

As in the proof of Theorem~\ref{thm:coproduct_filtration}, grouping ds-linear combinations of $2$-colored link diagrams into boxes gives $2$-colorings of link diagrams with boxes. More precisely, we require a \textbf{$2$-coloring} to assign the same color to the incoming (resp. outgoing) end of the detour and the outgoing (resp. incoming) end of the straight part at each box.
We say that $f$ is \textbf{degree-preserving} if each detour component contains at most one $f$-cutting vertex, and no $f$-cutting vertex is a self-crossing of a detour or a crossing involving a straight component.

\begin{figure}	\centering\includegraphics[width=70mm]{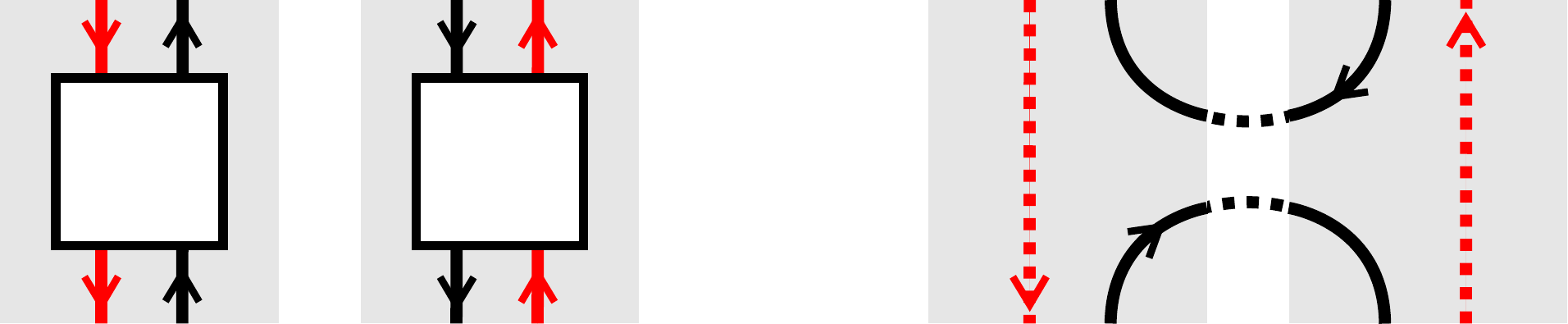}
\caption{Left: a local picture of two boxes.
Right: the link diagram with boxes obtained by deleting the detours and boxes in the left-hand diagram and reconnecting the strands. The deleted detours are shown by dotted lines.}
\label{fig:coproduct_connecting} 
\end{figure}

Let $b_1,\dots,b_k$ be the boxes of a link diagram $D$ with boxes.
For a degree-preserving $n$-coloring $f$ of $D$,
let $D_f$ be the link diagram with boxes obtained by deleting the detours $d_i,d_j$ meeting at each $f$-cutting vertex, together with the boxes $b_i,b_j$, and reconnecting as in Figure~\ref{fig:coproduct_connecting}. The resulting diagram $D_f$ consists of link diagrams with boxes colored by $1,\dots,n$. Denote the subdiagram colored by $i$ by $D_f^{(i)}$.
Let ${\rm Cut}(f)$ be the set of pairs $(i,j)$ with $i<j$ such that $d_i$ and $d_j$ meet at an $f$-cutting vertex.

For $2$-colorings of $D$, define $f\sim f'$ if they have the same pairs of detours meeting at cut vertices and assign the same labels to all $s$-components. This is an equivalence relation. Let $[f]$ denote the class of $f$.

The coproduct on $\Gr\widetilde{\cA}(\Sigma)$ is described by the degree-preserving $2$-colorings of link diagrams with boxes defined above.
By Theorem~\ref{thm:isom}, this gives a coproduct on the quantized necklace algebra. In homological terms, it can be expressed using the intersection form as follows.
\begin{theorem}\label{thm:coproduct_gr}
Suppose that $\rho=0$. Let $b_1,\dots,b_k$ be the boxes of a link diagram $D$ with boxes, and let $d_i$ be the detour of $b_i$. Assume that, whenever $i<j$, the detour $d_i$ passes over $d_j$ at every crossing between them.
Then the coproduct of $D$ in $\Gr\widetilde{\cA}(\Sigma)$ induced by Theorem~\ref{thm:coalg_Turaev} is given by
\begin{align}
    \Delta(D) = \sum_{[f]} \left(\prod_{(i,j)\in {\rm Cut}(f)}\mu(d_i,d_j)h^{\varepsilon_{ij}}\hbar^{1-\varepsilon_{ij}}\right) D_{f}^{(1)}\otimes D_{f}^{(2)},\label{eq:coproduct_necklace}
\end{align}
where $f$ ranges over degree-preserving $2$-colorings of $D$, $\mu$ is the intersection form on $H_1(\Sigma;\bZ)$, and $\varepsilon_{ij}\in\{0,1\}$ is chosen to preserve the component degree.
If there are no $f$-cutting vertices, so that ${\rm Cut}(f)=\emptyset$, the product in parentheses in \eqref{eq:coproduct_necklace} is understood to be $1$.
\end{theorem}

\begin{proof}For a link diagram $D$ with boxes,
consider the linear combination $\sum_{f'_i,i}c'_i\Delta(D'_i,f'_i)$ from the proof of Theorem~\ref{thm:coproduct_filtration}. The argument in the proof of Theorem~\ref{thm:coproduct_filtration} shows that the finite-type degree increases by at least one in each of the following cases:
\begin{itemize}
    \item A self-crossing of a detour portion is a cut vertex.
    \item A detour portion contains at least two cut vertices.
    \item A cut vertex involves a portion that is not a detour.
\end{itemize}
Thus, if a $2$-coloring $f$ satisfies any of these conditions, its contribution lies in strictly higher filtration degree and vanishes in the graded component of $\Delta(D)$ having the same total degree as $D$.
It therefore suffices to consider degree-preserving $2$-colorings.
If $f\sim f'$, then $D_f^{(1)}\otimes D_f^{(2)}=D_{f'}^{(1)}\otimes D_{f'}^{(2)}$ in $\Gr\widetilde{\cA}(\Sigma)\otimes\Gr\widetilde{\cA}(\Sigma)$.
For distinct detours $d_i,d_j$ with $i<j$, assume that $d_i$ passes over $d_j$ at every crossing between them. Let $\langle D'_l\mid f'_l\rangle$ denote the extension of \eqref{eq:coproduct_coeff} to link diagrams with boxes. The contribution of $d_i,d_j$ to $\langle D'_l\mid f'_l\rangle$ is the sign of their $f'_l$-cut vertex (Figure~\ref{fig:crossings_signs}) multiplied by $h$ or $\hbar$, with the parameter chosen so that the component degree is preserved.
Summing over the cut vertices between $d_i$ and $d_j$ therefore gives $\mu(d_i,d_j)h^{\varepsilon_{ij}}\hbar^{1-\varepsilon_{ij}}$.
This is the sum of the contributions over the representatives of the equivalence class $[f]$.
Taking the contribution from each pair of detours in ${\rm Cut}(f)$ gives \eqref{eq:coproduct_necklace}.
\end{proof}

\begin{remark}\label{rmk:bialg}
We explain how the coproduct in \cite{Sch05} is recovered from Theorem~\ref{thm:coproduct_gr}.
Set $\rho=0$ and $h=1$, and let $\Sigma=\Sigma_{g,1}$.
The condition $\rho=0$ is required in \cite{Tur91}, whereas the specialization $h=1$ is used in \cite{Sch05}.
Choose a symplectic basis $\{\alpha_i,\beta_i\mid i=1,\dots,g\}$ of $H_1(\Sigma)$,
and assume that the homology class of each detour in a link diagram with boxes is a member of this basis.
Through Theorem~\ref{thm:isom}, the coproduct of Theorem~\ref{thm:coproduct_gr} defines a coproduct on $\widetilde{\cN}(\Sigma)$.
The pairs $(i,j)$ in \eqref{eq:coproduct_necklace}, determined by a coloring of the link diagram with boxes, correspond to the pairs of beads in the coproduct of \cite{Sch05}. Thus this coproduct agrees with the coproduct of \cite{Sch05} for a bouquet quiver.
In particular, the resulting bialgebra structure on $\widetilde{\cN}(\Sigma)$ agrees with that of \cite{Sch05}. The biquantization in \cite{CHLZ26} refines \cite{Sch05} by retaining $h$ instead of specializing it to $1$. Thus Theorem~\ref{thm:coproduct_gr} recovers the biquantization of the Lie bialgebra structure for bouquet quivers in \cite{CHLZ26} from the associated graded of Turaev's skein biquantization \cite{Tur91}.
\end{remark}

%%%%%%%%%%%%%%%%%%%%%%%%%%%%%%%%%%%%%%%%%%%%%%%%%%%%%%%%%%%%%%%%%%%%%%%%%%%%%%%%
\section{Formality of HOMFLYPT skein algebras in genus zero}\label{sec:formality_genus0}
In this section, we prove the formality of the HOMFLYPT skein algebra of an oriented surface of genus zero.
This allows us to express the Dehn twist formula in the completed HOMFLYPT skein algebra in terms of the quantized necklace algebra.

Throughout this section, for a marked surface $(\Sigma,J)$, we assume $|J|=0$ or $2$ unless otherwise stated.

\subsection{Formality and homomorphic expansions}\label{subsec:formality_expansion}

Let $V$ be a filtered $\bQ$-algebra, and consider 
its the completed associated graded quotient
$$
\Gr V := \prod_{n=0}^\infty F_nV/F_{n+1}V. 
$$
We also consider the inverse limit
$$
\widehat{V} = \plim[n] V/F_nV
$$
The submodules $F_n\widehat{V}:=\Ker(\widehat{V}\to V/F_n V)$ define a filtration $\{F_n\widehat{V}\}_{n\in \bN}$ on $\widehat{V}$.

An isomorphism $\theta\colon \widehat{V}\to\Gr V$ of filtered modules induces a map on the associated graded quotient:
$$
\Gr  \theta: \Gr  V = \Gr  \widehat{V} \to \Gr (\Gr  V) = \Gr  V.
$$
We call $\theta$ an \textbf{expansion} of $V$ if $\Gr\theta=\mathrm{id}$.
If an expansion is a $\bQ$-linear map preserving an operation $\star$, we call it a \textbf{homomorphic expansion} with respect to $\star$. We say that the \textbf{formality} of $V$ with respect to $\star$ holds if such an expansion exists.

The isomorphism of Theorem~\ref{thm:isom} describes the associated graded quotient of the HOMFLYPT skein algebra in terms of the quantized necklace algebra.
However, as noted in Remark~\ref{rmk:HOMFLYPT}(3), when $\Sigma=\Sigma_{0,p+1}$, relations (A') and (B') in Definition~\ref{def:necklace_Lie_algebra} merely allow bead heights to be interchanged. Since the product on $\widetilde{\cA}(\Sigma)$ is generally noncommutative, except for certain surfaces, this definition is not suitable for studying the formality of $\widetilde{\cA}(\Sigma)$. To address this issue, we first modify the finite-type filtration of the HOMFLYPT skein algebra.

For a marked surface $(\Sigma,J)$ with $\Sigma=\Sigma_{0,p+1}$,
we modify the finite-type filtration $\{F^n\widetilde{\cA}(\Sigma,J)\}_{n\in\bN}$ of $\widetilde{\cA}(\Sigma,J)$ from Definition~\ref{def:finite-type_filtration}
to obtain a descending filtration $\{I^n\widetilde{\cA}(\Sigma,J)\}_{n\in\bN}$. We give its precise definition in Section~\ref{subsec:modified_filtration} and show in Proposition~\ref{prop:filtrations_relations} that it induces the same topology as $\{F^n\widetilde{\cA}(\Sigma,J)\}_{n\in\bN}$.
For a marked oriented surface $(\Sigma,J)$, let $\widehat{\widetilde{\cA}(\Sigma,J)}$ denote the inverse limit $\plim[n]\widetilde{\cA}(\Sigma,J)/I^n\widetilde{\cA}(\Sigma,J)$.

The main result of this subsection is the following.
\begin{theorem}\label{thm:formality}
Let $(\Sigma,J)$ be a marked surface with $\Sigma:=\Sigma_{0,p+1}$, where $p\geq0$.
There exists a homomorphic expansion $\widehat{\widetilde{\cA}(\Sigma,J)}\to\Gr_I\widehat{\widetilde{\cA}(\Sigma,J)}$ preserving the left and right actions of $\widehat{\widetilde{\cA}(\Sigma)}$ and composition. In other words, $\widehat{\widetilde{\cA}(\Sigma,J)}$ has formality with respect to these operations.
In particular, when $J=\emptyset$, formality also holds with respect to the product and the height involution on $\widetilde{\cA}(\Sigma)$.
\end{theorem}

By Theorem~\ref{thm:formality}, the multiplication by $h$ on $\widehat{\widetilde{\cA}(\Sigma,J)}$ is injective. Thus, we can define a Lie action $\sigma\colon\widehat{\widetilde{\cA}(\Sigma)}\to\End\widehat{\widetilde{\cA}(\Sigma,J)}$ by $\sigma(x)(y):=\dfrac{xy-yx}{h}$.
In particular, when $J=\emptyset$, we obtain a Lie bracket $[x,y]:=\sigma(x)(y)$.
Therefore, we obtain the following corollary.
\begin{cor}
Let $(\Sigma,J)$ be a marked surface with $\Sigma:=\Sigma_{0,p+1}$, where $p\geq0$.
Then, $\widehat{\widetilde{\cA}(\Sigma,J)}$ has formality with respect to the Lie action of $\widehat{\widetilde{\cA}(\Sigma)}$.
In particular, when $J=\emptyset$, formality holds with respect to the Lie bracket on $\widetilde{\cA}(\Sigma)$.
\end{cor}

\begin{remark}
\begin{enumerate}
    \item When $\rho=0$, one can define a coproduct on $\widetilde{\cA}(\Sigma)$ in the sense of \cite{Tur91}, but we do not establish formality for this coproduct. Nevertheless, as in \cite{BN+25}, formality with respect to the height involution implies formality with respect to the Turaev cobracket on the Goldman--Turaev Lie bialgebra $\bQ|\pi_1(\Sigma)|$.
    \item In \cite{BN+25}, the Kontsevich integral is used to prove the formality of the Goldman--Turaev Lie bialgebra. In this paper, we use a HOMFLYPT quotient obtained by imposing additional relations on the Conway quotient appearing in \cite{BN+25}, and our results include the corresponding version of the results of \cite{BN+25}. Although \cite{BN+25} works over $\bC$, the existence of associators over $\bQ$ allows us to work over $\bQ$.
    Although we give the proofs only over $\bQ$, the same arguments apply after extending scalars to $\bC$, and the corresponding results hold for any associator.
\end{enumerate}
\end{remark}

The following follows directly from the definition of formality.
\begin{prop}[{\cite[Proposition 2.7]{BND17}}]
Let $V$ be a complete filtered algebra, and let $\Gr V$ be its associated graded quotient. Suppose that $\mathcal{C}$ is a graded algebra equipped with a surjective graded map $\psi\colon\mathcal{C}\to\Gr V$. If there exists a filtered map $Z\colon V\to\mathcal{C}$ such that $\Gr Z\circ\psi=\id_{\mathcal{C}}$, then $\psi$ is an isomorphism and $\psi\circ Z$ is an expansion. If these maps also preserve an operation $\star$, then $\psi\circ Z$ is a homomorphic expansion of $V$ with respect to $\star$.
$$\begin{tikzcd}
V
\arrow[r, "Z"] \ar[rd, dashed, "\text{expansion}"'] 
&\cC\ar[d, two heads, "\psi"]\arrow[r, "{\rm id}_\cC"] &\cC\\
&\Gr\, V\ar[ru,  "{\Gr\,Z}"'] &\end{tikzcd}
$$
\end{prop}
To prove Theorem~\ref{thm:formality}, we use the Kontsevich integral as the algebra homomorphism $Z$ as in \cite{BN+25}. We then construct $\psi$ and prove that it is an isomorphism by showing that its inverse is $\Gr\,Z$. The key idea of this section is to use the quantized necklace structure of $\cC$ to prove that $\psi$ is an isomorphism.
More precisely, we show that the following commutative diagram satisfies property $(\ast)$ below.
\begin{equation}
\begin{tikzcd}
\widehat{\widetilde{\cA}(\Sigma,J)} 
\arrow[r, "Z_{\widetilde{\cA}}^J"] \ar[rd, dashed, "\text{expansion}"'] 
&\mathcal{D}_{\widetilde{\cA}}^J \ar[d, "\psi_{\widetilde{\cA}}^J","\cong"']\ar[rr, "\cong", "\psi_{\mathcal{D}}^J"'] &&\mathcal{D}_{\widetilde{\cN_I}}^J \ar[d, "\psi_{\widetilde{\cN_I}}^J", "\cong"']\\
&\widehat{\Gr}_I\widetilde{\cA}(\Sigma,J) %\ar[u, "\Gr_I Z_{\widetilde{\cA}}^J"]
&\widehat{\widetilde{\cN_I}(\Sigma,J)} \ar[l, "\Gr_I\Phi", "\cong"'] \ar[r, "\cong"] 
&\widehat{\widetilde{\cN_a}(\Sigma,J)} 
\end{tikzcd}\label{eq:commutative_diag_expansion}
\end{equation}
Here, $Z_{\widetilde{\cA}}^J$, $\psi_{\mathcal{D}}^J$, and $\psi_{\widetilde{\cN_I}}^J$ are $\bQ\llbracket a\abar,\rho_a\abar\rrbracket[a,\abar,\rho_a]$-module isomorphisms, and $\Gr\Phi$ is a $\bQ\llbracket h\hbar,\rho\hbar\rrbracket[h,\hbar,\rho]$-module isomorphism. The remaining maps are semilinear isomorphisms with respect to the coefficient-ring isomorphism $\bQ\llbracket h\hbar,\rho\hbar\rrbracket[h,\hbar,\rho]\to\bQ\llbracket a\abar,\rho_a\abar\rrbracket[a,\abar,\rho_a]$ induced by \eqref{eq:h_a} and \eqref{eq:chord_diag_rel}. When $J=\emptyset$, they are semilinear algebra isomorphisms.\footnote{A semilinear algebra isomorphism means a semilinear isomorphism preserving the algebraic structure.} Note that the grading also changes under these semilinear isomorphisms.
\paragraph{\textbf{Property $(\ast)$}}\label{aaa}
\begin{enumerate}
    \item If $|J|=0$, each module has an algebra structure, and the maps in the commutative diagram are $\bQ\llbracket a\abar,\rho_a\abar\rrbracket[a,\abar,\rho_a]$-module isomorphisms, $\bQ\llbracket\rho\hbar,h\hbar\rrbracket[\rho,h,\hbar]$-algebra isomorphisms or semilinear algebra isomorphisms.
    
    \item Let $J_1$ and $J_2$ be markings with $|J_1|=|J_2|=2$. If the incoming point of $J_1$ and the outgoing point of $J_2$ lie on the same boundary component, a composition is defined on the corresponding modules, and the maps in the commutative diagram preserve composition.
    
    \item The HOMFLYPT skein module $\widetilde{\cA}(\Sigma,J)$ is a $\widetilde{\cA}(\Sigma)$-bimodule.
    The other modules in the commutative diagram also carry bimodule structures, and the maps preserve the left and right actions.
\end{enumerate}

\begin{figure}	\centering\includegraphics[scale=0.15]{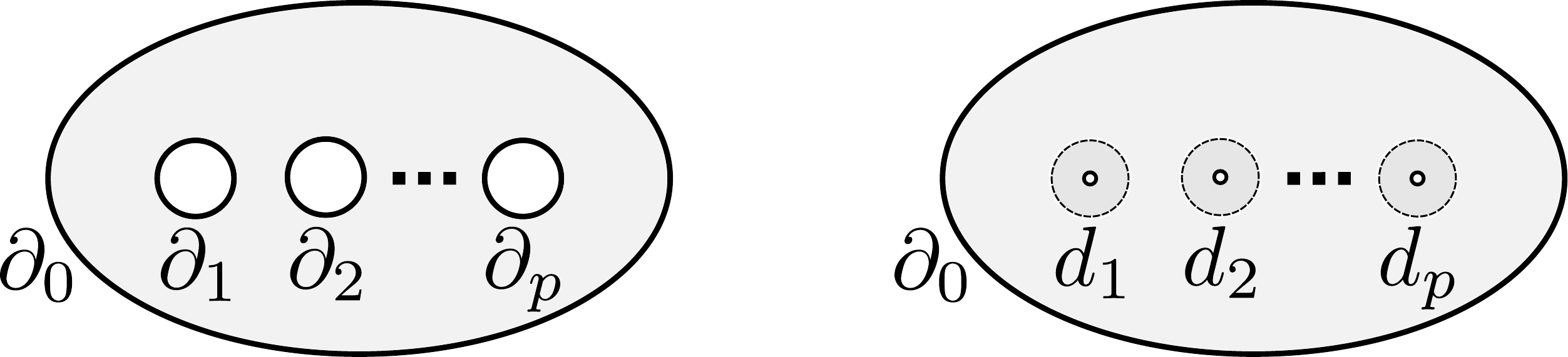}
    \caption{The surfaces $\Sigma_{0,p+1}$ (left) and $D_p$ (right).}
	\label{fig:disk_punctured}  
\end{figure}

Let $\partial_i (i=0,1,\dots,p)$ be the boundary components of $\Sigma_{0,p+1}$.
We depict this surface as a disk $D$ with boundary $\partial_0$ from which $p$ open disks have been removed, as in Figure~\ref{fig:disk_punctured}.
Let $D_p$ be the surface obtained by attaching a once-punctured disk $d_i$ along each boundary component $\partial_i$ for $i=1,\dots,p$; see Figure~\ref{fig:disk_punctured}.\footnote{In \cite{BN+25}, the notation $D_p$ denotes $\Sigma_{0,p+1}$.}
There is a natural embedding $e\colon\Sigma_{0,p+1}\hookrightarrow D_p$.
We call the puncture arising from the once-punctured disk $d_i$ attached along $\partial_i$ the \textbf{$i$-th puncture} of $D_p$.
We call the product of a puncture with $[0,1]$, namely each of the vertical lines removed from $D\times[0,1]$, a \textbf{pole}.

Consider $J$-tangles in $\Sigma_{0,p+1}\times[0,1]$, where $J\subset\partial\Sigma_{0,p+1}$.\footnote{In \cite{BN+25}, $J$-tangles have $J\subset\Sigma_{0,p+1}\times\{0,1\}$ and are simply called tangles. The two conventions are in bijection once a way of moving the endpoints is fixed.}
$$\{\text{$J$-tangles in $\Sigma_{0,p+1}\times[0,1]$}\}/\text{isotopy}$$
Let $\bQ\widetilde{\cT}^J$ denote the $\bQ$-module freely generated by this set.

\begin{figure}	\centering\includegraphics[scale=0.15]{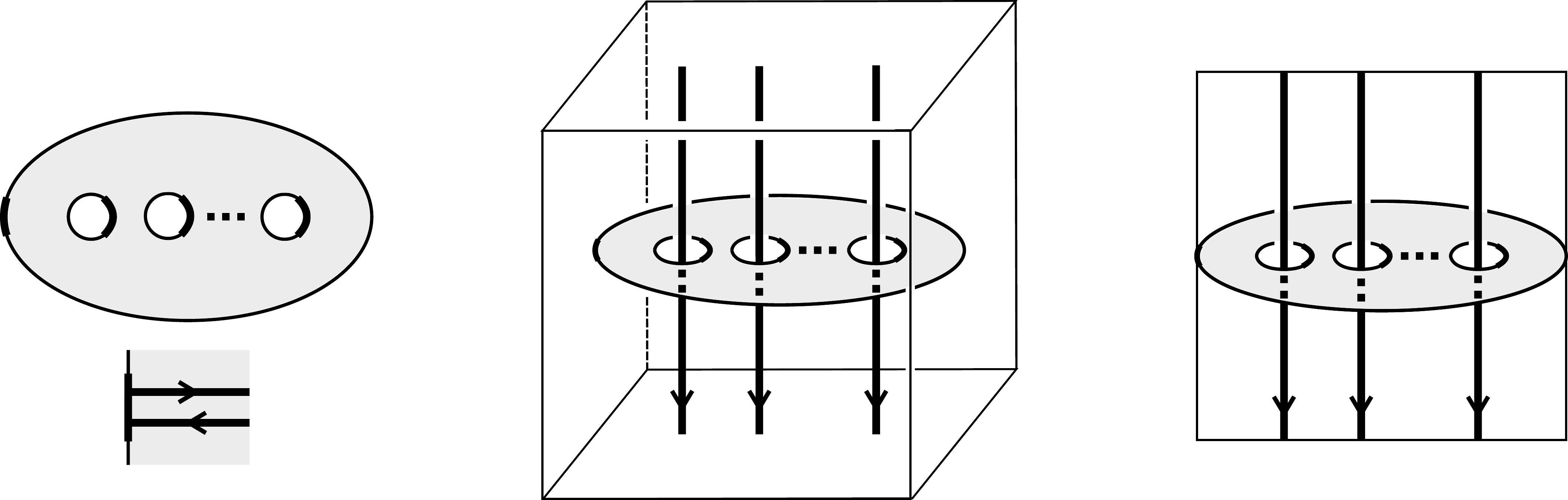}
    \caption{Left: $(\Sigma_{0,p+1},J)$ with $J$ contained in the thickened segments and the outgoing and incoming endpoints positioned as in the lower-left inset. Center: $\Sigma_{0,p+1}\times[0,1]$. Right: the horizontal projection and $\Sigma_{0,p+1}\times\{1/2\}$.}
	\label{fig:diagram_vertical_horizontal}  
\end{figure}

\begin{figure}	\centering\includegraphics[scale=0.15]{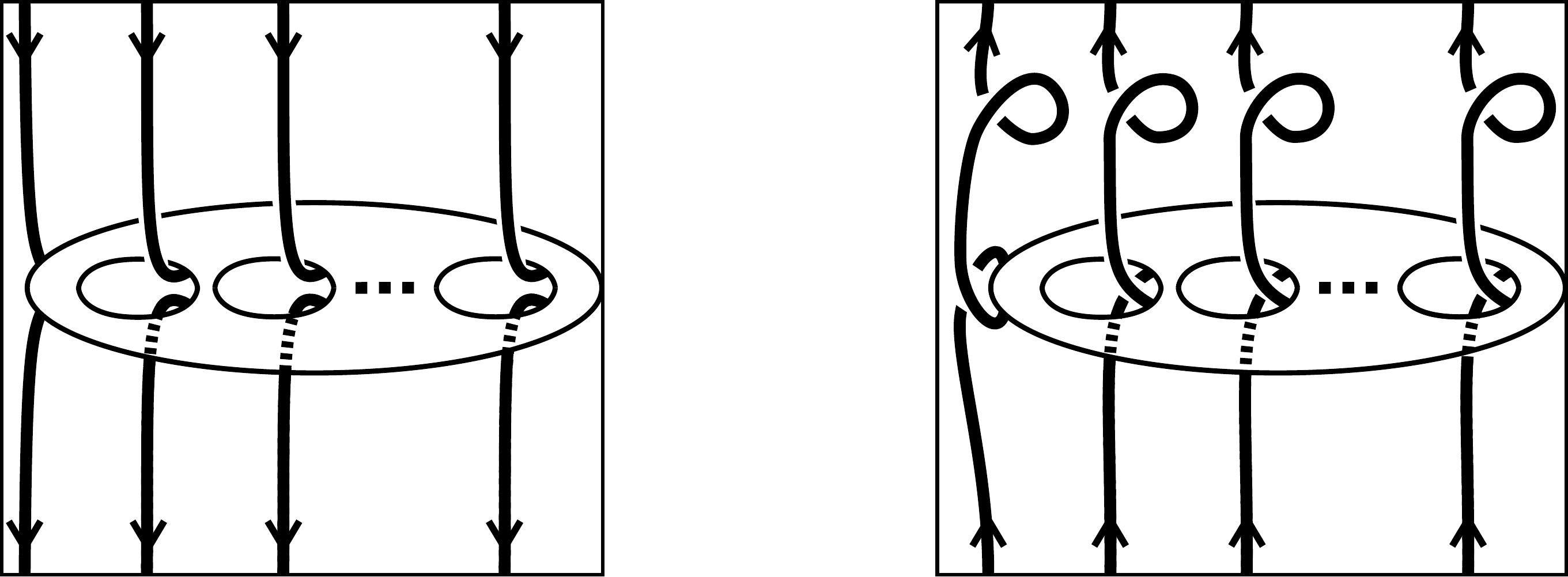}
    \caption{Two ways of turning a tangle diagram on $\Sigma$ into a horizontal $J$-tangle diagram: from top to bottom (left) and from bottom to top (right).}
	\label{fig:diagram_vertical_horizontal_updown}
\end{figure}

Let $|J|=0,2$. We reinterpret a $J$-tangle diagram on $(\Sigma_{0,p+1},J)$ as a horizontal projection as follows.
Draw $\Sigma_{0,p+1}$ so that $J$ lies in the thickened boundary segments shown on the left of Figure~\ref{fig:diagram_vertical_horizontal}.
Identify $\Sigma_{0,p+1}$ with $e(\Sigma_{0,p+1})\subset D_p\times\{1/2\}$, and draw the poles as downward-oriented lines.
Let $I_1=I_2=[0,1]$.
Identify $D_p$ with $(I_1\times I_2)\setminus\bigcup_{i=1}^p\{(i/(p+1),1/2)\}$ so that $(i/(p+1),1/2)$ corresponds to the $i$-th puncture of $D_p$.
We project $I_1\times I_2\times[0,1]$ along the $I_2$-direction.
the RHS of Figure~\ref{fig:diagram_vertical_horizontal} shows the horizontal projection after slightly tilting $e(\Sigma_{0,p+1})$.
When $|J|=2$, connect the outgoing (resp. incoming) endpoint to the top (resp. bottom) of the horizontal projection as on the left of Figure~\ref{fig:diagram_vertical_horizontal_updown}. This gives a $J$-tangle diagram whose endpoints lie at the top and bottom. We call the resulting diagram a \textbf{horizontal $J$-tangle diagram}. In such a diagram, the poles are drawn explicitly, and crossings with poles are included.
These horizontal $J$-tangle diagrams are used in \cite{BN+25}.
Let $T_1$ and $T_2$ be horizontal $J$-tangle diagrams such that the incoming endpoint of $T_1$ and the outgoing endpoint of $T_2$ have the same $I_1$-coordinate. Their \textbf{composition} $T_1\circ T_2$ is obtained by joining them with $T_1$ above $T_2$ in the $[0,1]$-direction.
This composition agrees with composition of $J$-tangles described immediately before Lemma~\ref{lem:sigma_Leibniz}.

\begin{remark}
The restriction on the position of $J$ on the left of Figure~\ref{fig:diagram_vertical_horizontal} ensures that composition of horizontal $J$-tangle diagrams agrees with composition of $J$-tangle diagrams.

The construction on the right of Figure~\ref{fig:diagram_vertical_horizontal_updown} also converts a $J$-tangle diagram on $\Sigma$ into a horizontal one. However, with this choice of endpoint connections, composition of horizontal $J$-tangle diagrams does not agree with composition induced by the embedding of $\pi_1(\Sigma)$ into $\widetilde{\cA}(\Sigma,J)$.
We use this second choice when we define the height involution in Section~\ref{subsec:height_involution}.
\end{remark}

\begin{remark}
For $|J|>2$, analogous generalizations of the statements in Section~\ref{sec:formality_genus0} can be proved by specifying how the endpoints are moved. To simplify the discussion, we give proofs only for $|J|=0,2$.
\end{remark}

For horizontal $J$-tangle diagrams, define two types of singularities by
\begin{align} 
\begin{array}{c}\includegraphics[scale=0.15]{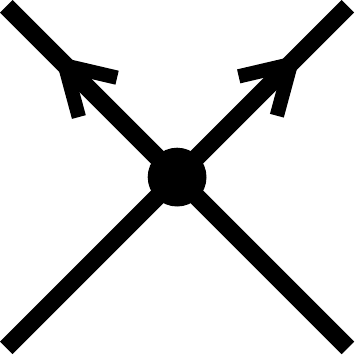}\end{array}:=\begin{array}{c}\includegraphics[scale=0.15]{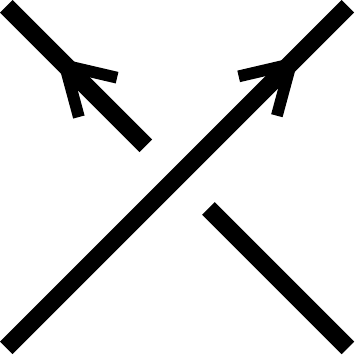}\end{array}-\begin{array}{c}\includegraphics[scale=0.15]{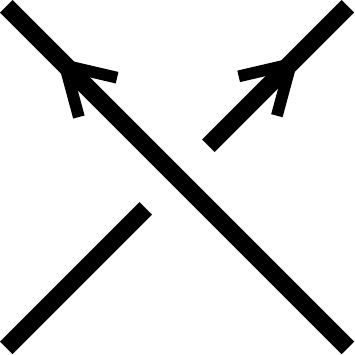}\end{array}\label{eq:singular_Vassiliev},\qquad
\begin{array}{c}\includegraphics[scale=0.14]{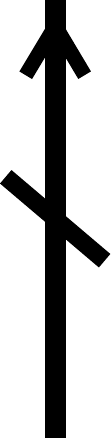}\end{array}=\begin{array}{c}\includegraphics[scale=0.14]{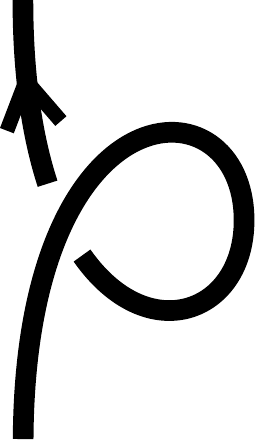}\end{array}-\begin{array}{c}\includegraphics[scale=0.14]{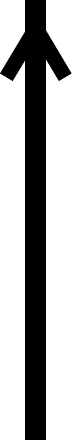}\end{array}
\end{align}
These formulas allow us to regard $J$-tangles with singularities as elements of $\bQ\widetilde{\cT}^J$. In the left-hand singularity of \eqref{eq:singular_Vassiliev}, one of the branches may be a pole, but crossings between two poles are excluded. Poles do not occur in the right-hand singularity.
We call a $J$-tangle with these singularities a \textbf{singular $J$-tangle}. Let $\widetilde{\cT}_t^J$ be the submodule of $\bQ\widetilde{\cT}^J$ generated by singular $J$-tangles with at least $t$ singularities.
The descending filtration $\{\widetilde{\cT}_t^J\}_{t=0}^\infty$ is called the \textbf{$t$-filtration} of $\bQ\widetilde{\cT}^J$.
Set $\Gr\bQ\widetilde{\cT}^J=\prod_{t\geq0}\widetilde{\cT}_t^J/\widetilde{\cT}_{t+1}^J$.
This $t$-filtration agrees with the one in \cite{BN+25}.

\begin{remark}
Whereas \cite{BN+25} proves formality only in degrees zero and one of the associated graded quotient for the $s$-filtration, we prove formality without introducing the $s$-filtration.
\end{remark}

Let $(\Sigma_{0,p+1},J)$ be a marked surface.
Recall that a $J$-tangle $T$ is a properly embedded oriented $1$-manifold in $\Sigma_{0,p+1}\times[0,1]$. Consider the following data:
\begin{enumerate}
    \item the disjoint union of the abstract $1$-manifold underlying $T$ and $p$ labeled oriented copies of $[0,1]$;
    \item the identification of the endpoints of the underlying $1$-manifold with $J$.
\end{enumerate}We call these data the \textbf{$J$-skeleton}, or simply the \textbf{skeleton}, of $T$. In (2), the orientations of the interval components are required to agree with $J_{\rm out}$ and $J_{\rm in}$. Each component of the underlying $1$-manifold in (1) is called a \textbf{strand}, and each additional oriented interval is called a \textbf{pole}.\footnote{These poles will later be identified with the poles of $D_p$.}
A \textbf{$J$-chord diagram}, or simply a chord diagram, is an orientation-preserving diffeomorphism class of a $J$-skeleton equipped with an even number of points and a perfect matching of these points. Each matched pair is joined by a dashed chord as in $\begin{array}{c}\includegraphics[scale=0.15]{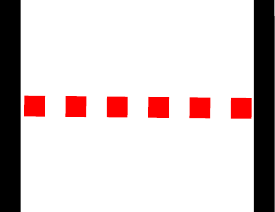}\end{array}$. We indicate the identification of the strand endpoints with $J$ and the labels of the poles, omitting them when they are clear or irrelevant. A circular strand with no chord endpoints is called a \textbf{trivial skeleton}.

A $J$-chord diagram is \textbf{admissible} if it has no chord joining two poles.
Let $\widetilde{\mathcal{D}}_t^J$ be the $\bQ$-vector space spanned by admissible $J$-chord diagrams with $t$ chords.
We call the $\bQ$-vector space $\widetilde{\mathcal{D}}^J:=\prod_{t\geq0}\widetilde{\mathcal{D}}_t^J$
the \textbf{universal space of admissible $J$-chord diagrams}.
The number of chords defines a grading on $\widetilde{\mathcal{D}}^J$, allowing infinite sums across degrees; its degree-$t$ part $\widetilde{\mathcal{D}}_t^J$ consists of finite sums.

In $\widetilde{\mathcal{D}}^J$,
the following relation is called the \textbf{$4T$ relation}:
\begin{align} 
\begin{array}{c}\includegraphics[scale=0.15]{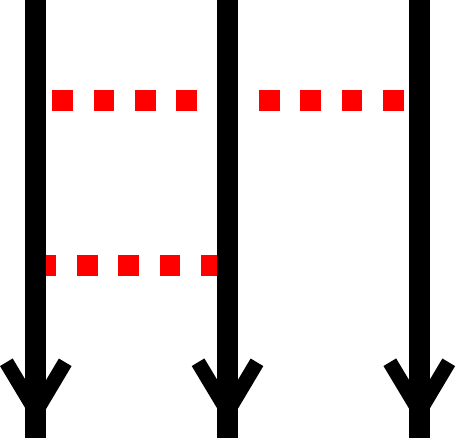}\end{array}-\begin{array}{c}\includegraphics[scale=0.15]{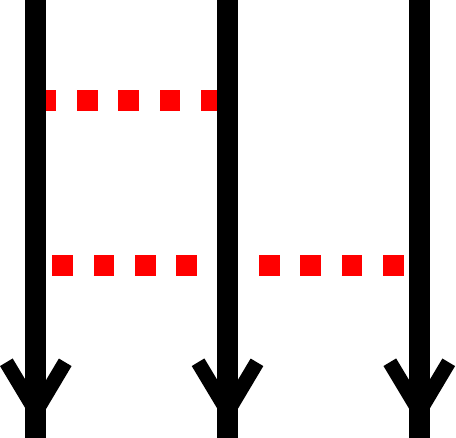}\end{array}
=\begin{array}{c}\includegraphics[scale=0.15]{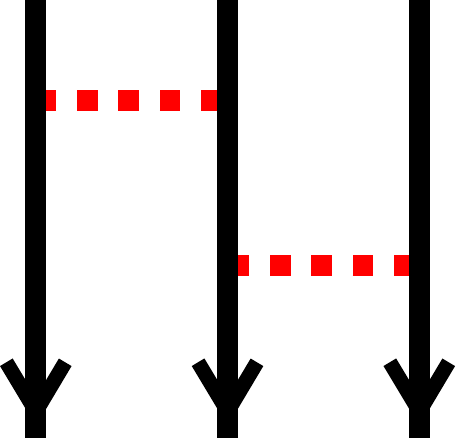}\end{array}
-\begin{array}{c}\includegraphics[scale=0.15]{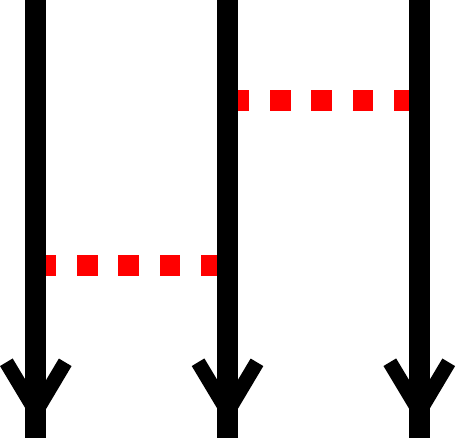}\end{array}\label{eq:4T}
\end{align}
A $4T$ relation is \textbf{admissible} if all four terms are admissible.
The quotient of $\widetilde{\mathcal{D}}^J$ by the admissible $4T$ relations is called the \textbf{space of admissible $J$-chord diagrams} and is denoted by $\mathcal{D}^J$.
Since \eqref{eq:4T} is homogeneous, the grading of $\widetilde{\mathcal{D}}^J$ descends to $\mathcal{D}^J$. Let $\mathcal{D}_t^J$ denote its degree-$t$ part.

To relate the HOMFLYPT quotient defined below to the modified quantized necklace algebra, introduce $a$ and $\abar$ satisfying
\begin{align}
h = \dfrac{a (\exp(\sqrt{a \abar}/2)-\exp(-\sqrt{a \abar}/2))}{\sqrt{a \abar}},\qquad
\hbar = \dfrac{\abar (\exp(\sqrt{a \abar}/2)-\exp(-\sqrt{a \abar}/2))}{\sqrt{a \abar}}\label{eq:h_a}
\end{align}
After expansion, neither RHS contains $\sqrt{a\abar}$.
Solving for $a$ and $\abar$ as power series in $h$ and $\hbar$, we obtain
\begin{align}
    a &= 2\sqrt{\dfrac{h}{\hbar}}\operatorname{arcsinh}\dfrac{\sqrt{h\hbar}}{2} 
    = \sum_{n=0}^{\infty} \frac{(-1)^n (2n)!}{16^n (n!)^2 (2n+1)} h^{n+1} \hbar^n =\sum_{n=0}^{\infty}\frac{(-1)^n (2n-1)!!}{4^n (2n)!! (2n+1)} h^{n+1} \hbar^n\label{eq:expansion_a_wrt_hhbar}\\
    \abar &= 2\sqrt{\dfrac{\hbar}{h}}\operatorname{arcsinh}\dfrac{\sqrt{h\hbar}}{2} 
    = \sum_{n=0}^{\infty} \frac{(-1)^n (2n)!}{16^n (n!)^2 (2n+1)} h^{n} \hbar^{n+1} =\sum_{n=0}^{\infty}\frac{(-1)^n (2n-1)!!}{4^n (2n)!! (2n+1)} h^{n} \hbar^{n+1}\label{eq:expansion_abar_wrt_hhbar}
\end{align}

Here $!!$ denotes the double factorial, with $(-1)!!=1$ and $0!!=1$.
In particular, if $(h,\hbar)$ denotes the ideal generated by $h$ and $\hbar$ in $\bQ\llbracket\rho\hbar,h\hbar\rrbracket[\rho,h,\hbar]$, then $a=h$ and $\abar=\hbar$ modulo $(h,\hbar)^2$.
Since $h$ and $\hbar$ have component degrees $1$ and $-1$, respectively, \eqref{eq:expansion_a_wrt_hhbar} and \eqref{eq:expansion_abar_wrt_hhbar} imply
that $a$ and $\abar$ have component degrees $1$ and $-1$, respectively.
Setting $a=\abar$ recovers the setting of \cite{BN+25}.

The \textbf{HOMFLYPT quotient} $\mathcal{D}_{\widetilde{\cA}}^J$ of $\widetilde{\mathcal{D}}^J$ is the quotient of $\widetilde{\mathcal{D}}^J\llbracket\rho\hbar,h\hbar\rrbracket[\rho,h,\hbar]$ by the admissible $4T$ relations and the following relations, graded using $a$, $\abar$, and $\rho_a$:
\begin{align}
&\begin{array}{c}\includegraphics[scale=0.15]{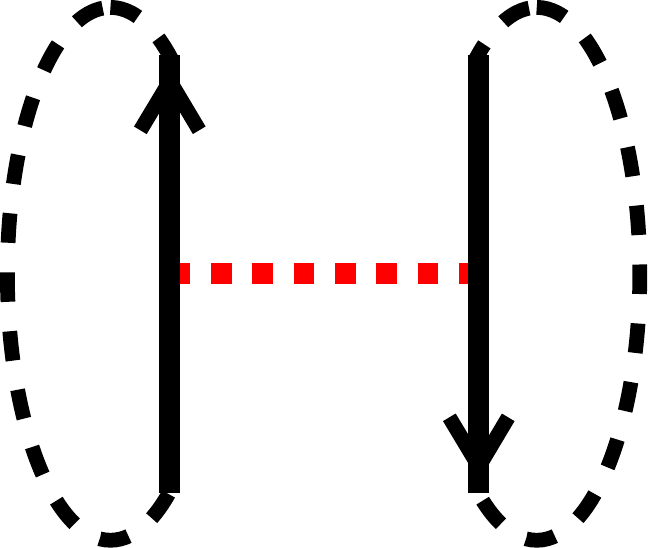}\end{array}=a\begin{array}{c}\includegraphics[scale=0.15]{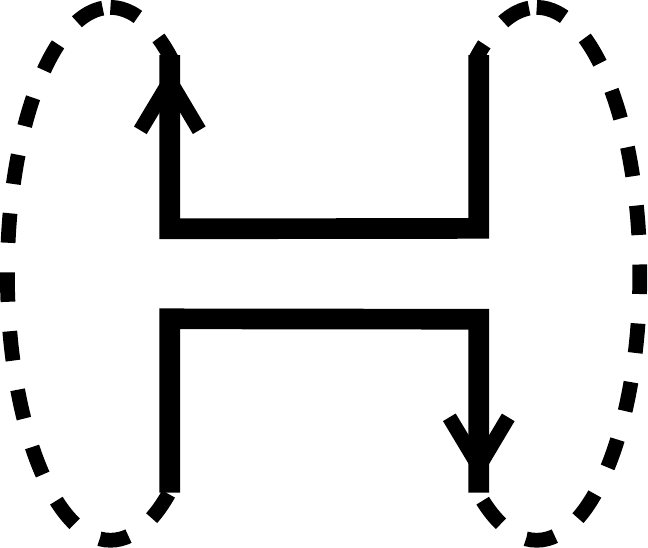}\end{array},%\label{eq:chord_a}\\
\qquad\begin{array}{c}\includegraphics[scale=0.15]{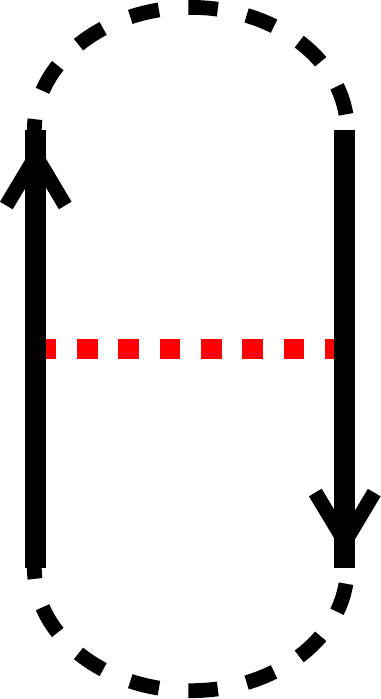}\end{array}=\abar\begin{array}{c}\includegraphics[scale=0.15]{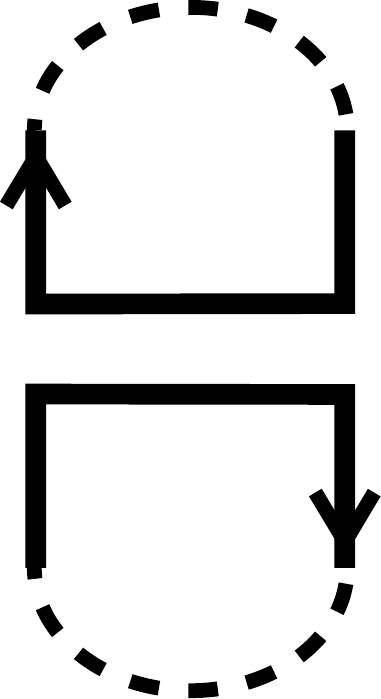}\end{array},%\label{eq:chord_abar}\\
\qquad\begin{array}{c}\includegraphics[scale=0.15]{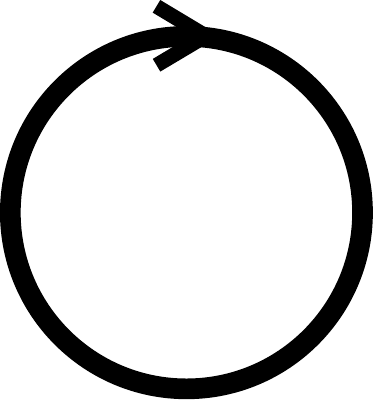}\end{array}=2\rho_a\label{eq:chord_diag_rel}
\end{align}
These relations apply only to chords joining strands. If a strand is an interval, the dotted connections indicate the identification of $J_{\rm out}$ with $J_{\rm in}$. To make the relations homogeneous, assign $t$-degree $1$ to both $a$ and $\abar$, and assign $t$-degree zero to
$\rho_a:=\nu^{-1}\dfrac{\sinh{(\rho\hbar)}}{\hbar}$
The invertible element $\nu=\dfrac{\cosh(\rho\hbar)}{\cosh(\sqrt{a\abar}/2)}=1-\dfrac{a\abar}{8}+\dfrac{\rho_a^2\abar^2}{2}+\dots$ has degree zero; see Remark~\ref{rmk:Kontsevich_well_defined}. The definition of $\nu$ is given as \eqref{eq:nu} and this identity is shown in Remark~\ref{rmk:Kontsevich_well_defined}. 
Until then, one may regard $\nu$ as the correction series at local maxima and minima that makes the associator-based Kontsevich integral invariant under isotopy.
Thus, $\mathcal{D}_{\widetilde{\cA}}^J$ is graded by the sum of the number of chords and the exponents of $a$ and $\abar$.

For example, if $\Sigma=\Sigma_{0,1}$ and $|J|=0,2$, then $\mathcal{D}_{\widetilde{\cA}}^J\cong\bQ\llbracket\rho\hbar,h\hbar\rrbracket[\rho,h,\hbar]$. We prove a more general statement in the proof of Lemma~\ref{lem:isom_diag_A_N}.

The \textbf{component degree} of a $J$-chord diagram is defined by
$$\#\{\text{strand components}\}+\text{(exponent of $a$)}-\text{(exponent of $\abar$)}+\text{
{(exponent of $\rho_a$)}}\footnote{For $|J|>2$, this must be defined using the operation corresponding to virtual closure.}.$$

\begin{figure}	\centering\includegraphics[scale=0.15]{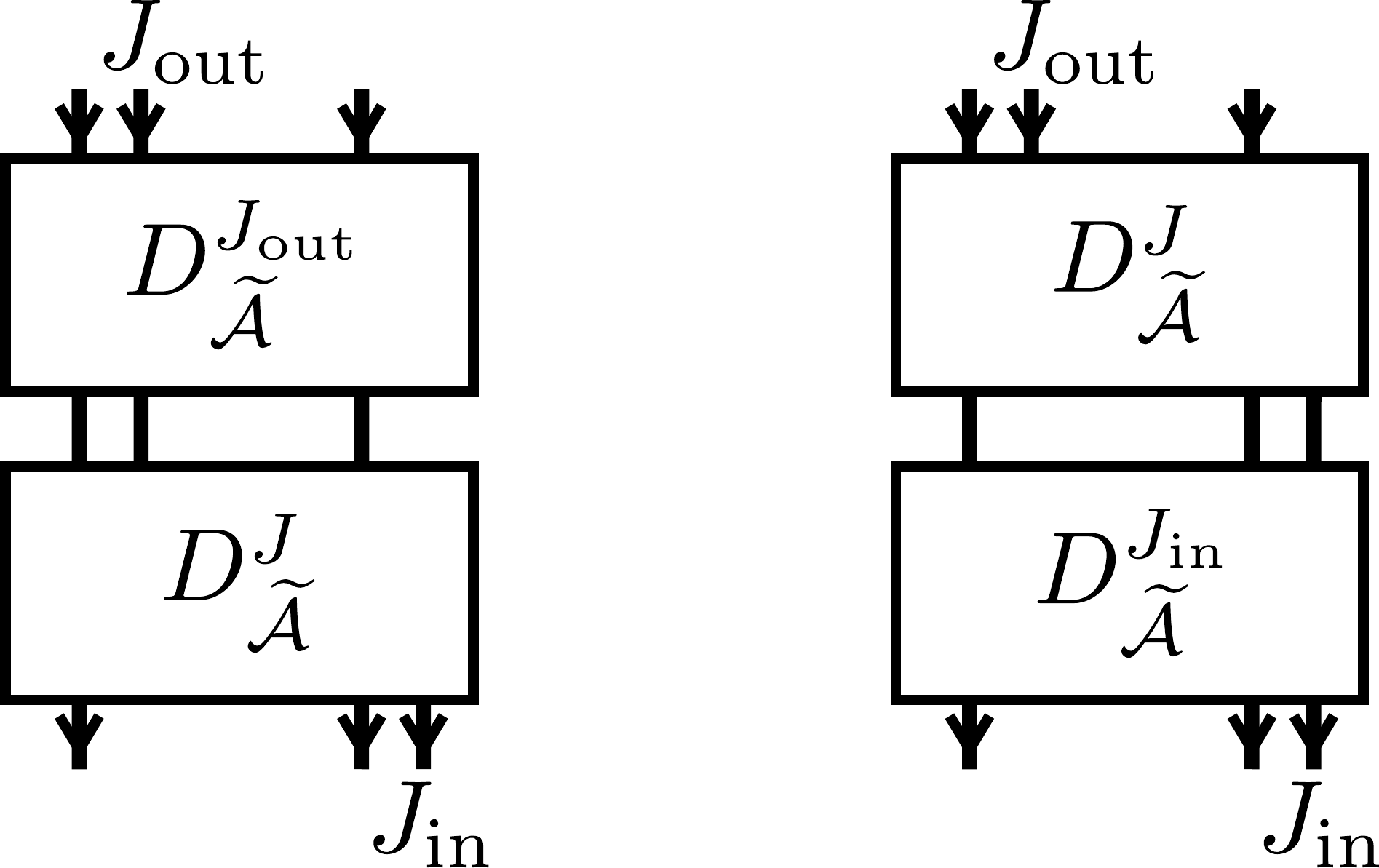}
\caption{Schematic pictures of the left action (left)
and the right action (right).}
\label{fig:action_D} 
\end{figure}

We now define the left action (resp. right action).
For $D_{\widetilde{\cA}}\in\mathcal{D}_{\widetilde{\cA}}:=\mathcal{D}_{\widetilde{\cA}}^\emptyset$, add a parallel strand near the pole adjacent to $J_{\rm out}$ (resp. $J_{\rm in}$). No parallel strand is added when the endpoint lies on the boundary corresponding to $\partial_0$.
At each chord endpoint on this pole, make the local replacement
$$\begin{array}{c}\includegraphics[scale=0.15]{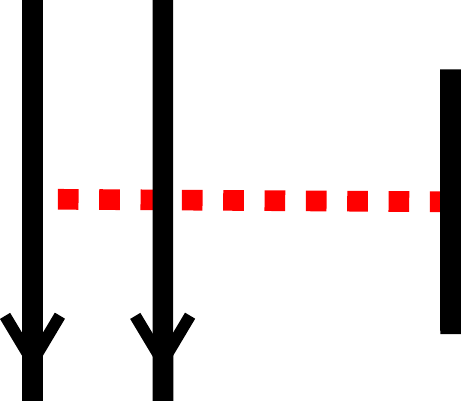}\end{array}\mapsto \begin{array}{c}\includegraphics[scale=0.15]{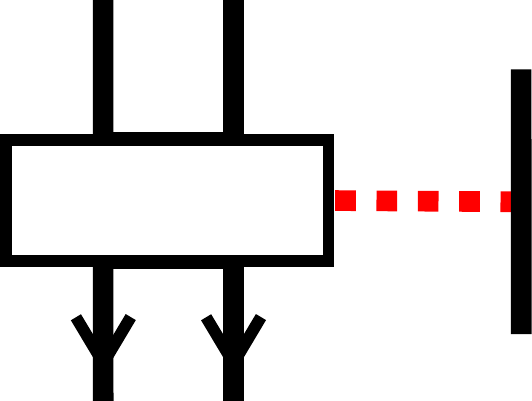}\end{array}$$
Denote the resulting element of $\mathcal{D}_{\widetilde{\cA}}^{J_{\rm out}}$ (resp. $\mathcal{D}_{\widetilde{\cA}}^{J_{\rm in}}$) by $D_{\widetilde{\cA}}^{J_{\rm out}}$ (resp. $D_{\widetilde{\cA}}^{J_{\rm in}}$).
Here $J_{\rm out}'$ (resp. $J_{\rm in}'$) is obtained by moving $J_{\rm out}$ (resp. $J_{\rm in}$) slightly along the oriented boundary of $\Sigma$, and $\mathcal{D}_{\widetilde{\cA}}^{J_{\rm out}}$ (resp. $\mathcal{D}_{\widetilde{\cA}}^{J_{\rm in}}$) denotes the submodule generated by $\{J_{\rm out},J_{\rm out}'\}$-tangles (resp. $\{J_{\rm in},J_{\rm in}'\}$-tangles). The coupon appearing in this operation denotes
\begin{align}
\begin{array}{c}\includegraphics[scale=0.15]{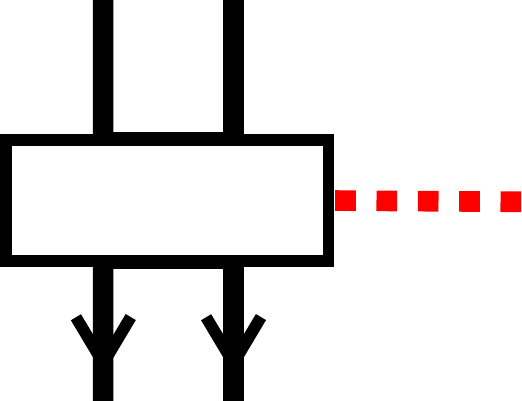}\end{array}
=\begin{array}{c}\includegraphics[scale=0.15]{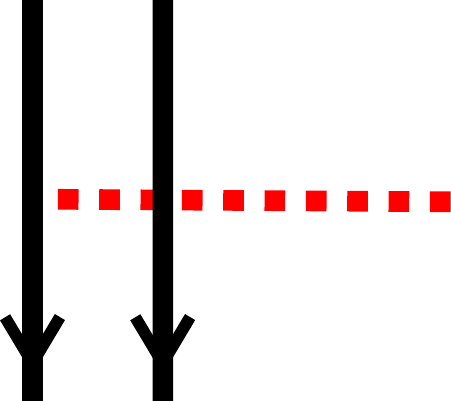}\end{array}
+\begin{array}{c}\includegraphics[scale=0.15]{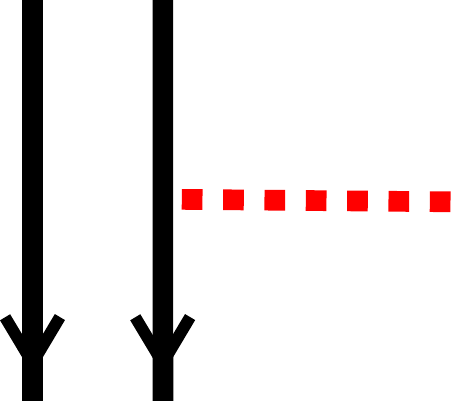}\end{array}.\label{eq:coupon}
\end{align}
Define the left action $D_{\widetilde{\cA}}^{J_{\rm out}}\cdot D_{\widetilde{\cA}}^J$ by composing $D_{\widetilde{\cA}}^{J_{\rm out}}$ with $D_{\widetilde{\cA}}^J$.
This is illustrated in Figure~\ref{fig:action_D}.
Vertical lines not labeled by $J_{\rm out}$ or $J_{\rm in}$ represent poles. The right action is defined similarly.

\begin{remark}
In \cite[Section~4.6]{BN+25}, the corresponding object is called the Conway quotient. There, $a=\abar$ is imposed, whereas the third relation in \eqref{eq:chord_diag_rel} is not; instead, trivial skeletons generate part of the coefficient ring.
\end{remark}

We now consider the framed Kontsevich integral
$Z\colon\bQ\widetilde{\cT}^J\to\mathcal{D}^J$.
Fixing a parenthesization of the pole endpoints and an associator determines a framed Kontsevich integral, which is an isotopy invariant of framed tangles. We use the normalization by the number of local maxima and minima for which, in the HOMFLYPT quotient, the image of a trivial framed loop is the product of $\nu$ from \eqref{eq:nu} and a trivial loop.
See \cite{CDM12, Oht02} for details. In the rest of this section, we use the same parenthesization of the pole endpoints at the top and bottom, so that the construction is compatible with composition.
We also assume that the associator has coefficients in $\bQ$.

For the HOMFLYPT skein module $\widetilde{\mathcal{A}}(\Sigma,J)$,
consider the $\bQ\llbracket\rho\hbar,h\hbar\rrbracket[\rho,h,\hbar]$-submodule
\begin{align}
\cI^n\widetilde{\cA}(\Sigma,J)={\rm span}\{h^i\hbar^j L\mid\text{$L$ is a singular $J$-tangle},\ i+j+\#\{\text{singularities}\}\geq n\}.
\end{align}
These submodules define
a descending filtration $\{\cI^n\widetilde{\mathcal{A}}(\Sigma,J)\}_{n\in\bN}$. Let $\widehat{\Gr_I}\widetilde{\mathcal{A}}(\Sigma,J)$ denote its completed associated graded quotient $\prod_{n=0}^\infty\cI^n\widetilde{\cA}(\Sigma,J)/\cI^{n+1}\widetilde{\cA}(\Sigma,J)$.

\begin{remark}
\begin{enumerate}
    \item Theorem~\ref{thm:two_filtrations} below shows that $\{\cI^n\widetilde{\mathcal{A}}(\Sigma,J)\}_{n\in\bN}$ agrees with the descending filtration defined using boxes.

    \item Proposition~\ref{prop:filtrations_relations} below shows that $\{\cI^n\widetilde{\mathcal{A}}(\Sigma,J)\}_{n\in\bN}$ induces the same topology as the filtration $\{F^n\widetilde{\mathcal{A}}(\Sigma,J)\}_{n\in\bN}$ of Definition~\ref{def:finite-type_filtration}. Hence, the two completions of the HOMFLYPT skein module coincide; this completion is denoted by $\widehat{\widetilde{\mathcal{A}}(\Sigma,J)}$.
\end{enumerate}
\end{remark}

\begin{lemma}[cf. \cite{BN+25}]
\label{lem:Kontsevich_descend}
\begin{enumerate}
    \item The Kontsevich integral induces a homomorphism $Z_{\widetilde{\cA}}^J\colon\widehat{\widetilde{\mathcal{A}}(\Sigma,J)}\to\mathcal{D}_{\widetilde{\cA}}^J$.
    Moreover, $Z_{\widetilde{\cA}}^J$ preserves the filtration $\{\cI^n\widehat{\widetilde{\mathcal{A}}(\Sigma,J)}\}_{n\in\bN}$ and the degree filtration on $\mathcal{D}_{\widetilde{\cA}}^J$, and satisfies property $(\ast)$.
    
    \item The map $\Gr Z_{\widetilde{\cA}}^J\colon\widehat{\Gr_I}\widetilde{\mathcal{A}}(\Sigma,J)\xrightarrow{\cong}\mathcal{D}_{\widetilde{\cA}}^J$ is given locally by
\begin{align}
&\begin{array}{c}\includegraphics[scale=0.15]{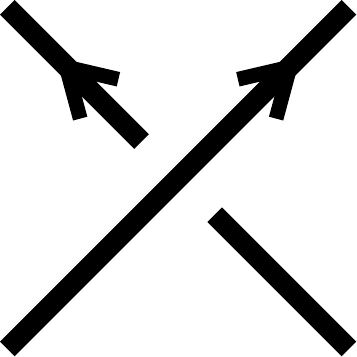}\end{array},\begin{array}{c}\includegraphics[scale=0.15]{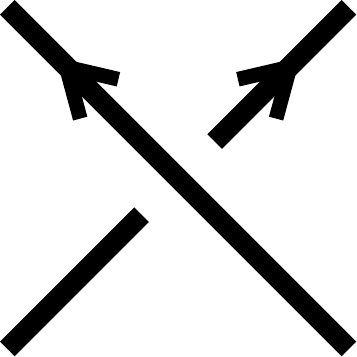}\end{array}\mapsto
    \begin{array}{c}\includegraphics[scale=0.15]{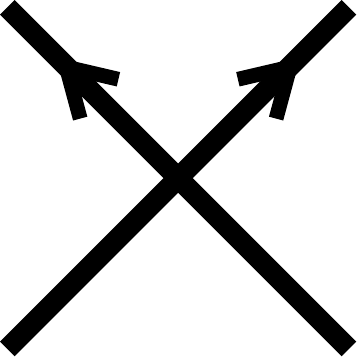}\end{array},\label{eq:GrZ_deg0}\\
    &\begin{array}{c}\includegraphics[scale=0.15]{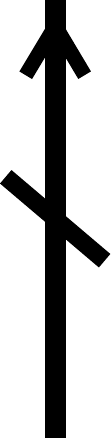}\end{array}
    \mapsto \rho_a\abar
    \begin{array}{c}\includegraphics[scale=0.15]{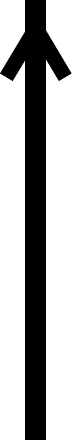}\end{array}\qquad
    \begin{array}{c}\includegraphics[scale=0.15]{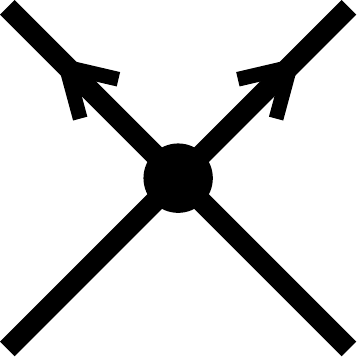}\end{array}
    \mapsto\eta
    \begin{array}{c}\includegraphics[scale=0.15]{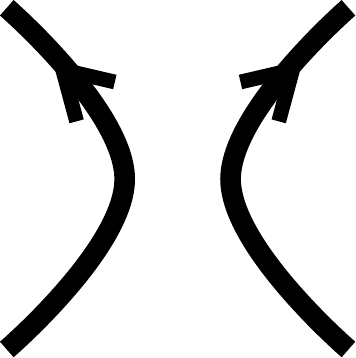}\end{array}, \qquad 
    \begin{array}{c}\includegraphics[scale=0.15]{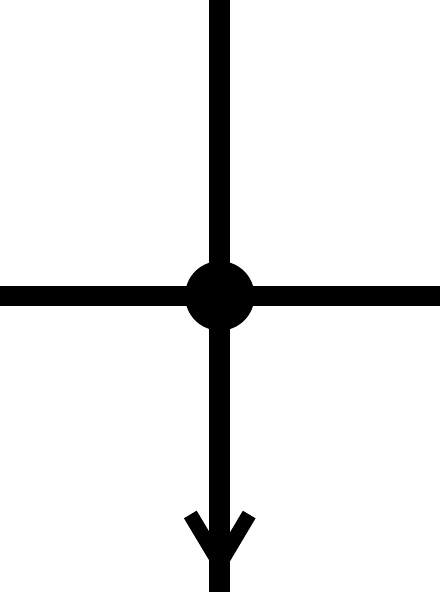}\end{array}\mapsto
    \begin{array}{c}\includegraphics[scale=0.15]{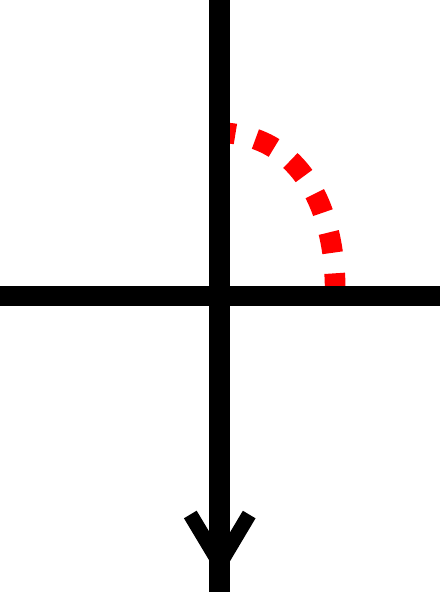}\end{array}.\label{eq:GrZ_deg1}
    \end{align}
    No pole appears in the third correspondence, whereas the vertical line in the fourth correspondence represents a pole.
    The parameter $\eta\in\{a,\abar\}$ in \eqref{eq:singular} is determined by preservation of component degree. Moreover, $\Gr Z_{\widetilde{\cA}}^J$ satisfies property $(\ast)$.
\end{enumerate}
\end{lemma}
We will show below that $Z_{\widetilde{\cA}}^J$ is an expansion.
\begin{proof}
\textbf{Preservation of the skein relations by $Z$.}
In this proof, we abbreviate $Z_{\widetilde{\cA}}^J$ to $Z$.
\begin{align*}
&Z\Big(\begin{array}{c}\includegraphics[scale=0.15]{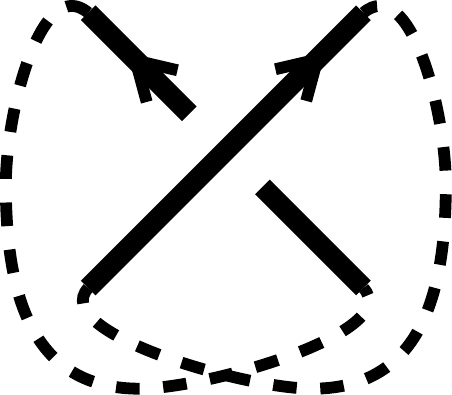}\end{array}\Big)=
    \begin{array}{c}\includegraphics[scale=0.15]{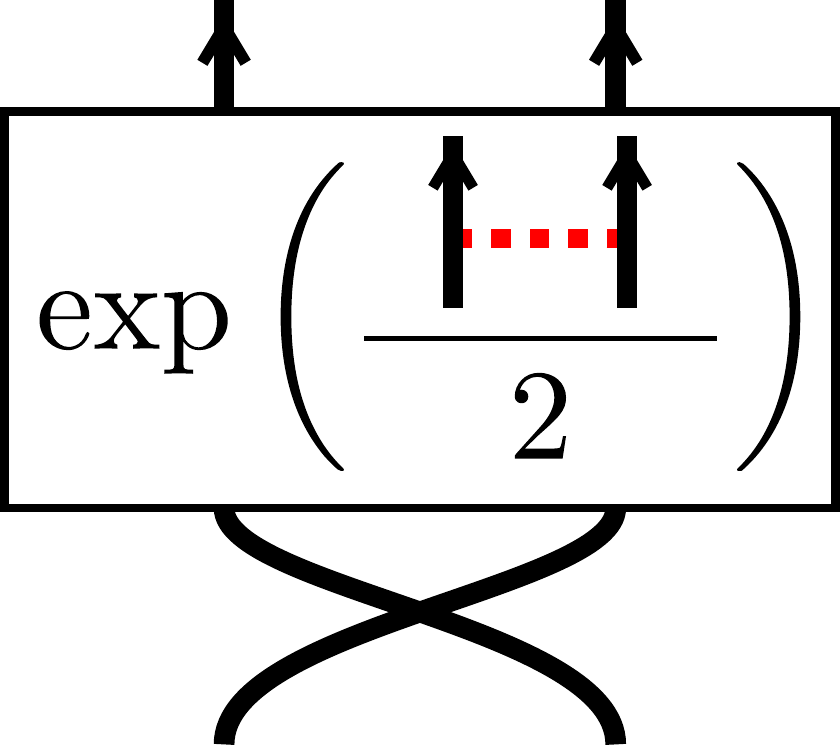}\end{array} 
    =\dfrac{a (\exp(\frac{\sqrt{a\abar}}{2})-\exp(-\frac{\sqrt{a\abar}}{2}))}{2\sqrt{a\abar}}\begin{array}{c}\includegraphics[scale=0.15]{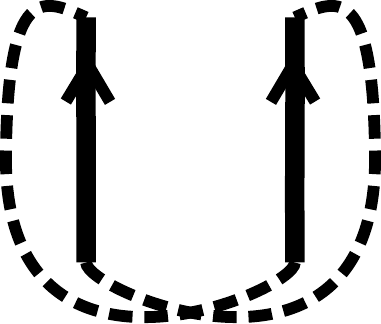}\end{array}+
    \dfrac{(\exp(\frac{\sqrt{a\abar}}{2})+\exp(-\frac{\sqrt{a\abar}}{2}))}{2}\begin{array}{c}\includegraphics[scale=0.15]{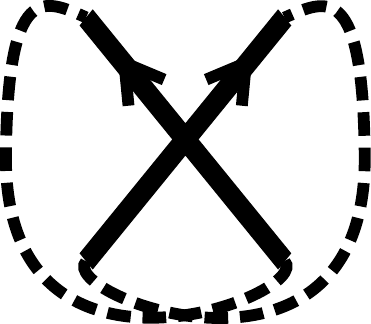}\end{array},\medskip\\
&Z\Big(\begin{array}{c}\includegraphics[scale=0.15]{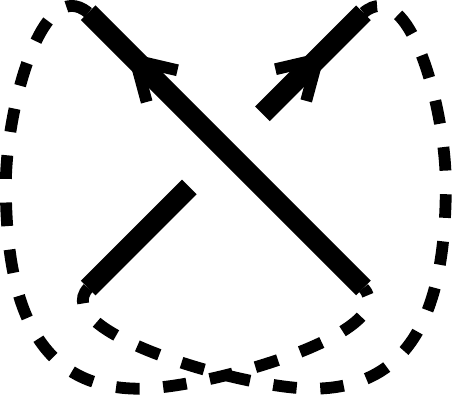}\end{array}\Big)
=\begin{array}{c}\includegraphics[scale=0.15]{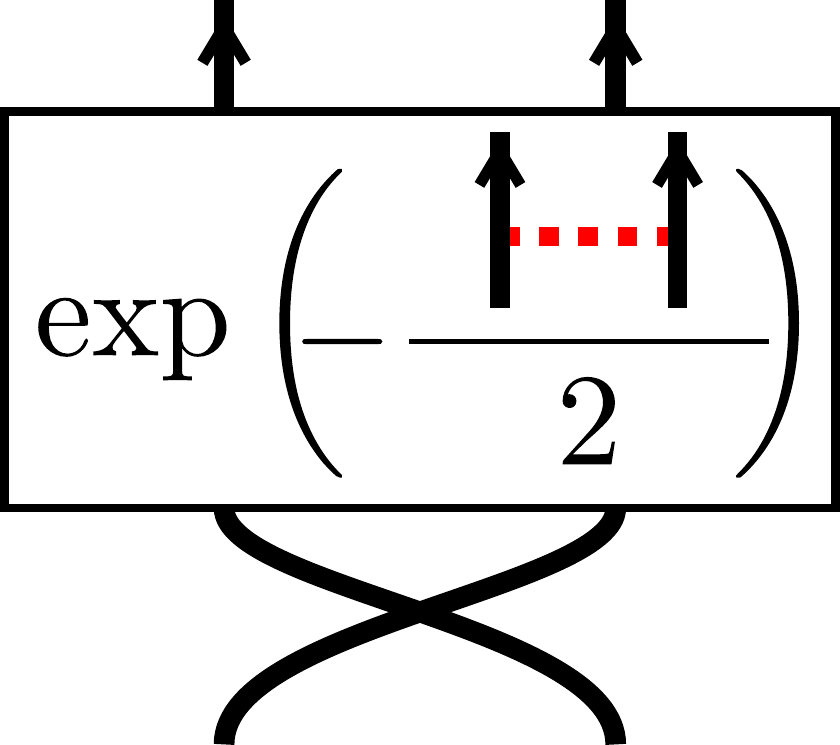}\end{array}=
    \dfrac{-a (\exp(\frac{\sqrt{a\abar}}{2})-\exp(-\frac{\sqrt{a\abar}}{2}))}{2\sqrt{a\abar}}\begin{array}{c}\includegraphics[scale=0.15]{draws/chord_2strands_0crossing02.pdf}\end{array}+
    \dfrac{(\exp(\frac{\sqrt{a\abar}}{2})+\exp(-\frac{\sqrt{a\abar}}{2}))}{2}\begin{array}{c}\includegraphics[scale=0.15]{draws/chord_2strands_1crossing02.pdf}\end{array},\medskip\\
&Z\Big(\begin{array}{c}\includegraphics[scale=0.15]{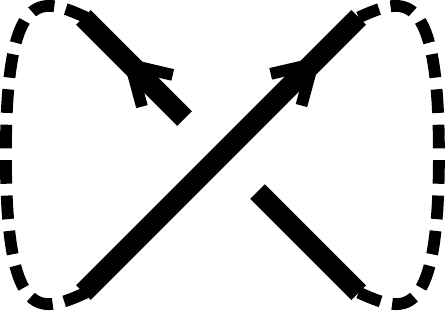}\end{array}\Big)
=\begin{array}{c}\includegraphics[scale=0.15]{draws/chord_2strands_exp_plus.pdf}\end{array}=
    \dfrac{\abar (\exp(\frac{\sqrt{a\abar}}{2})-\exp(-\frac{\sqrt{a\abar}}{2}))}{2\sqrt{a\abar}}\begin{array}{c}\includegraphics[scale=0.15]{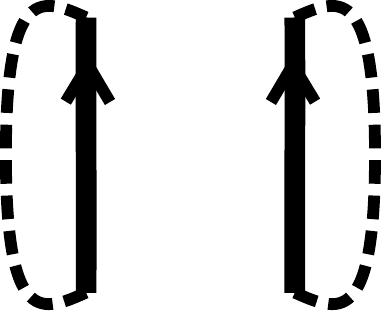}\end{array}+
    \dfrac{(\exp(\frac{\sqrt{a\abar}}{2})+\exp(-\frac{\sqrt{a\abar}}{2}))}{2}\begin{array}{c}\includegraphics[scale=0.15]{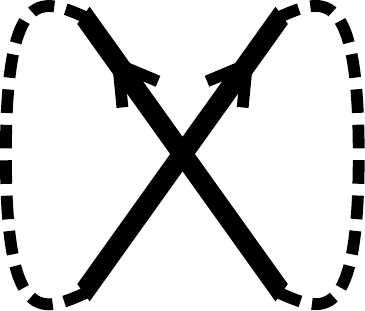}\end{array},\medskip\\
&Z\Big(\begin{array}{c}\includegraphics[scale=0.15]{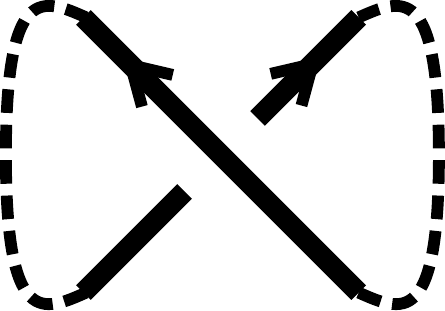}\end{array}\Big)
    =\begin{array}{c}\includegraphics[scale=0.15]{draws/chord_2strands_exp_minus.pdf}\end{array}=
    \dfrac{-\abar (\exp(\frac{\sqrt{a\abar}}{2})-\exp(-\frac{\sqrt{a\abar}}{2}))}{2\sqrt{a\abar}}\begin{array}{c}\includegraphics[scale=0.15]{draws/chord_2strands_0crossing.pdf}\end{array}+
    \dfrac{(\exp(\frac{\sqrt{a\abar}}{2})+\exp(-\frac{\sqrt{a\abar}}{2}))}{2}\begin{array}{c}\includegraphics[scale=0.15]{draws/chord_2strands_1crossing.pdf}\end{array}.
\end{align*}
Hence
\begin{align*}
&Z\Big(\begin{array}{c}\includegraphics[scale=0.15]{draws/crossing_posi_diff.pdf}\end{array}\Big)-Z\Big(\begin{array}{c}\includegraphics[scale=0.15]{draws/crossing_nega_diff.pdf}\end{array}\Big)=\dfrac{a (\exp(\frac{\sqrt{a\abar}}{2})-\exp(-\frac{\sqrt{a\abar}}{2}))}{\sqrt{a\abar}}\begin{array}{c}\includegraphics[scale=0.15]{draws/chord_2strands_0crossing02.pdf}\end{array}=h \cdot Z\Big(\begin{array}{c}\includegraphics[scale=0.15]{draws/chord_2strands_0crossing02.pdf}\end{array}\Big)\\
&Z\Big(\begin{array}{c}\includegraphics[scale=0.15]{draws/crossing_posi_same.pdf}\end{array}\Big)-Z\Big(\begin{array}{c}\includegraphics[scale=0.15]{draws/crossing_nega_same.pdf}\end{array}\Big)=\dfrac{\abar (\exp(\frac{\sqrt{a\abar}}{2})-\exp(-\frac{\sqrt{a\abar}}{2}))}{\sqrt{a\abar}}\begin{array}{c}\includegraphics[scale=0.15]{draws/chord_2strands_0crossing.pdf}\end{array}=\hbar  \cdot Z\Big(\begin{array}{c}\includegraphics[scale=0.15]{draws/chord_2strands_0crossing.pdf}\end{array}\Big)
\end{align*}
Thus, $Z$ preserves the skein relations (A) and (B).

For the values at local maxima and minima, set
\begin{align}
\nu:=\left(\begin{array}{c}\includegraphics[scale=0.15]{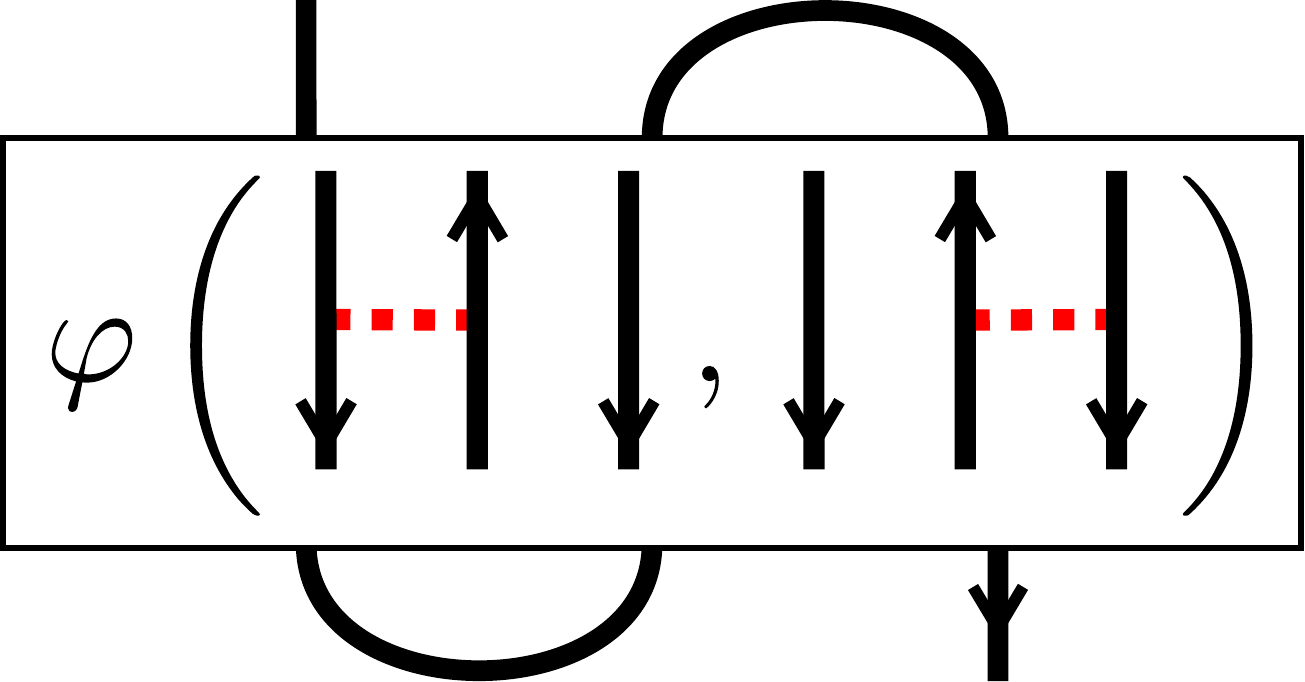}\end{array}\right)^{-1}\label{eq:nu}
\end{align}
and define
$Z\left(\begin{array}{c}\includegraphics[scale=0.15]{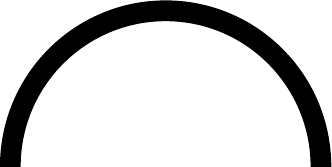}\end{array}\right)=Z\left(\begin{array}{c}\includegraphics[scale=0.15]{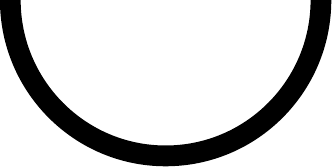}\end{array}\right)=\nu^{1/2}$.

The relations in the HOMFLYPT quotient give
$$Z\left(\begin{array}{c}\includegraphics[scale=0.15]{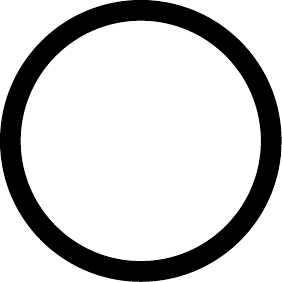}\end{array}\right)=\nu \cdot 2\rho_a Z\left(\emptyset\right)
=\nu \cdot \nu^{-1}\dfrac{2\sinh{(\rho \hbar)}}{\hbar}\cdot Z\left(\emptyset\right)
=
\dfrac{2\sinh{(\rho \hbar)}}{\hbar}
\cdot Z\left(\emptyset\right)$$
so $Z$ preserves skein relation (C).

Preservation of the skein relations (B) and (C) by $Z$ gives
\begin{align}
Z\left(\begin{array}{c}\includegraphics[scale=0.15]{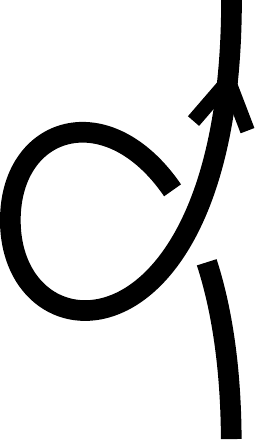}\end{array}\right)-Z\left(\begin{array}{c}\includegraphics[scale=0.15]{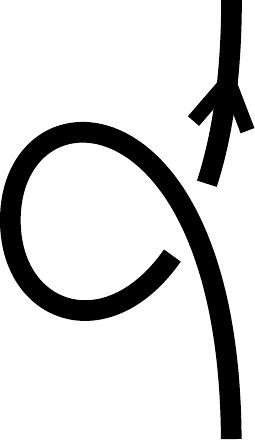}\end{array}\right)=
Z\left(\begin{array}{c}\includegraphics[scale=0.15]{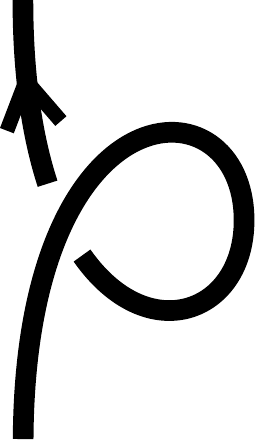}\end{array}\right)-Z\left(\begin{array}{c}\includegraphics[scale=0.15]{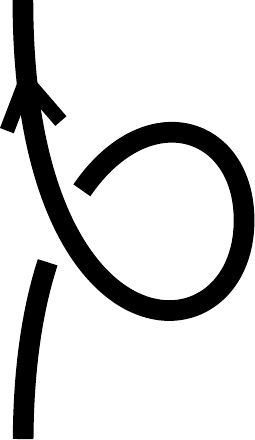}\end{array}\right)=\hbar\cdot Z\left(\begin{array}{c}\includegraphics[scale=0.15]{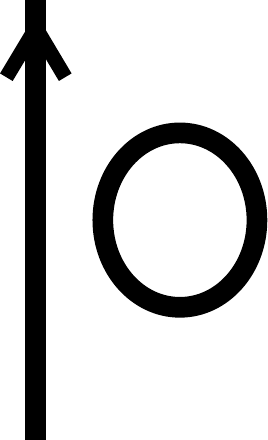}\end{array}\right)\nonumber\\
=\nu \cdot \nu^{-1} 2\sinh{\rho\hbar}\cdot
Z\left(\begin{array}{c}\includegraphics[scale=0.15]{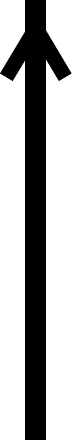}\end{array}\right)=2\sinh{\rho\hbar}\cdot
Z\left(\begin{array}{c}\includegraphics[scale=0.15]{draws/kink_noloop.pdf}\end{array}\right)\label{eq:HOMFLYPT_quotient_relAB}
\end{align}

Moreover,
using $Z\Big(\begin{array}{c}\includegraphics[scale=0.15]{draws/crossing_posi_same.pdf}\end{array}\Big)$ and the associator, we obtain $$Z\left(\begin{array}{c}\includegraphics[scale=0.15]{draws/kink_left_posi.pdf}\end{array}\right)=(1+(\text{terms of positive $\hbar$-degree}))
\cdot Z\left(\begin{array}{c}\includegraphics[scale=0.15]{draws/kink_noloop.pdf}\end{array}\right)$$
The analogous formulas hold for the other kinks.
In particular, there exist invertible elements $x,y\in\bQ\llbracket\rho\hbar,h\hbar\rrbracket[\rho,h,\hbar]$ such that
\begin{align}
Z\left(\begin{array}{c}\includegraphics[scale=0.15]{draws/kink_left_posi.pdf}\end{array}\right)=x\cdot Z\left(\begin{array}{c}\includegraphics[scale=0.15]{draws/kink_noloop.pdf}\end{array}\right),\qquad Z\left(\begin{array}{c}\includegraphics[scale=0.15]{draws/kink_left_nega.pdf}\end{array}\right)
=y\cdot Z\left(\begin{array}{c}\includegraphics[scale=0.15]{draws/kink_noloop.pdf}\end{array}\right)\label{eq:Kontsevich_RI}
\end{align}

Invariance of the Kontsevich integral under Reidemeister II and III moves gives
\begin{align}
Z\left(\begin{array}{c}\includegraphics[scale=0.15]{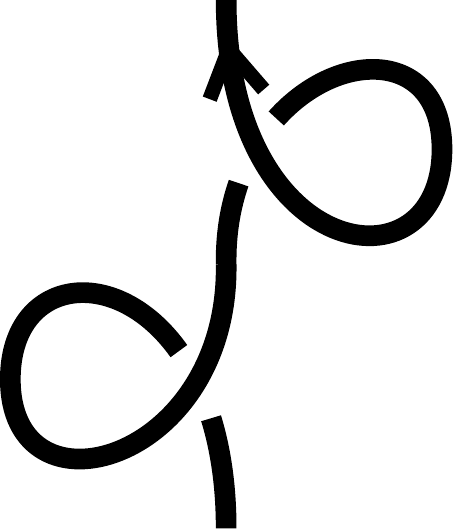}\end{array}\right)=Z\left(\begin{array}{c}\includegraphics[scale=0.15]{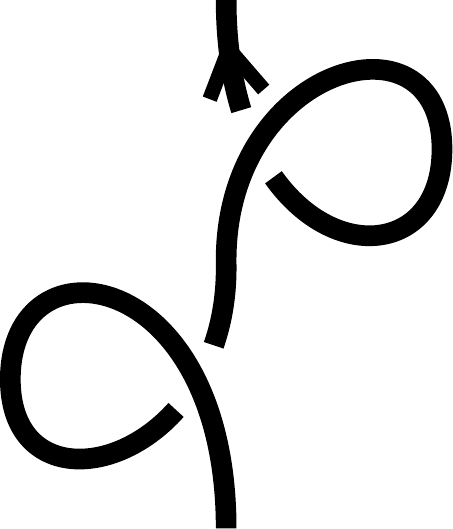}\end{array}\right)=Z\left(\begin{array}{c}\includegraphics[scale=0.15]{draws/kink_noloop.pdf}\end{array}\right)\label{eq:Kontsevich_RII_RIII}
\end{align}

Applying \eqref{eq:Kontsevich_RI} to the left and middle terms of \eqref{eq:Kontsevich_RII_RIII}, and using the invertibility of $x$ and $y$, gives
$$Z\left(\begin{array}{c}\includegraphics[scale=0.15]{draws/kink_right_posi.pdf}\end{array}\right)=y^{-1}\cdot Z\left(\begin{array}{c}\includegraphics[scale=0.15]{draws/kink_noloop.pdf}\end{array}\right),\qquad 
Z\left(\begin{array}{c}\includegraphics[scale=0.15]{draws/kink_right_nega.pdf}\end{array}\right)=x^{-1} \cdot Z\left(\begin{array}{c}\includegraphics[scale=0.15]{draws/kink_noloop.pdf}\end{array}\right)$$

By \eqref{eq:HOMFLYPT_quotient_relAB}, we have
$$x-y=y^{-1}-x^{-1}=2\sinh{\rho\hbar}=\exp(\rho\hbar)-\exp(-\rho\hbar)$$

Since $y^{-1}-x^{-1}=x^{-1}y^{-1}(x-y)$ and the coefficient ring is an integral domain,
it follows that $1-x^{-1}y^{-1}=0$, or equivalently, $y=x^{-1}$.
Hence
$$x-x^{-1}=\exp(\rho\hbar)-\exp(-\rho\hbar)$$
The constant term of $x$ is $1$, so
$x=\exp(\rho\hbar)$.
Therefore, $Z$ preserves skein relation (D).

For a singularity between a pole and a $J$-tangle, its image is locally an infinite sum of chord diagrams, each containing at least one chord.
It follows that $Z_{\widetilde{\cA}}^J$ preserves the filtration $\{\cI^n\widehat{\widetilde{\mathcal{A}}(\Sigma,J)}\}_{n\in\bN}$ and the degree filtration on $\mathcal{D}_{\widetilde{\cA}}^J$.
The doubling operation for the Kontsevich integral corresponds to \eqref{eq:coupon}. Since the products, the left and right actions, and compositions of horizontal $J$-tangle diagrams agree with those of $J$-tangle diagrams, property $(\ast)$ follows.

By the definition of $Z_{\widetilde{\cA}}^J$ and the leading terms of the identities above,
the map $\Gr_I Z_{\widetilde{\cA}}^J\colon\Gr\widetilde{\mathcal{A}}(\Sigma,J)\xrightarrow{\cong}\mathcal{D}_{\widetilde{\cA}}^J$ is given locally by \eqref{eq:GrZ_deg0} and \eqref{eq:GrZ_deg1}.
Moreover, property $(\ast)$ for $Z_{\widetilde{\cA}}^J$ implies the same property for $\Gr Z_{\widetilde{\cA}}^J$.
\end{proof}

\begin{remark}\label{rmk:Kontsevich_well_defined}
\begin{enumerate}
    \item 
    $\begin{array}{c}\includegraphics[scale=0.15]{draws/kink_right_posi.pdf}\end{array}+\begin{array}{c}\includegraphics[scale=0.15]{draws/kink_right_nega.pdf}\end{array}$
    Applying the skein relations to this and then applying $Z$ gives
    $(\exp(\rho\hbar)+\exp(-\rho\hbar))Z\left(\begin{array}{c}\includegraphics[scale=0.15]{draws/kink_noloop.pdf}\end{array}\right)=2\cosh(\rho\hbar)Z\left(\begin{array}{c}\includegraphics[scale=0.15]{draws/kink_noloop.pdf}\end{array}\right)$.
    On the other hand, $Z\left(\begin{array}{c}\includegraphics[scale=0.15]{draws/kink_right_posi.pdf}\end{array}\right)+Z\left(\begin{array}{c}\includegraphics[scale=0.15]{draws/kink_right_nega.pdf}\end{array}\right)=\nu\cdot\dfrac{2\abar}{\sqrt{a\abar}}\cosh(\sqrt{a\abar}/2)\begin{array}{c}\includegraphics[scale=0.15]{draws/kink_noloop.pdf}\end{array}$.
    Hence, $\nu=\dfrac{\cosh(\rho\hbar)}{\cosh(\sqrt{a\abar}/2)}$.

    \item The case $a=\abar$ of Lemma~\ref{lem:Kontsevich_descend} is proved in \cite{BN+25}. In this paper, we also show that $\Gr\widetilde{\mathcal{A}}(\Sigma,J)$ is isomorphic to the modified quantized necklace algebra.

    \item The third relation in \eqref{eq:chord_diag_rel} corresponds to the relation (C) 
    %relations (C) and (D) 
    of the HOMFLYPT skein algebra. It is imposed to study formality for the HOMFLYPT skein algebra. Without this relation, one still obtains formality for the Conway quotient.
\end{enumerate}
\end{remark}

\subsection{The necklace-like HOMFLYPT quotient and the modified quantized necklace algebra}
We use the necklace description to construct $\psi_{\widetilde{\cA}}^J$. For this purpose, we introduce the necklace-like HOMFLYPT quotient $\mathcal{D}_{\widetilde{\cN_I}}^J$.
In the preceding subsection, we introduced a modified finite-type filtration on the HOMFLYPT skein algebra. To define a corresponding degree on the quantized necklace algebra, we modify its defining relations.

The following relation is called the \textbf{necklace-like $4T$ relation}:
\begin{align} 
\begin{array}{c}\includegraphics[scale=0.15]{draws/4T_13_12_d.pdf}\end{array}-\begin{array}{c}\includegraphics[scale=0.15]{draws/4T_12_13_d.pdf}\end{array}
=\eta\left(\begin{array}{c}\includegraphics[scale=0.15]{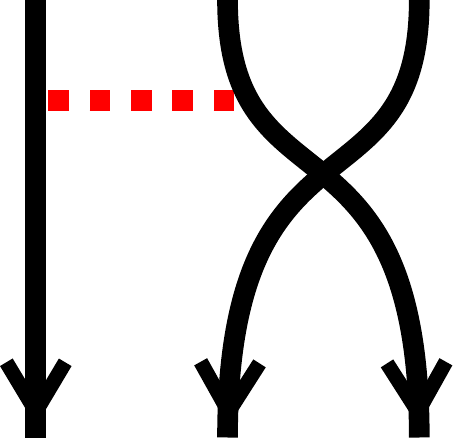}\end{array}
-\begin{array}{c}\includegraphics[scale=0.15]{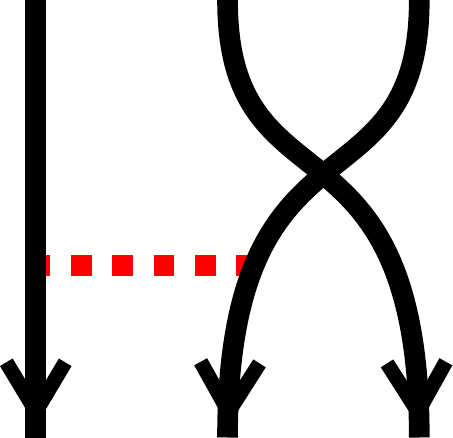}\end{array}\right)\label{eq:4T_necklace}
\end{align}
In each term, the leftmost of the three lines is a pole, and the other two lines are parts of strands. The parameter $\eta$ is $a$ if the middle and right strands belong to different components, and $\abar$ if they belong to the same component.

A chord diagram is \textbf{necklace-like} if it is admissible and contains no chord joining two strands. Thus, only pole--strand chords are allowed.
Consider the $\bQ\llbracket\rho\hbar,h\hbar\rrbracket[\rho,h,\hbar]$-submodule of $\mathcal{D}^J\llbracket\rho\hbar,h\hbar\rrbracket[\rho,h,\hbar]$ spanned by necklace-like $J$-chord diagrams. Its quotient by the third relation in \eqref{eq:chord_diag_rel}, concerning trivial skeletons, and by \eqref{eq:4T_necklace} is called the \textbf{necklace-like HOMFLYPT quotient} and is denoted by $\mathcal{D}_{\widetilde{\cN_I}}^J$.

The following observation is central to our approach.
\begin{lemma}\label{lem:isom_diag_A_N}
The natural inclusion of diagrams induces a well-defined $\bQ\llbracket\rho\hbar,h\hbar\rrbracket[\rho,h,\hbar]$-module isomorphism $\mathcal{D}_{\widetilde{\cN_I}}^J\to\mathcal{D}_{\widetilde{\cA}}^J$ satisfying property $(\ast)$.
\end{lemma}
\begin{proof}
Let $\widetilde{\mathcal{D}}_{\widetilde{\cA}}^J$ be the quotient of $\widetilde{\mathcal{D}}^J$ by \eqref{eq:chord_diag_rel}. Imposing the admissible $4T$ relations then gives $\mathcal{D}_{\widetilde{\cA}}^J$. Similarly, take the subspace of $\mathcal{D}^J$ spanned by necklace-like $J$-chord diagrams and quotient it by the third relation in \eqref{eq:chord_diag_rel}, concerning trivial skeletons.
Denote this quotient by $\widetilde{\mathcal{D}}_{\widetilde{\cN_I}}^J$. We first prove that $\widetilde{\mathcal{D}}_{\widetilde{\cA}}^J$ and $\widetilde{\mathcal{D}}_{\widetilde{\cN_I}}^J$ are isomorphic as $\bQ\llbracket\rho\hbar,h\hbar\rrbracket[\rho,h,\hbar]$-modules, and then deduce the corresponding statement for $\mathcal{D}_{\widetilde{\cA}}^J$ and $\mathcal{D}_{\widetilde{\cN_I}}^J$.

Define the complexity of a chord diagram to be the pair
$$\text{(number of strand--strand chords, number of trivial skeletons)},$$
ordered lexicographically. For a necklace-like chord diagram, the first entry is zero.
The choice of $a$ or $\abar$ is determined by preservation of component degree. Consequently, the defining relations of $\widetilde{\mathcal{D}}_{\widetilde{\cA}}^J$ and $\widetilde{\mathcal{D}}_{\widetilde{\cN_I}}^J$ commute with one another. In particular, the reduction system is locally confluent. By the diamond lemma, the resulting necklace-like chord diagram is independent of the order in which the relations are applied \cite{SW07}. The diagrams of complexity zero form a basis. Since
$$\{\text{complexity-zero diagrams in $\widetilde{\mathcal{D}}_{\widetilde{\cN_I}}^J$}\}=\{\text{complexity-zero diagrams in $\widetilde{\mathcal{D}}_{\widetilde{\cA}}^J$}\},$$
the two bases coincide.
The module homomorphism induced by the inclusion of diagrams $\widetilde{\mathcal{D}}_{\widetilde{\cN_I}}^J\to\widetilde{\mathcal{D}}_{\widetilde{\cA}}^J$ sends each basis diagram to itself and is therefore bijective.

Next, impose \eqref{eq:4T} and \eqref{eq:4T_necklace} on $\widetilde{\mathcal{D}}_{\widetilde{\cA}}^J$ and $\widetilde{\mathcal{D}}_{\widetilde{\cN_I}}^J$, respectively,
to obtain $\mathcal{D}_{\widetilde{\cA}}^J$ and $\mathcal{D}_{\widetilde{\cN_I}}^J$.
When one of the lines in \eqref{eq:4T} is a pole, applying the first two relations of \eqref{eq:chord_diag_rel} gives \eqref{eq:4T_necklace}.
When all lines in \eqref{eq:4T} are strands, the two sides become equal after applying \eqref{eq:chord_diag_rel}.
Thus, \eqref{eq:4T} and \eqref{eq:4T_necklace} impose equivalent relations on $\widetilde{\mathcal{D}}_{\widetilde{\cA}}^J$ and $\widetilde{\mathcal{D}}_{\widetilde{\cN_I}}^J$.
Therefore, $\mathcal{D}_{\widetilde{\cA}}^J$ and $\mathcal{D}_{\widetilde{\cN_I}}^J$ are isomorphic as $\bQ\llbracket\rho\hbar,h\hbar\rrbracket[\rho,h,\hbar]$-modules.

Property $(\ast)$ follows from the fact that the map is induced by inclusion of diagrams.
Although the left and right actions on $\mathcal{D}_{\widetilde{\cN_I}}^J$ are induced from $\mathcal{D}_{\widetilde{\cA}}^J$, note that rewriting their values as necklace-like diagrams without coupons may increase the number of terms.
\end{proof}

Denote the inverse of the isomorphism in Lemma~\ref{lem:isom_diag_A_N} by $$\psi_{\mathcal{D}}^J\colon\mathcal{D}_{\widetilde{\cA}}^J\to\mathcal{D}_{\widetilde{\cN_I}}^J.$$
The map $\psi_{\mathcal{D}}^J$ applies \eqref{eq:chord_diag_rel} to an admissible chord diagram
and expresses it as a linear combination of necklace-like chord diagrams.

As in Figure~\ref{fig:disk_punctured},
let $\mathsf{B}$ be the set of classes $\partial_1,\dots,\partial_p\in H_1(\Sigma_{0,p+1};\bQ)$ represented by clockwise $i$-th peripheral loops. This set is a basis of $H_1(\Sigma_{0,p+1};\bQ)$.

To discuss formality in the presence of endpoints, we introduce $J$-necklace ensembles.
Let $(\Sigma,J)$ be a marked surface.
In the definition of an $H_1(\Sigma)$-necklace in Section~\ref{subsec:quantized_necklace}, allow compact oriented $1$-manifolds in place of closed oriented $1$-manifolds. A \textbf{$J$-necklace ensemble} is a disjoint union of such $H_1(\Sigma)$-necklaces, equipped with a total order on its beads, called their height order, and a bijection between the set of wire endpoints and $J$. Each interval wire is oriented from the endpoint corresponding to $J_{\rm out}$ to the endpoint corresponding to $J_{\rm in}$.

The relations defining the following modified quantized necklace algebra are chosen to match those of $\mathcal{D}_{\widetilde{\cN_I}}^J$.
\begin{definition}\label{def:modified_necklace}
Let $(\Sigma,J)$ be a marked oriented surface.
Consider the $\bQ\llbracket h\hbar,\rho\hbar\rrbracket[\rho,h,\hbar]$-module generated by $J$-necklace ensembles modulo relation (H) below.
Suppose that $J=\emptyset$.
For $J$-necklace ensembles $\alpha_i$ with $m_i$ beads, $i=1,2$, define $\alpha_1\alpha_2$ to be their abstract disjoint union $\alpha_1\sqcup\alpha_2$, retaining the heights of the beads of $\alpha_2$ and replacing each height $k$ in $\alpha_1$ by $m_2+k$. 
Here, we represent the total order by consecutive integers. Equipped with this product, the resulting $\bQ\llbracket h\hbar,\rho\hbar\rrbracket[\rho,h,\hbar]$-algebra is called the
\textbf{universal modified quantized necklace module} of $(\Sigma,J)$ and is denoted by $\widetilde{\mathrm{N}}(\Sigma,J)$.
The quotient by (H), (A')--(C'), (A''), and (B'') is called the \textbf{modified quantized necklace Lie module} of $(\Sigma,J)$ and is denoted by $\widetilde{\cN_I}(\Sigma,J)$.
When $J=\emptyset$, the product is defined as in Definition~\ref{def:necklace_Lie_algebra}. We call the resulting algebra the \textbf{modified quantized necklace Lie algebra} of $\Sigma$, denoted by  $\widetilde{\cN_I}(\Sigma)$.
\begin{align*}
&{\rm (H)} 
\begin{array}{c}\includegraphics[scale=0.26]{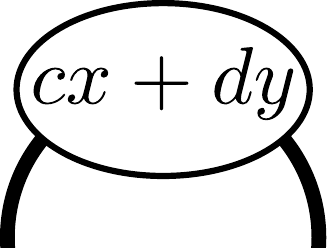}\end{array}=c\begin{array}{c}\includegraphics[scale=0.26]{draws/ax.pdf}\end{array}+d\begin{array}{c}\includegraphics[scale=0.26]{draws/by.pdf}\end{array}\quad (x,y\in H_1(\Sigma),\ c,d\in \bQ),\\
&{\rm (A')} 
\begin{array}{c}\includegraphics[scale=0.26]{draws/xy_i+1i.pdf}\end{array}-\begin{array}{c}\includegraphics[scale=0.26]{draws/xy_ii+1.pdf}\end{array}=0,\\
&{\rm (A'')} \begin{array}{c}\includegraphics[scale=0.26]{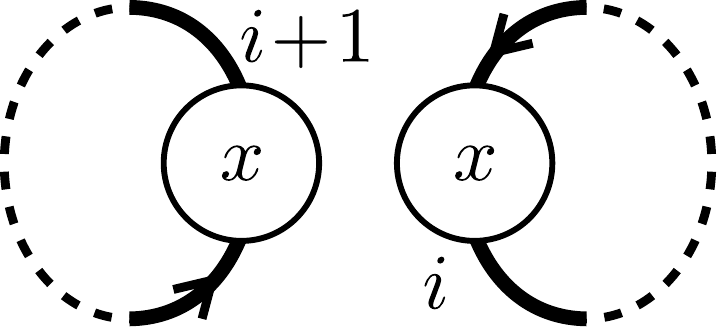}\end{array}-\begin{array}{c}\includegraphics[scale=0.26]{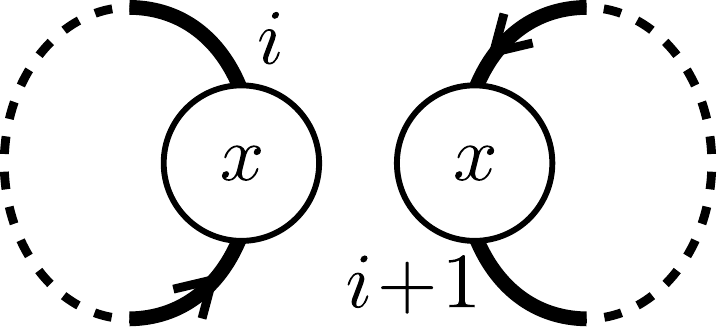}\end{array}=h\left(\begin{array}{c}\includegraphics[scale=0.26]{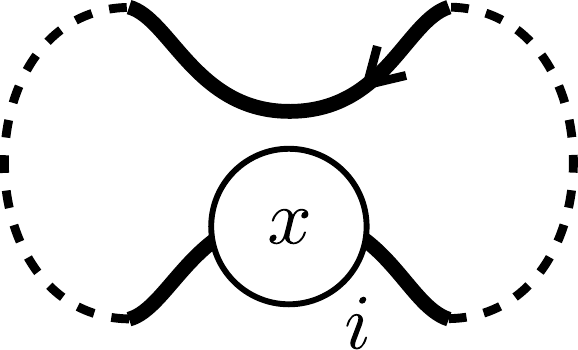}\end{array}
-\begin{array}{c}\includegraphics[scale=0.26]{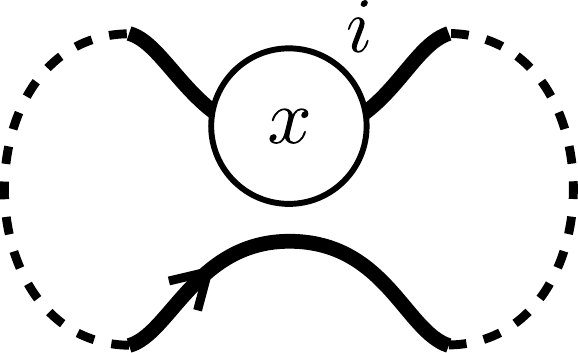}\end{array}\right),\\
&{\rm (B')} 
\begin{array}{c}\includegraphics[scale=0.26]{draws/xy_i+1i_ou.pdf}\end{array}- \begin{array}{c}\includegraphics[scale=0.26]{draws/xy_ii+1_ou.pdf}\end{array}=0,\\
&{\rm (B'')} \begin{array}{c}\includegraphics[scale=0.26]{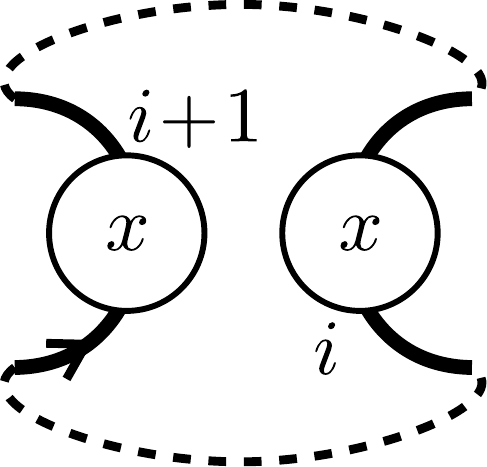}\end{array}-\begin{array}{c}\includegraphics[scale=0.26]{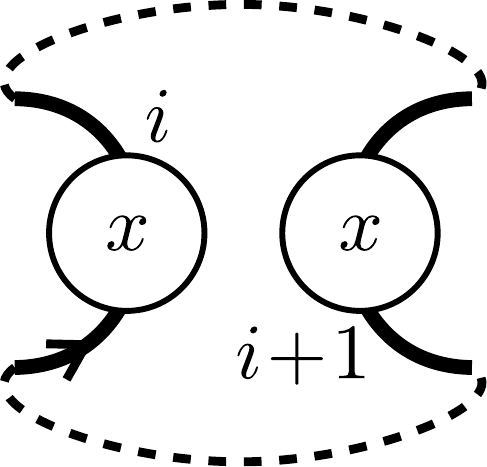}\end{array}=\hbar\left(\begin{array}{c}\includegraphics[scale=0.26]{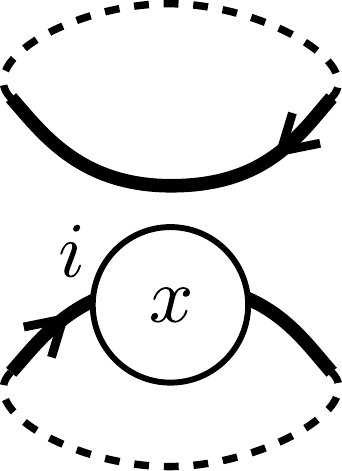}\end{array}
-\begin{array}{c}\includegraphics[scale=0.26]{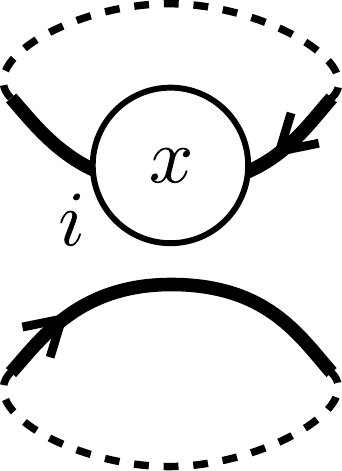}\end{array}\right)\\
&{\rm (C')} 
\begin{array}{c}\includegraphics[scale=0.25]{draws/bracelet_no_beads.pdf}\end{array}=2\rho.
\end{align*}
In (A'), (A''), (B'), and (B''), we require $x,y\in\mathsf{B}$; in (A') and (B'), we further require $x\neq y$. When a wire is an interval, the dotted connections indicate the identification of $J_{\rm out}$ with $J_{\rm in}$.
\end{definition}
The restrictions on the homology labels differ from those in the unmodified definition, so these relations cannot be extended linearly in the same way.

\begin{figure}	\centering\includegraphics[scale=0.15]{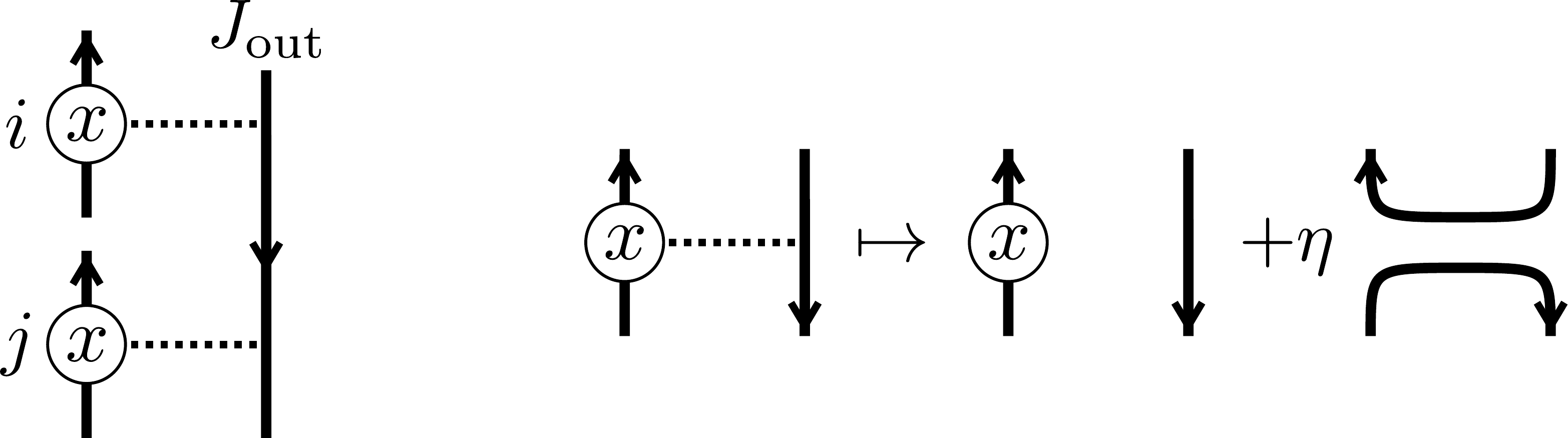}
\caption{Left: a matching compatible with heights, with $i>j$.
Right: the local replacement defining the left action.}
\label{fig:action_N} 
\end{figure}

Assume that $|J|=2$. Define the left action of $\widetilde{\cN_I}(\Sigma)$ on $\widetilde{\cN_I}(\Sigma,J)$ as follows.
By (H), it suffices to consider necklace ensembles whose beads are labeled by elements of $\mathsf{B}$. Suppose that $J_{\rm out}\subset d_k$ for some $k$, and let $x\in\mathsf{B}$ be the class represented by $\partial_k$.
Match the beads labeled by $x$ with points near the endpoint of the wire labeled by $J_{\rm out}$, preserving the height order of the beads; see the dotted lines on the left of Figure~\ref{fig:action_N}.
Define the left action by the linear combination obtained from the local replacements on the right of Figure~\ref{fig:action_N}. The fact that this defines an action follows from the product structure of $\widetilde{\cN_I}(\Sigma,J)$. The parameter $\eta\in\{h,\hbar\}$ is determined by preservation of component degree. If there are $n$ beads labeled by $x$, the resulting linear combination has $2^n$ terms.
The right action is defined by the same procedure with $J_{\rm out}$ replaced by $J_{\rm in}$.

Let $\widetilde{\cN_a}(\Sigma,J)$ be the $\bQ\llbracket h\hbar,\rho\hbar\rrbracket[\rho,h,\hbar]$-module obtained from Definition~\ref{def:modified_necklace} by replacing $h$ and $\hbar$ on the RHS's of (A'') and (B'') by $a$ and $\abar$, respectively, and replacing the RHS of (C') by $\dfrac{2\sinh(\rho\hbar)}{\hbar}$.
Replacing $h$ and $\hbar$ by $a$ and $\abar$ in the definitions of the left and right actions of $\widetilde{\mathrm{N}}(\Sigma,J)$ defines the corresponding actions on $\widetilde{\cN_a}(\Sigma,J)$. Let $\Gr^n\widetilde{\cN_I}(\Sigma,J)$ (resp. $\Gr^n\widetilde{\cN_a}(\Sigma,J)$) denote the $\bQ[\rho]$-submodule generated by
\begin{align*}
    &\{h^i\hbar^j N\mid N:\text{$J$-necklace ensemble},\ \#\{\text{beads}\}+\deg(h)+\deg(\hbar)=n\}\\
    &(\text{resp. }\{a^i\abar^j N\mid N:\text{$J$-necklace ensemble},\ \#\{\text{beads}\}+\deg(a)+\deg(\abar)=n\})
\end{align*}
in $\widetilde{\cN_I}(\Sigma,J)$ (resp. $\widetilde{\cN_a}(\Sigma,J)$). We call this degree the \textbf{modified finite-type degree}.

The map $\bQ\llbracket h\hbar,\rho\hbar\rrbracket[\rho,h,\hbar]\to\bQ\llbracket h\hbar,\rho\hbar\rrbracket[\rho,h,\hbar]$ extending
\begin{align}
h\mapsto a,\qquad \hbar\mapsto \abar,\qquad \rho\mapsto \rho_a:=\dfrac{\sinh{(\rho\hbar)}}{\hbar}\label{eq:h_to_a}
\end{align}
is a $\bQ$-algebra isomorphism. For later convenience, we describe the induced isomorphism in the direction $\widetilde{\cN_a}(\Sigma,J)\to\widetilde{\cN_I}(\Sigma,J)$.
\begin{lemma}\label{lem:a_h}
The map $\widetilde{\cN_a}(\Sigma,J)\to\widetilde{\cN_I}(\Sigma,J)$ extending
$$a \mapsto h,\qquad \abar \mapsto \hbar,\qquad \rho_a\mapsto \rho,
\qquad L\mapsto L$$ is an isomorphism of graded modules, where $L$ is a $J$-necklace ensemble. Moreover, it satisfies property $(\ast)$.
\end{lemma}

\begin{remark}
The module $\widetilde{\cN_I}(\Sigma,J)$ is introduced to connect with the preceding sections, whereas $\widetilde{\cN_a}(\Sigma,J)$ is used to study formality.
\end{remark}

Define $\psi_{\widetilde{\cN_I}}^J\colon\mathcal{D}_{\widetilde{\cN_I}}^J\to\widetilde{\cN_a}(\Sigma,J)$ as follows.
Given a necklace-like $J$-chord diagram, replace each strand by a wire and retain the $J$-labels at the endpoints of interval components.
Replace each chord by a bead labeled by the basis element of $\mathsf{B}$ corresponding to its pole. The order of the beads along each oriented wire agrees with the order of the chord endpoints along the corresponding strand.
Choose a total order to the beads so that the restriction to the beads obtained from the same pole is compatible with the opposite orientation of the pole.
%on the chord endpoints on the poles that is compatible with their heights in the $[0,1]$-factor of $\Sigma\times[0,1]$.
Thus, among endpoints on the same pole, a higher endpoint is larger in the total order.
The relative order of endpoints on different poles is arbitrary.
The endpoints of the wires inherit the $J$-labels of the corresponding strand endpoints.
This associates a necklace ensemble to the necklace-like $J$-chord diagram.
Extending $\bQ\llbracket h\hbar,\rho\hbar\rrbracket[\rho,h,\hbar]$-linearly gives $\psi_{\widetilde{\cN_I}}^J\colon\mathcal{D}_{\widetilde{\cN_I}}^J\to\widetilde{\cN_a}(\Sigma,J)$.

\begin{theorem}\label{thm:psi_N_isom}
The map $\psi_{\widetilde{\cN_I}}^J$ is a well-defined $\bQ\llbracket h\hbar,\rho\hbar\rrbracket[\rho,h,\hbar]$-module isomorphism satisfying property $(\ast)$.
\end{theorem}
\begin{proof}
Applying the necklace-like $4T$ relation gives
$$\begin{array}{c}\includegraphics[scale=0.15]{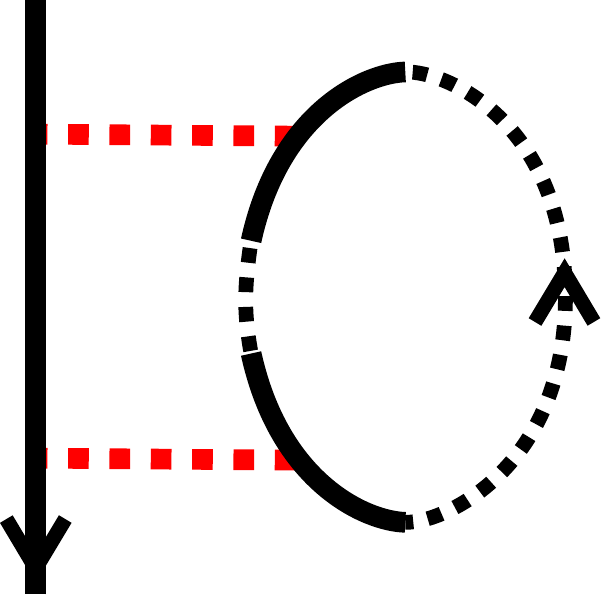}\end{array}-\begin{array}{c}\includegraphics[scale=0.15]{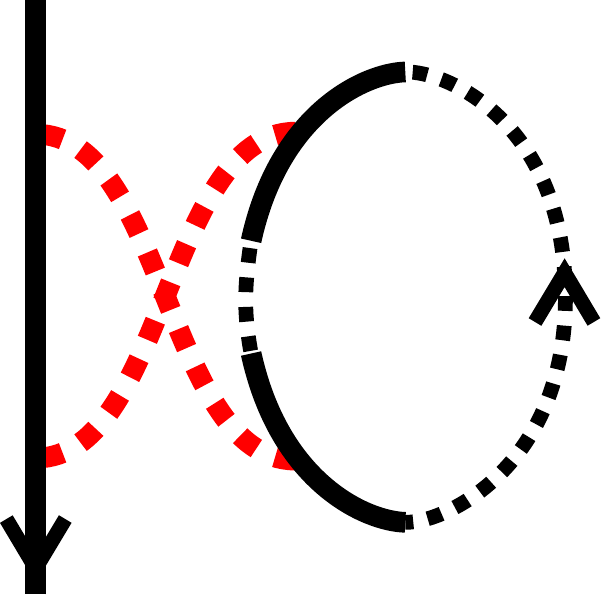}\end{array}
=\abar \begin{array}{c}\includegraphics[scale=0.15]{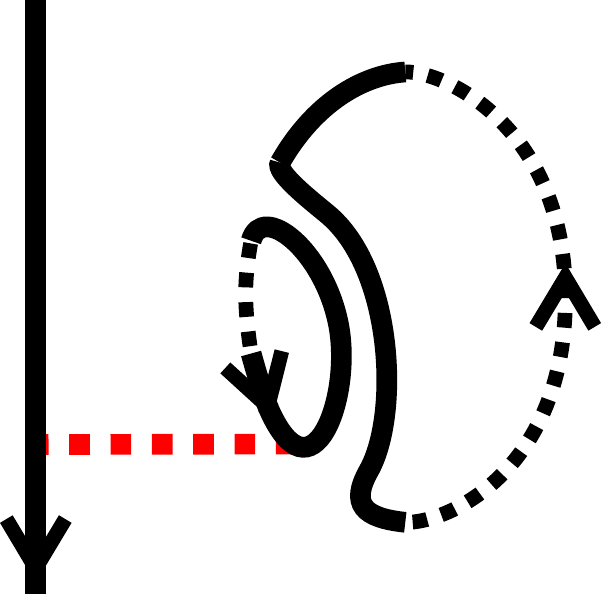}\end{array}
-\abar \begin{array}{c}\includegraphics[scale=0.15]{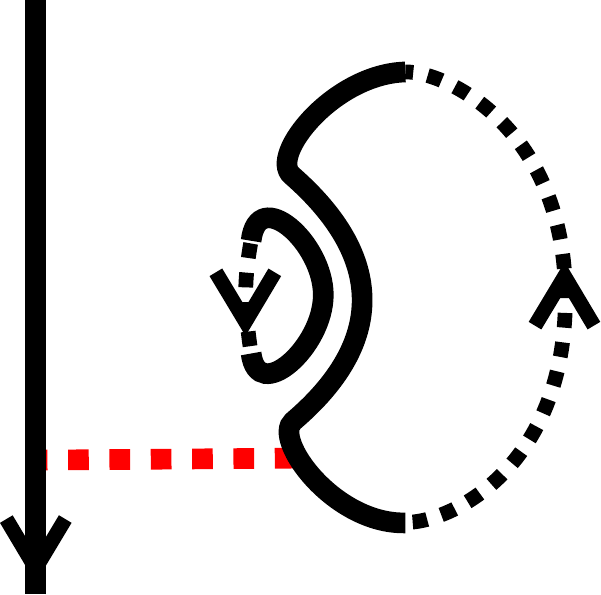}\end{array}$$
This corresponds to the relation (B'') in $\widehat{\widetilde{\cN_I}(\Sigma)}$.
The relation corresponding to (A'') follows in the same way.
Independence of the ordering of endpoints on different poles follows from (A') and (B'). The third relation in \eqref{eq:chord_diag_rel} corresponds to (C').
Hence, $\psi_{\widetilde{\cN_I}}^J$ is well-defined.

By construction, $\psi_{\widetilde{\cN_I}}^J$ is a surjective $\bQ\llbracket h\hbar,\rho\hbar\rrbracket[\rho,h,\hbar]$-module homomorphism.

We construct its inverse as follows.
Apply (H) to a $J$-necklace ensemble to expand each bead label in the basis $\mathsf{B}$.
In each term, replace every wire by a strand and every bead by a chord joining that strand to the pole corresponding to its label in $\mathsf{B}$. The heights of chord endpoints on poles are chosen to agree with the bead heights. The strand endpoints inherit the $J$-labels of the corresponding wire endpoints. Extend this construction $\bQ\llbracket h\hbar,\rho\hbar\rrbracket[\rho,h,\hbar]$-linearly.
This map is well-defined for the following reasons:
\begin{enumerate}
    \item[(C')] trivial wires correspond to trivial strands;
    \item[(A'), (B')] the corresponding chord diagrams coincide;
    \item[(A''), (B'')] the corresponding identities are the necklace-like $4T$ relations.
\end{enumerate}

By construction, this map is inverse to $\psi_{\widetilde{\cN_I}}^J$.
Therefore, $\psi_{\widetilde{\cN_I}}^J$ is a $\bQ\llbracket\rho\hbar,h\hbar\rrbracket[\rho,h,\hbar]$-module isomorphism.
Since the number of chords equals the number of beads, $\psi_{\widetilde{\cN_I}}^J$ preserves the modified finite-type degree in $a$ and $\abar$.
Property $(\ast)$ follows from the construction.
\end{proof}

\subsection{Modified degrees on the HOMFLYPT skein algebra and the modified quantized necklace algebra}\label{subsec:modified_filtration}
To show that the map $Z_{\widetilde{\cA}}^J$ induced by the Kontsevich integral preserves the filtration, we introduce two filtrations on $\widetilde{\cA}(\Sigma_{0,p+1},J)$ and prove that they agree.

Let $D$ be a disk, and identify
$\Sigma_{0,p+1}=D\setminus\text{($p$ open disks)}$.
Recall that $|J|=0,2$.

For $J=\{\bullet,\ast\}\subset\Sigma$, let $\pi_J^+$ be the group generated by regular homotopy classes of immersions, as defined in Section~\ref{subsec:module_str_HOMFLYPT}.
Let $\varepsilon\colon\bQ\pi_J^+\to\bQ$ be the augmentation map.

Let $\iota\colon\pi_J^+\to\widetilde{\cA}(\Sigma_{0,p+1},J)$ be the homomorphism that assigns to a regular homotopy class its descending lift to $\Sigma\times[0,1]$.

The detours of a tangle with boxes are \textbf{separable}
if they are contained in pairwise disjoint embedded $3$-balls in $D\times[0,1]$, with exactly one detour in each ball.
These balls may intersect poles and $s$-components.
We also call such a $J$-tangle with boxes \textbf{separable}.
The following $\bQ\llbracket\rho\hbar,h\hbar\rrbracket[\rho,h,\hbar]$-submodules define two descending filtrations:
\begin{align}
    &I^n\widetilde{\cA}(\Sigma_{0,p+1},J)={\rm span}\{h^i\hbar^j L\mid\text{$L$ is a separable tangle diagram with boxes},\ i+j+\#\{\text{boxes}\}\geq n\}\label{eq:filtration_separable}\\
    &\cI^n\widetilde{\cA}(\Sigma_{0,p+1},J)={\rm span}\{h^i\hbar^j L\mid\text{$L$ is a singular $J$-tangle},\ i+j+\#\{\text{singularities}\}\geq n\}\nonumber.
\end{align}
We abbreviate these submodules to $I^n$ and $\cI^n$, omitting $\widetilde{\cA}(\Sigma_{0,p+1},J)$.
We call $\{I^n\}_{n\in\bN}$ the
\textbf{modified finite-type filtration}. Its associated graded quotient $\Gr_I(\widetilde{\cA}(\Sigma,J))$ (resp. completed associated graded quotient $\widehat{\Gr_I}(\widehat{\widetilde{\cA}(\Sigma,J)})=\widehat{\Gr_I}(\widetilde{\cA}(\Sigma,J))=\prod_{n\in\bN}I^n/I^{n+1}$) is called the \textbf{modified associated graded quotient} of $\widetilde{\cA}(\Sigma,J)$ (resp. $\widehat{\widetilde{\cA}(\Sigma,J)}$).
\begin{theorem}\label{thm:two_filtrations}
The filtrations $\{I^n\}_{n\in\bN}$ and $\{\cI^n\}_{n\in\bN}$ coincide.
\end{theorem}
\begin{proof}
We first prove $I^n\supset\cI^n$.
Convert singular $J$-tangles into $J$-tangles with boxes using the following local correspondences:
\begin{align}
\begin{array}{c}\includegraphics[scale=0.15]{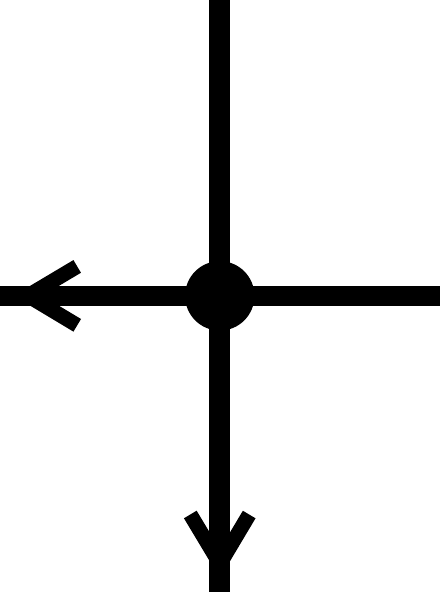}\end{array}
&=\begin{array}{c}\includegraphics[scale=0.15]{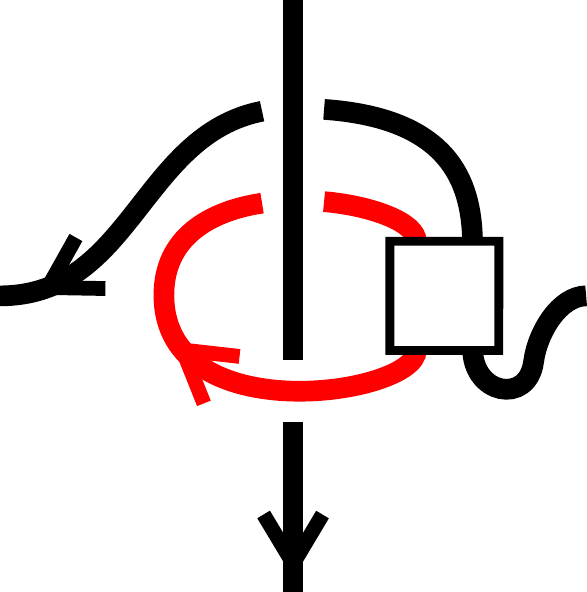}\end{array}\label{eq:pole}\\
\begin{array}{c}\includegraphics[scale=0.15]{draws/filtration_equiv_crossing03.pdf}\end{array}
&=(\exp(\rho\hbar)-1)
\begin{array}{c}\includegraphics[scale=0.15]{draws/filtration_equiv_crossing04.pdf}\end{array}\label{eq:framing}\\
\begin{array}{c}\includegraphics[scale=0.15]{draws/filtration_equiv_crossing05.pdf}\end{array}
&=\eta
\begin{array}{c}\includegraphics[scale=0.15]{draws/filtration_equiv_crossing06.pdf}\end{array}\qquad (\eta\in \{h,\hbar\})\label{eq:singular}
\end{align}
The vertical line in \eqref{eq:pole} represents a pole, and the orientation of the straight component agrees with that of the singular $J$-tangle.
The parameter $\eta\in\{h,\hbar\}$ in \eqref{eq:singular} is determined by preservation of component degree.
By construction, the resulting $J$-tangles with boxes lie in $I^n$.
Thus, $I^n\supset\cI^n$.

Next, we prove $I^n\subset\cI^n$.
Draw a $J$-tangle with boxes horizontally, as in the middle of Figure~\ref{fig:diagram_vertical_horizontal}.
Following a detour from its box back to the same box, label by $c_1,\dots,c_k$ the crossings at which it passes over a straight component or a pole.
Then the detour can be expressed using singularities as follows:
\begin{align}
\begin{array}{c}\includegraphics[scale=0.15]{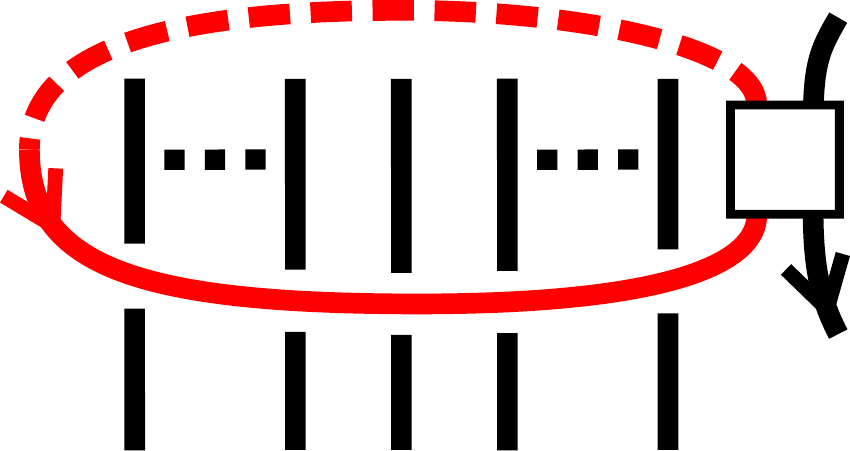}\end{array}=\sum_{i=1}^{k}\varepsilon(c_i)
\begin{array}{c}\includegraphics[scale=0.15]{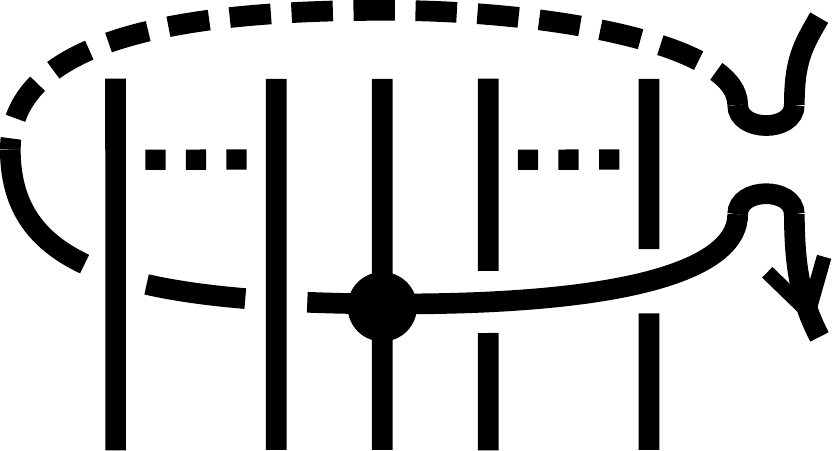}\end{array}+\begin{array}{c}\includegraphics[scale=0.15]{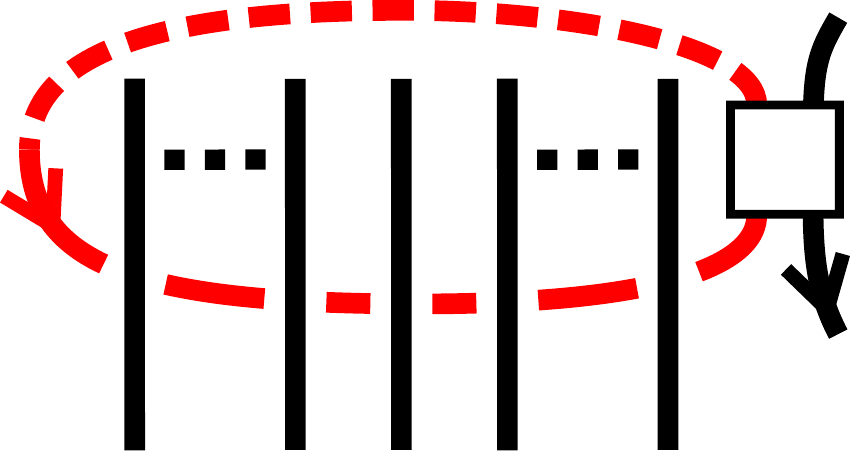}\end{array}\label{eq:filtration_equiv_alternative}
\end{align}
On the LHS of \eqref{eq:filtration_equiv_alternative}, the crossings appear from left to right in the order $c_1,\dots,c_k$, and $\varepsilon(c_i)\in\{-1,1\}$ is the sign of $c_i$.
Although \eqref{eq:filtration_equiv_alternative} depicts a detour that becomes trivial after the poles are removed, the same argument applies to a detour with more complicated entanglements.
The detour in the last term on the RHS of \eqref{eq:filtration_equiv_alternative} lies in an embedded $3$-ball in $\Sigma\times[0,1]$ containing no pole.
Consequently, this term belongs to $\hbar\widetilde{\mathcal{A}}(\Sigma,J)$. 
Therefore, the RHS lies in $\cI^n$, proving $I^n\subset\cI^n$.
\end{proof}

\begin{remark}The filtration $\{\mathscr{I}^n\widetilde{\cA}(\Sigma_{0,p+1},J)\}_{n\in \bN}$ defined by 
\begin{align*}
&\mathscr{I}^n\widetilde{\cA}(\Sigma_{0,p+1},J)\\
&={\rm span}\{h^i\hbar^j |\iota(x_1)|\dots|\iota(x_k)|y\mid x_i\in (\ker\varepsilon)^{\ell_i},\ y\in (\ker\varepsilon)^{\ell_0}\pi^+(\Sigma,\bullet,\ast),\ i+j+\ell_0+\dots+\ell_k\geq n\}
\end{align*}
can be shown to agree with $\{I^n\widetilde{\cA}(\Sigma_{0,p+1},J)\}_{n\in\bN}$. Replacing $i+j$ by $2i+2j$ in the definition of $\mathscr{I}^n\widetilde{\cA}(\Sigma_{0,p+1},J)$ yields the filtration $\{F^n\widetilde{\cA}(\Sigma_{0,p+1},J)\}_{n\in\bN}$.
Thus, $\{I^n\widetilde{\cA}(\Sigma_{0,p+1},J)\}_{n\in\bN}$ and $\{F^n\widetilde{\cA}(\Sigma_{0,p+1},J)\}_{n\in\bN}$ can be viewed as arising from different degree assignments to a bifiltration, 
which is given by 
\begin{align*}
&F^{n,m}\widetilde{\cA}(\Sigma_{0,p+1},J)=\\
&{\rm span}\{h^i\hbar^j |\iota(x_1)|\dots|\iota(x_k)|y\mid x_i\in (\ker\varepsilon)^{\ell_i},\ y\in (\ker\varepsilon)^{\ell_0}\pi^+(\Sigma,\bullet,\ast),\ i+j\geq n, \ell_0+\dots+\ell_k\geq m\}.
\end{align*}
From this viewpoint, Proposition~\ref{prop:filtrations_relations} below is immediate, and one obtains another proof of Corollary~\ref{cor:h_times_genus0}, concerning injectivity of multiplication by $h$ and $\hbar$.
\end{remark}

\begin{figure}	\centering\includegraphics[scale=0.15]{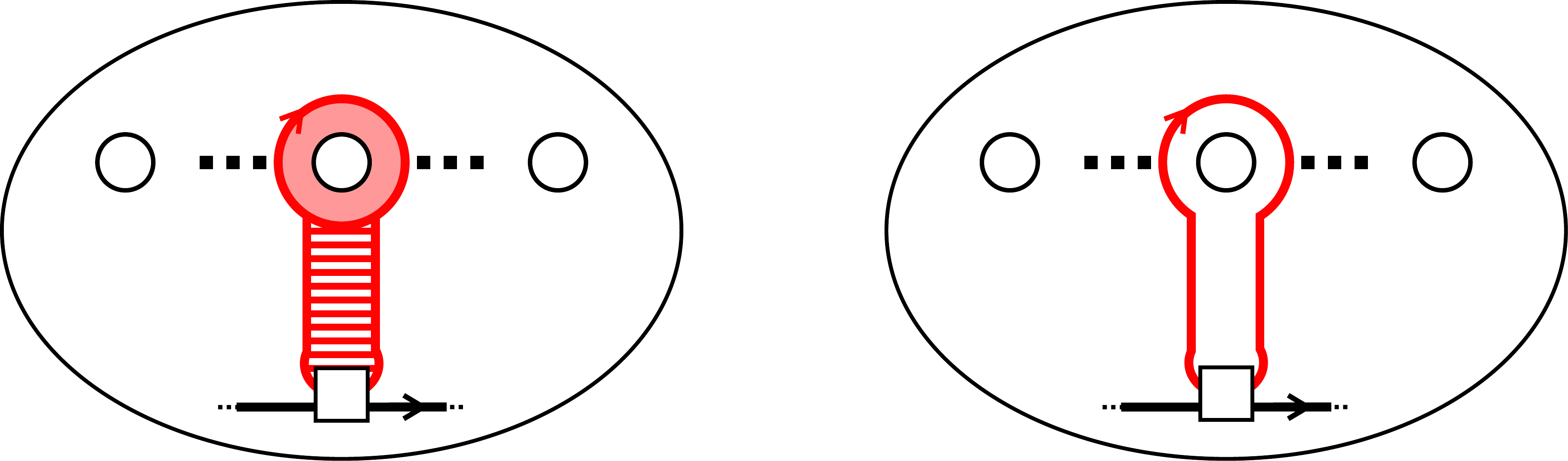}
\caption{Left: surfaces bounded by a peripheral loop and a trivial detour, together with a band joining them.
Right: the $d$-component forming the boundary of the resulting surface.}
\label{fig:bandsum_disk} 
\end{figure}

\begin{figure}	\centering\includegraphics[width=40mm]{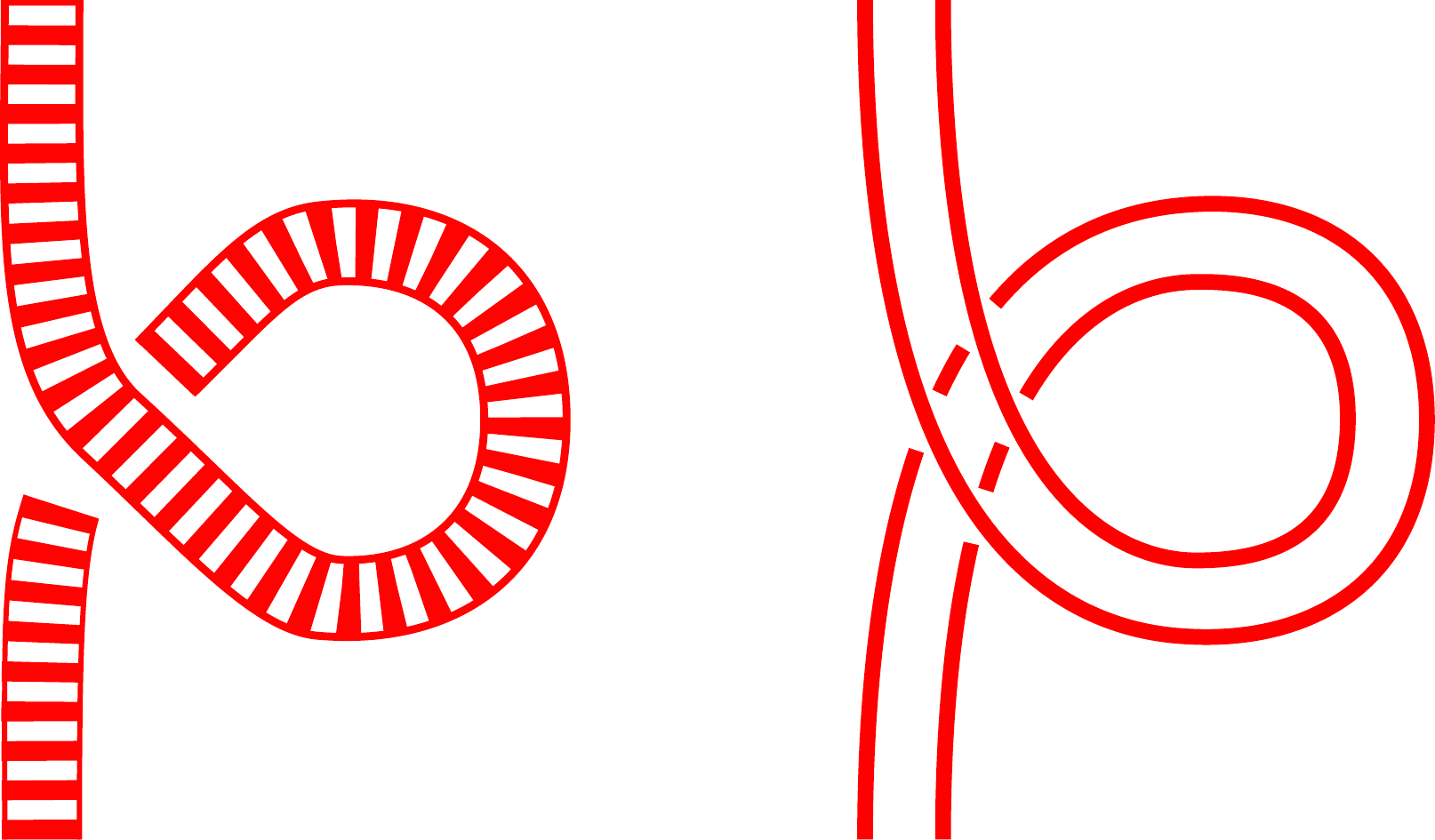}
\caption{A twist in the band and the corresponding twist in part of the $d$-component.}
\label{fig:twist_band_handrail} 
\end{figure}

The following proposition states that, for $\Sigma_{0,p+1}$, the finite-type filtration of Section~\ref{sec;HOMFLYPT} and the modified finite-type filtration induce the same topology.
\begin{prop}\label{prop:filtrations_relations}
The filtrations $\{I^n\widetilde{\cA}(\Sigma_{0,p+1},J)\}_{n\in\bN}$ and $\{F^n\widetilde{\cA}(\Sigma_{0,p+1},J)\}_{n\in\bN}$ satisfy
$$F^{2n}\widetilde{\cA}(\Sigma_{0,p+1},J)\subset I^n\widetilde{\cA}(\Sigma_{0,p+1},J)\subset F^n\widetilde{\cA}(\Sigma_{0,p+1},J)$$

In particular, $\{I^n\widetilde{\cA}(\Sigma_{0,p+1},J)\}_{n\in\bN}$ and $\{F^n\widetilde{\cA}(\Sigma_{0,p+1},J)\}_{n\in\bN}$ induce the same topology.
\end{prop}
\begin{proof}
The inclusion $I^n\widetilde{\cA}(\Sigma_{0,p+1},J)\subset F^n\widetilde{\cA}(\Sigma_{0,p+1},J)$ follows since separability is an additional condition on tangles with boxes, while the degrees of $h$ and $\hbar$ are doubled in the latter filtration.

We prove $F^{2n}\widetilde{\cA}(\Sigma_{0,p+1},J)\subset I^n\widetilde{\cA}(\Sigma_{0,p+1},J)$
by induction on $n$.
For $n=0$, this follows from
$$F^{0}\widetilde{\cA}(\Sigma_{0,p+1},J)= I^0\widetilde{\cA}(\Sigma_{0,p+1},J)=\widetilde{\cA}(\Sigma_{0,p+1},J)$$

Fix $k\in\bN$, and assume that
$F^{2n}\widetilde{\cA}(\Sigma_{0,p+1},J)\subset I^n\widetilde{\cA}(\Sigma_{0,p+1},J)$
holds for all $n\leq k$.
If a generator of $F^{2(k+1)}\widetilde{\cA}(\Sigma_{0,p+1},J)$ contains a factor $\eta\in\{h,\hbar\}$, it belongs to $\eta F^{2k}\widetilde{\cA}(\Sigma_{0,p+1},J)$. By the induction hypothesis,
\begin{align}
    \eta F^{2k}\widetilde{\cA}(\Sigma_{0,p+1},J)\subset 
    \eta I^{k}\widetilde{\cA}(\Sigma_{0,p+1},J)\subset I^{k+1}\widetilde{\cA}(\Sigma_{0,p+1},J)\label{eq:two_filtrations}
\end{align}
Thus, it suffices to consider generators of $F^{2(k+1)}\widetilde{\cA}(\Sigma_{0,p+1},J)$ containing no factors of $h$ or $\hbar$, namely $J$-tangle diagrams with $2(k+1)$ boxes.
As in \eqref{eq:surj_crossing}, applying (A) or (B) to switch a crossing involving a detour produces a correction term
in $\eta F^{2k}\widetilde{\cA}(\Sigma_{0,p+1},J)$.
Repeating these crossing changes expresses the original diagram as a linear combination of separable $J$-tangle diagrams with boxes and elements of $hF^{2k}\widetilde{\cA}(\Sigma_{0,p+1},J)+\hbar F^{2k}\widetilde{\cA}(\Sigma_{0,p+1},J)$.
The latter terms lie in $I^{k+1}\widetilde{\cA}(\Sigma_{0,p+1},J)$ by \eqref{eq:two_filtrations}.
The separable diagrams also lie in $I^{k+1}\widetilde{\cA}(\Sigma_{0,p+1},J)$, since they have at least $k+1$ boxes.
Hence, the original diagram lies in $I^{k+1}\widetilde{\cA}(\Sigma_{0,p+1},J)$. 
This proves the inclusion.
\end{proof}

For the modified quantized necklace algebra, $\Gr_I\Phi$, defined later, does not satisfy the properties of Lemma~\ref{lem:transparency}, but an analogous statement holds.
We first introduce some terminology.
\begin{definition}\label{def:detour_band_leaf}
Suppose a tangle diagram with boxes has a $d$-component obtained by taking a band sum of a trivial $d$-component with the peripheral loop $\partial_i$ around the $i$th puncture of $D_p$. Assume that the resulting $d$-component is isotopic in $D_p$ to a peripheral loop around the $i$th boundary component.
We call
the portion of the detour arising from the boundary of the band the \textbf{band part} of the $d$-component. We draw the band explicitly when we wish to emphasize this part.
The $d$-component bounds a once-punctured disk; see Figure~\ref{fig:bandsum_disk}.
A tangle diagram with boxes is \textbf{clasper-like} if every $d$-component has such a band part.
\end{definition}

For a $J$-tangle with boxes in $\Sigma\times[0,1]$, the clasper-like condition can be expressed as follows:
\begin{align}
&\text{Each detour bounds a disk embedded in $D\times[0,1]$,}\nonumber\\
&\text{meeting a pole in exactly one point and disjoint from all other parts.}\label{eq:condition}
\end{align}
Recall that $D$ denotes the closed disk to obtain $D_p$.
Under this condition, the $J$-tangle with boxes is separable, and each detour can be viewed as the boundary of a band part and a leaf part, as in Definition~\ref{def:detour_band_leaf}.

The following lemma will be used to construct $\widetilde{\cN_I}(\Sigma,J)\to\Gr\widetilde{\cA}(\Sigma,J)$ and prove that it is well-defined.
We impose an additional condition beyond separability, which will hold for the images of necklace ensembles.
\begin{lemma}\label{lem:transparency_modified}
Let $\Sigma=\Sigma_{0,p+1}$ and $|J|=0,2$. For clasper-like $J$-tangles with boxes satisfying \eqref{eq:condition}, the following operations preserve the clasper-like condition and give equal elements in the modified associated graded quotient $\Gr^n\widetilde{\cA}(\Sigma,J)$.
\begin{enumerate}
    \item Crossings between $s$-components may be switched, whether the branches belong to the same component or to different components. Such changes leave the element of $\Gr^n\widetilde{\cA}(\Sigma,J)$ unchanged. The framings of $s$-components may also be changed without changing the element of $\Gr^n\widetilde{\cA}(\Sigma,J)$.
    
    \item A crossing between the band part of a $d$-component and an $s$-component, the band part of another $d$-component, or a pole may be switched. These changes leave the element of $\Gr^n\widetilde{\cA}(\Sigma,J)$ unchanged.
    
    \item The following identity holds:
    \begin{align*}
        \begin{array}{c}\includegraphics[scale=0.15]{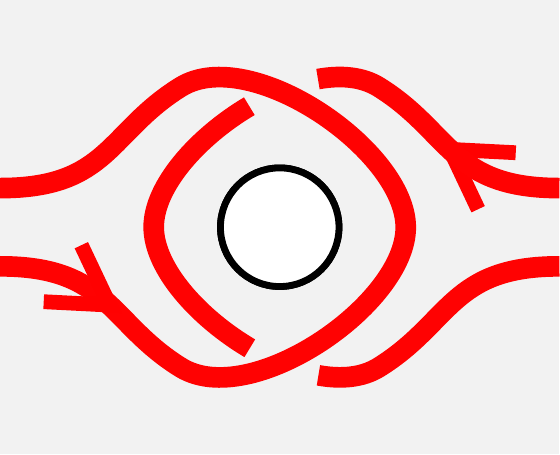}\end{array}
        -\begin{array}{c}\includegraphics[scale=0.15]{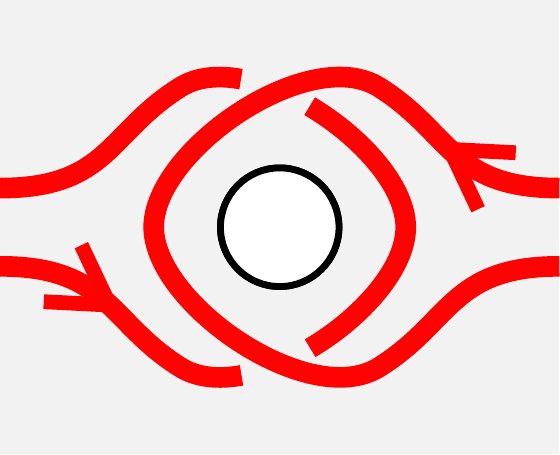}\end{array}
        =\eta\left(\begin{array}{c}\includegraphics[scale=0.15]{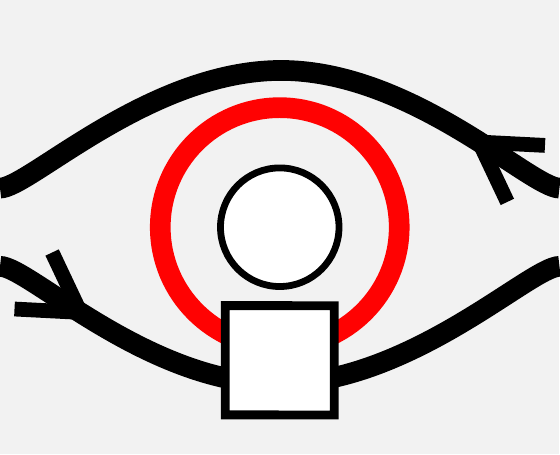}\end{array}
        -\begin{array}{c}\includegraphics[scale=0.15]{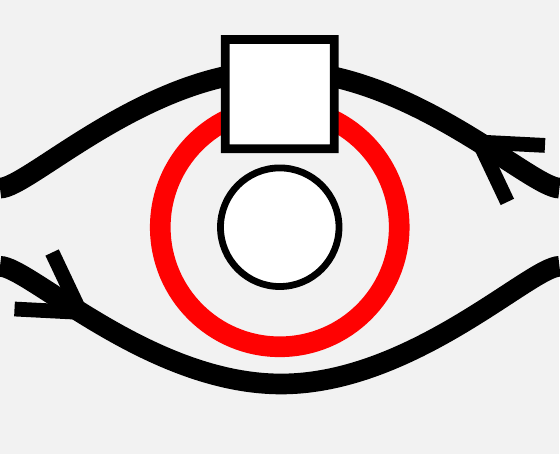}\end{array}\right).
    \end{align*}
    Here, $\eta\in\{h,\hbar\}$ is determined by preservation of component degree.
    
    \item The element is independent of twists in the band part. Thus, adding or removing twists in the band part of a $d$-component does not change the element of $\Gr^n\widetilde{\cA}(\Sigma,J)$.
    
    \item Crossings between an $s$-component and a pole may be switched.
\end{enumerate}
\end{lemma}
This differs substantially from Lemma~\ref{lem:transparency}, so the proof of Theorem~\ref{thm:welldefined} does not apply directly. In particular,
one cannot switch a crossing between an $s$-component and a part of a $d$-component outside its band part.

\begin{proof}
Although the filtration and degree differ from those in Section~\ref{sec;isom}, the proof of (1) is analogous to the proof of Lemma~\ref{lem:transparency}.

Part (2) follows from
\begin{align*}
    \begin{array}{c}\includegraphics[scale=0.15]{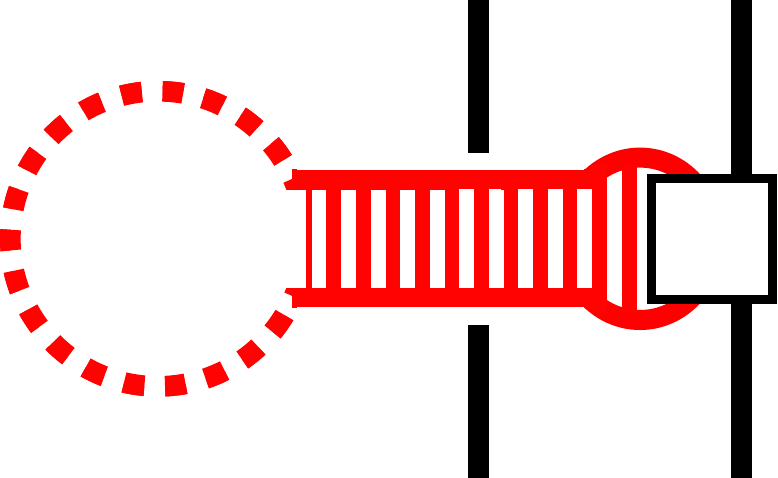}\end{array}-\begin{array}{c}\includegraphics[scale=0.15]{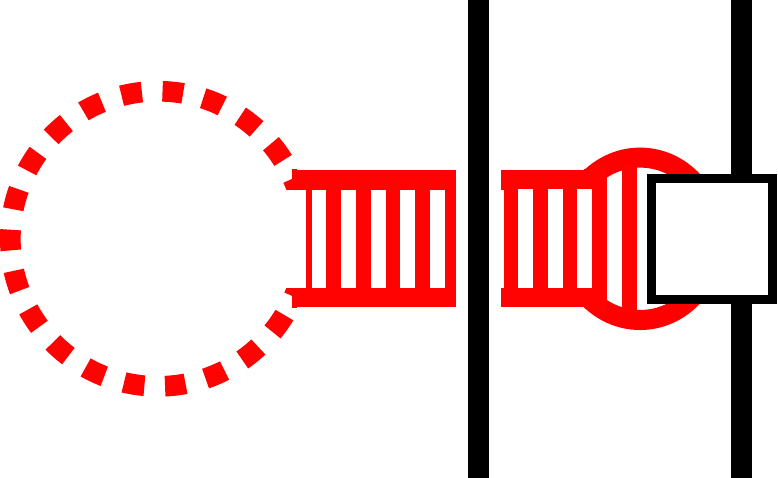}\end{array}
    &=\begin{array}{c}\includegraphics[scale=0.15]{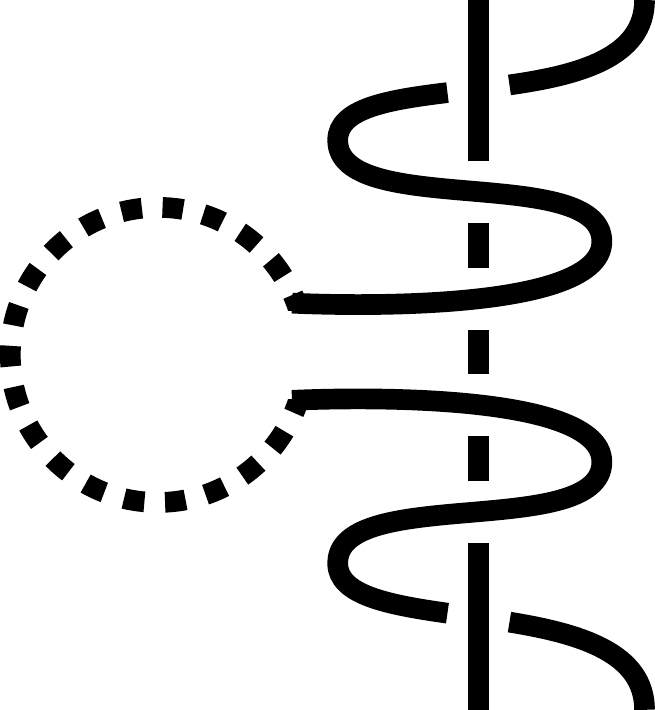}\end{array}-\begin{array}{c}\includegraphics[scale=0.15]{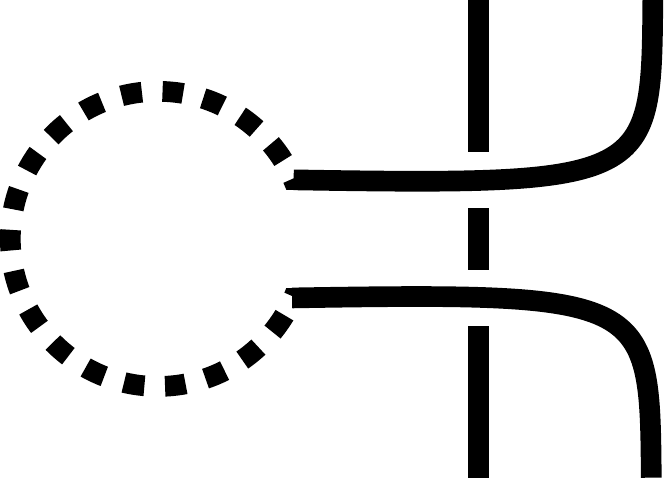}\end{array} \\
    &=\begin{array}{c}\includegraphics[scale=0.15]{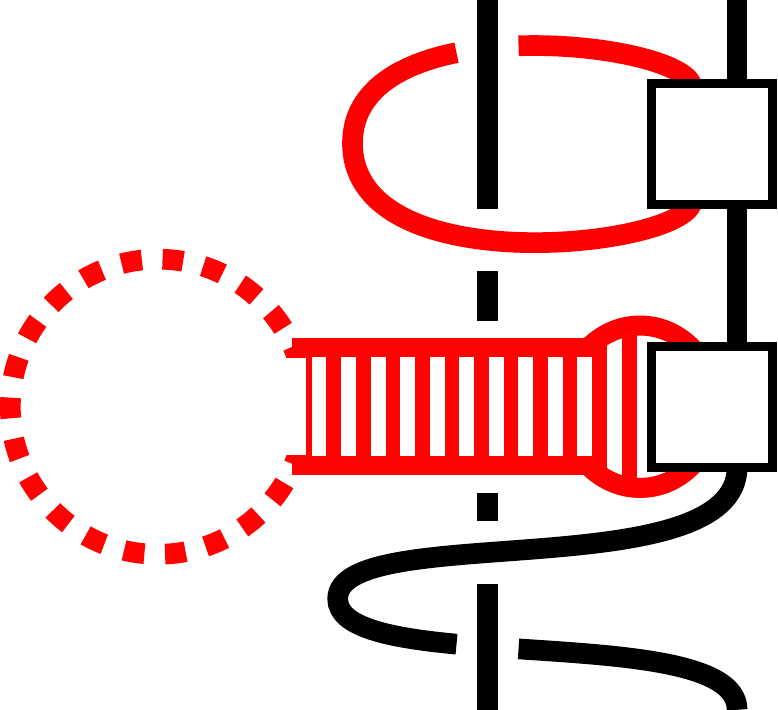}\end{array}-\begin{array}{c}\includegraphics[scale=0.15]{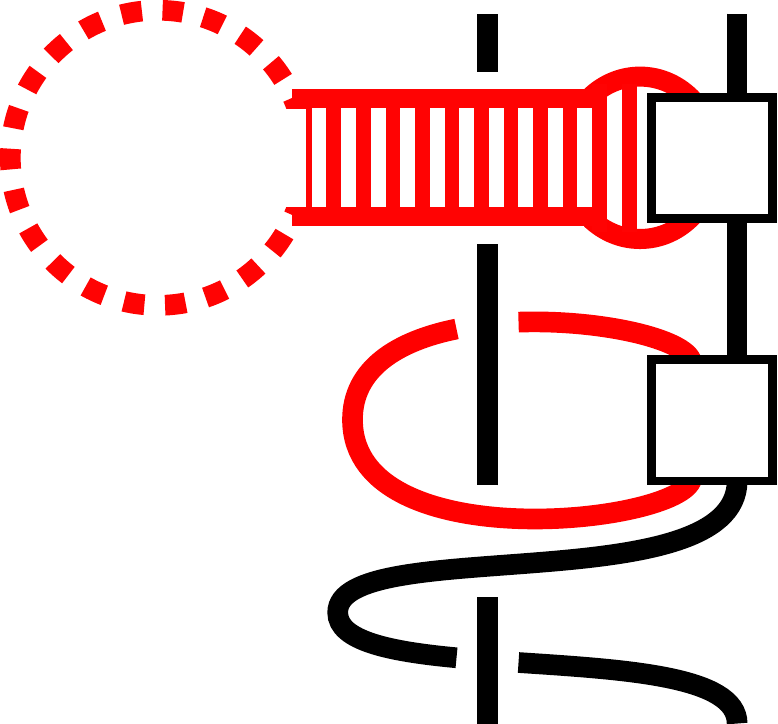}\end{array}=0 \mod I^2
\end{align*}
In each term, the left vertical line represents an $s$-component, the band part of a $d$-component, or a pole.

Part (3) follows from
    \begin{align*}
        \begin{array}{c}\includegraphics[scale=0.15]{draws/modified_filtration_lem3_01.pdf}\end{array}
        -\begin{array}{c}\includegraphics[scale=0.15]{draws/modified_filtration_lem3_02.pdf}\end{array}
        &=\begin{array}{c}\includegraphics[scale=0.15]{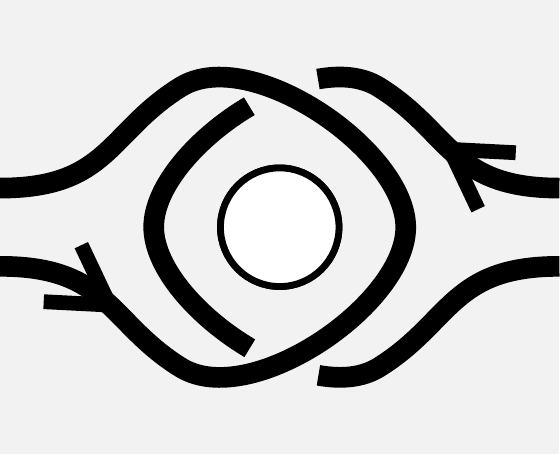}\end{array}
        -\begin{array}{c}\includegraphics[scale=0.15]{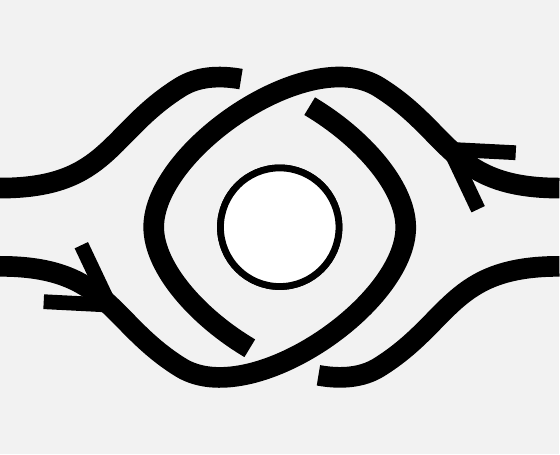}\end{array}
        +\begin{array}{c}\includegraphics[scale=0.15]{draws/modified_filtration_lem3_07.pdf}\end{array}
        -\begin{array}{c}\includegraphics[scale=0.15]{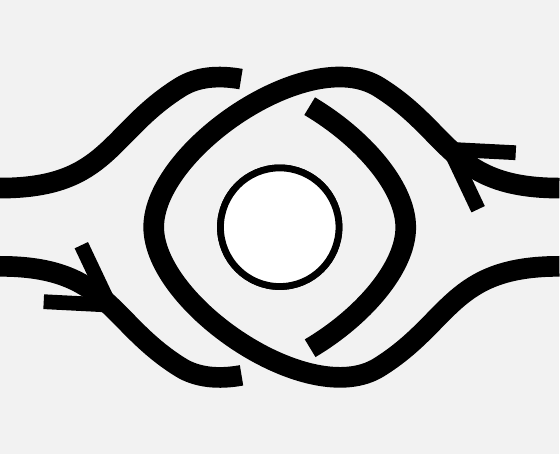}\end{array}\\
        &=\eta\left(\begin{array}{c}\includegraphics[scale=0.15]{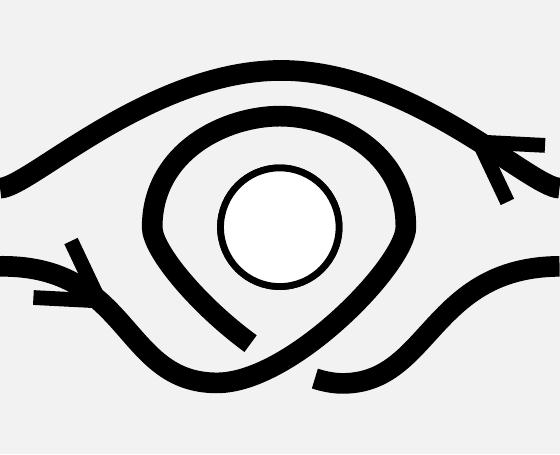}\end{array}
        -\begin{array}{c}\includegraphics[scale=0.15]{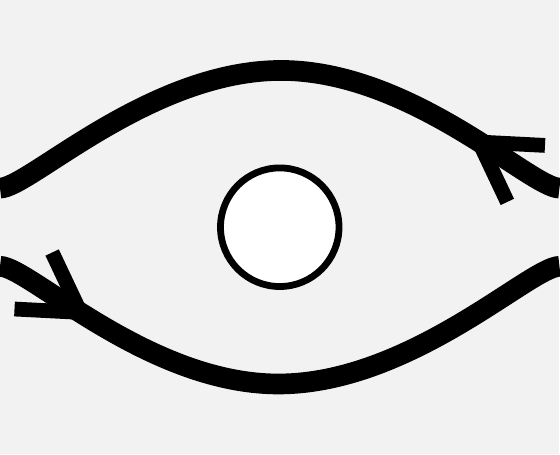}\end{array}
        +\begin{array}{c}\includegraphics[scale=0.15]{draws/modified_filtration_lem3_08.pdf}\end{array}
        -\begin{array}{c}\includegraphics[scale=0.15]{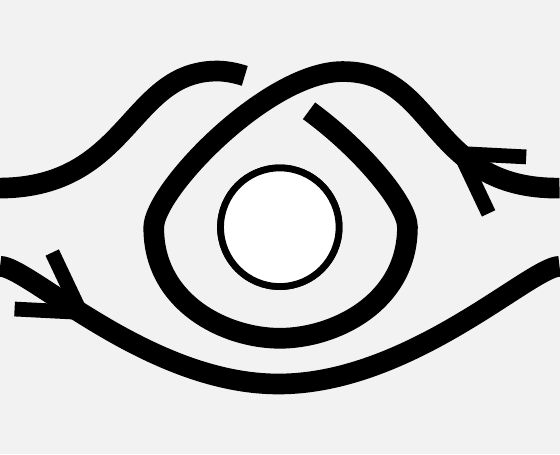}\end{array}\right)\\
        &=\eta\left(\begin{array}{c}\includegraphics[scale=0.15]{draws/modified_filtration_lem3_03.pdf}\end{array}
        -\begin{array}{c}\includegraphics[scale=0.15]{draws/modified_filtration_lem3_04.pdf}\end{array}\right).
    \end{align*}
    The second equality follows from the relation (A) or (B).
    The parameter $\eta\in\{h,\hbar\}$ is $h$ if the boxes of the detours on the left lie on different components of the virtual closure, and $\hbar$ if they lie on the same component.
    The third equality uses the framing independence in (1).

Part (4) follows from
\begin{align*}
\begin{array}{c}\includegraphics[scale=0.22]{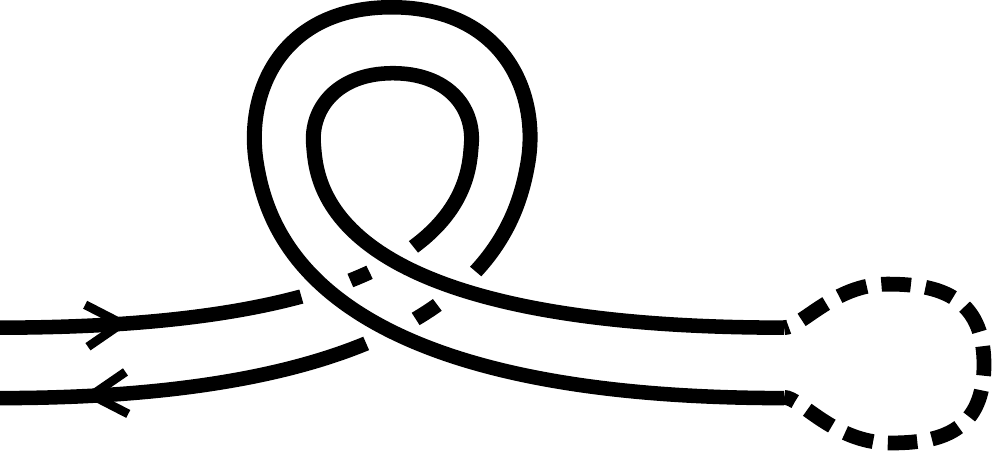}\end{array}
&=\begin{array}{c}\includegraphics[scale=0.22]{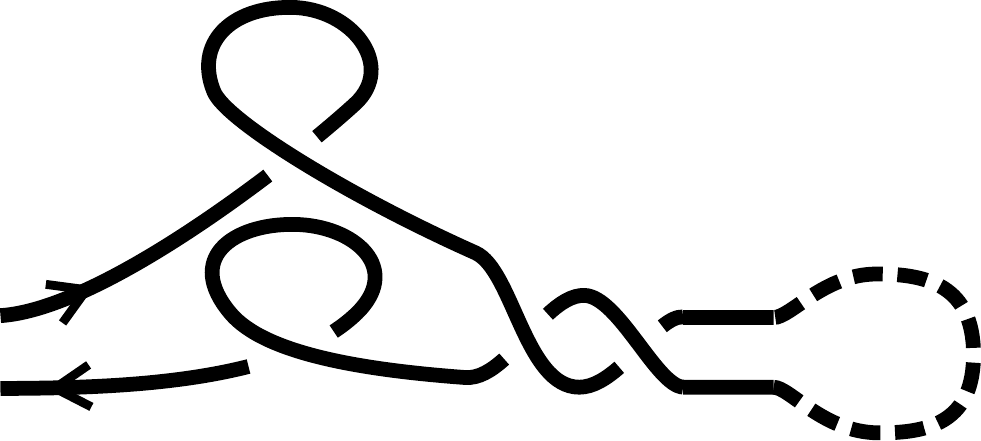}\end{array}
=\exp(2\rho\hbar)\begin{array}{c}\includegraphics[scale=0.22]{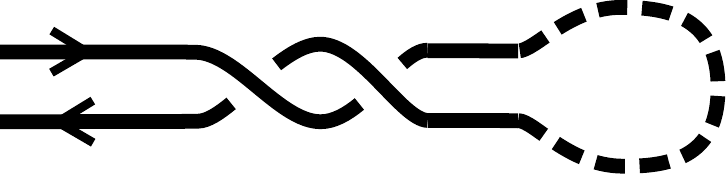}\end{array}\\
&=\exp(2\rho\hbar)\left(\begin{array}{c}\includegraphics[scale=0.22]{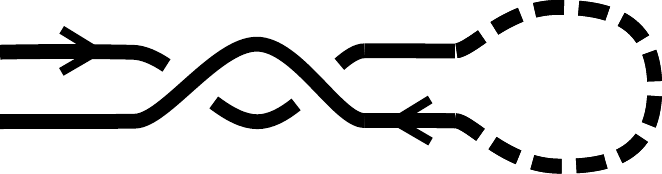}\end{array}+\hbar \begin{array}{c}\includegraphics[scale=0.22]{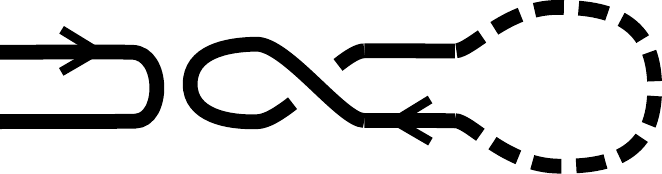}\end{array}\right)\\
&=\exp(2\rho\hbar)\left(\begin{array}{c}\includegraphics[scale=0.22]{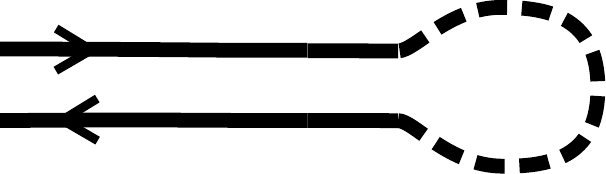}\end{array}-\hbar \exp(-\rho\hbar)\begin{array}{c}\includegraphics[scale=0.22]{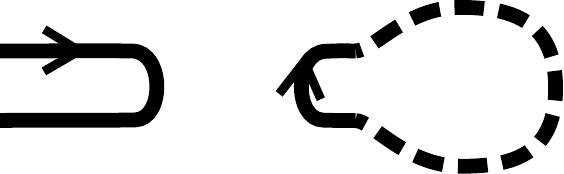}\end{array}\right)\\
&\equiv (1+2\rho \hbar)\begin{array}{c}\includegraphics[scale=0.22]{draws/genus0_self_twist06.pdf}\end{array}-\hbar \begin{array}{c}\includegraphics[scale=0.22]{draws/genus0_self_twist07.pdf}\end{array}\mod \hbar^2\\
&\equiv \begin{array}{c}\includegraphics[scale=0.22]{draws/genus0_self_twist06.pdf}\end{array} 
+ \hbar\left(\begin{array}{c}\includegraphics[scale=0.22]{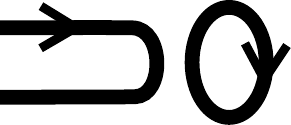}\end{array}-\begin{array}{c}\includegraphics[scale=0.22]{draws/genus0_self_twist07.pdf}\end{array}\right) \mod I^2\\
&= \begin{array}{c}\includegraphics[scale=0.22]{draws/genus0_self_twist06.pdf}\end{array} 
+ \hbar \begin{array}{c}\includegraphics[scale=0.22]{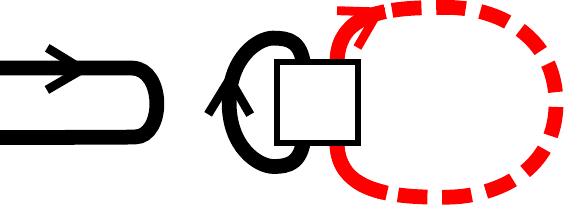}\end{array}
\equiv  
\begin{array}{c}\includegraphics[scale=0.22]{draws/genus0_self_twist06.pdf}\end{array} \mod I^2.
\end{align*}

Part (5) follows from
$$\begin{array}{c}\includegraphics[scale=0.15]{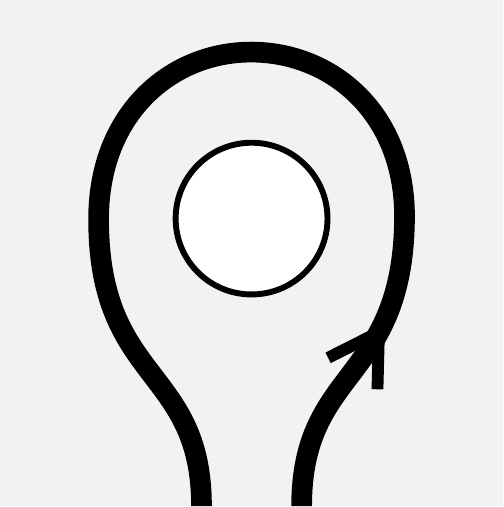}\end{array} -\begin{array}{c}\includegraphics[scale=0.15]{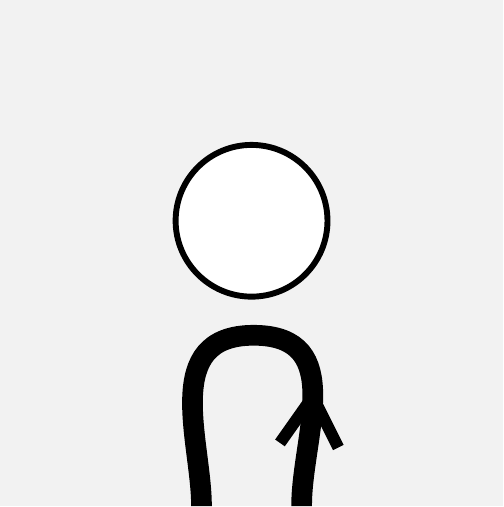}\end{array} =\begin{array}{c}\includegraphics[scale=0.15]{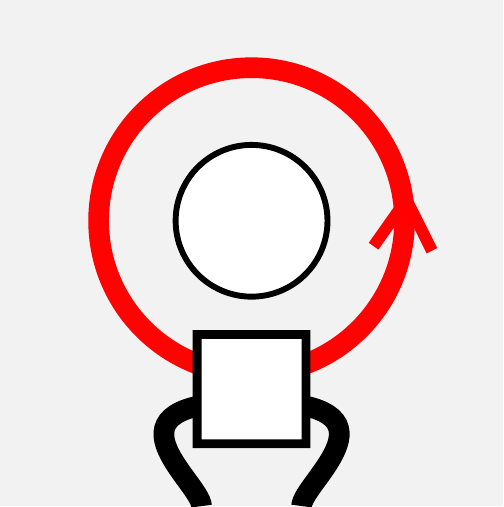}\end{array}.$$
\end{proof}

Define $\Gr_I\Phi\colon\widetilde{\cN_I}(\Sigma,J)\to\Gr_I\widetilde{\cA}(\Sigma,J)$ as follows. Let $\mathsf{B}$ be the basis of $H_1(\Sigma_{0,p+1})$ represented by the peripheral loops $\partial_i\ (i=1,\dots,p)$ on the right of Figure~\ref{fig:disk_punctured}. By applying (H), a necklace ensemble can be expanded in terms of $\mathsf{B}$. Therefore, it suffices to define the map on ensembles whose bead labels belong to $\mathsf{B}$ and extend linearly. Assume that every bead is labeled by an element of $\mathsf{B}$.
Draw a $J$-tangle diagram on $\Sigma$ with one component for each wire, preserving the orientations and endpoint labels of interval wires.
The diagram may have crossings, and its overcrossing and undercrossing data may be chosen arbitrarily.
Choose points $\ast$ corresponding to the beads so that their order along each oriented tangle component agrees with the order of the beads along the corresponding oriented wire.
For a bead labeled by $\partial_i$, replace a neighborhood of the corresponding point $\ast$ by a box and draw a detour representing $\partial_i$ so that the resulting diagram is clasper-like.
Do this for all beads, choosing crossings between detours according to bead heights: the detour corresponding to the higher bead passes over the other detour. All remaining crossing data may be chosen arbitrarily.
This gives a clasper-like $J$-tangle diagram. Extending linearly yields $\Gr_I\Phi$, which is a well-defined module homomorphism by Theorem~\ref{thm:isom_modified} below.

Compose $\psi_{\widetilde{\cN_I}}^J\circ\psi_{\mathcal{D}}^J\circ\Gr Z_{\widetilde{\cA}}^J\colon\widehat{\Gr}\widetilde{\cA}(\Sigma,J)\to\widehat{\widetilde{\cN_a}(\Sigma,J)}$ with $\widehat{\widetilde{\cN_a}(\Sigma,J)}\to\widehat{\widetilde{\cN_I}(\Sigma,J)}$, and denote the resulting map by
$$
\Psi\colon  \widehat{\Gr_I}\widetilde{\cA}(\Sigma,J)\to \mathcal{D}_{\widetilde{\cA}}^J\to \mathcal{D}_{\widetilde{\cN_I}}^J \to\widehat{\widetilde{\cN_a}(\Sigma,J)}\to \widehat{\widetilde{\cN_I}(\Sigma,J)}
$$

In the following theorem, we show that $\Psi$ is inverse to $\Gr_I\Phi$.

\begin{theorem}\label{thm:isom_modified}
Let $\Sigma=\Sigma_{0,p+1}$, and let $J\subset\partial\Sigma_{0,p+1}$ with $|J|=0,2$.
The modified quantized necklace Lie module $\widetilde{\cN_I}(\Sigma,J)$ and the modified associated graded quotient $\Gr_I\widetilde{\cA}(\Sigma,J)$ of the HOMFLYPT skein module satisfy the following.
\begin{enumerate}
    \item The map $\Gr_I\Phi\colon\widehat{\widetilde{\cN_I}(\Sigma,J)}\to\widehat{\Gr_I}\widetilde{\cA}(\Sigma,J)$ is a well-defined $\bQ\llbracket\rho\hbar,h\hbar\rrbracket[\rho,h,\hbar]$-module isomorphism.
Thus, the modified quantized necklace module $\widehat{\widetilde{\cN_I}(\Sigma,J)}$ and $\widehat{\Gr_I}\widetilde{\cA}(\Sigma,J)$ are isomorphic as $\bQ[\rho,h,\hbar]$-modules. Moreover, this isomorphism satisfies property $(\ast)$.

\item The map $\Gr Z_{\widetilde{\cA}}^J$ described in Lemma~\ref{lem:Kontsevich_descend} is an isomorphism. Its composition with the diagram identifications and coefficient change above is $\Psi=(\Gr_I\Phi)^{-1}$.
\end{enumerate}
\end{theorem}
\begin{proof}
Extending $\bQ\llbracket\rho\hbar,h\hbar\rrbracket[\rho,h,\hbar]$-linearly and applying Lemma~\ref{lem:transparency_modified} gives a well-defined module homomorphism. In particular, Lemma~\ref{lem:transparency_modified}(3) gives the relation for interchanging heights.
The argument of Lemma~\ref{lem:alghom} shows that this is an algebra homomorphism, and the argument of Lemma~\ref{lem:surj} proves surjectivity. Relations analogous to \eqref{eq:homologous} hold in $\Gr_I\widetilde{\cA}(\Sigma,J)$, from which preservation of the left and right actions follows.

To prove injectivity of $\Gr_I\Phi$,
we show that $\Psi\circ\Gr_I\Phi$ is the identity.
The image of any necklace ensemble under $\Gr_I\Phi$ is a linear combination of clasper-like tangle diagrams with boxes.
At the pole encircled by a clasper-like detour, use the following correspondence to redraw the diagram as a singular $J$-tangle diagram:
$$
\begin{array}{c}\includegraphics[scale=0.15]{draws/filtration_equiv_crossing02.pdf}\end{array}
=\begin{array}{c}\includegraphics[scale=0.15]{draws/filtration_equiv_crossing01.pdf}\end{array}. 
$$
Each box is replaced by one singularity, and all resulting singularities occur between a pole and a strand.
By Lemma~\ref{lem:transparency_modified}, the image depends only on the orders in which these singularities occur along the poles and along the strands.
Applying $\Gr Z_{\widetilde{\cA}}^J$ to the singular $J$-tangle diagram gives the local correspondences of Lemma~\ref{lem:Kontsevich_descend}(2).
The image is necklace-like, Hence, belongs to $\mathcal{D}_{\widetilde{N}}$.
Applying $\psi_{\widetilde{N}}^J$ to these necklace-like chord diagrams recovers the original necklace ensemble.
Although the coefficient ring changes during this construction, the coefficient of the ensemble remains $1$. Indeed, $\widehat{\widetilde{\cN_a}(\Sigma,J)}\to\widehat{\widetilde{\cN_I}(\Sigma,J)}$ changes coefficients and sends $1$ to $1$, so it leaves the necklace ensemble unchanged.

Now consider the product of a necklace ensemble with an element of $\bQ[\rho,h,\hbar]$. The image of the coefficient under $\Gr Z_{\widetilde{\cA}}^J$ is given by \eqref{eq:h_to_a}.
The maps $\psi_{\mathcal{D}}^J$ and $\psi_{\widetilde{N}}^J$ leave coefficients unchanged, while $\widehat{\widetilde{\cN_a}(\Sigma,J)}\to\widehat{\widetilde{\cN_I}(\Sigma,J)}$ changes them as in Lemma~\ref{lem:a_h}.
Thus, the original coefficient is recovered.
This proves that the $\bQ[\rho,h,\hbar]$-module homomorphism $\Gr_I\Phi$ is injective.

In particular, $\Psi$ is inverse to $\Gr_I\Phi$.
\end{proof}

Theorem~\ref{thm:isom_modified} implies the following.
\begin{cor}\label{cor:h_times_genus0}
Let $(\Sigma,J)$ be a marked surface with $\Sigma=\Sigma_{0,p+1}$. Then the following hold.
\begin{enumerate}
    \item The associated graded quotient $\Gr_I\widetilde{\cA}(\Sigma,J)$ and the modified quantized necklace Lie module $\widetilde{\cN_I}(\Sigma,J)$ are isomorphic as graded $\bQ[\rho,h,\hbar]$-modules, and the isomorphism satisfies property $(\ast)$.
    
    \item The following $\bQ[\rho]$-module homomorphisms on the modified associated graded quotient $\Gr_I^n\widetilde{\cA}(\Sigma,J)$ are injective:% multiplication by $h$ and by $\hbar$
\begin{align*}
    \Gr_{I}^n\widetilde{\cA}(\Sigma,J)\to \Gr^{n+1}_{I}\widetilde{\cA}(\Sigma,J),\quad x\mapsto hx\\ 
    \Gr_{I}^n\widetilde{\cA}(\Sigma,J)\to \Gr_{I}^{n+1}\widetilde{\cA}(\Sigma,J),\quad x\mapsto \hbar x
\end{align*}

\end{enumerate}
\end{cor}

In the commutative diagram \eqref{eq:commutative_diag_expansion}, the map $\psi_{\widetilde{\cA}}^J:=(\Gr Z_{\widetilde{\cA}}^J)^{-1}$ is an isomorphism. Hence, we obtain the following.
\begin{theorem}\label{thm:formality02}
For the filtration $\{I^n\}_{n\in\bN}$ on $\widetilde{\cA}(\Sigma,J)$ defined in \eqref{eq:filtration_separable},
formality holds in $\widehat{\Gr_I}\widetilde{\cA}(\Sigma,J)$ with respect to the product, composition, and left and right actions. In particular, it also holds with respect to the Lie action $\sigma$.
\end{theorem}

\begin{remark}\label{rmk:notorsion_hhbar_completion_genus0}
\begin{enumerate}
    \item As in Remark~\ref{rmk:notorsion_hhbar_completion}, the completed modified skein algebra $\widehat{\widetilde{\cA}(\Sigma)}=\plim[n]\widetilde{\cA}(\Sigma)/I^n\widetilde{\cA}(\Sigma)$ has no zero divisors. In particular, neither $h$ nor $\hbar$ is a zero divisor.

    \item For $\Sigma_{0,b}$, \cite[Proposition 4.20]{BN+25} proves injectivity of multiplication by $h$ and $\hbar$ only for the maps $\bQ\cT/\bQ\cT^1\to\bQ\cT/\bQ\cT^2$. Remark~\ref{rmk:notorsion_hhbar_completion_genus0} and Corollary~\ref{cor:isom_completions} give the analogous statement for $\bQ\cT/\bQ\cT^k\to\bQ\cT/\bQ\cT^{k+1}$ for all $k\geq1$.
\end{enumerate}
\end{remark}

An isomorphism between $\widehat{\widetilde{\cA}(\Sigma)}$, before passing to the associated graded quotient, and the target of the Kontsevich integral was constructed in \cite{BN+25}. Combining this with the preceding results gives the following corollary.
\begin{cor}\label{cor:isom_completions}
For $\Sigma=\Sigma_{0,p+1}$ with $p\in\bN$, the completion $\widehat{\widetilde{\cN_I}(\Sigma)}$ of the modified quantized necklace algebra, the associated graded quotient $\widehat{\Gr_I}\widetilde{\cA}(\Sigma)$ of the HOMFLYPT skein algebra, and its completion $\widehat{\widetilde{\cA}(\Sigma)}$ are related by semilinear $\bQ\llbracket h,\hbar\rrbracket[\rho]$-algebra isomorphisms.
\end{cor}

\subsection{The height involution and formality}\label{subsec:height_involution}
Define the \textbf{height involution} $\widetilde{\cA}(\Sigma,J)\to\widetilde{\cA}(\Sigma,J)$ by
$$\begin{array}{c}\includegraphics[scale=0.15]{draws/crossing_positive.pdf}\end{array}\mapsto \begin{array}{c}\includegraphics[scale=0.15]{draws/crossing_negative.pdf}\end{array},\qquad 
\begin{array}{c}\includegraphics[scale=0.15]{draws/crossing_negative.pdf}\end{array}\mapsto \begin{array}{c}\includegraphics[scale=0.15]{draws/crossing_positive.pdf}\end{array},\qquad
h\mapsto -h,\qquad \hbar\mapsto -\hbar,\qquad \rho\mapsto\rho$$

This operation preserves the modified finite-type filtration and therefore induces a height involution $\Gr_I\widetilde{\cA}(\Sigma,J)\to\Gr_I\widetilde{\cA}(\Sigma,J)$.
Under the isomorphism of Theorem~\ref{thm:isom_modified}, it induces the \textbf{height involution} $\widetilde{\cN_I}(\Sigma,J)\to\widetilde{\cN_I}(\Sigma,J)$ given by
\begin{align*}
&\text{necklace ensemble}\mapsto\text{the same ensemble with the bead height order reversed},\\
&h\mapsto -h,\qquad \hbar\mapsto -\hbar,\qquad \rho\mapsto\rho
\end{align*}
This operation leaves the beads and the connectivity of the necklace ensemble unchanged.
It also induces a height involution on $\mathcal{D}_{\widetilde{\cA}}^J$, which multiplies a $J$-chord diagram by $-1$ for each strand--strand chord and reverses the orientations of the poles.
By abuse of notation, let $\tau$ also denote the height involution on each of these modules and algebras.

To study formality, we work with the completion $\widehat{\widetilde{\cA}(\Sigma,J)}$. Assume that $|J|=2$ and that $J$ lies on the boundary component $\partial_0$ of $\Sigma$; see the LHS of Figure~\ref{fig:disk_punctured}. Fix a one-component $J$-tangle diagram $x$.
By the skein relations (A) and (B), every term of $x-\tau(x)$ has positive $\hbar$-degree. Therefore, theorem~\ref{thm:h_times} and Remark~\ref{rmk:notorsion_hhbar_completion} allow us to consider $\frac{1}{\hbar}(x-\tau(x))$.
Each term in the constant part of $\frac{1}{\hbar}(x-\tau(x))$ has two components: one loop and one one-component $J$-tangle.
Note that the operation $x\mapsto\frac{1}{\hbar}(x-\tau(x))$ commutes with closure.

For $\ell\in\pi:=\pi_1(\Sigma,\ast)$, let $x_\ell$ be the $J$-tangle diagram obtained by choosing crossing data at each double point so that the diagram is descending along its orientation. Here, $\Jout=\ast$ and $\Jin$ is obtained by moving $\ast$ slightly along the orientation of the boundary.
Forgetting the crossing data in $\frac{1}{\hbar}(x_\ell-\tau(x_\ell))$ and substituting $\rho=|\textbf{1}|/2$ allows us to regard it as an element of $\bQ|\pi|\otimes\bQ\pi$.
One checks that this agrees with the Kawazumi--Kuno coaction $\mu^+\colon\bQ\pi^+\to\bQ|\pi|\otimes\bQ\pi$. Our convention for $\mu^+$ differs by a minus sign from \cite[Eq.(36)]{AKKN18}.
Since $\frac{1}{\hbar}(x_\ell-\tau(x_\ell))$ commutes with closure, closing the $J$-tangle component, antisymmetrizing, and forgetting crossings recovers the Turaev cobracket of the Goldman--Turaev Lie bialgebra \cite[Section 5.2]{BN+25}. More precisely, if $x_1$ is the component obtained by closure and $x_2$ is the other component, the antisymmetrization is $x_1\otimes x_2-x_2\otimes x_1$.
This viewpoint is used in \cite{BN+25} to prove formality of the Turaev cobracket.
Here, we retain the full quantum expression $\frac{1}{\hbar}(x-\tau(x))$, including its higher-order terms. Thus, formality for the height involution can be viewed as a quantization of the formality result in \cite{BN+25}.

To determine how the height involution behaves under
the Kontsevich integral $Z_{\widetilde{\cA}}^J\colon\widehat{\widetilde{\mathcal{A}}(\Sigma,J)}\xrightarrow{\cong}\mathcal{D}_{\widetilde{\cA}}^J$, we describe it on horizontal $J$-tangle diagrams.

The \textbf{horizontal mirror image} of a horizontal $J$-tangle is obtained by reversing the orientation of the second factor of $\Sigma\times[0,1]$, while retaining the orientations of the poles.
As in \cite[Lemma~5.19]{BN+25}, $Z_{\widetilde{\cA}}$ has formality with respect to the horizontal mirror image.\footnote{\cite{BN+25} uses the term ``(vertical) flip.''}

\begin{figure}	\centering\includegraphics[scale=0.15]{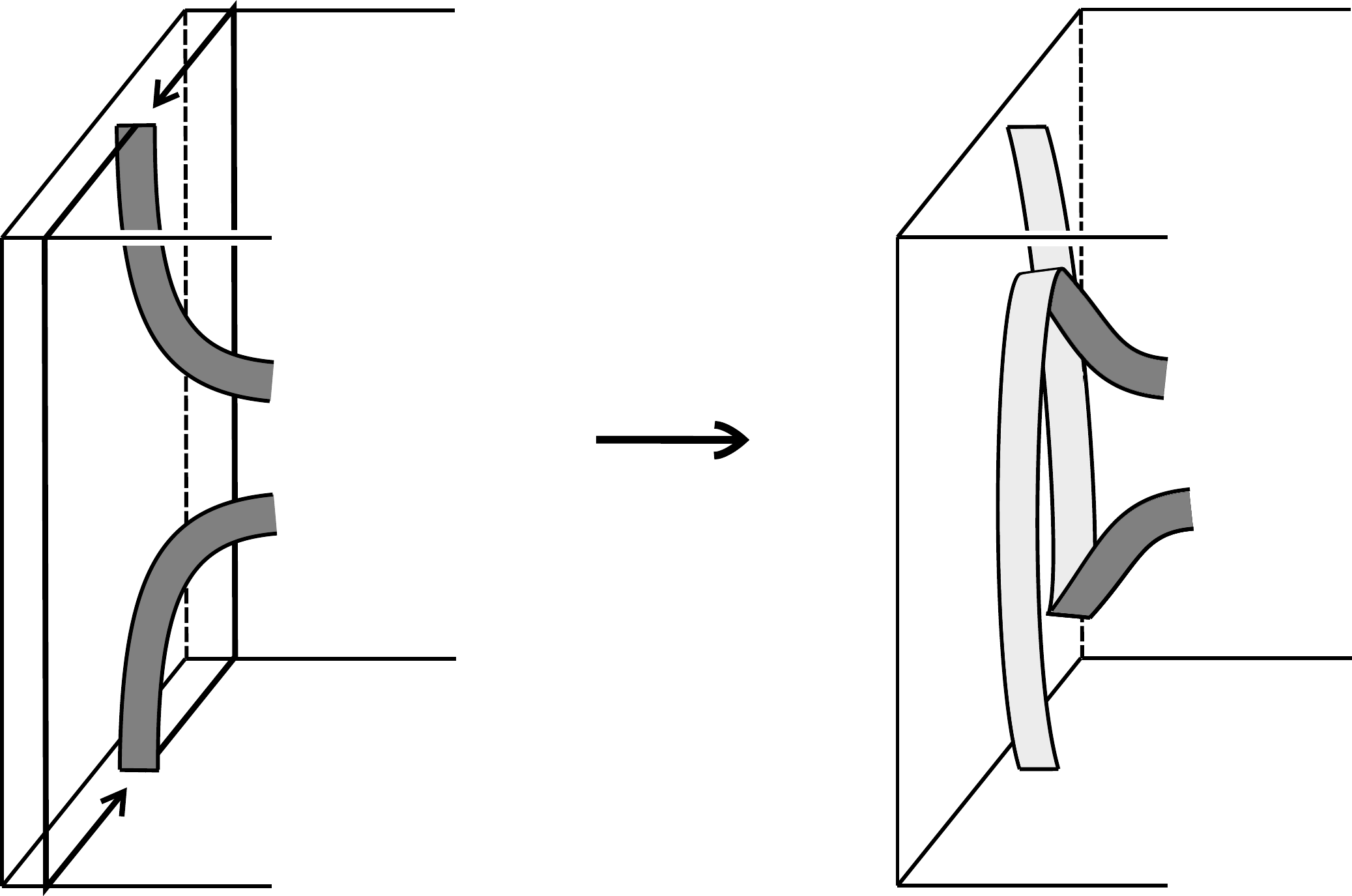}
\caption{An example of sliding.}
\label{fig:sliding_J_tangle} 
\end{figure}

\begin{figure}
\centering\includegraphics[scale=0.15]{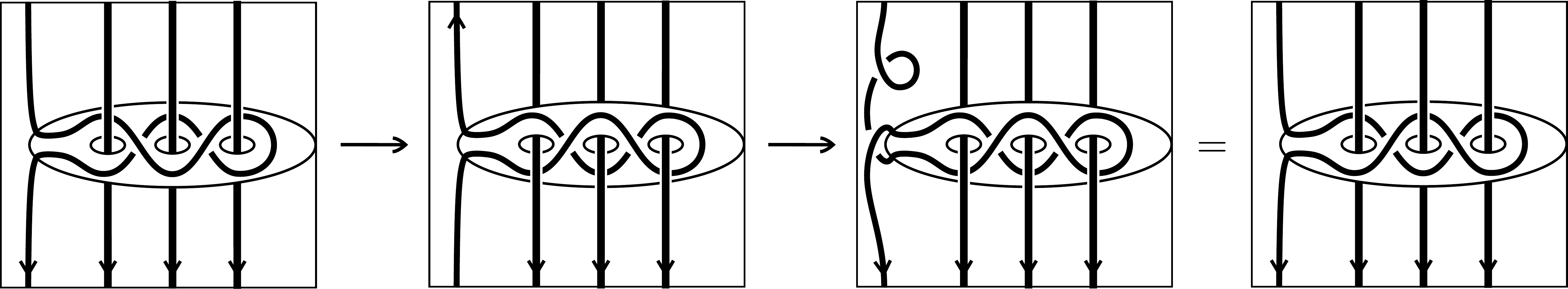}
\caption{The height involution of a horizontal $J$-tangle diagram. From left to right, we apply horizontal mirroring and sliding. The last equality follows from isotopy. The first and last diagrams differ only in their crossing data on the surface.}
\label{fig:height_involution_horizontal_ex} 
\end{figure}

Define sliding of a horizontal $J$-tangle around $\partial_0$ as in Figure~\ref{fig:sliding_J_tangle}. This operation may introduce new kinks.
Sliding near a pole is defined similarly, with the operation supported in a neighborhood of that pole.
Define the \textbf{height involution} of a horizontal $J$-tangle as composition of horizontal mirroring and sliding. It agrees with the height involution on $\widetilde{\cA}(\Sigma,J)$; see Figure~\ref{fig:height_involution_horizontal_ex}.

We now use the height involution $\tau$ and $Z_{\widetilde{\cA}}$ to prove formality of the Goldman--Turaev Lie bialgebra with respect to the framed Turaev cobracket $\delta^{\rm fr}$.
Viewing $\Sigma_{0,p+1}$ as a subspace of the Euclidean plane gives a natural framing. We take clockwise rotation to have rotation number $+1$.
For $\gamma\in\pi^+$, let $\rot(\gamma)$ denote its rotation number.
Define $i\colon\pi:=\pi_1(\Sigma)\to\pi^+$
by choosing $i(\gamma)$ so that its closure $\operatorname{cl}(i(\gamma))$ has rotation number $1$.
Since
$\rot(\gamma_1\gamma_2)=\rot(\gamma_1)+\rot(\gamma_2)-1$, the map $i$ is a group homomorphism.
In \cite{AKKN18}, $i$ is defined so that closure has rotation number zero. Their different convention for the product on $\pi^+$ makes the map a group homomorphism as well. We abuse $i$ for the induced homomorphism $\bQ\pi\to\bQ\pi^+$.
With our convention, the descending lift of $i(\gamma)$ has writhe zero. It can be drawn as a monotonically descending horizontal $J$-tangle diagram.

We define the framed Turaev cobracket $\delta^{\rm fr}$ by
\begin{align}
\delta^{\rm fr}\colon \bQ|\pi|\to \bQ|\pi|\otimes \bQ|\pi|,
\qquad\gamma
\mapsto \operatorname{Alt}(\id\otimes |\cdot |)(\mu^+(i(\Tilde{\gamma})))+\gamma\otimes \textbf{1}_\circ-\textbf{1}_\circ\otimes \gamma,\label{eq:Turaev_cobracket}
\end{align}
where, $\Tilde{\gamma}$ is a lift of $\gamma$ to $\bQ\pi$, and $\operatorname{Alt}\colon\bQ|\pi|\otimes\bQ|\pi|\to\bQ|\pi|\otimes\bQ|\pi|$ sends $x\otimes y$ to $x\otimes y-y\otimes x$. The symbol $\textbf{1}_\circ$ denotes the trivial loop,
and $\textbf{1}_\circ\equiv2\rho$ in the classical limit.
The map $\delta^{\rm fr}$ is well-defined and is the framed version of the original Turaev cobracket \cite[Proposition~5.10]{AKKN18}.

The height involutions on $\widetilde{\cA}(\Sigma,J)$ and $\Gr_I\widetilde{\cA}(\Sigma,J)$ induce two height involutions on $\mathcal{D}_{\widetilde{\cA}}$.
We first describe their difference. Since formality holds for horizontal mirroring, we separate the horizontal mirror and sliding operations. Thus, we compare the images under $Z_{\widetilde{\cA}}$ of the first and third diagrams from the left in Figure~\ref{fig:height_involution_horizontal_ex}.

Recall that $\text{horizontal mirror}=\Gr\tau=Z^{-1}\tau Z$ and $\Gr\mu^+=(Z(T_\gamma)-Z(\tau(T_\gamma)))/\hbar$.
The terms of $\hbar$-degree one in $(\mu^+-\Gr\mu^+)(\gamma)$ are relevant to the Kawazumi--Kuno coaction. Therefore, we work modulo $h\hbar$ and $\hbar^2$.
Writing $\equiv$ for equality modulo $h\hbar$ and $\hbar^2$, the degree-one part of $(\mu^+-\Gr\mu^+)(\gamma)$ is given by
\begin{align}\label{eq:quantum_height_involution}
&\nu\exp(-\rho\hbar)\begin{array}{c}\includegraphics[scale=0.15]{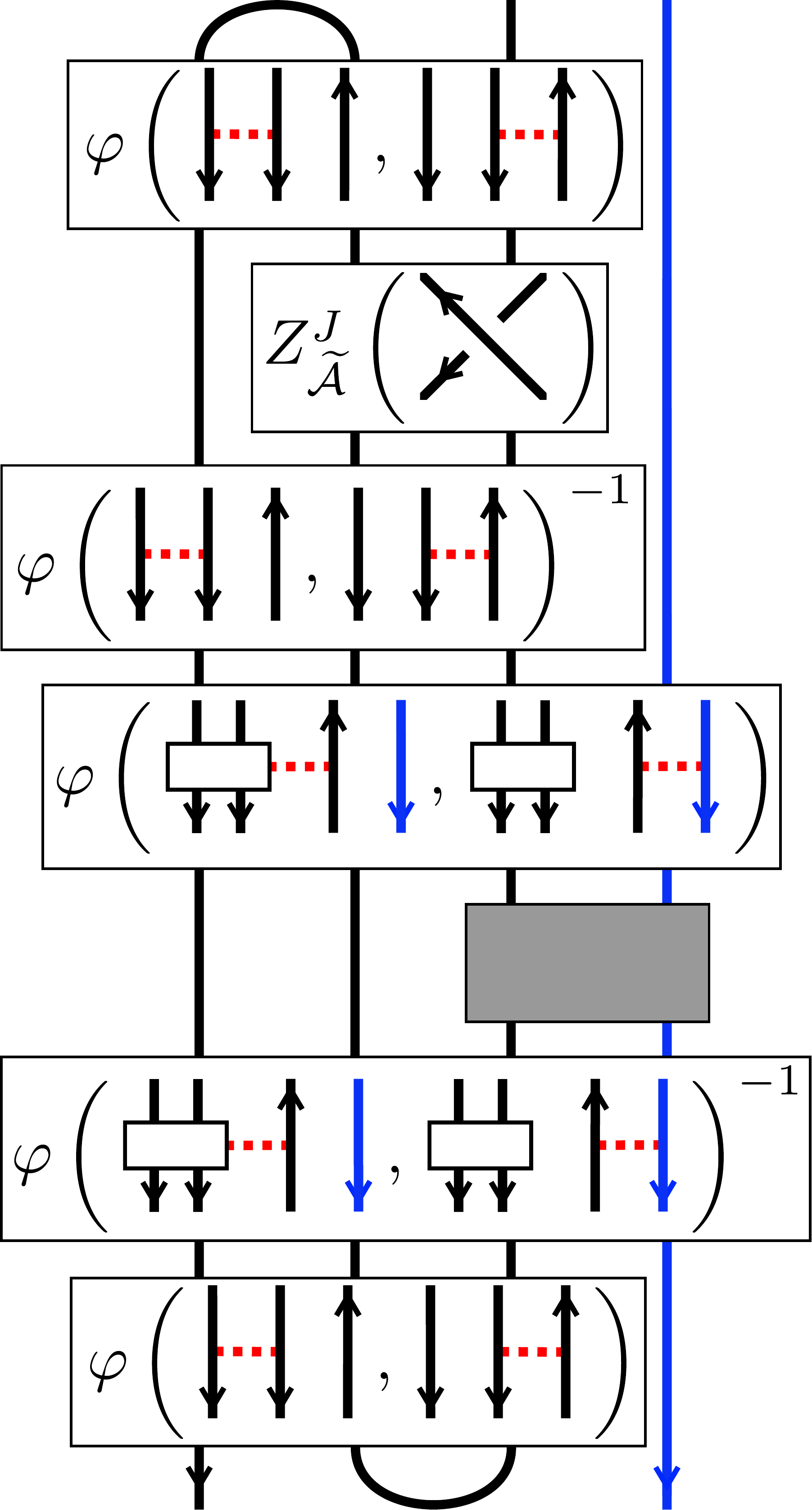}\end{array}-\begin{array}{c}\includegraphics[scale=0.15]{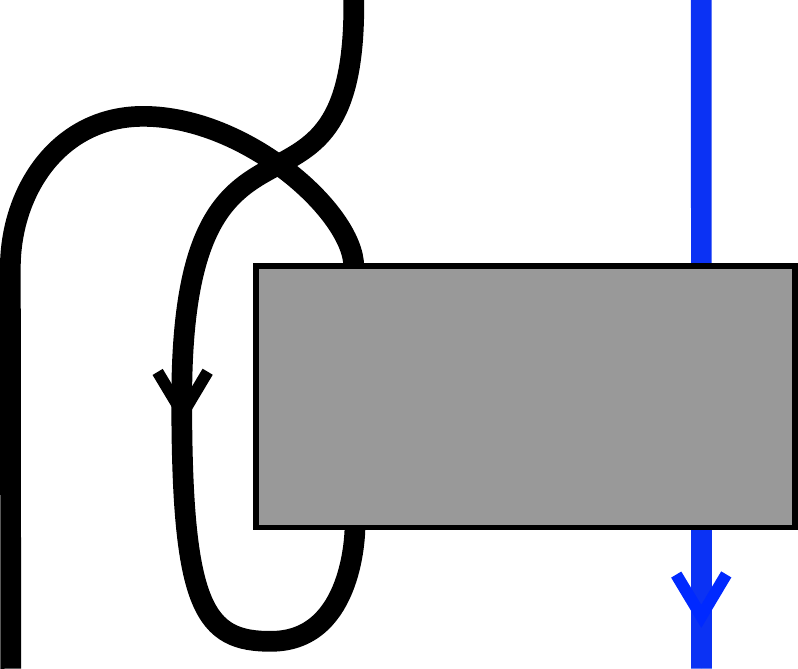}\end{array}\\
&\equiv\exp(-\rho\hbar)\begin{array}{c}\includegraphics[scale=0.15]{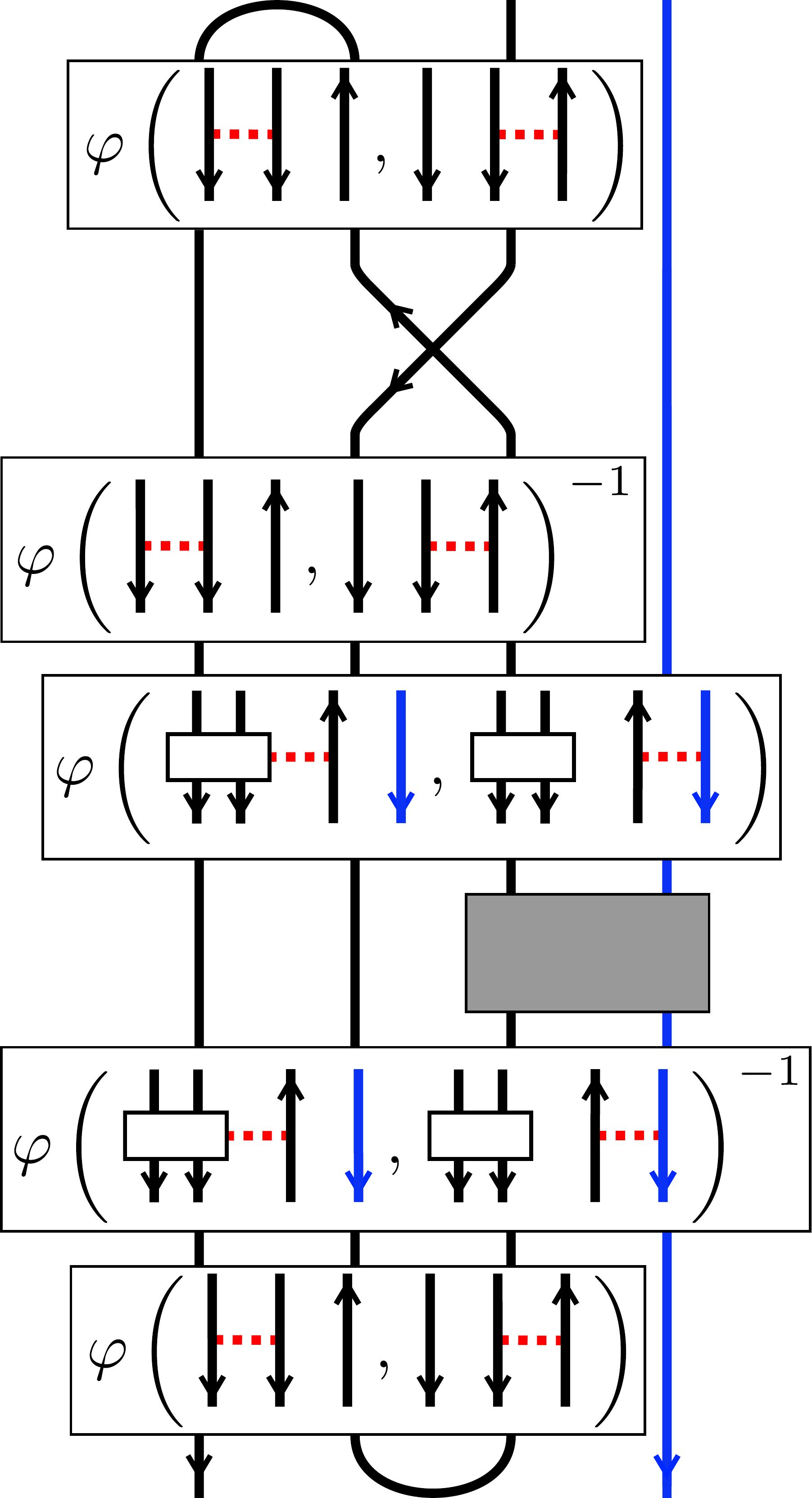}\end{array}
+\exp(-\rho\hbar)\frac{\abar}{2}\begin{array}{c}\includegraphics[scale=0.15]{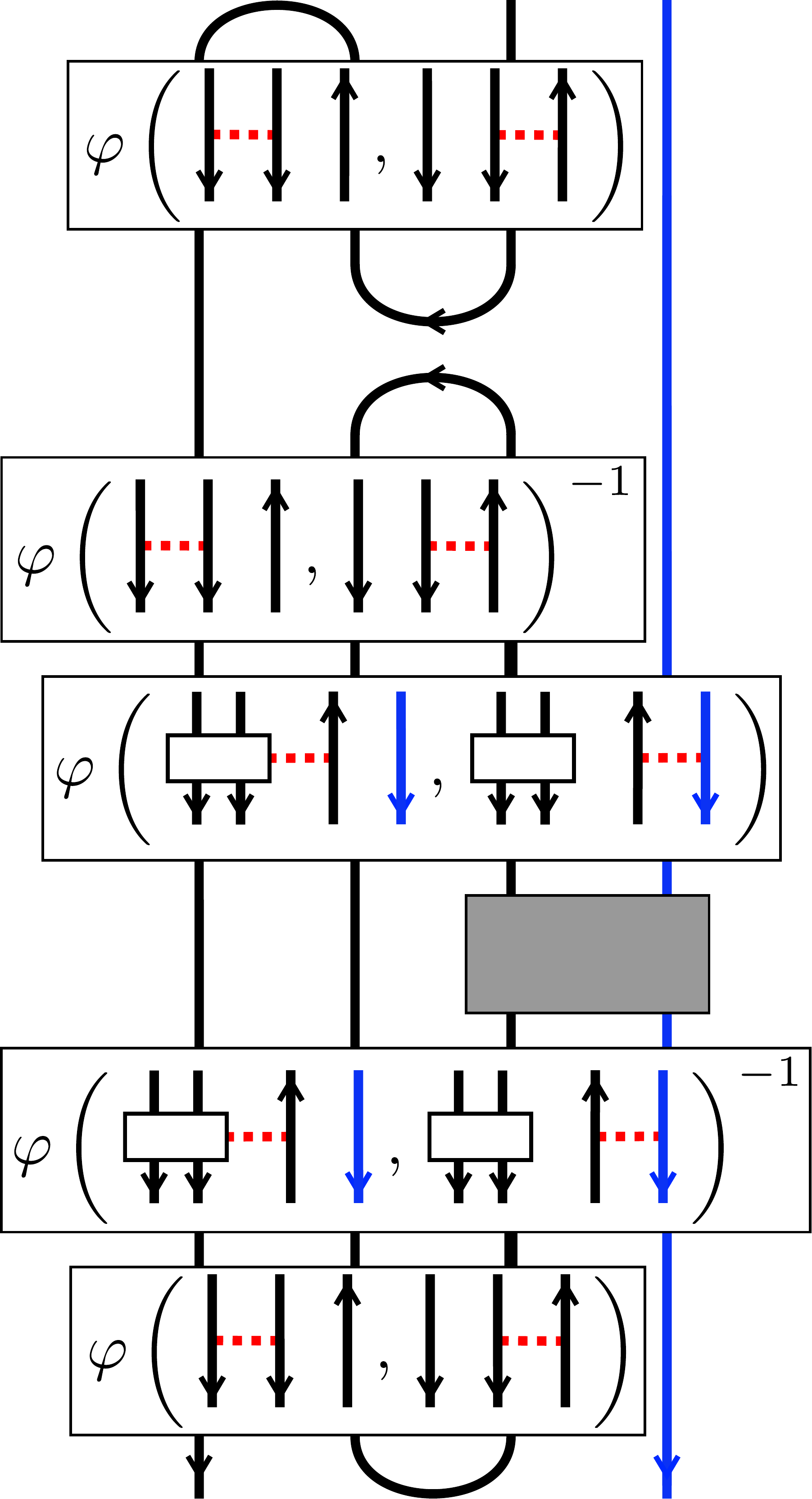}\end{array}
-\begin{array}{c}\includegraphics[scale=0.15]{draws/associator_coaction00.pdf}\end{array}\nonumber\\
&\equiv\begin{array}{c}\includegraphics[scale=0.15]{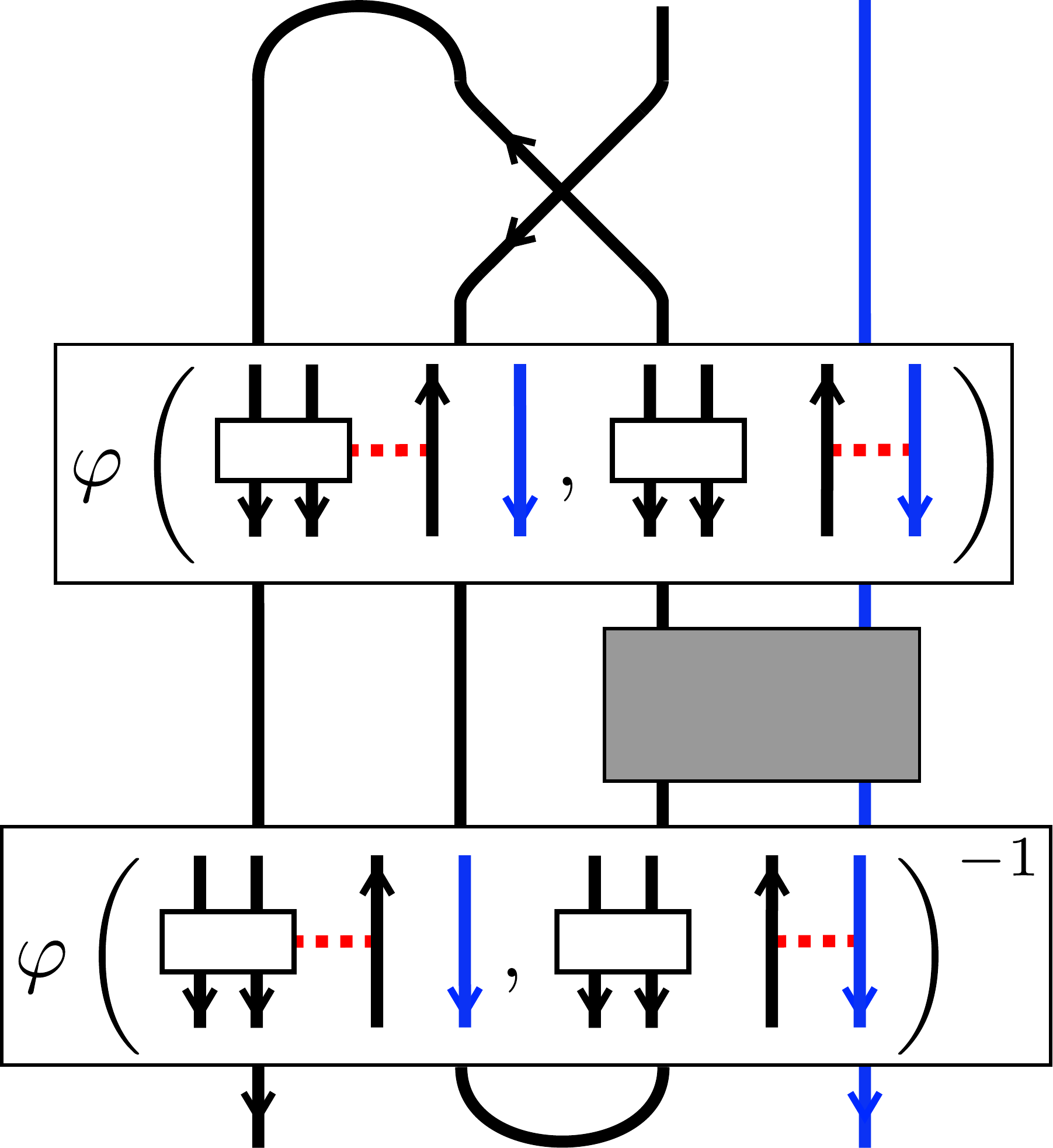}\end{array}
-\begin{array}{c}\includegraphics[scale=0.15]{draws/associator_coaction00.pdf}\end{array}
-\rho \hbar\begin{array}{c}\includegraphics[scale=0.15]{draws/associator_coaction00.pdf}\end{array}
+\frac{\hbar}{2} \begin{array}{c}\includegraphics[scale=0.15]{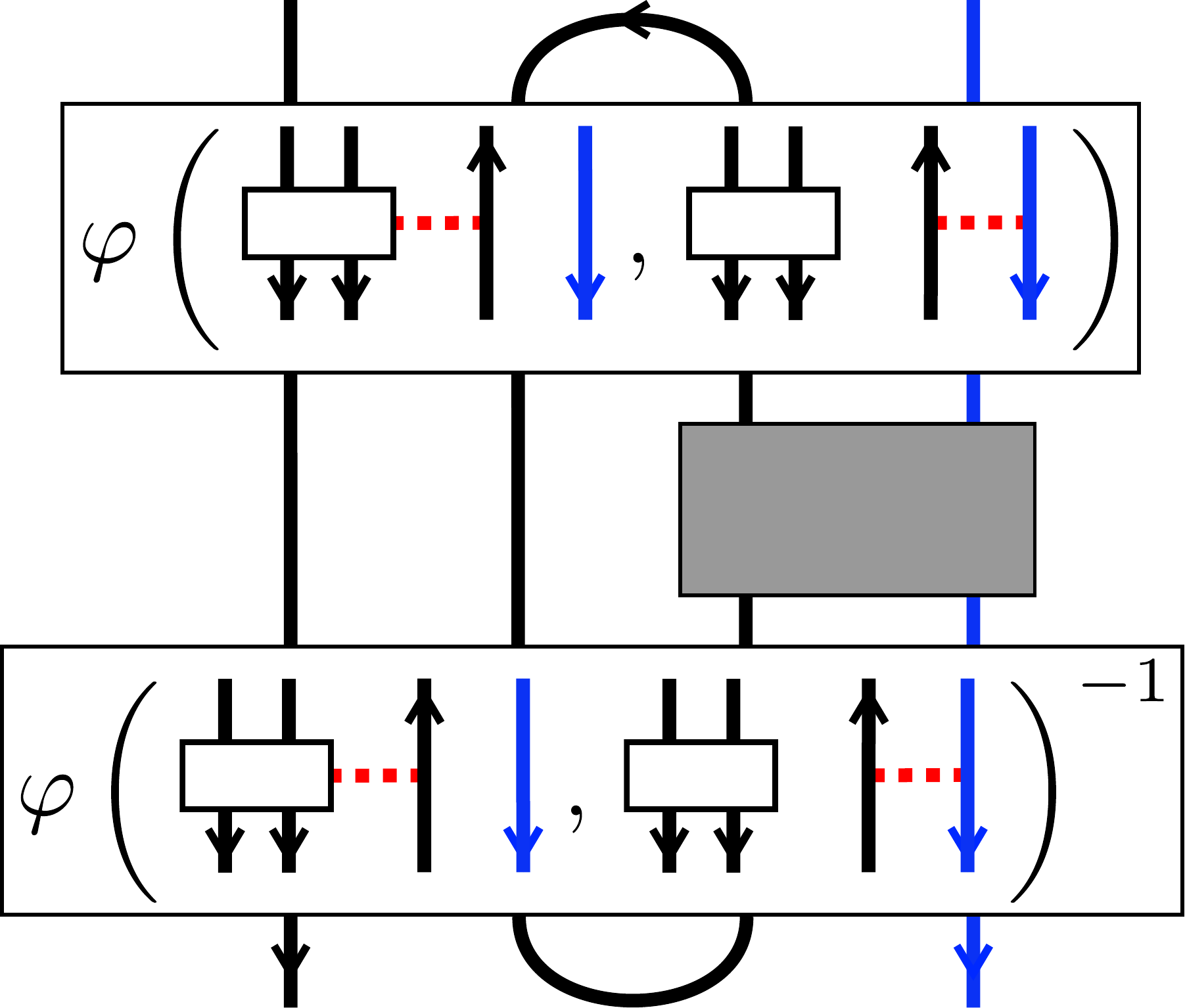}\end{array}\nonumber
\end{align}
The blue vertical line denotes a bundle of poles, and $\varphi$ denotes the associator.
The gray region indicates the part corresponding to horizontal mirroring.
As explained in Section~\ref{subsec:formality_expansion}, $\abar\equiv\hbar$.
In these computations, associators involving only strands may be omitted modulo $h\hbar$ and $\hbar^2$, since their nonconstant terms have degree at least two.

Rewriting part of the preceding expression gives
\begin{align*}
&\begin{array}{c}\includegraphics[scale=0.15]{draws/associator_coaction03.pdf}\end{array}
-\begin{array}{c}\includegraphics[scale=0.15]{draws/associator_coaction00.pdf}\end{array}\\
&\equiv
\begin{array}{c}\includegraphics[scale=0.15]{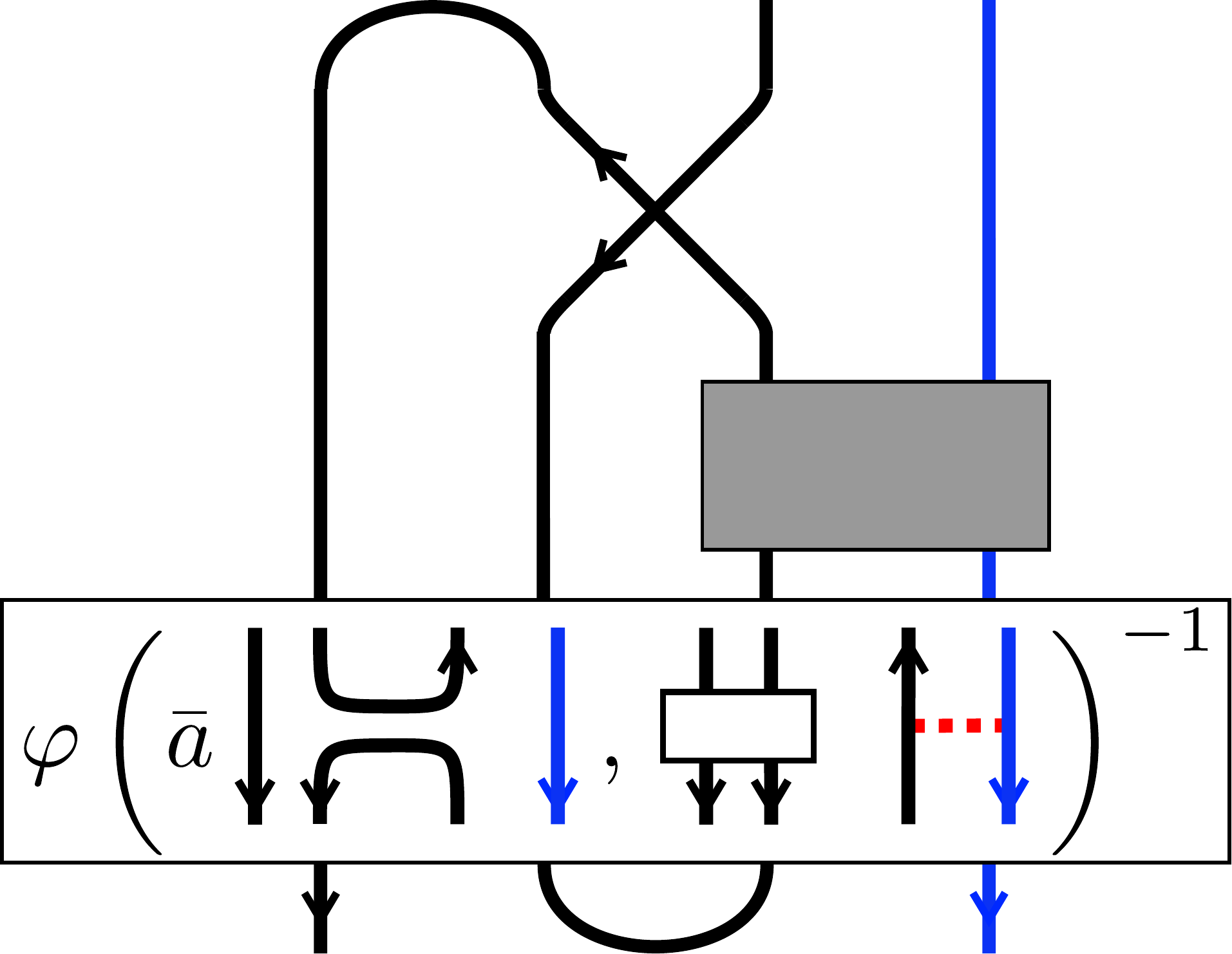}\end{array}
+\begin{array}{c}\includegraphics[scale=0.15]{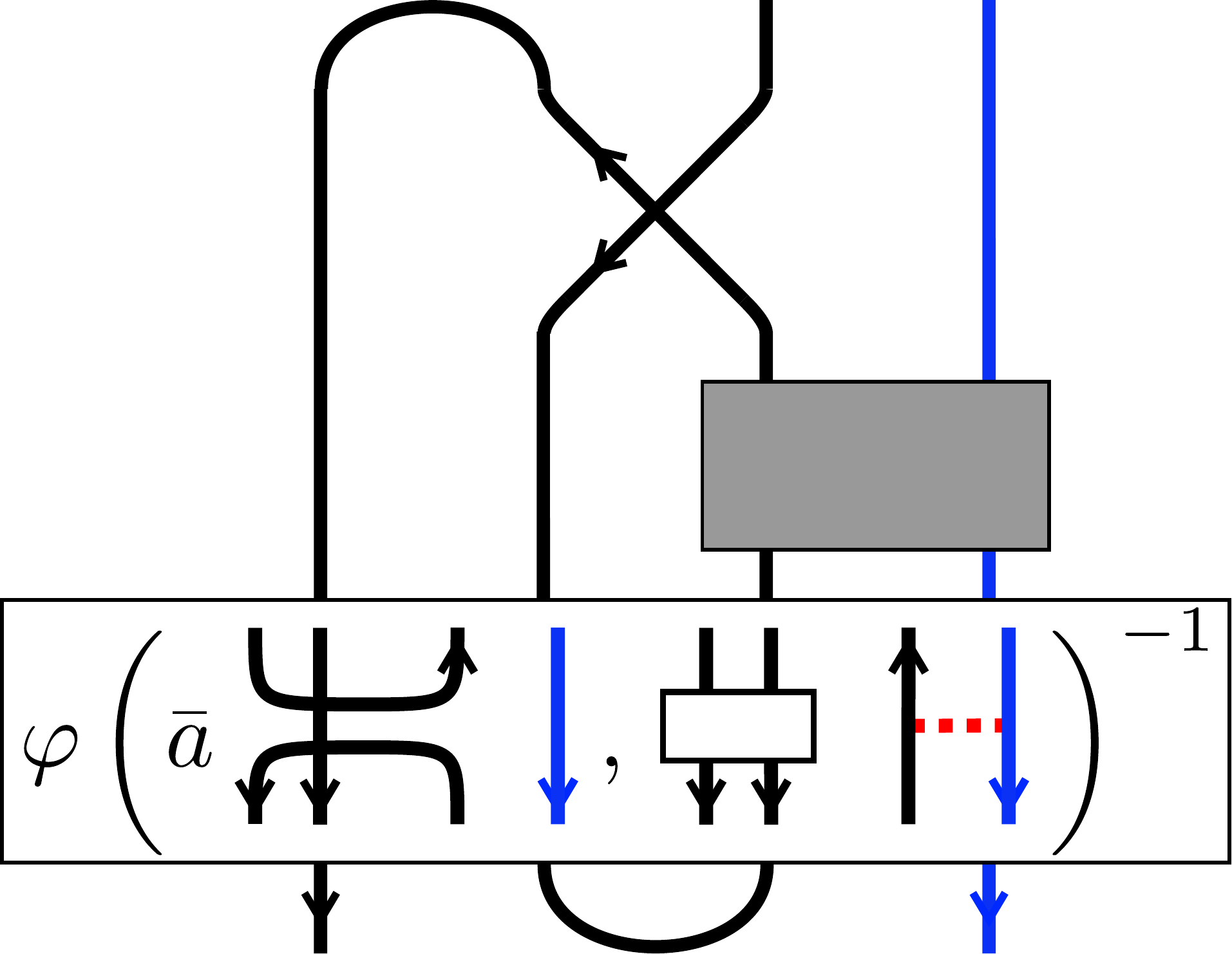}\end{array}\\
&+\begin{array}{c}\includegraphics[scale=0.15]{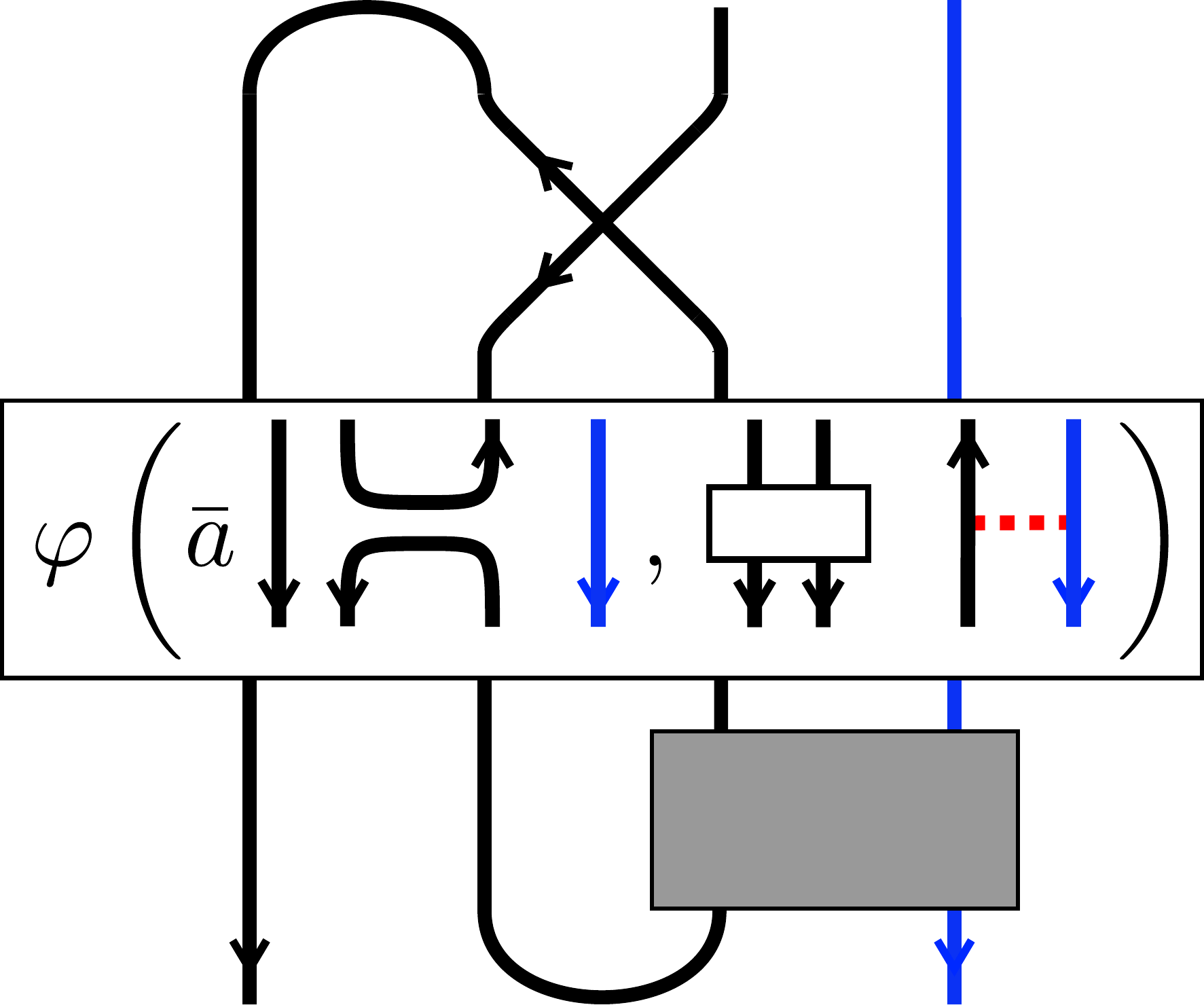}\end{array}
+\begin{array}{c}\includegraphics[scale=0.15]{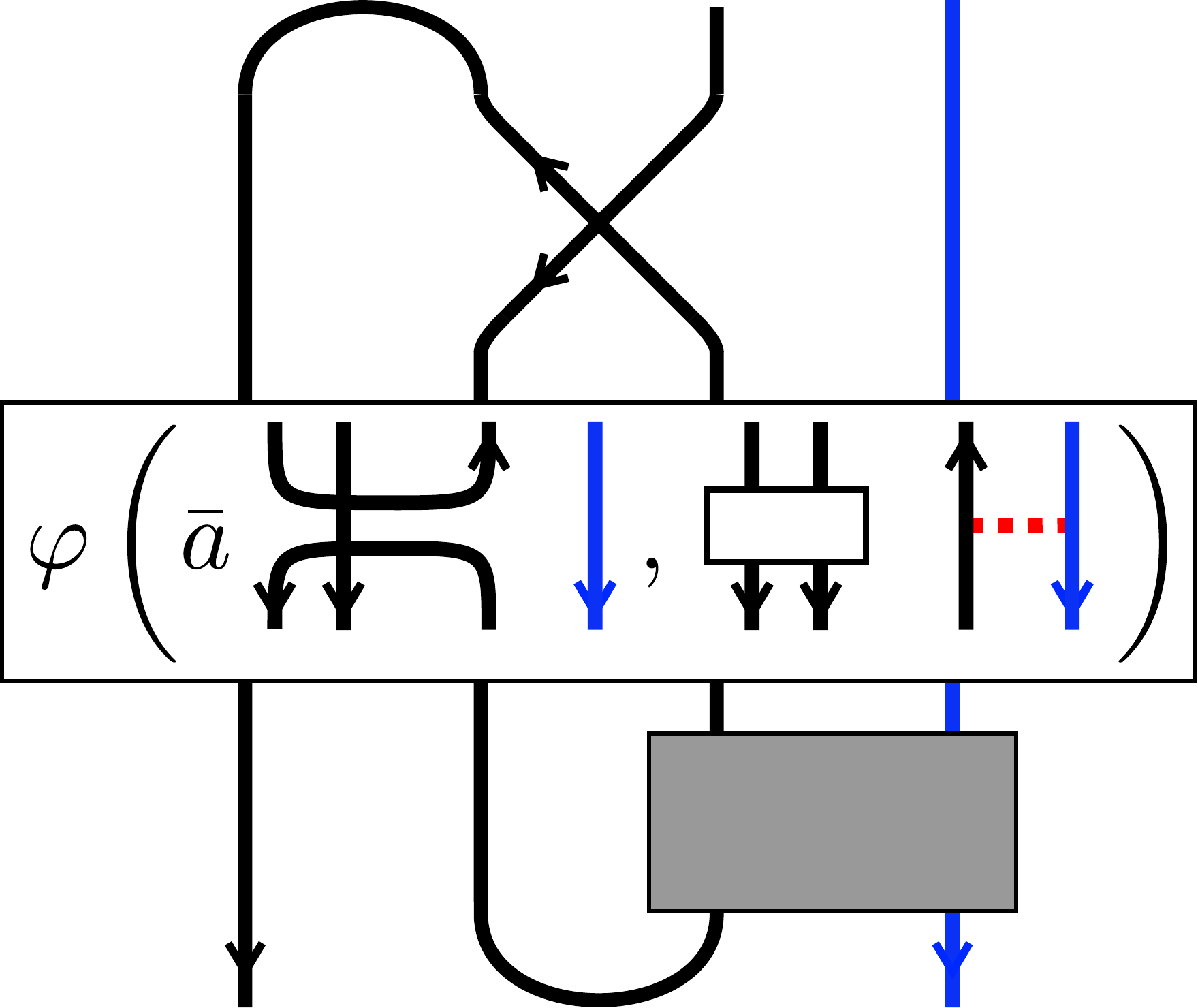}\end{array}
-4\begin{array}{c}\includegraphics[scale=0.15]{draws/associator_coaction00.pdf}\end{array}\\
&\equiv\abar\sum_{\ell}q_{\ell}\sum_{i=0}^\ell\binom{\ell}{i}\Bigg(-(-1)^{i}
\begin{array}{c}\includegraphics[scale=0.15]{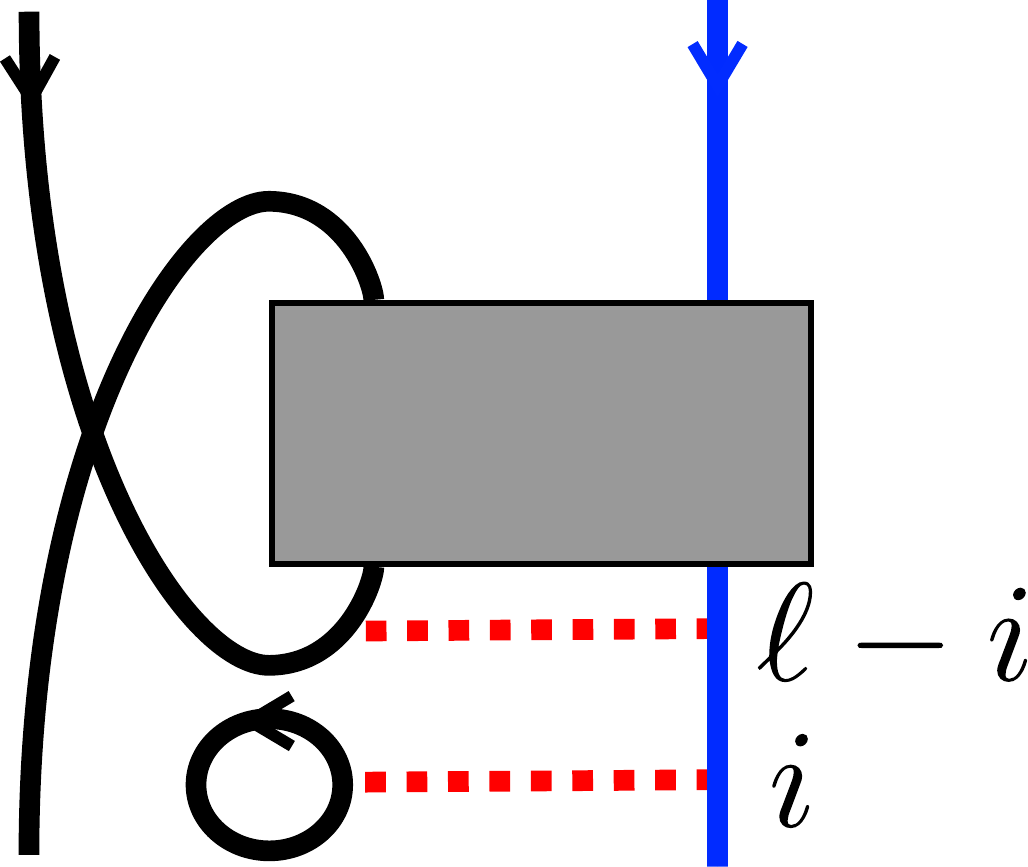}\end{array}
-(-1)^{i}\begin{array}{c}\includegraphics[scale=0.15]{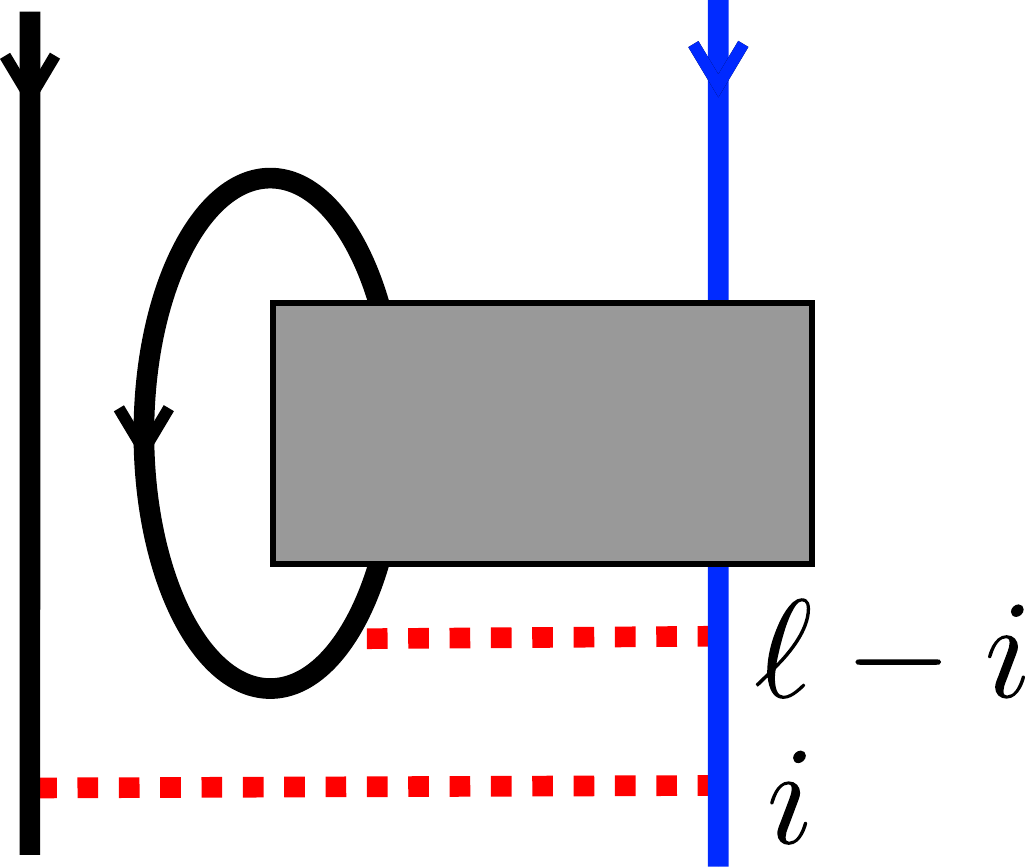}\end{array}\\
&+(-1)^{i}\begin{array}{c}\includegraphics[scale=0.15]{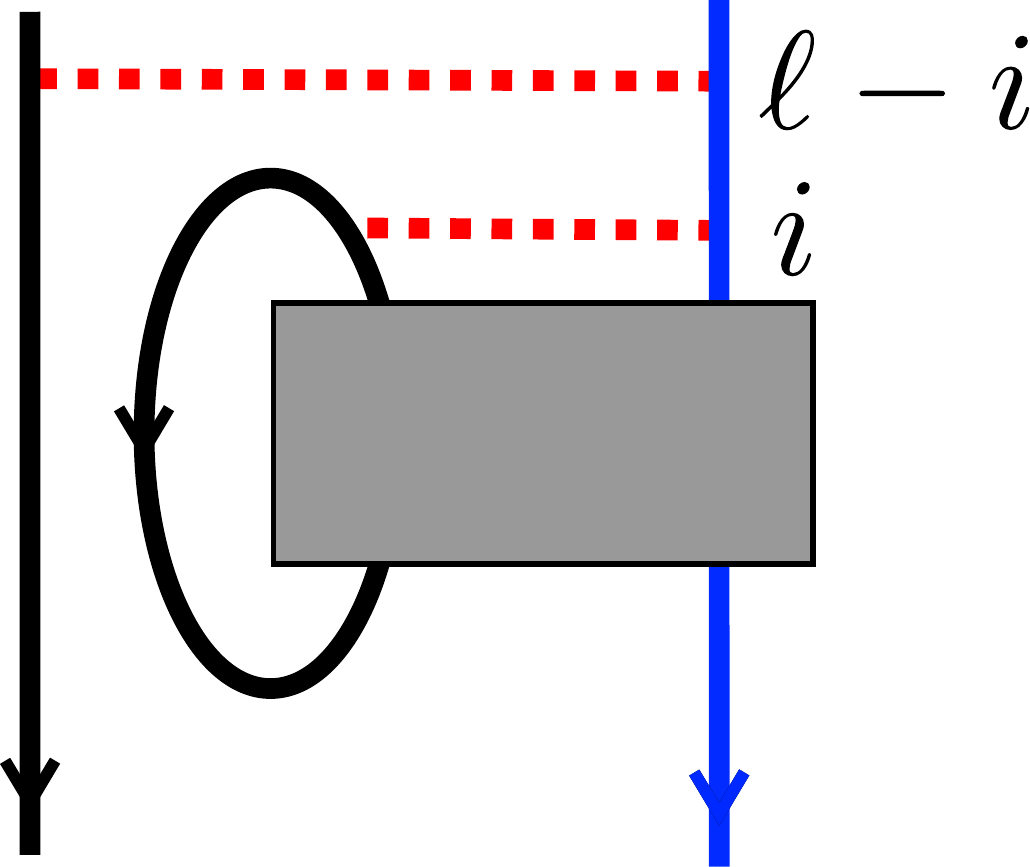}\end{array}
+(-1)^{i}\begin{array}{c}\includegraphics[scale=0.15]{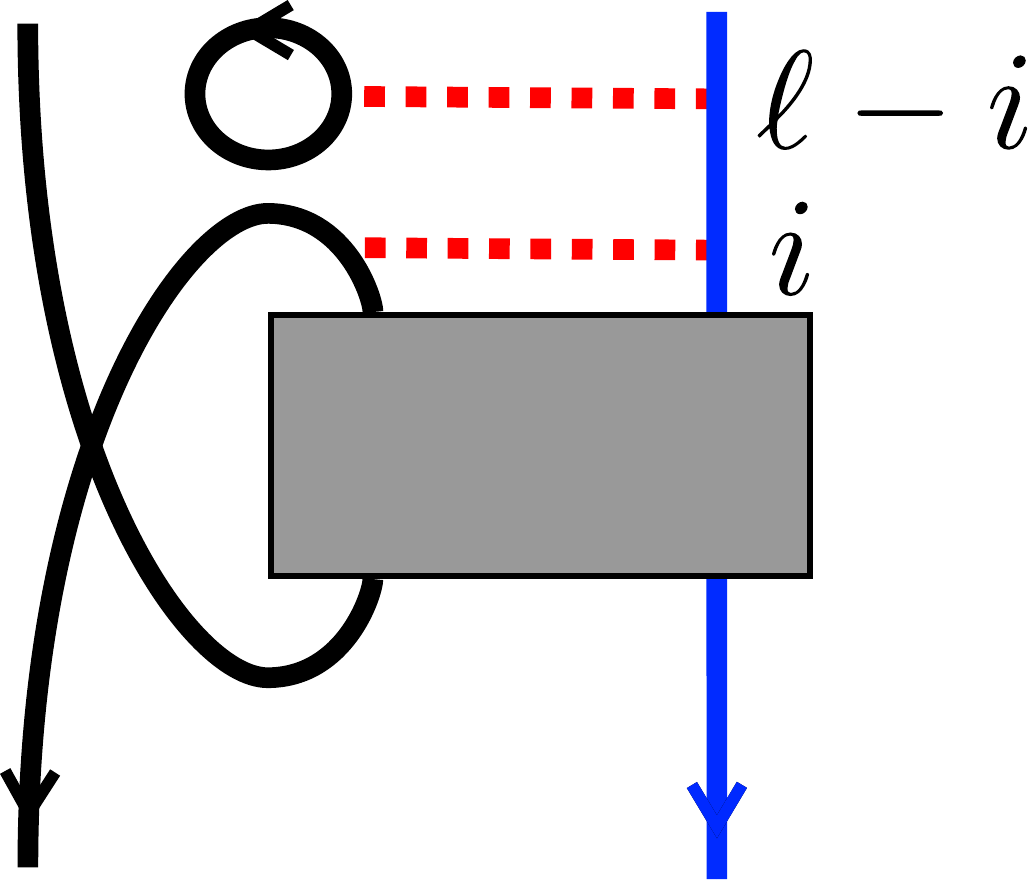}\end{array}\Bigg).
\end{align*}
Present the associator $\varphi(x,y)$ in the form
\begin{align}
    \varphi(x,y)=\exp(\sum_{1\leq\ell}q_\ell \mathrm{ad}_y^\ell(x)+O(x^2))\label{eq:associator}
\end{align}
where $\mathrm{ad}_y$ denotes the adjoint action of $y$ and $O(x^2)$ denotes terms containing at least two occurrences of $x$.
Using the isomorphism between $\mathcal{D}_{\widetilde{\cA}}^J$ and the quantized necklace algebra,
we regard elements of the classical limit $(\mathcal{D}_{\widetilde{\cA}}^J)^{(2)}/(h(\mathcal{D}_{\widetilde{\cA}}^J)^{(1)}+\hbar(\mathcal{D}_{\widetilde{\cA}}^J)^{(3)})$
as elements of $|TH|\otimes TH$. Here, $TH$ is the tensor algebra of $H_1(\Sigma_{0,p+1})$, and $|TH|$ is its cyclic quotient.

We will apply the previous computations modulo $h$ and $\hbar$.
Choose $\gamma\in|\pi|$
and a lift $\Tilde{\gamma}\in\pi$. Let
$T_\gamma\in \widetilde{\cA}
(\Sigma_{0,p+1},J)$ be the $J$-tangle diagram obtained by choosing over/under information so that it is descending along its orientation.
Abbreviate $Z_{\cA}^J$ to $Z$ and apply $\rho=|\textbf{1}|/2$. By abuse of notation, identify the classical limits of $\mathcal{D}_{\widetilde{\cA}}^J$ and $\widetilde{\cN}_I(\Sigma,J)$ with $|TH|\otimes TH$. Let $[\partial_0]$ be the homology class of the peripheral loop $\partial_0$. Then
\begin{align}
&Z(\mu^+(\Tilde{\gamma}))-\Gr\mu^+(Z(\Tilde{\gamma}))
\equiv
\frac{1}{\hbar}\left(Z(T_{\gamma})- Z(\tau(T_{\gamma}))-(Z(T_{\gamma})-
\tau  Z(T_{\gamma}))\right)
\equiv
-\frac{1}{\hbar}(Z(\tau(T_{\gamma}))-
\tau  Z(T_{\gamma}))
\nonumber\\
&\equiv 
-(\frac{1}{2}(-|\textbf{1}|\otimes Z(\Tilde{\gamma}) + Z(\gamma) \otimes \textbf{1})\nonumber\\
&+
\sum_{\ell}q_{\ell}\sum_{i=0}^\ell\binom{\ell}{i}(-1)^{i}(-|[\partial_0]^{\ell-i}Z(\tilde{\gamma})|\otimes [\partial_0]^i-|[\partial_0]^{i}|\otimes Z(\tilde{\gamma}) [\partial_0]^{\ell-i}
\nonumber\\
&+|[\partial_0]^{i}Z(\Tilde{\gamma})|\otimes  [\partial_0]^{\ell-i}
+|[\partial_0]^{\ell-i}|\otimes Z(\Tilde{\gamma}) [\partial_0]^{i}))\label{eq:difference_mu+}
\end{align}
Here, $\textbf{1}$ denotes the trivial path, or equivalently the unit $\textbf{1}\in H^{\otimes0}\subset TH$, while
$|\textbf{1}|=\textbf{1}_\circ$ denotes the trivial loop,
so $|\textbf{1}|\in|H^{\otimes0}|\subset|TH|$. 
Since
$$T_\gamma\equiv\tilde{\gamma}\mod h,\hbar,\qquad
|T_\gamma|\equiv\gamma\mod h,\hbar,$$ we use the notation
$Z(\Tilde{\gamma}):=Z(T_\gamma)$,
$Z(\gamma):=Z(|T_\gamma|)$
for the corresponding images.
After applying $\id\otimes|\cdot|$, the summation term in \eqref{eq:difference_mu+} is invariant under interchanging the two tensor factors. Therefore, applying $\operatorname{Alt}$ annihilates it.

The previous computation expresses, in terms of $\Gr\delta^{\mathrm{fr}}$, the cobracket on the quantized necklace algebra transported from $\delta^{\mathrm{fr}}$ by the module isomorphism $Z$. The term $\gamma\otimes\textbf{1}_\circ-\textbf{1}_\circ\otimes\gamma$ in the definition of $\delta^{\rm fr}$ vanishes in the associated graded quotient.
We show that $Z$ is a homomorphic expansion preserving $\delta^{\mathrm{fr}}$, namely,
$Z\circ\delta^{\rm fr}=\Gr\delta^{\rm fr}\circ Z$.
Using \eqref{eq:difference_mu+},
we obtain, for $\Tilde{\gamma}\in\pi$,
\begin{align*}
&(Z\circ \delta^{\rm{fr}} -\Gr \delta^{\rm{fr}}\circ Z )(|\Tilde{\gamma}|) \\
&=\operatorname{Alt}(\id\otimes |\cdot |)
Z(\mu^+ (i(\Tilde{\gamma})))
+|Z(\Tilde{\gamma})|\otimes \textbf{1}_\circ-\textbf{1}_\circ\otimes |Z(\Tilde{\gamma})|-\operatorname{Alt}(\id\otimes |\cdot |)
\Gr\mu^+(Z(i(\Tilde{\gamma})))
\\
&=\operatorname{Alt}(\id\otimes |\cdot |)
(Z(\mu^+ (i(\Tilde{\gamma})))-\Gr\mu^+(Z(i(\Tilde{\gamma}))))
+|Z(\Tilde{\gamma})|\otimes \textbf{1}_\circ-\textbf{1}_\circ\otimes |Z(\Tilde{\gamma}) |\\
&=-(-|\textbf{1}|\otimes|Z(\Tilde{\gamma})| + |Z(\Tilde{\gamma})|\otimes |\textbf{1}|)+|Z(\Tilde{\gamma})|\otimes \textbf{1}_\circ-\textbf{1}_\circ\otimes |Z(\Tilde{\gamma} )|
= 0
\end{align*}

Combining this with Remark~\ref{rmk:Goldman_bracket_quantum}(1), which gives a quantum interpretation of the Goldman Lie bracket, and Theorem~\ref{thm:formality02}, we obtain the following.
\begin{theorem}\label{thm:formality_Goldman_Turaev}
[cf. {\cite[Theorems~5.8, 5.24]{BN+25}}]
Let $\Sigma=\Sigma_{0,p+1}$ and set $h=\hbar=0$. Under \eqref{eq:commutative_diag_expansion}, $Z_{\widetilde{\cA}}$ induces a module homomorphism from the Goldman--Turaev Lie bialgebra to the necklace Lie bialgebra.
This induced map is a homomorphic expansion preserving the Goldman bracket and the Turaev cobracket. It is also a homomorphic expansion with respect to the Kawazumi--Kuno action, which is the classical limit of $\sigma$.
\end{theorem}

\begin{remark}\label{rmk:quantum_formality}
\begin{enumerate}
    \item One checks that \eqref{eq:difference_mu+} agrees with Eq.~(47) of \cite{AKKN18}.
    Our $\Gr\mu^+$ corresponds to $\mu^{\rm alg}$ in \cite{AKKN18}. The Duflo function is expressed in terms of the coefficients $q_\ell$ above and is determined by the associator.
    The Goldman--Turaev Lie bialgebra is known not to have formality with respect to the Kawazumi--Kuno coaction \cite[Theorem~7.6 (i)]{AKKN18}. We do not give a new proof of this nonformality statement. Rather, we explicitly describe the obstruction for our construction and show that it agrees with the expression in \cite{AKKN18}.
    
    \item Our rotation-number correction differs from that in \cite{AKKN18} because composition conventions on $\pi^+$ differ. Inserting a kink identifies our composition with the one in \cite{AKKN18}. Under this identification, the homomorphic expansion preserving the framed Turaev cobracket is essentially the same as in \cite{AKKN18}.
    In contrast, the definition of the framed Turaev cobracket $\delta^{\rm fr}$ in \cite{BN+25} does not require the last two terms of \eqref{eq:Turaev_cobracket}; see Remark~3.5 of \cite{BN+25}.

    \item In \cite{AKKN18,BN+25}, the coaction is described formally by focusing on first-order terms. By viewing the coaction as $\frac{\id-\tau}{\hbar}$, we obtain formal computations in all degrees. From the viewpoint of the height involution, \eqref{eq:quantum_height_involution} can be regarded as a quantization of the Kawazumi--Kuno coaction, extending \cite{BN+25}.
    
    \item Equation~\eqref{eq:difference_mu+} refines Proposition~5.2 of \cite{BN+25}.
\end{enumerate}
\end{remark}

\subsection{Application: formality for the annulus and Dehn twists}
The commutative diagram \eqref{eq:commutative_diag_expansion} gives two descriptions of formality: one using the modified quantized necklace Lie module and the other using the HOMFLYPT quotient of chord diagrams.
In this subsection, we give both descriptions for the annulus and use them to give a simpler proof of a result of \cite{Tsu18a}.

Assume that $|J|=2$ and that the incoming and outgoing endpoints lie on distinct boundary components of $S^1\times[0,1]$. Define
\begin{align*}
\gamma^0=\begin{array}{c}\includegraphics[scale=0.15]{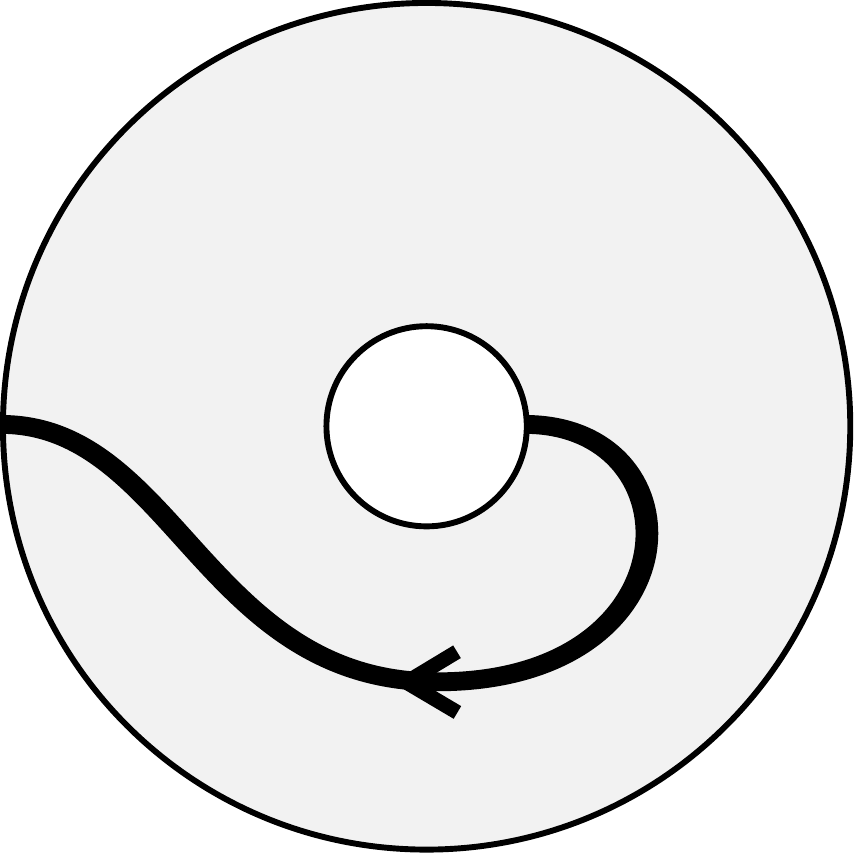}\end{array},
\quad\gamma^k=
\begin{array}{c}\includegraphics[scale=0.15]{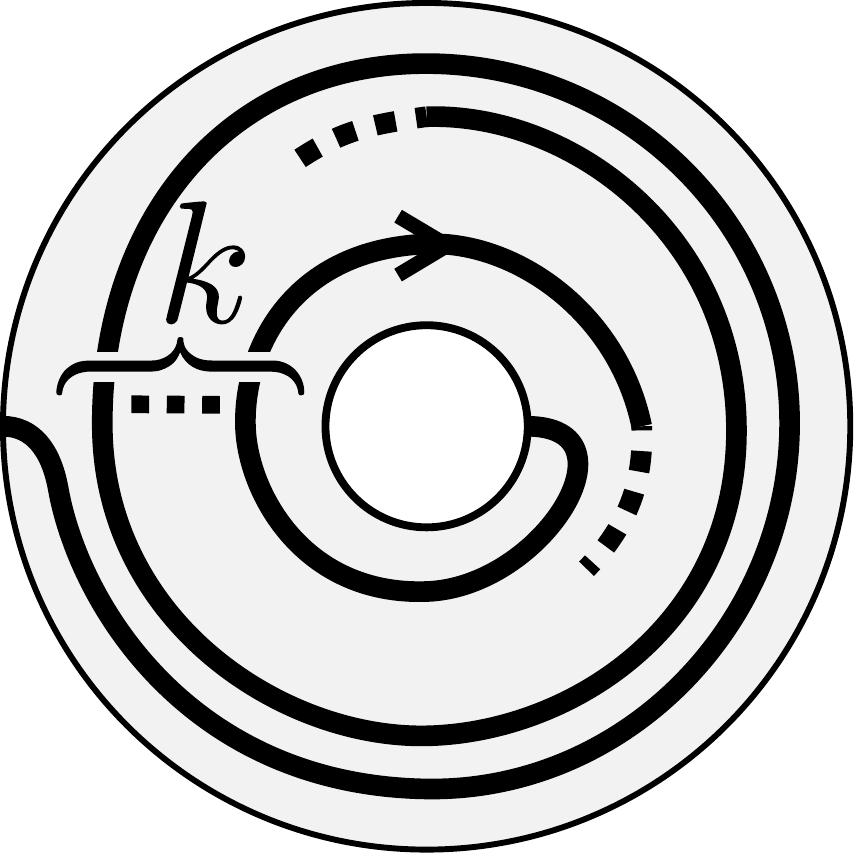}\end{array}
\in \widehat{\widetilde{\mathcal{A}}(S^1\times [0,1],J)},
\end{align*}
where, $\gamma^k$ is obtained from $\gamma^0$ by applying $k$ right-handed Dehn twists about the core of the annulus.

In what follows,
for a marked surface $(\Sigma,J')$ with $|J'|=0,2$ and $\Sigma=S^1\times[0,1]$,
we use the inverse of the expansion
$(Z_{\widetilde{\cA}}^{J'})^{-1}\colon\mathcal{D}_{\widetilde{\cA}}^{J'}(\Sigma,J')\xrightarrow{\cong}\widehat{\widetilde{\cA}(S^1\times[0,1],J')}$.
\begin{lemma}[cf. {\cite[Lemma 5.8]{Tsu18a}}]\label{lem:Dehn_twist_annulus}
Let $(S^1\times[0,1],J)$ be a marked annulus, and set
$$ L_{\widetilde{\cA}} := \frac{h}{2a} \begin{array}{c}\includegraphics[scale=0.15]{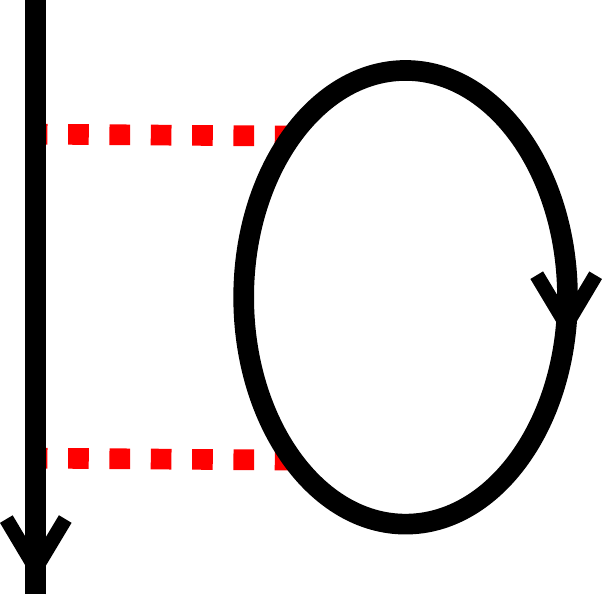}\end{array}
\in\widehat{\widetilde{\mathcal{A}}(S^1\times[0,1])}.$$
For $C=S^1\times\{1/2\}$, the right-handed Dehn twist $t_C$ about $C$ is given by
\begin{align}
    t_C=\exp(\sigma(L_{\widetilde{\cA}}))\in\Aut(\widehat{\widetilde{\mathcal{A}(\Sigma,J)})}\label{eq:Dehn_twist_C}
\end{align} 

The vertical line in the chord diagram represents a pole.
Moreover, $L_{\widetilde{\cA}}$ is uniquely determined by \eqref{eq:Dehn_twist_C} together with the condition $\iota^\ast(L_{\widetilde{\cA}})=0$, where $\iota^\ast$ is induced by the natural embedding $\iota\colon S^1\times[0,1]\to D^2$ into a disk.
\end{lemma}
For the element $L_{\widetilde{\cA}}$ displayed in the lemma, $\iota^\ast(L_{\widetilde{\cA}})=0$ follows since $\iota^\ast$ forgets the pole.
Also, \eqref{eq:h_a}, \eqref{eq:expansion_a_wrt_hhbar}, and \eqref{eq:expansion_abar_wrt_hhbar} show that $h/a$ is a power series in $h\hbar$.
\begin{proof}
First, assume that $|J|=2$ and that the incoming and outgoing endpoints lie on distinct boundary components of $S^1\times[0,1]$. We show that
$$ \sigma(L_{\widetilde{\cA}})(\gamma^{0}) = \sum_{n=1}^\infty \frac{(-1)^{n-1}}{n} \sum_{i=0}^n \binom{n}{i} (-1)^{n-i} \gamma^i 
=\sum_{n=1}^\infty\frac{(-1)^{n-1}}{n}(\gamma-\gamma^0)^n=\log\gamma,$$
where $\gamma^i$ is the $J$-tangle diagram obtained by joining $i$ copies of $\gamma$ with $i-1$ copies of
$\begin{array}{c}\includegraphics[scale=0.15]{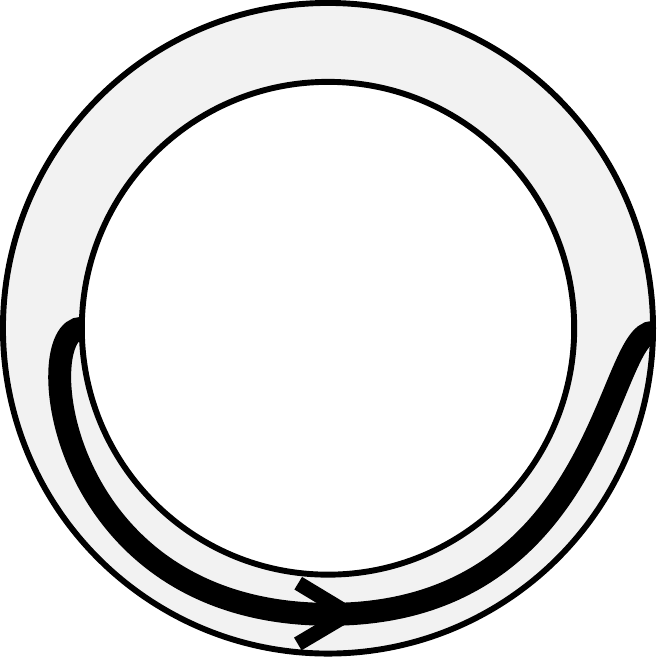}\end{array}$ inserted between them. After moving the endpoints, this agrees with the product on $\widetilde{\cA}(\Sigma,J)$ defined along the boundary of the annulus.
For example, joining two copies of $\gamma^0$ gives the following diagram, so $\gamma^0$ is the unit for this product:
$$\begin{array}{c}\includegraphics[scale=0.15]{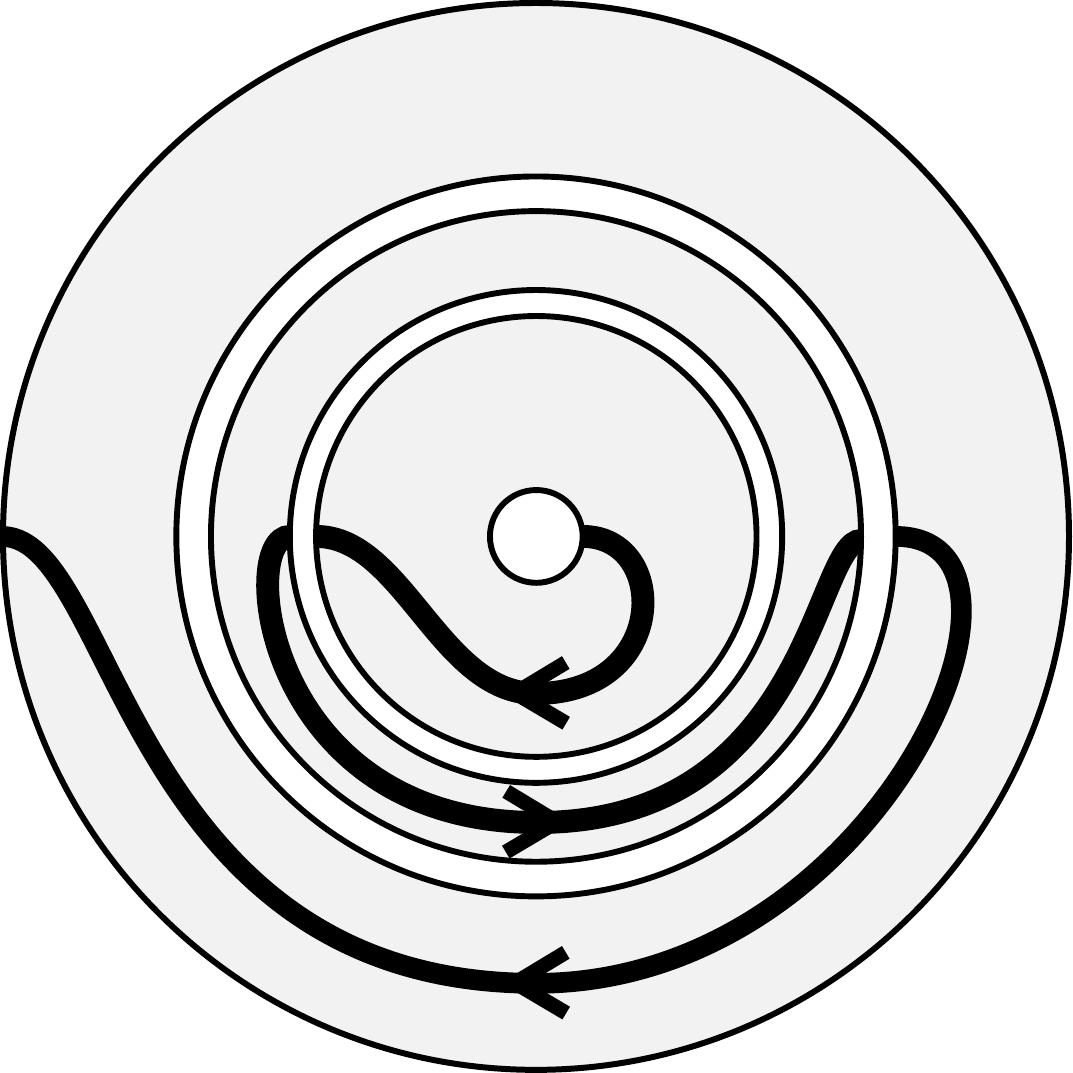}\end{array}=\gamma^0.$$

Both $L_{\widetilde{\cA}}$ and $\gamma^0$ are invariant under the height involution, which interchanges the left and right actions.
Consequently,
the height involution $\tau$ satisfies
$\tau(L_{\widetilde{\cA}}\gamma^0)=\tau(\gamma^0)\tau(L_{\widetilde{\cA}})=\gamma^0L_{\widetilde{\cA}}$.
Hence
\begin{align}
\sigma(L_{\widetilde{\cA}})(\gamma^{0}) = \frac{1}{h}(L_{\widetilde{\cA}} \gamma^{0} - \gamma^{0} L_{\widetilde{\cA}})
=\frac{1}{h}(L_{\widetilde{\cA}} \gamma^{0} - \tau(L_{\widetilde{\cA}} \gamma^{0}))
= \frac{1}{h}(\mathrm{id} - \tau)(L_{\widetilde{\cA}} \gamma^{0})
\end{align}

Therefore, to compute $\sigma(L_{\widetilde{\cA}})$, it suffices to compute $L_{\widetilde{\cA}}\gamma^0$ and separate its $\tau$-invariant and $\tau$-anti-invariant parts. Using the formality isomorphism, the computation near the pole gives
\begin{align*}
L_{\widetilde{\cA}}\gamma^0= \frac{h}{2a}\begin{array}{c}\includegraphics[scale=0.15]{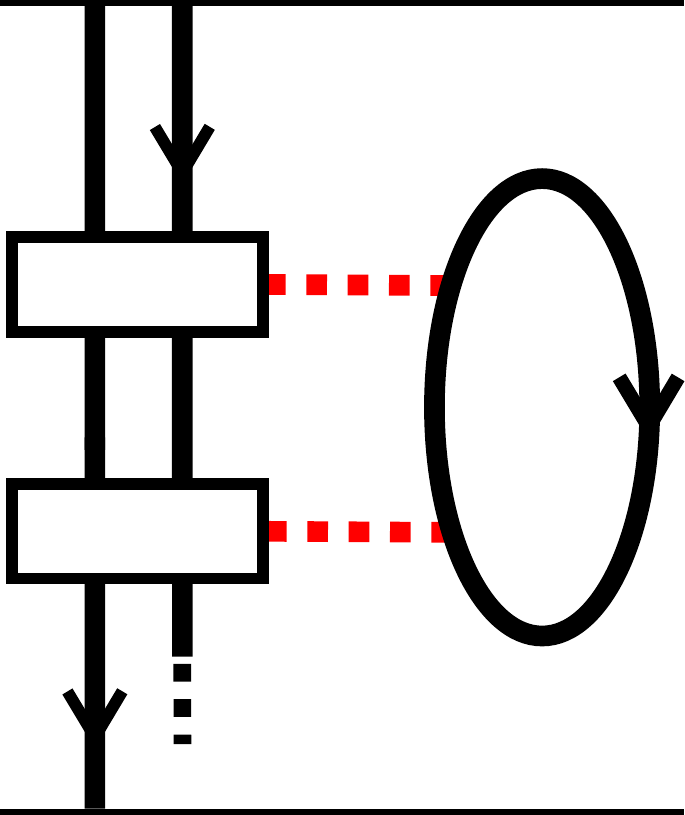}\end{array}
= \frac{h}{2a}\left(\begin{array}{c}\includegraphics[scale=0.15]{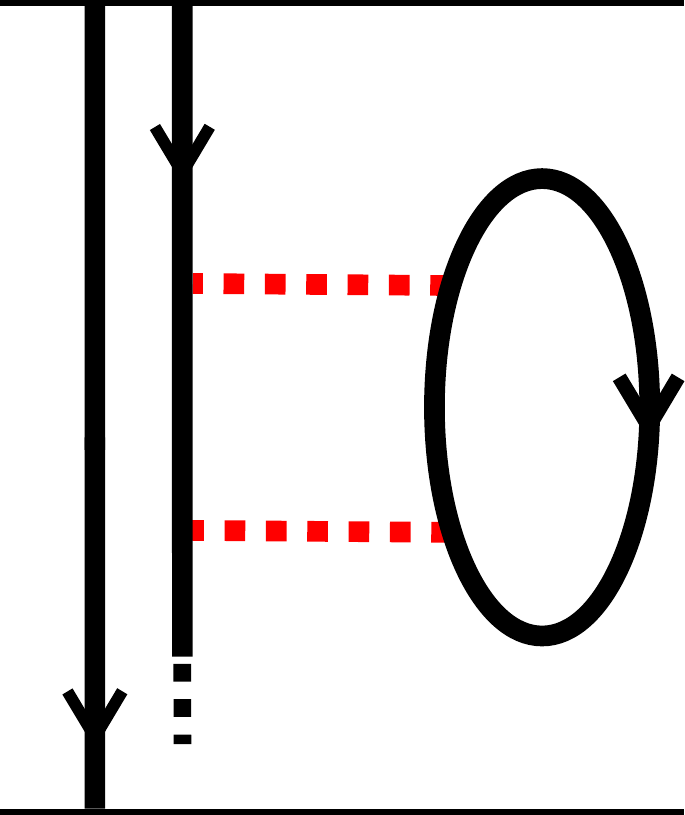}\end{array}
+\begin{array}{c}\includegraphics[scale=0.15]{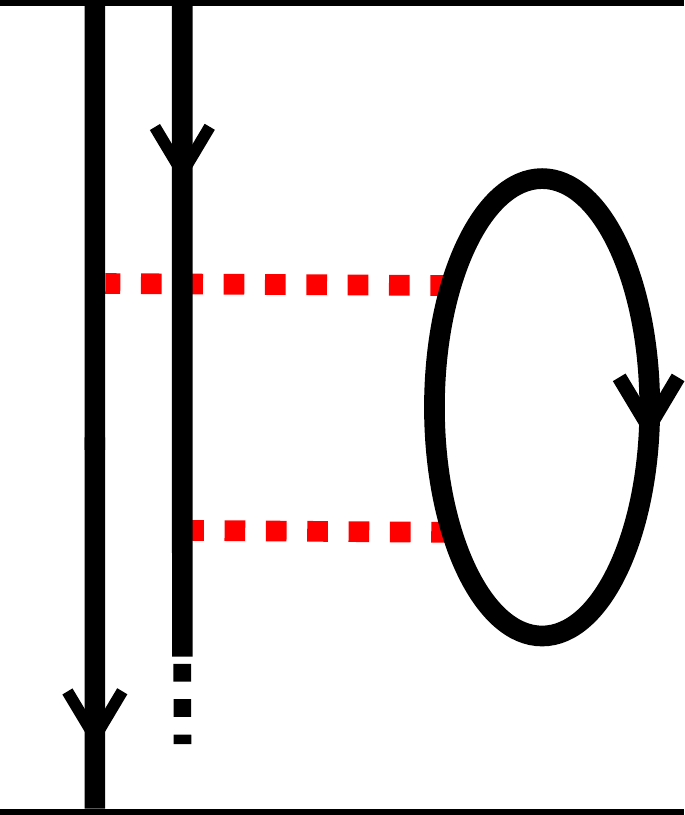}\end{array}
+\begin{array}{c}\includegraphics[scale=0.15]{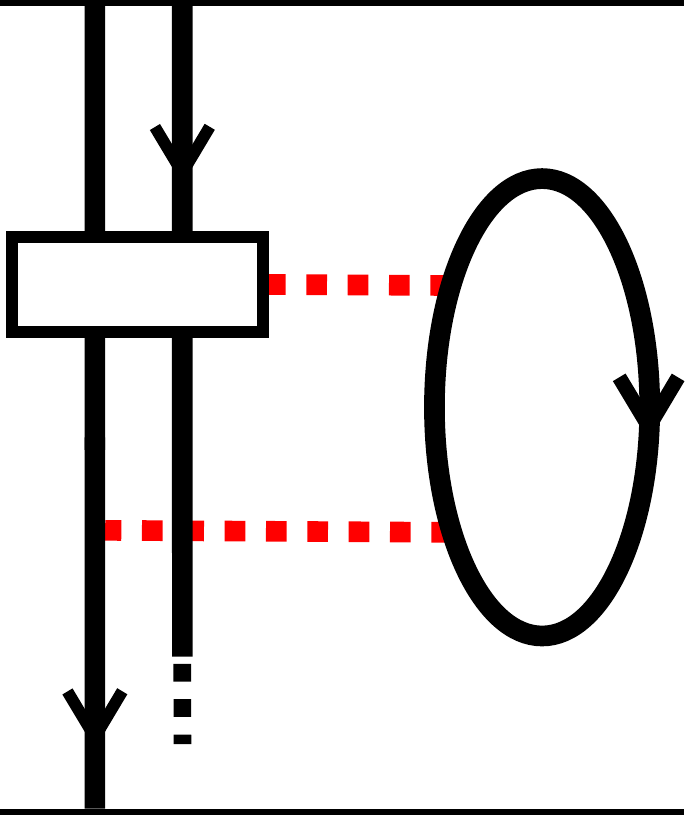}\end{array}\right)
\end{align*}

By the definition of the coupons, the first and third terms on the right are invariant under the height involution. The definitions of the Kontsevich integral and the HOMFLYPT quotient give
\begin{align}
\frac{h}{2a}\begin{array}{c}\includegraphics[scale=0.15]{draws/Dehn_twist02.pdf}\end{array}&=
\frac{h}{2}\begin{array}{c}\includegraphics[scale=0.15]{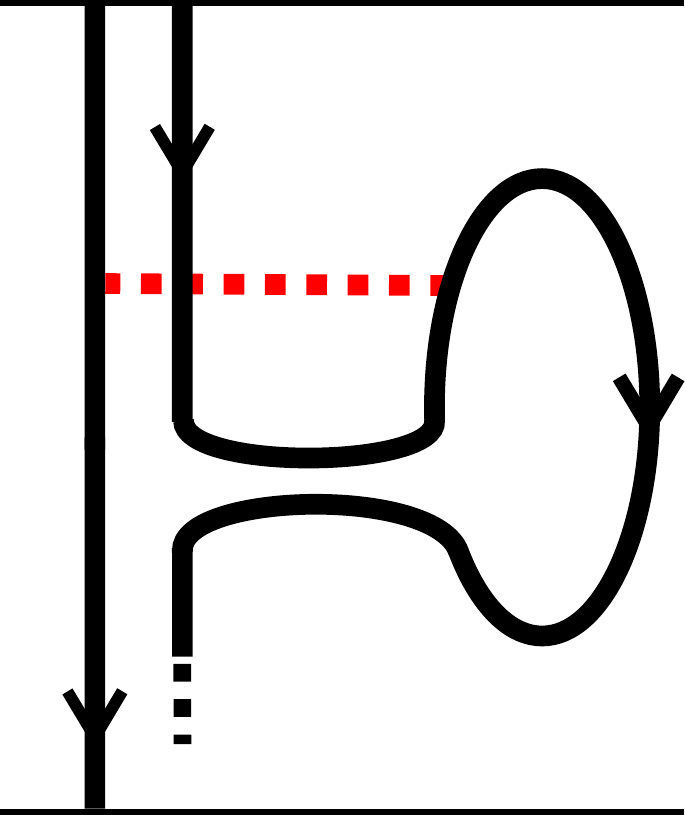}\end{array}
=\frac{h}{2}\begin{array}{c}\includegraphics[scale=0.15]{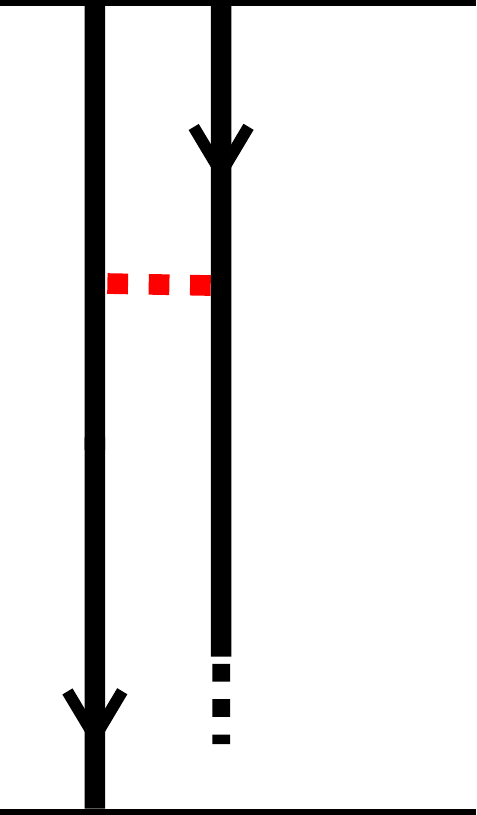}\end{array}=\frac{h}{2}Z_{\widetilde{\cA}}^J\left(\begin{array}{c}\includegraphics[scale=0.15]{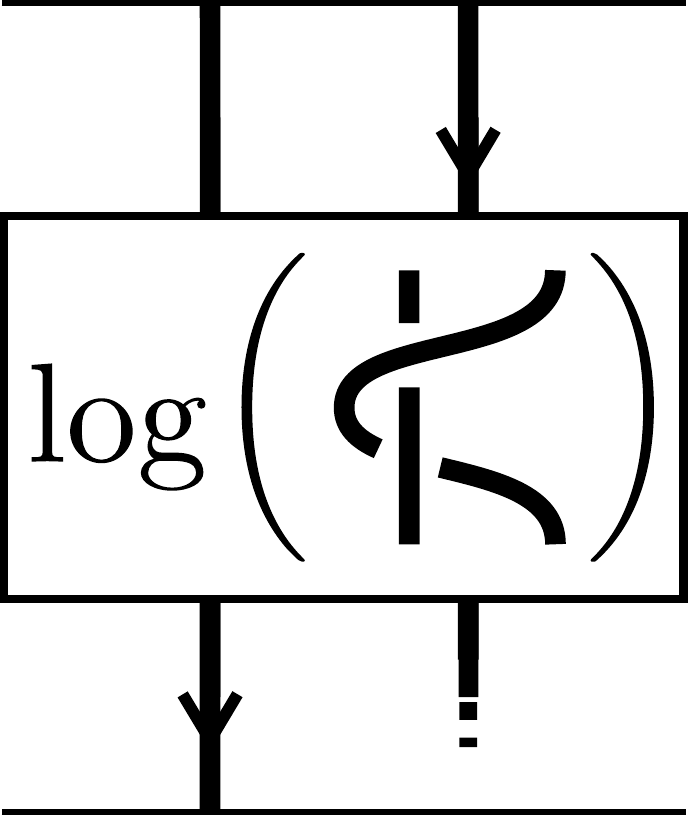}\end{array}\right)\label{eq:Dehn_twist01}\\
&= \frac{h}{2}  \sum_{n=1}^\infty \frac{(-1)^{n-1}}{n} (\exp(-\rho\hbar)\gamma - \gamma^0)^n 
= \frac{h}{2}  \log(\gamma^0 + \exp(-\rho \hbar)\gamma - \gamma^0) \label{eq:Dehn_twist02}\\
&= \frac{h}{2}  \log(\exp(-\rho\hbar)\gamma) 
= -\frac{h}{2}  \rho \hbar \gamma^0 + \frac{h}{2}  \log \gamma\nonumber
\end{align}
Equation~\eqref{eq:Dehn_twist01} is interpreted in $\mathcal{D}_{\widetilde{\cA}}^J$. From \eqref{eq:Dehn_twist02} onward, we use the isomorphism $Z_{\widetilde{\cA}}^J$ to work in $\widehat{\widetilde{\cA}(\Sigma,J)}$.
The term $-\frac{h}{2}\rho\hbar\gamma^0$ is $\tau$-invariant, while $\frac{h}{2}\log\gamma$ is $\tau$-anti-invariant.
Thus,
\begin{align}
\sigma(L_{\widetilde{\cA}})(\gamma^{0}) &= \frac{1}{h} \left( \frac{h}{2}  \log \gamma - \left(-\frac{h}{2}  \log \gamma\right) \right) = \log \gamma
\end{align}
This proves the Dehn twist formula on $\gamma^0$.

Next, we show that $\exp(\sigma(L_{\widetilde{\cA}}))$ agrees with the Dehn twist for general $J$.
The logarithm $\log t_C$ in \eqref{eq:Dehn_twist_C} satisfies the properties of Lemma~\ref{lem:sigma_Leibniz}. As in the proof of Lemma~\ref{lem:sigma_generators_gr}, it suffices to compare the actions on $\gamma^0$ and on tangles along the boundary of $S^1\times[0,1]$.
Both $\sigma(L_{\widetilde{\cA}})$ and $\log t_C$
vanish on tangles along the boundary of $S^1\times[0,1]$, and their actions on $\gamma^0$ agree by the computation above.
Hence, $\sigma(L_{\widetilde{\cA}})=\log t_C$ for general $J$.
Equivalently, $t_C=\exp(\sigma(L_{\widetilde{\cA}}))$.

It remains to prove uniqueness.
Let $L_{\widetilde{\cA}}'$ represent the same Dehn twist as $L_{\widetilde{\cA}}$ and satisfy $\iota^\ast(L_{\widetilde{\cA}}')=0$. 
In the rest of the proof, we abbreviate $F^{n}\widehat{\widetilde{\mathcal{A}}(S^1\times[0,1])}$ to $F^{n}$. 
We prove by induction that $L_{\widetilde{\cA}}-L_{\widetilde{\cA}}'\in F^i$ for every $i\in\bN$. The case $i=0$ is immediate.
Suppose that $L_{\widetilde{\cA}}-L_{\widetilde{\cA}}'\in F^n$.
The first component of its image under \eqref{eq:GrA_into_Hom} vanishes since $L_{\widetilde{\cA}}$ and $L_{\widetilde{\cA}}'$ represent the same Dehn twist.
The second component of its image under \eqref{eq:GrA_into_Hom} vanishes since $\iota^\ast(L_{\widetilde{\cA}})=\iota^\ast(L_{\widetilde{\cA}}')=0$.
Injectivity of \eqref{eq:GrA_into_Hom} implies $L_{\widetilde{\cA}}-L_{\widetilde{\cA}}'\in F^{n+1}$. Thus, by induction, $L_{\widetilde{\cA}}-L_{\widetilde{\cA}}'\in F^i$ for every $i\in\bN$.
This means that $L_{\widetilde{\cA}}$ and $L_{\widetilde{\cA}}'$ are equal in $\widehat{\widetilde{\mathcal{A}}(S^1\times[0,1])}$.
\end{proof}

\begin{remark}\label{rmk:independence_orientation_annulus}
\begin{enumerate}
    \item There are two choices for the orientation of the $S^1$-factor in this embedding. Since a Dehn twist about a simple closed curve does not depend on its orientation, $\sigma(L_{\widetilde{\cA}})$ is unchanged when the orientations of both $S^1$ and $[0,1]$ are reversed. The condition $\iota^\ast(L_{\widetilde{\cA}})=0$ also remains unchanged. Uniqueness therefore shows that the two embeddings give the same element.
    
    \item For the annulus $S^1\times[0,1]$, Turaev's result \cite{Tur88}, or \cite[Theorem 5.2]{Tur91}, identifies, after a suitable adjustment, a generating set of the completion $\widehat{\widetilde{\cA}(S^1\times[0,1])}$ without endpoints with a generating set of $\widehat{\widetilde{\cN_I}(S^1\times[0,1])}$.
\end{enumerate}
\end{remark}

\begin{remark}
By the uniqueness statement of Lemma~\ref{lem:Dehn_twist_annulus} and the result of \cite{Tsu18a}, we have
\begin{align}
L_{\widetilde{\cA}} = \frac{h}{2a} \begin{array}{c}\includegraphics[scale=0.15]{draws/bracelet_Dehn_twist.pdf}\end{array}= \left(\frac{\frac{\sqrt{h\hbar}}{2}}{\operatorname{arcsinh} \frac{\sqrt{h\hbar}}{2}}\right)^2 \left( \sum_{n=2}^{\infty} u_n \ell_n^{**} - \frac{1}{3} \rho^3 \hbar^2 \right)
=\left(\frac{h}{a}\right)^2 \left( \sum_{n=2}^{\infty} u_n \ell_n^{**} - \frac{1}{3} \rho^3 \hbar^2 \right),\label{eq:Dehn_twist_LA}
\end{align}
where $\ell_n^{**}$ and $u_n$ are defined as follows. In \cite[Theorem 5.2]{Tsu18a}, let $\ell_n$ be the link represented by a diagram winding $n$ times counterclockwise around the annulus with $n-1$ positive crossings, and set
$$
\ell_n^* = \sum_{k=1}^n \frac{(-\hbar)^{k-1}}{k} \sum_{i_1 + \dots + i_k = n} \ell_{i_1} \cdots \ell_{i_k},\quad \ell_0^*=2\rho
$$
Then, define
$$
\ell_n^{**} := \sum_{k=0}^n \frac{n!}{k!(n-k)!} (-1)^k \ell_{n-k}^*
$$
The paper \cite{Tsu18a} deals with the case $h=\hbar$; our $\ell_n^\ast$ is denoted by $\ell'(n)$ in \cite{Tsu18a}.
The coefficients $u_n$ are determined by
$$
\frac{1}{2}(\log x)^2 = \sum_{n=2}^\infty u_n (x-1)^n \in \mathbb{Q}\llbracket x-1\rrbracket
$$
where $x$ is a formal variable.
In \cite{Tsu18a}, convergence with respect to the filtration is established using an embedding of $\widehat{\bQ\pi^+}$ into $\widehat{\widetilde{\cA}(S^1\times[0,1],J)}$. The same straightforward argument shows that $\ell_n^{**}\in F^n$ when $h\neq\hbar$.

Although \cite{Tsu18a} considers only $h=\hbar$, preservation of component degree uniquely determines the case $h\neq\hbar$ from the result of \cite{Tsu18a}.
\end{remark}

The annulus case gives a corresponding result for a general surface $\Sigma$.
The following theorem can be viewed as a skein-module extension of the theories of Kawazumi--Kuno \cite{KK14,KK15} and Massuyeau--Turaev \cite{MT13,MT14}.
\begin{theorem}[cf. {\cite[Theorem 5.1]{Tsu18a}}]\label{thm:Dehn_twist_general}
Let $C\subset\Sigma$ be a simple closed curve.
Define $L_{\widetilde{\cA}}(C)$ as the image of $L_{\widetilde{\cA}}$ under the map induced by an embedding of $S^1\times[0,1]$ onto a neighborhood of $C$.
Then $L_{\widetilde{\cA}}(C)\in\widehat{\widetilde{\mathcal{A}}(\Sigma)}$ satisfies
$$t_C=\exp(\sigma(L_{\widetilde{\cA}}(C)))\in\Aut(\widehat{\widetilde{\mathcal{A}}(\Sigma,J)}).$$

\end{theorem}
\begin{proof}
Both the Dehn twist and $\sigma(L_{\widetilde{\cA}}(C))$ are determined by their actions in a neighborhood of $C$. This neighborhood is homeomorphic to an annulus, where the formula holds by Lemma~\ref{lem:Dehn_twist_annulus}.
\end{proof}

\begin{remark}
Independence on the choice of the orientation of the surface $\Sigma$ follows from Remark~\ref{rmk:independence_orientation_annulus}. 
\end{remark}

\end{document}